\documentclass[10pt,reqno]{smfart}
\usepackage{amsmath,amssymb,smfthm,stmaryrd,epic}
\usepackage[frenchb]{babel}
\usepackage[latin1]{inputenc}
\usepackage[all]{xy}
\usepackage{a4wide}

\usepackage{xypic}

\NumberTheoremsIn{subsection}\SwapTheoremNumbers

\begin{document}

\newtheorem{prop-defi}[smfthm]{Proposition-Définition}
\newtheorem{theo-defi}[smfthm]{Théorème-définition}
\newtheorem{lem-defi}[smfthm]{Lemme-définition}
\newtheorem{notas}[smfthm]{Notations}
\newtheorem{nota}[smfthm]{Notation}
\newtheorem{defis}[smfthm]{Définitions}
\newtheorem{remas}[smfthm]{Remarques}

\newtheorem{theob}{Théorème}[section]
\def\thetheob{\arabic{section}.\arabic{theob}}
\newtheorem{propb}[theob]{Proposition}
\newtheorem{lemb}[theob]{Lemme}
\newtheorem{corob}[theob]{Corollaire}
\newtheorem{defib}[theob]{Définition}
\newtheorem{defisb}[theob]{Définitions}
\newtheorem{remab}[theob]{Remarque}

\renewcommand{\theequation}{\Roman{part}.\arabic{section}.\arabic{subsection}.\arabic{equation}}

\def\Am{{\mathbb A}}
\def\Fm{{\mathbb F}}
\def\Mm{{\mathbb M}}
\def\Nm{{\mathbb N}}
\def\Pm{{\mathbb P}}
\def\Qm{{\mathbb Q}}
\def\Zm{{\mathbb Z}}
\def\Dm{{\mathbb D}}
\def\Cm{{\mathbb C}}
\def\Rm{{\mathbb R}}
\def\Gm{{\mathbb G}}

\def\AC{{\mathcal A}}
\def\CC{{\mathcal C}}
\def\DC{{\mathcal D}}
\def\EC{{\mathcal E}}
\def\FC{{\mathcal F}}
\def\GC{{\mathcal G}}
\def\HC{{\mathcal H}}
\def\IC{{\mathcal I}}
\def\JC{{\mathcal J}}
\def\KC{{\mathcal K}}
\def\LC{{\mathcal L}}
\def\MC{{\mathcal M}}
\def\NC{{\mathcal N}}
\def\OC{{\mathcal O}}
\def\PC{{\mathcal P}}
\def\UC{{\mathcal U}}
\def\VC{{\mathcal V}}
\def\XC{{\mathcal X}}
\def\YC{{\mathcal Y}}

\def\AF{{\mathfrak A}}
\def\GF{{\mathfrak G}}
\def\EF{{\mathfrak E}}
\def\CF{{\mathfrak C}}
\def\DF{{\mathfrak D}}
\def\JF{{\mathfrak J}}
\def\LF{{\mathfrak L}}
\def\MF{{\mathfrak M}}
\def\NF{{\mathfrak N}}
\def\XF{{\mathfrak X}}
\def\UF{{\mathfrak U}}
\def\KF{{\mathfrak K}}

\def \longmapright#1{\smash{\mathop{\longrightarrow}\limits^{#1}}}
\def \mapright#1{\smash{\mathop{\rightarrow}\limits^{#1}}}
\def \lexp#1#2{\kern \scriptspace \vphantom{#2}^{#1}\kern-\scriptspace#2}
\def \linf#1#2{\kern \scriptspace \vphantom{#2}_{#1}\kern-\scriptspace#2}
\def \linexp#1#2#3 {\kern \scriptspace{#3}_{#1}^{#2} \kern-\scriptspace #3}

\def \a{\alpha}
\def \b{\beta}
\def \d{\delta}
\def \e{\epsilon}
\def \g{\gamma}
\def \k{\kappa}
\def \l{\lambda}
\def \m{\mu}
\def \n{\nu}
\def \o{\omega}
\def \r{\rho}
\def \s{\sigma}
\def \t{\tau}
\def \th{\theta}
\def \u {\upsilon}
\def \x{\chi}
\def \vphi {\varphi}

\let \leq=\leqslant
\let \geq=\geqslant
\def \lefto{\longleftarrow}
\def \fin{\hfill $\square$}
\let \DS=\displaystyle
\let \SS=\scriptstyle
\let \longto=\longrightarrow
\let \oo=\infty

\def \FH{\mathop{\mathrm{FH}}\nolimits}
\def \FPH{\mathop{\mathrm{FPH}}\nolimits}
\def \coh{\mathop{\mathrm{Coh}}\nolimits}
\def \res{\mathop{\mathrm{res}}\nolimits}
\def \op{\mathop{\mathrm{op}}\nolimits}
\def \rec {\mathop{\mathrm{rec}}\nolimits}
\def \art{\mathop{\mathrm{Art}}\nolimits}
\def \hyp {\mathop{\mathrm{Hyp}}\nolimits}
\def \cusp {\mathop{\mathrm{Cusp}}\nolimits}
\def \Iw {\mathop{\mathrm{Iw}}\nolimits}
\def \JL {\mathop{\mathrm{JL}}\nolimits}
\def \speh {\mathop{\mathrm{Speh}}\nolimits}
\def \isom {\mathop{\mathrm{Isom}}\nolimits}
\def \Vect {\mathop{\mathrm{Vect}}\nolimits}
\def \groth {\mathop{\mathrm{Groth}}\nolimits}
\def \lef {\mathop{\mathrm{Lef}}\nolimits}
\def \fix {\mathop{\mathrm{Fix}}\nolimits}
\def \hom {\mathop{\mathrm{Hom}}\nolimits}
\def \deg {\mathop{\mathrm{deg}}\nolimits}
\def \val {\mathop{\mathrm{val}}\nolimits}
\def \det {\mathop{\mathrm{det}}\nolimits}
\def \rep {\mathop{\mathrm{Rep}}\nolimits}
\def \spec {\mathop{\mathrm{Spec}}\nolimits}
\def \fr {\mathop{\mathrm{Fr}}\nolimits}
\def \frob {\mathop{\mathrm{Frob}}\nolimits}
\def \ker {\mathop{\mathrm{Ker}}\nolimits}
\def \im {\mathop{\mathrm{Im}}\nolimits}
\def \Red {\mathop{\mathrm{Red}}\nolimits}
\def \red {\mathop{\mathrm{red}}\nolimits}
\def \aut {\mathop{\mathrm{Aut}}\nolimits}
\def \diag {\mathop{\mathrm{diag}}\nolimits}
\def \spf {\mathop{\mathrm{Spf}}\nolimits}
\def \Def {\mathop{\mathrm{Def}}\nolimits}
\def \twist {\mathop{\mathrm{Twist}}\nolimits}
\def \supp {\mathop{\mathrm{Supp}}\nolimits}
\def \Id {{\mathop{\mathrm{Id}}\nolimits}}
\def \bar {\overline}
\def \ind {\mathop{\mathrm{Ind}}\nolimits}
\def \mod {\mathop{\mathrm{mod}}\nolimits}
\def \ker {\mathop{\mathrm{Ker}}\nolimits}
\def \coker {\mathop{\mathrm{Coker}}\nolimits}
\def \mult {\mathop{\mathrm{mult}}\nolimits}
\def \vide{\emptyset}
\def \bad {{\mathop{\mathrm{Bad}}\nolimits}}
\def \gal {{\mathop{\mathrm{Gal}}\nolimits}}
\def \Nr {{\mathop{\mathrm{Nr}}\nolimits}}
\def \rn {{\mathop{\mathrm{rn}}\nolimits}}
\def \vol {{\mathop{\mathrm{vol}}\nolimits}}
\def \ad {{\mathop{\mathrm{ad}}\nolimits}}
\def \tr {{\mathop{\mathrm{Tr~}}\nolimits}}
\def \Sp {{\mathop{\mathrm{Sp}}\nolimits}}
\def \lie {{\mathop{\mathrm{Lie}}\nolimits}}
\def \st {{\mathop{\mathrm{Sp}}\nolimits}}
\def \sp{{\mathop{\mathrm{Sp}}\nolimits}}
\def \card{{\mathop{\mathrm{card}}\nolimits}}
\def \sym{{\mathop{\mathrm{Sym}}\nolimits}}
\def \perv{\mathop{\mathrm{Perv}}\nolimits}
\def \sh {{\mathop{\mathrm{Sh}}\nolimits}}
\def \const {{\mathop{\mathrm{Const}}\nolimits}}

\def \ele{élément }
\def \eles{éléments }
\def \cad{c'est à dire }
\def \rem{{\noindent\textit{Remarque:~}}}
\def \exem{{\noindent \textit{Exemple:~}}}
\def \ssi{~si et seulement si~}
\def \cl {{\mathop{\mathrm{cl}}\nolimits}}
\def \Tw {{\mathop{\mathrm{Tw}}\nolimits}}
\def \ob {{\mathop{\mathrm{Ob}}\nolimits}}
\def \ext {{\mathop{\mathrm{Ext}}\nolimits}}
\def \End {{\mathop{\mathrm{End}}\nolimits}}
\def \inv {{\mathop{\mathrm{inv}}\nolimits}}
\def \fix {{\mathop{\mathrm{Fix}}\nolimits}}

\def\semi{\mathrel{>\!\!\!\triangleleft}}


\setcounter{secnumdepth}{3} \setcounter{tocdepth}{3}

\newcommand{\marque}{\addtocounter{smfthm}{1}
{\smallskip \noindent \textit{\thesmfthm}~---~}}

\renewcommand\atop[2]{\ensuremath{\genfrac..{0pt}{1}{#1}{#2}}}

\title[Monodromie du faisceau pervers des cycles évanescents]{Monodromie du faisceau pervers des cycles évanescents
de
quelques variétés de Shimura simples et applications}

\alttitle{Monodromy of the perverse sheaf of vanishing cycles of some simple Shimura varieties and
applications}

\author[Boyer Pascal]{Boyer Pascal \\ et avec un appendice de Laurent Fargues \\ arxiv:math.AG/0511448}
\email{boyer@math.jussieu.fr}
\address{Institut de mathématiques de Jussieu \\ UMR 7586, université Paris 6 \\
175 rue du Chevaleret Paris 13}

\frontmatter

\begin{abstract} Dans la situation géométrique des variétés de Shimura simples de Kottwitz étudiées dans le livre
d'Harris et Taylor \cite{h-t},
on décrit la filtration de monodromie du complexe des cycles évanescents ainsi que la suite
spectrale correspondante. On prouve en particulier que cette filtration coincide avec celle donnée par les poids à un
décalage près. D'après le théorème de comparaison de Berkovich-Fargues, on en déduit une description de la
filtration de monodromie-locale du modèle de Deligne-Carayol. En application on obtient la description des
composantes locales des représentations automorphes cohomologiques de certains groupes unitaires, une
correspondance de Jacquet-Langlands globale entre deux tels groupes et une preuve de la conjecture de
monodromie-poids pour la cohomologie de ces variétés de Shimura.
\end{abstract}

\begin{altabstract} In the geometric situation of the simple Shimura varieties of Kottwitz studied in Harris and
Taylor's book \cite{h-t}, we describe
the monodromy filtration of the vanishing cycles complex and the spectral sequence associated to it. We prove in
particular
that this filtration coincides with the weight one up to shift. Thanks to the Berkovich-Fargues'
theorem, we deduce the description of the local monodromy filtration of the Deligne-Carayol model. In
application, we obtain the description of the local components of cohomological automorphic representations of
certains unitary groups, a global Jacquet-Langlands correspondence between two such groups and a proof of the
weight-monodromy conjecture for the cohomology of these Shimura varieties.
\end{altabstract}



\altkeywords{Shimura varieties, formal modules, Langlands correspondences, Jacquet-Langlands correspondences,
monodromy filtration, weight-monodromy conjecture, perverse sheaves, vanishing cycles}

\maketitle

\tableofcontents

\pagestyle{headings} \pagenumbering{arabic}

\section*{Introduction}

\renewcommand{\theequation}{\arabic{equation}}

\backmatter

\noindent \textbf{0.1.} --- Soit $K$ une extension finie de $\Qm_p$, d'anneau des entiers $\OC_K$. Pour un
entier $d$ strictement positif fixé, on considère le groupe $D_{K,d}^\times$ (resp. $W_K$) des éléments
inversibles de ``l'''algèbre à division centrale sur $K$ d'invariant $1/d$ (resp. le groupe de Weil de $K$).
Pour un nombre premier $l \neq p$, Langlands (resp. Jacquet-Langlands) a (resp. ont) conjecturé l'existence
d'une bijection $\rec_{K}$ (resp. $\JL$) entre les $\bar \Qm_l$-représentations irréductibles admissibles de
$GL_d(K)$ et les représentations $l$-adiques indécomposables de $W_K$ (resp. entre les représentations
admissibles irréductibles de $D_{K,d}^\times$ et les représentations essentiellement de carré intégrable de
$GL_d(K)$) qui sont compatibles à la formation des fonctions $L$; on renvoie à \cite{he} pour des énoncés
précis.

A l'aide de la cohomologie étale, Deligne a construit une série de représentations $\UC_{K,l,d}^{i}$ du
produit de ces trois groupes. Pour $d=2$ et pour $\rho$ une représentation irréductible de $D_{K,2}^\times$
tel que $\pi:=\JL(\rho)$ soit une représentation cuspidale de $GL_2(K)$, Carayol dans \cite{ca}, montre que
la composante $\rho$-isotypique $\UC_{K,l,2}^{1}(\rho)$ de $\UC_{K,l,2}^{1}$ réalise les correspondances de
Langlands et de Jacquet-Langlands, i.e.
$$\UC_{K,l,2}^{1}(\JL^{-1}(\pi^\vee)) \simeq \pi \otimes \rec_K(\pi^\vee)(-\frac{1}{2}).$$
Le cas $d$ quelconque est traité dans \cite{h-t} et dans \cite{boy} en
égale caractéristique. En outre pour $d=2$, Carayol décrit également ce qui se passe pour les autres
représentations.

Le but premier de ce travail est de faire de même pour $d$ quelconque, i.e. calculer complètement les
$\UC_{K,l,d}^{i}(\rho)$ pour $\rho$ une représentation irréductible quelconque de $D_{K,d}^\times$. Dans le
cas où $\rho$ est la représentation triviale, le résultat se formule comme suit.

\medskip

\noindent \textbf{Théorème 1} \textit{Pour $0 \leq i \leq d-1$, on a $\UC_{K,l,d}^{i}(1^\vee)=\pi_{i} \otimes
1(-i)$ où $\pi_{i}$ est l'unique quotient irréductible de l'induite parabolique
$$\ind_{P_{i,d}(K)}^{GL_{sg}(K)} \st_{i} \otimes 1 $$
où $P_{i,d}$ est le parabolique standard associée aux $i$ premiers vecteurs et $\st_i$ est la représentation
de Steinberg de $GL_i(K)$.}

\medskip

Le cas général, cf. le théorème (\ref{theo-ripsi-local}), s'énonce de manière similaire en faisant intervenir
les correspondances de Langlands et Jacquet-Langlands. Une autre formulation de notre résultat revient à dire
qu'il n'y a pas d'annulation dans l'expression, calculée dans \cite{h-t}, de la représentation virtuelle
$\sum_{i=0}^{d-1} (-1)^i [\UC_{K,l,d}^{i}(\rho)]$ où $\UC_{K,l,d}^{d-i}(\rho)$, pour $1 \leq i \leq d$, y est
alors donnée par le $i$-ème terme de plus haut poids.

\medskip

\noindent \textbf{0.2.} --- La preuve du théorème 1 procède par globalisation via l'étude des variétés de
Shimura étudiées dans \cite{h-t}, qui sont de type PEL, définies sur un corps CM, $F$, et associées à un
groupe de similitudes $G/\Qm$ tel que le groupe unitaire correspondant sur $\Rm$ soit de la forme $U(1,d-1)
\times U(0,d) \times \cdots \times U(0,d)$. Une place $v$ au dessus de $p$ est fixée de sorte que le complété
$F_v$ du localisé de $F$ en $v$ soit isomorphe au corps local précédemment noté $K$, d'anneau des entiers
$\OC_v$ et tel que $G(F_v) \simeq GL_d(F_v)$. A certains sous-groupes compacts $U^p(m)$ de $G(\Am^\oo)$, on
associe alors une variété de Shimura $X_{U^p,m} \to \spec \OC_v$ munie d'une action par correspondances de
$G(\Am^\oo)$, classifiant des variétés abéliennes munies de structures additionnelles; on renvoie au \S
\ref{rappels-globaux} pour les détails.

\medskip

\noindent \textbf{0.3.} --- La fibre spéciale $\bar X_{U^p,m}$ de $X_{U^p,m}$ est stratifiée par des
sous-schémas localement fermés $\bar X_{U^p,m}^{(d-h)}$ pour $1 \leq h \leq d$, de pure dimension $d-h$ tels
que le théorème de Serre-Tate donne un isomorphisme entre le complété de l'hensélisé strict de l'anneau local
de $X_{U^p,m}$ en un point géométrique de $\bar X_{U^p,m}^{(d-h)}$ et l'anneau $R_{F_v,d-h,m_1}$ qui
représente le foncteur des déformations de niveau $m_1$ d'un $\OC_v$-module de Barsotti-Tate sur $\bar \k(v)$
de hauteur $h$. Par ailleurs pour $1 \leq h < d$, il existe un sous-schéma fermé $\bar
X_{U^p,m,M_{d-h}}^{(d-h)}$ de $\bar X_{U^p,m}^{(d-h)}$ stable sous les correspondances associées aux \eles du
sous-groupe parabolique $P_{h,d}(F_v)$ de $GL_d(F_v)$ (cf. la définition (\ref{defi-parab})) et tel que
$$\bar X_{U^p,m}^{(d-h)}=\bar X_{U^p,m,M_{d-h}}^{(d-h)} \times_{P_{h,d}(\OC_v/\MC_v^{m_1})}
GL_d(\OC_v/\MC_v^{m_1}):$$
on dit que \textit{les strates non supersingulières sont géométriquement induites.}

\medskip

\noindent \textbf{0.4.} --- Pour un nombre premier $l$ distinct de $p$, les auteurs de \cite{h-t}
introduisent sur chacune des strates $\bar X_{U^p,m,M_{d-h}}^{(d-h)}$, des systèmes locaux $\FC_{\tau_v}$
associés aux représentations irréductibles $\tau_v$ de $D_{v,h}^\times$. Ils décrivent alors la restriction
des cycles évanescents $R^i \Psi_v:=R^i\Psi_{\eta_v}(\bar \Qm_l)$ à la strate $\bar X_{U^p,m}^{(d-h)}
\times_{\Fm_q} \bar \Fm_q$ en fonction des $\FC_{\t_v}$ et des cycles évanescents locaux
$\Psi_{F_v,l,h,n}^{i}$, cf. (\ref{iso1}). D'après le théorème de comparaison de Berkovich, le théorème local
se déduit alors de la connaissance de la fibre en un point supersingulier des $R^i\Psi_v$.

\medskip

\noindent \textbf{0.5.} --- Le complexe $R\Psi_v[d-1]$ est vu comme un faisceau pervers muni d'une filtration
de monodromie dont on notera $gr_k$ les gradués. Notre deuxième résultat concerne la description des $gr_k$
dans la catégorie des faisceaux pervers sur $\bar X_{U^p,m}$ munis d'une action compatible de $G(\Am^\oo)
\times W_v$. Pour $1 \leq tg \leq d$, et pour $\pi_v$ (resp. $\Pi_t$) une représentation irréductible
cuspidale de $GL_g(F_v)$ (resp. quelconque de $GL_{tg}(F_v)$), on introduit le faisceau $HT(g,t,\pi_v,\Pi_t)$
sur la strate $\bar X_{U^p,m}^{(d-tg)}$, "induit" à partir du système local
$\FC_{\JL^{-1}(\st_t(\pi_v)^\vee)} \otimes \Pi_t$ sur la composante $\bar X_{U^p,m,M_{d-tg}}^{(d-tg)}$. Les
composantes $\pi_v$-isotypiques $gr_{k,\pi_v}$ des $gr_k$, cf. la proposition (\ref{prop-so}), se décrivent
alors, théorème (\ref{theo-global1}), au moyen des faisceaux pervers $j^{\geq tg}_{!*}
HT(g,t,\pi_v,\st_t(\pi_v))[d-tg]$ où $j^{\geq tg}$ désigne l'injection de la strate $\bar
X_{U^p,m}^{(d-tg)}$; en ce qui concerne le cas $\pi_v$ triviale, l'énoncé est le suivant:

\medskip

\noindent\textbf{Théorème 2} \textit{Les faisceaux pervers $gr_{k,1_v}$ sont nuls pour $k \geq d$ et sinon on
a}
$$gr_{k,1_v}=\bigoplus_{\genfrac{}{}{0pt}{}{|k| < t \leq d}{t \equiv k-1 \mod 2}} j^{\geq t}_{!*} HT(1,t,1_v,\st_t)
(-\frac{t+k-1}{2})[d-t]$$

\medskip

L'énoncé pour $\pi_v$ quelconque est similaire et fait intervenir les correspondances de Langlands et
Jacquet-Langlands.

\medskip

\noindent \textbf{0.6.} --- En utilisant la conjecture de monodromie-poids, la preuve du théorème 2
découlerait de la connaissance du théorème 1 pour toutes les hauteurs $h < d$ ainsi que de la description des
restrictions aux strates des faisceaux des cycles évanescents en fonction des systèmes locaux $\FC_{\tau_v}$
comme rappelé en (0.4.). Par ailleurs le théorème 1 en hauteur $d$, découle d'après le théorème de
comparaison de Berkovich-Fargues, du calcul des germes aux points supersinguliers des faisceaux de
cohomologie des $gr_k$.

On raisonne alors par récurrence en supposant connus\footnote{En fait on suppose plutôt connu les gradués
locaux $gr_{h,k,loc}$ de la filtration de monodromie-locale du complexe des cycles évanescents
$\Psi_{F_v,l,h}^\bullet$.} les $\UC_{F_v,l,h}^{i}$ du modèle de Deligne-Carayol de hauteur $h$ pour tout $1
\leq h <d$. On en déduit alors le théorème 2 sauf en ce qui concerne les faisceaux pervers supportés aux
points supersinguliers i.e. on ne sait pas à quel $k$ associer chacun des $j^{\geq d}_! HT(1,d,1_v,\st_d)(-
\frac{d+k'-1}{2})$ pour $|k'| < d$ et $k' \equiv d-1 \mod 2$. De même on sait déterminer tous les faisceaux
de cohomologie $h^i gr_k$ des $gr_k$ en dehors des points supersinguliers. La technique repose sur l'étude de
la suite spectrale associée à la filtration de monodromie:

\begin{equation} \label{ss-intro}
E_1^{i,j}=h^{i+j} gr_{-i} \Rightarrow R^{i+j+d-1}\Psi_v
\end{equation}
en essayant d'en deviner les termes initiaux, l'aboutissement étant connu. En outre en utilisant la
perversité des $gr_k$ ainsi que la compatibilité de $R\Psi_v$ à la dualité de Verdier, on obtient un contrôle
sur les germes aux points supersinguliers des $h^i gr_k$. On étudie ensuite la suite spectrale des cycles
évanescents dont à nouveau on essaie de deviner les termes initiaux alors que l'aboutissement est connu, en
utilisant en particulier le théorème de Lefschetz difficile. Le contrôle obtenu précédemment nous permet
alors de prouver le théorème 1. On prouve ensuite la conjecture de monodromie-poids faisceautique en
"plaçant" chacun des $j^{\geq d}_! HT(1,d,1_v,\st_d)(- \frac{d+k-1}{2})$ sur $gr_k$ et on explicite la suite
spectrale de monodromie (\ref{ss-intro}).

\medskip

\noindent \textbf{0.7.} --- En ce qui concerne les résultats globaux que l'on obtient, citons

\begin{itemize}
\item la description ``explicite'' des gradués pour la filtration de monodromie du faisceau pervers des cycles
évanescents
ainsi que de la suite spectrale associée;

\item la pureté de la filtration de monodromie du complexe des cycles évanescents;

\item la détermination des extensions intermédiaires des systèmes locaux d'Harris-Taylor;

\item le calcul de tous les groupes de cohomologie des différents faisceaux ou complexes de faisceaux qui sont
introduits.
\end{itemize}

\medskip

\noindent \textbf{0.8.} --- Comme conséquence directe de ces calculs on obtient une description des
composantes locales des représentations automorphes cohomologiques de $G_\tau(\Am)$ ainsi qu'une
correspondance de Jacquet-Langlands globale entre deux formes intérieures de groupes unitaires qui complète
\cite{H-L} \S 3. A partir de cette description, on prouve alors la conjecture de monodromie-poids, i.e. la
pureté de la filtration de monodromie des groupes de cohomologie de la fibre générique des variétés de
Shimura étudiées.

\medskip

\noindent Signalons que récemment Taylor et Yoshida ont indépendamment obtenu dans \cite{y-t}, la pureté de
la filtration de monodromie de la cohomologie globale dans le cas tempéré, i.e. dans le cas où la cohomologie
est concentré en degré médian. Pour cela ils se ramènent par changement de base au cas Iwahori et utilisent
la suite spectrale de Rapoport-Zink. Dans notre cas, nous avons un dévissage qui permet, en utilisant aussi
un changement de base, de se ramener au cas de $GL_2$ traité par Carayol dans \cite{ca}. Par ailleurs en ce
qui concerne la conjecture de monodromie-poids cohomologique, le cas des variétés de Shimura uniformisées par
les demi-espaces de Drinfeld est traité dans \cite{ito} à partir de la version faisceautique qui découle de
l'existence d'un modèle semi-stable.

\bigskip

Décrivons succinctement le contenu des divers paragraphes.

\noindent \textbf{0.9.} --- \textit{La première partie} est consacrée à des rappels sur les données
géométriques locales et globales tirées de \cite{h-t}. On commence par des rappels de \cite{ze} sur les
induites paraboliques, \S \ref{defi-induite}, et les foncteurs de Jacquet \S \ref{jacquet}.

Le modèle local de Deligne-Carayol est introduit au paragraphe \ref{rapel-DC}. En ce qui concerne les données
globales, \S \ref{rappels-globaux}, outre le fait que les strates non supersingulières soient géométriquement
induites, la donnée fondamentale est celle des systèmes locaux d'Harris-Taylor $\FC(g,t,\pi_v)$ sur la strate
$tg$, attachés aux représentations $\JL^{-1}(\st_t(\pi_v)^\vee)$ des inversibles $D_{v,tg}^\times$ de
l'algèbre à division centrale sur $F_v$ d'invariant $1/tg$. La propriété essentielle que l'on utilisera sur
ces systèmes locaux est l'isomorphisme (\ref{iso1}).

On redonne en outre la description de l'ensemble des points supersinguliers et on définit pour tout diviseur
$g$ de $d=sg$, le faisceau $\FC(g,s,\pi_v)$ à support sur les points supersinguliers. On notera par ailleurs
que les $\FC(g,t,\pi_v)$ ne sont pas à priori irréductibles car ils proviennent de la restriction à
$\DC_{v,tg}^\times$ de $\JL^{-1}(\st_t(\pi_v)^\vee)$.

\medskip

\noindent \textbf{0.10.} --- Le paragraphe \ref{rappels-coho} rappelle les résultats d'ordre cohomologique de
\cite{h-t}. Au niveau local, \S \ref{rap-loc}, on explicite le lien entre $\Psi_{F_v,l,h}^{i}$ et
$\UC_{F_v,l,h}^{i}$ et on introduit l'entier $e_{\pi_v}$ qui est le cardinal de la classe d'équivalence
inertielle, définition (\ref{defi-inert}), de $\pi_v$. En particulier $e_{\pi_v}$ est égal au nombre de
sous-représentations irréductibles de la restriction à $\DC_{o,tg}^\times$ de $\JL^{-1}(\st_t(\pi_v)^\vee)$.

De \cite{h-t}, on retient principalement la détermination des parties cuspidales des groupes de cohomologie
des modèles locaux de Deligne-Carayol ainsi que le calcul, cf. le théorème (VII.1.5) de \cite{h-t}, de la
somme alternée, dans le groupe de Grothendieck des $GL_d(F_v) \times W_v$-modules
$$\sum_i (-1)^i \UC_{F_v,l,d}^{i}(\JL^{-1}(\st_s(\pi_v)^\vee)).$$
En global, \S \ref{ht-prop}, on exploite l'isomorphisme (\ref{iso1}) et notamment la description de l'action,
par la proposition (\ref{prop-hic}). On introduit selon \cite{lrs}, pour toute représentation irréductible
algébrique $\xi$ de $G(\Qm)$, le système local $\LC_{\xi}$ sur les variétés $X_{U^p,m}$ et on considère les
groupes de cohomologie du produit tensoriel de ce dernier avec un système local d'Harris-Taylor.

Le résultat fondamental qui résulte des arguments de comptage de points est la proposition
(\ref{prop-somme-alternee}) qui calcule la somme alternée des groupes de cohomologie à support compact des
systèmes locaux d'Harris-Taylor, dans le groupe de Grothendieck des $G^{(h)}(\Am^\oo) \times \Zm$-modules, où
pour tout $h$, $\Zm$ est identifié au quotient $D_{v,h}^\times/ \DC_{v,h}^\times$ via la valuation de la
norme réduite. On introduit pour cela le caractère $\Xi$ de $\Zm \longto \bar \Qm_l^\times$ défini par
$\Xi(1)=\frac{1}{q}$. Ainsi les représentations automorphes $\Pi$ qui interviennent dans cette description
ainsi que dans celle de la cohomologie de la fibre générique à valeurs dans $\LC_{\xi}$, vérifient des
conditions $\hyp(\oo)$ aux place infinies que l'on donne à la définition (\ref{defi-hyp}).

\medskip

\noindent \textbf{0.11.} --- \textit{La deuxième partie} donne l'énoncé des théorèmes locaux et globaux. Le
théorème (\ref{theo-ripsi-local}) est la motivation initiale de ce travail, c'est à dire décrire, pour tout
diviseur $g$ de $d$ et toute représentation irréductible cuspidale $\pi_v$ de $GL_g(F_v)$, chacun des
$\UC_{F_v,l,d}^{i}(\JL^{-1}(\st_s(\pi_v)^\vee))$. Finalement on obtient un résultat plus précis, théorème
(\ref{theo-local-fil}), qui est la description des gradués de la filtration de monodromie-locale définie par
Fargues en appendice et de la suite spectrale correspondante.

\medskip

\noindent \textbf{0.12.} --- On énonce ensuite, \S \ref{enonces-glob}, les résultats globaux. On commence,
après avoir fait quelques rappels sur les faisceaux pervers \S \ref{intro}, par découper, proposition
(\ref{prop-so}), nos faisceaux pervers selon leur composantes isotypiques pour le sous-groupe d'inertie $I_v$
et on note ainsi $gr_{k,\pi_v}$ (resp. $gr_{k,\xi,\pi_v}$) la composante
$\rec_{F_v}^\vee(\pi_v)_{|I_v}$-isotypiques du gradué $gr_k$ (resp. $gr_{k,\xi}$) pour la filtration de
monodromie du faisceau pervers $R\Psi_v[d-1]$ (resp. $(R\Psi_v \otimes \LC_{\xi})[d-1]$ pour $\xi$ une
représentation irréductible algébrique de $G(\Qm)$), où $\pi_v$ est une représentation irréductible cuspidale
de $GL_g(F_v)$ avec $1 \leq g \leq d$. La représentation $\rec_{F_v}^\vee(\pi_v)_{|I_v}$ n'étant pas
irréductible, pour un $I_v$-module $V$, $V_{\pi_v}$ est le facteur direct de $V$ sur lequel l'action de $I_v$
se décompose en une somme de sous-représentations irréductibles qui sont aussi des sous-représentations de
$\rec_{F_v}^\vee(\pi_v)_{|I_v}$.

On introduit ensuite \S \ref{defi-fph}, certaines catégories de faisceaux pervers de Hecke qui fourniront le
cadre catégoriel des différents complexe de faisceaux que nous considérerons dans la suite et on donne,
définition (\ref{defi-type}), un certain nombre de notations attachées aux systèmes locaux d'Harris-Taylor
$\FC(g,t,\pi_v)_1$ et aux faisceaux induits $HT(g,t,\pi_v,\Pi_t)$ qui leurs sont associés sur la strate $\bar
X_{U^p,m}^{(d-tg)}$ où $\Pi_t$ est une représentation de $GL_{tg}(F_v)$ qui sera le plus souvent elliptique
de type $\pi_v$. On introduit alors le système projectif $\GF$ des groupes de Grothendieck des faisceaux
pervers de Hecke sur la tour des $\bar X_{U^p,m}$ munis d'une action compatible de $G(\Am^\oo) \times W_v$.

On énonce alors les théorèmes globaux qui précisent, théorème (\ref{theo-global0}) (resp.
(\ref{theo-global1})), les faisceaux pervers simples de $(R\Psi_v \otimes \LC_{\xi})[d-1]$ (resp. des
$gr_{k,\xi,\pi_v}$) en termes des faisceaux pervers

$$j^{\geq tg}_{!*} HT(g,t,\pi_v,\st_t(\pi_v))[d-tg] \otimes \rec_{F_v}^\vee(\pi_v) (-\frac{tg-1+k}{2})$$
avec $1 \leq t \leq d/g$
(resp. $|k| < t \leq s_g$ et $t \equiv k-1 \mod 2$). D'après le théorème de comparaison de Berkovich-Fargues,
le théorème local (\ref{theo-ripsi-local}) se déduit alors du calcul, théorème (\ref{theo-global2}), des
faisceaux de cohomologie des $j^{\geq tg}_{!*} \FC(g,t,\pi_v)$ et de la détermination, théorème
(\ref{theo-ss}), de la suite spectrale de monodromie associée.

On donne ensuite le schéma de la preuve qui procède par récurrence en supposant connu (\ref{theo-local-fil})
pour tout $d'<d$. On renvoie le lecteur à \S \ref{schema} pour un apercu des implications logiques entre les
divers énoncés locaux et globaux. Avant de se lancer dans la preuve proprement dite, on donne des
conséquences de l'hypothèse de récurrence, sur les représentations obtenues en prenant la cohomologie des
systèmes locaux d'Harris-Taylor.

\medskip

\noindent \textbf{0.13.} --- \textit{La troisième partie} s'occupe des preuves des théorèmes locaux et
globaux. On démontre en premier lieu, \S \ref{glob0}, le théorème (\ref{theo-global0}), i.e. on donne l'image
de $R\Psi_v[d-1]$ dans le système projectif $\GF$ des groupes de Grothendieck des faisceaux pervers de Hecke
sur la tour des $\bar X_{U^p,m}$.

La preuve procède en plusieurs étapes. Tout d'abord de la description (\ref{iso1}) des restrictions aux
strates $\bar X_{U^p,m}^{(d-h)}$ des $R^i\Psi_v$ et du calcul (\ref{somme-alternee}) de $\sum_i (-1)^i
[\UC_{F_v,l,d}^{i}]$, on en déduit, proposition (\ref{prop-libre}), l'égalité suivante où on a posé
$s_g=\lfloor \frac{d}{g} \rfloor$, la partie entière de $d/g$:

\begin{multline} \label{egaliteR}
[R\Psi_v[d-1]]=\sum_{g=1}^d \sum_{\pi_v \in \cusp_v(g)} \frac{1}{e_{\pi_v}}
\sum_{i=1}^{s_g} \sum_{t=i}^{s_g} (-1)^{t-i} \\
[j^{\geq tg}_! HT(g,t,\pi_v, [\overleftarrow{i-1},\overrightarrow{t-i}]_{\pi_v})[d-tg] \otimes
\rec_{F_v}^\vee(\pi_v) (-\frac{tg-2+2i-t}{2})]
\end{multline}

L'étape suivante consiste alors à exprimer l'image des faisceaux pervers qui interviennent dans le membre de
droite de l'égalité ci-dessus soit:

\begin{multline} \label{egalitej}
[j^{\geq tg}_! HT(g,t,\pi_v,[\overleftarrow{t-1}]_{\pi_v})[d-tg]]= \\ \sum_{i=0}^{s_g-t} j^{\geq (t+i)g}_{!*}
HT(g,t+i,\pi_v,[\overleftarrow{t-1}]_{\pi_v} \overrightarrow{\times} [\overleftarrow{i-1}]_{\pi_v})[d-(t+i)g]
\otimes \Xi^{i(g-1)/2}
\end{multline}
Pour prouver cette dernière égalité, on raisonne par récurrence descendante sur la dimension des supports des
faisceaux pervers simples dans $j^{\geq tg}_! HT(g,t,\pi_v,[\overleftarrow{t-1}]_{\pi_v})[d-tg]$ en
réinjectant (\ref{egalitej}) pour les différents $t$, avec $g$ et $\pi_v$ fixés, dans (\ref{egaliteR}). On
argumente tout d'abord sur le fait que le membre de droite de (\ref{egaliteR}) doit s'écrire comme une somme
à coefficients positifs de faisceaux pervers simples, autoduale pour la dualité de Verdier. On invoque
ensuite le théorème de comparaison de Berkovich-Fargues afin d'utiliser l'hypothèse de récurrence sur les
modèles locaux en hauteur $h<d$, afin d'obtenir des renseignements sur les germes aux points non
supersinguliers de ces faisceaux pervers simples.

On démontre alors, proposition (\ref{prop-p}), le résultat hors des points supersinguliers au sens où les
égalités du théorème précédent et de (\ref{egalitej}) sont vraies si on rajoute une somme alternée de
faisceaux pervers à support sur les points supersinguliers. Le fait est qu'on utilise vraiment (\ref{iso1})
et pas seulement un calcul de somme alternée des restrictions des faisceaux des cycles évanescents, ce qui
explique l'indétermination au niveau des points supersinguliers.

\medskip

\noindent \textbf{0.14.} --- Le paragraphe (\ref{fp-ponctuel}) est consacré à la détermination de ces
faisceaux pervers ponctuels qui nous manquent. Une idée naive est que pour connaître un faisceau ponctuel, on
peut commencer par calculer son groupe de cohomologie $H^0$. Étant donné un tel faisceau pervers ponctuel
$\PC$ à déterminer, nous verrons en fait que la connaissance de son $H^0$ suffit à le déterminer
complètement: en effet grâce à (\ref{egalite-groth}), nous montrerons que $\PC$ contient $\FC(g,s,\pi_v)
\otimes (\Pi_l \overrightarrow{\times} [\overleftarrow{s-t-1}]_{\pi_v}) \otimes \Xi^{(s-t)(g-1)/2}$ avec
$$H^0(\FC(g,s,\pi_v)\otimes\Pi_t \overrightarrow{\times} [\overleftarrow{s-t-1}]_{\pi_v}) \otimes \Xi^{(s-t)(g-1)/2}
=H^0(\PC)$$ de sorte que $\PC =\FC(g,s,\pi_v) \otimes \Pi_t \overrightarrow{\times}
[\overleftarrow{s-t-1}]_{\pi_v} \otimes \Xi^{(s-t)(g-1)/2}$.

On commence alors par calculer les groupes de cohomologie des faisceaux pervers d'Harris-Taylor, ou tout du
moins leur $\Pi^{\oo,v}$-partie, pour $\Pi$ automorphe vérifiant $\hyp(\oo)$ avec $\Pi_v=\st_s(\pi_v)$. On
montre à la proposition (\ref{prop-coho1}) que ceux-ci sont alors tous nuls de sorte que, d'après la
proposition (\ref{prop-p}), l'égalité (\ref{somme-alternee}) fournit, corollaire (\ref{coro-pnul}), le $H^0$
des faisceaux ponctuels manquant ainsi que leur détermination. Le théorème (\ref{theo-global0}) découle alors
directement de ces résultats, cf. le corollaire (\ref{global0-ok}).

\medskip

\noindent \textbf{0.15.} --- Le huitième paragraphe est consacré, proposition (\ref{prop-hij}), au calcul des
faisceaux de cohomologies des faisceaux pervers $j^{\geq tg}_{!*} HT(g,t,\pi_v,\Pi_l)[d-tg]$. D'après la
proposition (\ref{prop-p}), on peut procéder par récurrence en utilisant la suite spectrale
(\ref{suite-spectrale}) dont l'aboutissement est connu sauf au niveau des points supersinguliers. Ainsi ces
faisceaux de cohomologie ne sont pas complètement déterminés aux points supersinguliers mais il ne reste
qu'un nombre réduit de possibilités, cf. le lemme (\ref{lem-hij3}), que l'on peut obtenir, d'après
(\ref{prop-p}), par récurrence en utilisant par exemple la suite spectrale (\ref{ss-dualite}).

\medskip

\noindent \textbf{0.16.} --- Au neuvième paragraphe, on prouve les théorèmes globaux sous la proposition
(\ref{cas-m}). D'après Berkovich-Fargues la connaissance du théorème local (\ref{theo-local-fil}) revient à
celle des germes aux points supersinguliers des faisceaux de cohomologie des faisceaux pervers
d'Harris-Taylor ainsi que les flèches correspondantes dans la suite spectrale (\ref{suite-spectrale}). Si
nous disposions du théorème (\ref{theo-local-fil}) pour $d$, le raisonnement de la preuve de la proposition
(\ref{prop-hij1}) nous permettrait de déterminer complètement les faisceaux de cohomologie des $gr_k$. Par
souci d'efficacité nous montrons, proposition (\ref{prop-hij-ss-g}), qu'il suffit, en utilisant l'opérateur
de monodromie $N$, en fait de connaître les parties de poids $s(g-1)$ de l'aboutissement de la suite
spectrale de monodromie-locale de (\ref{theo-local-fil}). On est alors ramené à prouver la proposition
(\ref{cas-m}), i.e. à déterminer les parties de poids $s(g-1)$ des $\UC_{F_v,l,d}^{i}$.

On donne ensuite, proposition (\ref{prop-grk-final}), une preuve de la conjecture de monodromie-poids version
faisceautique. On remarque tout d'abord que pour $s>2$, le raisonnement par récurrence implique le résultat
de sorte qu'il reste à initialiser la récurrence et donc à prouver que pour $s=2$, la monodromie n'est pas
triviale, en particulier sur la cohomologie. Le résultat est connu pour $s=2$ et $g=1$ d'après \cite{ca}. Le
principe, qui nous a été suggéré par M. Harris, est de s'y ramener via un changement de base adéquat.

\medskip

\noindent \textbf{0.17.} --- Le dixième paragraphe est consacré à la preuve de la proposition (\ref{cas-m}).
Celle-ci repose sur l'étude de la suite spectrale des cycles évanescents. Pour $\Pi$ une représentation
irréductible automorphe de $D_\Am^\times$ telle que $\Pi_v \simeq \st_s(\pi_v)$ avec $\pi_v$ cuspidale, la
$\Pi^{\oo,v}$-partie de son aboutissement est connue d'après \cite{h-t} et \cite{y-t}, ou peut-être
rapidement recalculée à partir de la proposition (\ref{prop-coho1}). La détermination de cet aboutissement
nous restreint, corollaire (\ref{coro-1}), alors le nombre de possibilités pour les
$\UC_{F_v,l,d}^{i}(\JL^{-1}(\st_s(\pi_v)^\vee))$ ce qui nous permet de conclure en utilisant une propriété
d'invariance des $\UC_{F_v,l,d}^\bullet$ sous l'involution de Zelevinski, cf. le théorème (\ref{involution}).

Cependant la preuve de celle-ci repose sur tout ou partie du théorème de comparaison de Faltings à partir
d'un énoncé similaire du coté espace de Drinfeld, ce qui dépasse le cadre de ce texte. Par ailleurs la preuve
de la conjecture de monodromie-poids version cohomologique demande une étude des $\Pi^{\oo}$-parties des
divers groupes de cohomologie, pour $\Pi$ automorphe vérifiant $\hyp(\oo)$ et tel que $\Pi_v \simeq
\speh_s(\pi_v)$. La $\Pi^{\oo,v}$-partie de la cohomologie de la fibre générique n'est pas à priori connue.
On se propose dans un premier temps de la calculer, ou tout du moins les bouts de poids $s(g-1)$, le cas
général étant traité, de manière indépendante, à la proposition (\ref{prop-lrs2}).

On calcule tout d'abord, proposition (\ref{prop-not}), les $\Pi^{\oo,v}$-parties des groupes de cohomologie
des faisceaux pervers d'Harris-Taylor. On en déduit alors, corollaire (\ref{coro-ssce-min}), la connaissance
des $\Pi^{\oo,v}$-parties de poids $s(g-1)$ des groupes de cohomologie de la fibre générique. Pour ce faire
on utilise le théorème de Lefschetz difficile. On étudie alors, proposition (\ref{prop-ssce-poids}), les
$\Pi^{\oo,v}$-parties des termes $E_2^{p,q}$ de la suite spectrale des cycles évanescents. On montre, qu'à
travers les suites spectrales associées à la stratification, qu'en ce qui concerne les bouts de poids
$s(g-1)$, seule la strate supersingulière contribue de sorte que l'on se retrouve dans une situation
similaire à celle de \cite{boy} dans le cas cuspidal, où les bouts de poids $s(g-1)$ des
$E_2^{p,q}[\Pi^{\oo,v}]$ sont nuls pour $p \neq 0$. Cette constatation découle du contrôle, lemme
(\ref{lem-rj-combi}), des contributions des strates non supersingulières. La proposition (\ref{cas-m}),
découle alors directement de la connaissance des $\Pi^{\oo,v}$-parties de l'aboutissement de la suite
spectrale des cycles évanescents et du fait qu'en ce qui concerne celles de poids $s(g-1)$, celles-ci ne
proviennent que des points supersinguliers.

\medskip

\noindent \textbf{0.18.} --- \textit{La quatrième partie} donne des compléments et applications des résultats
précédents. On commence, proposition (\ref{prop-lrs2}), par calculer les $\Pi^{\oo,v}$-parties des groupes de
cohomologie de la fibre générique pour $\Pi$ automorphe tel que $\Pi_v \simeq \speh_s(\pi_v)$, alors que nous
n'avions traité, corollaire (\ref{coro-ssce-min}), que celles de poids $s(g-1)$.

On donne ensuite une correspondance de Jacquet-Langlands globale qui complète les résultats de \cite{H-L} \S
3. On décrit ensuite, proposition (\ref{prop-compo-locale}), les composantes locales des représentations
automorphes et cohomologiques de $G_\tau(\Am)$.

On conclut enfin, théorème (\ref{theo-mono-glob}), par la preuve de la conjecture de monodromie-poids
globale, i.e. les gradués de la filtration de monodromie des groupes de cohomologie de la fibre générique
sont purs.

\medskip

\noindent \textbf{0.19.} --- Enfin dans \textit{cinquième partie}, on récapitule les différents résultats
obtenus tout au long du texte et dans \textit{sixième} on illustre par des figues les différentes suites
spectrales étudiées.


\bigskip

\noindent \textit{Remerciements} Je tiens à exprimer ma profonde gratitude envers Gérard Laumon pour son
soutien constant tout au long de ces années ainsi que pour sa perspicacité à dépister les fausses bonnes
idées et à mettre en avant les autres. Merci aussi à Laurent Fargues de m'avoir fourni une bonne notion de
monodromie locale qui m'a permis de simplifier grandement les preuves; un deuxième merci pour ses nombreuses 
explications
sur les variétés de Shimura et les groupes unitaires. J'adresse des remerciements à
Jean-François Dat pour m'avoir fourni l'habillage catégoriel des faisceaux de Hecke ainsi que pour ses
nombreux conseils. Enfin je remercie Michael Harris pour son aide ainsi que l'IHES pour son accueil.


\mainmatter

\renewcommand{\theequation}{\arabic{section}.\arabic{subsection}.\arabic{equation}}

\part{Notations et rappels}

\section{Rappels sur les représentations de $GL_d(F_v)$}
\label{rappel-ze-0}

Tous les énoncés ainsi que leur preuve peuvent être trouvés dans \cite{ze}.

\subsection{Induites paraboliques}
\label{defi-induite}

\begin{defi} \label{defi-parab}
Pour une suite $0< r_1 < r_2 < \cdots < r_k=d$, on note $P_{r_1,r_2,\cdots,r_k}$ le
sous-groupe parabolique
de $GL_d$ standard associé au sous-groupe de Levi $GL_{r_1}(F_v) \times
GL_{r_2-r_1}(F_v) \times \cdots
\times GL_{r_k-r_{k-1}}(F_v)$ et on note $N_{r_1,\cdots,r_k}$ son radical unipotent.
\end{defi}

\marque Soient $\pi_1$ et $\pi_2$ des représentations de respectivement $GL_{n_1}(F_v)$
et $GL_{n_2}(F_v)$;
on note selon la coutume, $\pi_1 \times \pi_2$ l'induite parabolique
$\ind_{P_{n_1,n_1+n_2}(F_v)}^{GL_{n_1+n_2}(F_v)} \pi_1(n_2/2) \otimes \pi_2(-n_1/2)$.

\rem Le symbole $\times$ est associatif, i.e. $\pi_1 \times (\pi_2 \times \pi_3)=(\pi_1
\times \pi_2) \times
\pi_3$ que l'on notera donc $\pi_1 \times \pi_2 \times \pi_3$.

\begin{defis} Soit $g$ un diviseur de $d=sg$ et $\pi_v$ une représentation cuspidale
irréductible de $GL_g(F_v)$:

- les sous-quotients irréductibles de
$V(\pi_v,s):=\pi_v(\frac{1-s}{2}) \times \pi_v(\frac{3-s}{2}) \times \cdots \times
\pi_v(\frac{s-1}{2})$
seront dits elliptiques de type $\pi_v$;

- $V(\pi_v,s)$ possède un unique quotient (resp. sous-espace) irréductible que l'on
notera
$[\overleftarrow{s-1}]_{\pi_v}$ (resp.  $[\overrightarrow{s-1}]_{\pi_v})$); c'est une
représentation de
Steinberg (resp. de Speh) généralisée notée habituellement $\st_s(\pi_v)$ (resp.
$\speh_s(\pi_v)$);

- pour $\pi_1$ et $\pi_2$ des représentations respectivement de $GL_{t_1g}(F_v)$ et
$GL_{t_2g}(F_v)$, on
notera
$\pi_1 \overrightarrow{\times} \pi_2$ (resp. $\pi_1 \overleftarrow{\times} \pi_2$)
l'induite parabolique $\pi_1(-t_2/2) \times \pi_2(t_1/2)$ (resp. $\pi_1(t_2/2) \times
\pi_2(-t_1/2)$).
\end{defis}

La propriété habituelle de transitivité des induites paraboliques donne le lemme
suivant.

\begin{lemm} Pour $(\pi_i)_{1 \leq i \leq 3}$ des représentations de $GL_{t_ig}(F_v)$,
on a les égalités suivantes:

$$\begin{array}{rl}
(\pi_1 \overrightarrow{\times} \pi_2)^\vee \simeq & \pi_1^\vee \overleftarrow{\times}
\pi_2^\vee \\
(\pi_1 \overrightarrow{\times} \pi_2) \overrightarrow{\times_g} \pi_3 = & \pi_1
\overrightarrow{\times}
(\pi_2
\overrightarrow{\times} \pi_3) \\
(\pi_1 \overleftarrow{\times} \pi_2) \overleftarrow{\times_g} \pi_3= & \pi_1
\overleftarrow{\times} (\pi_2
\overleftarrow{\times} \pi_3)
\end{array}$$
En outre si $\pi_1$ et $\pi_2$ sont elliptiques de type $\pi_v$, il en est de même de
$\pi_1
\overrightarrow{\times} \pi_2$ et donc de $\pi_1 \overleftarrow{\times} \pi_2$.
\end{lemm}

\rem Si $\pi_1$, $\pi_2$ et $\pi_3$ sont des représentations elliptiques de type $\pi_v$
on a en outre
$$\pi_1 \overrightarrow{\times} (\pi_2 \overleftarrow{\times} \pi_3)= \pi_2
\overleftarrow{\times} (\pi_1
\overrightarrow{\times} \pi_3)$$ En effet d'après \cite{ze}, on sait que $\pi_1 \times
\pi_2 \simeq \pi_2
\times \pi_1$ si les supports cuspidaux de $\pi_1$ et $\pi_2$ sont disjoints
\footnote{En fait loc. cit.
donne un critère plus précis; dans le cas où $\pi_1$ et $\pi_2$ sont irréductibles, il
faut et il suffit que
leurs supports cuspidaux ne soient pas liés.}.

\begin{prop-defi} Pour $g$ un diviseur de $d=sg$ et $\pi_v$ une représentation
irréductible cuspidale de
$GL_g(F_v)$, on a
pour $1 \leq t \leq s$:

\begin{itemize}
\item[-] $[\overleftarrow{t-1}]_{\pi_v} \overrightarrow{\times}
[\overrightarrow{s-t-1}]_{\pi_v}$ (resp.
$[\overrightarrow{t-1}]_{\pi_v} \overrightarrow{\times}
[\overleftarrow{s-t-1}]_{\pi_v}$) est de longueur
$2$; on notera
$[\overleftarrow{t-1},\overrightarrow{s-t}]_{\pi_v}$, (resp.
$[\overrightarrow{t},\overleftarrow{s-t-1}]_{\pi_v}$) son
unique sous-espace irréductible, et
$[\overleftarrow{t},\overrightarrow{s-t-1}]_{\pi_v}$ (resp.
$[\overrightarrow{t-1},\overleftarrow{s-t}]_{\pi_v}$) son
unique quotient irréductible;

\item[-] par dualité $[\overleftarrow{t-1}]_{\pi_v} \overleftarrow{\times}
[\overrightarrow{s-t-1}]_{\pi_v}$ (resp.
$[\overrightarrow{t-1}]_{\pi_v} \overleftarrow{\times} [\overleftarrow{s-t-1}]_{\pi_v}$)
est de longueur $2$
avec $[\overrightarrow{s-t-1},\overleftarrow{t}]_{\pi_v}$ (resp.
$[\overleftarrow{s-t},\overrightarrow{t-1}]_{\pi_v}$)
pour unique sous-espace irréductible et
$[\overrightarrow{s-t},\overleftarrow{t-1}]_{\pi_v}$ (resp.
$[\overleftarrow{s-t-1},\overrightarrow{t}]_{\pi_v}$) pour
unique quotient irréductible;

\item[-] on notera $\lfloor \pi \rfloor$ (resp. $\lceil \pi \rceil$) l'unique, s'il
existe, sous-espace (resp.
quotient)
irréductible de $\pi$. Pour $\pi_1$ et $\pi_2$ des représentations irréductibles
elliptiques de type $\pi_v$,
$\pi_1 \overrightarrow{\times} \pi_2$ et $\pi_1 \overleftarrow{\times} \pi_2$ ont un
unique sous-espace et un
unique quotient irréductible et on a:
$$\lfloor \pi_1 \overrightarrow{\times} \pi_2 \rfloor = \lceil \pi_2^\vee
\overleftarrow{\times} \pi_1^\vee
\rceil^\vee \qquad
\lceil \pi_1 \overrightarrow{\times} \pi_2 \rceil = \lfloor \pi_2^\vee
\overleftarrow{\times} \pi_1^\vee
\rfloor^\vee$$

Pour $\pi_3$ une troisième représentation irréductible elliptique de type $\pi_v$, on a
$$\lfloor \lfloor \pi_1 \overrightarrow{\times} \pi_2 \rfloor \overrightarrow{\times}
\pi_3 \rfloor=\lfloor
\pi_1 \overrightarrow{\times} \lfloor \pi_2 \overrightarrow{\times} \pi_3 \rfloor
\rfloor=\lfloor \pi_1
\overrightarrow{\times} \pi_2 \overrightarrow{\times} \pi_3 \rfloor$$
$$\lceil \lceil \pi_1 \overrightarrow{\times} \pi_2 \rceil \overrightarrow{\times}
\pi_3 \rceil=\lceil \pi_1 \overrightarrow{\times} \lceil \pi_2 \overrightarrow{\times}
\pi_3 \rceil
\rceil=\lceil \pi_1 \overrightarrow{\times} \pi_2 \overrightarrow{\times} \pi_3
\rceil$$

\item[-] Soient $r \geq 1$ et $\Gamma^s=(a_i,\e_i)_{1 \leq i \leq r}$ tel que les $a_i$
sont des entiers strictement
positifs
avec $s-1=a_1+\cdots + a_r$ et $\e_i= \pm 1$; on notera $\Gamma^s$ sous la forme
$(\overleftarrow{a_1},
\cdots , \overrightarrow{a_r})$ où pour tout $i$ la flèche au dessus de $a_i$ est
$\overleftarrow{a_i}$
(resp. $\overrightarrow{a_i}$) si $\e_i=-1$ (resp. $\e_i=1$). On associe à $\Gamma^s$ un
sous-quotient
irréductible $[\Gamma^s]$ de $V(s,\pi_v)$ que l'on note aussi sous la forme
$[\overleftarrow{a_1}, \cdots ,
\overrightarrow{a_r}]_{\pi_v}$. On convient par ailleurs des égalités suivantes:
$$[\cdots,\overleftarrow{a},\overleftarrow{b},\cdots]_{\pi_v}=[\cdots,\overleftarrow{a+b},\cdots]_{\pi_v}
\qquad
[\cdots,\overrightarrow{a},\overrightarrow{b},\cdots]_{\pi_v}=[\cdots,\overrightarrow{a+b},\cdots]_{\pi_v}$$
D'après \cite{ze}, on obtient alors une bijection, modulo les identifications ci-dessus,
entre les
sous-quotients irréductibles de $V(s,\pi_v)$ et l'ensemble des $[\Gamma^s]$ telle que
l'on ait les relations
suivantes, où $\overleftrightarrow{c}$ désigne arbitrairement $\overleftarrow{c}$ ou
$\overrightarrow{c}$:
$$\lfloor [\cdots,\overleftrightarrow{a}]_{\pi_v} \overrightarrow{\times}
[\overleftrightarrow{b},\cdots]_{\pi_v}\rfloor
=
[\cdots,\overleftrightarrow{a},\overrightarrow{1},\overleftrightarrow{b},\cdots]_{\pi_v}$$
$$\lceil [\cdots,\overleftrightarrow{a}]_{\pi_v} \overrightarrow{\times}
[\overleftrightarrow{b},\cdots]_{\pi_v}\rceil
=
[\cdots,\overleftrightarrow{a},\overleftarrow{1},\overleftrightarrow{b},\cdots]_{\pi_v}$$
$$\lfloor [\cdots,\overleftrightarrow{a}]_{\pi_v} \overleftarrow{\times}
[\overleftrightarrow{b},\cdots]_{\pi_v}\rfloor
=
[\cdots,\overleftrightarrow{b},\overleftarrow{1},\overleftrightarrow{a},\cdots]_{\pi_v}$$
$$\lceil [\cdots,\overleftrightarrow{a}]_{\pi_v} \overleftarrow{\times}
[\overleftrightarrow{b},\cdots]_{\pi_v}\rceil
=
[\cdots,\overleftrightarrow{b},\overrightarrow{1},\overleftrightarrow{a},\cdots]_{\pi_v}$$

\end{itemize}
\end{prop-defi}

\marque Dans \cite{ze}, $(\overleftarrow{a_1}, \cdots , \overrightarrow{a_r})$ est
représenté sous la forme
$\circ^{\frac{1-s}{2}} \lefto \circ \cdots \circ \lefto \cdots \circ \to
\circ^{\frac{s-1}{2}}$
où les $a_1$ premières flèches vont dans le sens de la flèche au dessus de $a_1$, les
$a_2$ suivantes dans le
sens de la flèche au dessus de $a_2$ ... Ainsi graphiquement on a
$$[\Gamma^{s_1}_1] \overrightarrow{\times} [\Gamma^{s_2}_2]=\overbrace{\circ \lefto
\circ \cdots \to
\circ}^{\Gamma^{s_1}_1} \longleftrightarrow \overbrace{\circ \to \circ \cdots
\circ}^{\Gamma^{s_2}_2}, \qquad
[\Gamma^{s_1}_1] \overleftarrow{\times} [\Gamma_2^{s_2}]= \overbrace{\circ \lefto \circ
\cdots \to
\circ}^{\Gamma_2^{s_2}} \longleftrightarrow \overbrace{\circ \to \circ \cdots
\circ}^{\Gamma_1^{s_1}}$$

\marque \noindent \textit{Exemples}: la représentation cuspidale irréductible $\pi_v$ de
$GL_g(F_v)$ étant
fixé, avec nos notations, on a:
$$\begin{array}{rl}
{[}\overleftarrow{s-1}{]}_{\pi_v} &= \lfloor \overbrace{[\overrightarrow{0}]_{\pi_v}
\overleftarrow{\times}
\cdots \overleftarrow{\times} [\overrightarrow{0}]_{\pi_v}}^{s} \rfloor = \lceil
\overbrace{[\overleftarrow{0}]_{\pi_v} \overrightarrow{\times} \cdots
\overrightarrow{\times}
[\overleftarrow{0}]_{\pi_v}}^s
\rceil \\
{[}\overrightarrow{s-1}{]}_{\pi_v} &= \lfloor \overbrace{[\overleftarrow{0}]_{\pi_v}
\overrightarrow{\times}
\cdots \overrightarrow{\times} [\overleftarrow{0}]_{\pi_v}}^s \rfloor =\lceil
\overbrace{[\overrightarrow{0}]_{\pi_v} \overleftarrow{\times} \cdots
\overleftarrow{\times}
[\overrightarrow{0}]_{\pi_v}}^s
 \rceil \\
{[} \overleftarrow{t-1},\overrightarrow{s-t} {]}_{\pi_v} &= \lfloor
[\overleftarrow{t-1}]_{\pi_v}
\overrightarrow{\times} [\overrightarrow{s-t-1}]_{\pi_v} \rfloor =\lceil
[\overleftarrow{t-2}]_{\pi_v}
\overrightarrow{\times}
[\overrightarrow{s-t}]_{\pi_v} \rceil \\
&= \lceil [\overrightarrow{s-t-1}]_{\pi_v} \overleftarrow{\times}
[\overleftarrow{t-1}]_{\pi_v} \rceil=
\lfloor [\overrightarrow{s-t}]_{\pi_v} \overleftarrow{\times}
[\overleftarrow{t-2}]_{\pi_v} \rfloor
\end{array}$$

\noindent \textbf{Notation}: \textit{Dans la suite $[\overleftrightarrow{s-1}]_{\pi_v}$
désignera une
représentation elliptique quelconque de type $\pi_v$ de $GL_{sg}(F_v)$.}

\subsection{Foncteur de Jacquet}
\label{jacquet}

Soit $P=MN$ un parabolique de $GL_d$ de Lévi $M$ et de radical unipotent $N$.

\begin{defi}
Pour $\pi$ une représentation admissible de $GL_d(F_v)$, l'espace des vecteurs
$N(F_v)$-coinvariants est
stable sous l'action de $M(F_v) \simeq P(F_v)/N(F_v)$. On notera $J_N(\pi)$ cette
représentation tordue par
$\d_P^{-1/2}$.
\end{defi}

\begin{lemm} Soit $g$ un diviseur de $d=sg$ et $\pi_v$ une représentation irréductible
cuspidale de
$GL_g(F_v)$. Pour $1 \leq h \leq d$, le foncteur de Jacquet vérifie les propriétés
suivantes:
\begin{itemize}
\item si $g$ ne divise pas $h$, alors
$$J_{N_{h,d}^{op}}([\overleftarrow{s-1}]_{\pi_v})=J_{N_{h,d}^{op}}([\overrightarrow{s-1}]_{\pi_v})=
J_{N_{h,d}}([\overleftarrow{s-1}]_{\pi_v})=J_{N_{h,d}}([\overrightarrow{s-1}]_{\pi_v})
=(0)$$

\item si $h=tg$ alors
$$J_{N_{tg,sg}}([\overleftarrow{s-1}]_{\pi_v})=[\overleftarrow{t-1}]_{\pi_v((s-t)/2)}
\otimes
[\overleftarrow{s-t-1}]_{\pi_v(-t/2)}$$
$$J_{N_{tg,sg}}([\overrightarrow{s-1}]_{\pi_v})=[\overrightarrow{t-1}]_{\pi_v((t-s)/2)}
\otimes
[\overrightarrow{s-t-1}]_{\pi_v(t/2)}$$
$$J_{N_{tg,sg}^{op}}([\overleftarrow{s-1}]_{\pi_v})=[\overleftarrow{t-1}]_{\pi_v((t-s)/2)}
\otimes
[\overleftarrow{s-t-1}]_{\pi_v(t/2)}$$
$$J_{N_{tg,sg}^{op}}([\overrightarrow{s-1}]_{\pi_v})=[\overrightarrow{t-1}]_{\pi_v((s-t)/2)}
\otimes
[\overrightarrow{s-t-1}]_{\pi_v(-t/2)}$$
\end{itemize}
\end{lemm}

\begin{proof} soit $\Gamma=(a_i,\e_i)_{1 \leq i \leq r}$, les résultats découlent alors
des propriétés suivantes
que l'on trouve dans \cite{ze}:

- $J_{N_{g,2g,\cdots,sg}^{op}}([\Gamma^s])$ est de la forme $\pi_v(\frac{1-s}{2}+\s(0))
\otimes \cdots
\pi_v(\frac{1-s}{2}+\s(s-1))$ où $\s$ est une permutation de l'ensemble $\{ 0,\cdots,s-1
\}$ soumise à la
règle suivante: soit $1 \leq i \leq r$ avec $\e_i=1$ (resp. $\e_i=-1$): pour tout
$a_1+\cdots + a_{i-1} \leq
r < r' \leq a_1+ \cdots +a_i$ alors $\s^{-1}(r) < \s^{-1}(r')$ (resp. $\s^{-1}(r) >
\s^{-1}(r')$).\footnote{Autrement dit $\s$ est compatible aux orientations des
flèches.}

- en ce qui concerne $J_{N_{g,2g,\cdots,sg}}([\Gamma^s])$ la règle est inversée, i.e.
$\s^{-1}(r)>
\s^{-1}(r')$ (resp. $\s^{-1}(r) < \s^{-1}(r')$).

\end{proof}

\section{Rappels des données géométriques}

Soit $K$ une extension finie de $\Qm_p$, on note $\OC_K$ son anneau des entiers, $\PC_K$
l'idéal maximal et
$\kappa=\OC_K/\PC_K$ le corps résiduel de cardinal $q=p^f$. On notera aussi $\varpi_K$
une uniformisante.
L'extension maximal non ramifiée de $K$ sera notée $K^{nr}$ de complété $\hat K^{nr}$,
d'anneau des entiers
respectifs $\OC_{K^{nr}}$ et $\OC_{\hat K^{nr}}$.

\begin{defisb}
- Soit $\art_K^{-1}:W_K \longto K^\times$ le morphisme de la théorie du corps de classe
qui envoie les
frobenius géométriques $\fr$ du groupe de Weil $W_K$ de $K$ sur les uniformisantes, i.e.
$\val(\art_K^{-1}
(\fr))=-1$. Pour $c \in W_K$, on notera $\deg(c):=\val(\art^{-1}(c))$.

- Étant donnés une représentation $\s$ (resp. $\pi$) de $W_K$ (resp. de $GL_d(K)$) et un
entier $r$, on
notera $\s(r)$ (resp. $\pi(r)$) la représentation $\s \otimes |\art_K^{-1}|^{r}$ (resp.
$\pi \otimes |\det
|^{r}$).
\end{defisb}

\subsection{Le modèle local de Deligne-Carayol}
\label{rapel-DC}

On considère dans la suite la définition restrictive suivante de $\OC_K$-module de
Barsotti-Tate, cf. par exemple
\cite{h-t} II.1.

\begin{defi} Soit $S$ un $\OC_K$-schéma dans lequel $p$ est localement nilpotent. Un
$\OC_K$-module de Barsotti-Tate
sur $S$, est un groupe de Barsotti-Tate $H/S$ muni d'une injection $i:\OC_K
\hookrightarrow \End(H)$
vérifiant la propriété suivante. L'algèbre de Lie de $H$, $\lie(H)$, est un faisceau
localement libre sur $S$
muni de deux actions de $\OC_K$ définies soit par $i$ soit par le morphisme structural
$\OC_K \to \OC_S$; on
demande alors que ces deux actions coincident et que $\lie(H)$ soit de dimension $1$.
\end{defi}

\marque On rappelle alors les propriétés suivantes issues de la conjonction des
propriétés générales sur les
groupes de Barsotti-Tate et de la théorie des modules elliptiques de Drinfeld:
\begin{itemize}
\item pour tout $d \geq 1$, il existe un unique $\OC_K$-module de Barsotti-Tate formel
$\Sigma_{K,d}$ sur
$\bar \k$ de hauteur $d$ caractérisé par les propriétés suivantes sur une loi de groupe
$F$ et sa structure
de $\OC_K$-module $x \in \OC_K \mapsto \vphi_x(X) \in \bar \k [[X]]$: $F(X,Y) \equiv X+Y
\mod (X,Y)^{q^d}$,
$\vphi_{\varpi_K}(X)=X^{q^d}$ et $\vphi_{\{ a \} }(X) \equiv aX \mod X^{q^d}$, pour tout
$a \in \k$ où $\{ a
\} \in \OC_K$ désigne le représentant de Teichmuller de $a$.

\item tout $\OC_K$-module de Barsotti-Tate sur $\bar \k$ est de la forme $\Sigma_{K,g}
\times (K/\OC_K)^h$.
\end{itemize}

\begin{prop} \label{defi-p-strate} (cf. \cite{h-t})
Soit $S$ un schéma localement noethérien sur lequel $p$ est localement nilpotent et soit
$H/S$
un $\OC_K$-module de Barsotti-Tate. Alors pour tout $h \geq 0$, il existe un sous-schéma
fermé réduit
$S^{[h]}$ de $S$ tel que:

\begin{itemize}
\item $S^{[h-1]} \subset S^{[h]}$ et la codimension de toute composante de $S^{[h-1]}$
dans $S^{[h]}$
est au plus $1$;

\item pour tout point géométrique $s$ de $S$, $s$ appartient à $S^{[h]}$ si et seulement
si
$\sharp H[p](\k^{sep}(s)) \leq p^{[K:\Qm_p]h};$

\item sur $S^{(h)}:=S^{[h]}-S^{[h-1]}$, on a une suite exacte courte de $\OC_K$-modules
de Barsotti-Tate
$0 \to H^0 \longto H \longto H^{et} \to 0$ où $H^0$ est formel et $H^{et}$ est ind-étale
de hauteur $h$.
\end{itemize}
\end{prop}

\begin{defi} Soit $C$ la catégorie des $\OC_K$-algèbres locales, artiniennes, de corps
résiduel $\bar \k$.
Pour un objet $R$ de $C$, une déformation sur $R$ d'un $\OC_K$-module de Barsotti-Tate
$H_0/\k$ est un
couple $(H,f)$ où $H$ est un $\OC_K$-module de Barsotti-Tate sur $R$ et $f:H_0
\longmapright{\sim} H \times_R
\bar \k$.
\end{defi}

\begin{prop} (cf. \cite{dr1}) Le foncteur qui à un objet $R$ de $C$ associe l'ensemble
des classes d'isomorphismes
des
 déformations sur $R$ de $\Sigma_{K,d}$ est pro-représentable par un anneau local
complet noethérien
 $R_{K,d}$ de corps résiduel $\bar \k$. En outre il existe un isomorphisme
 $$R_{K,d} \simeq \OC_{\hat K^{nr}} [[T_1,\cdots,T_{d-1}]]$$
 tel que la déformation universelle $(\widetilde{\Sigma_{K,d}},\tilde f)$ sur $R_{K,d}$
est caractérisée par
 les propriétés suivantes sur une loi de groupe formel $F$ et sa structure de
$\OC_K$-module $x \mapsto
 \vphi_x$:

\begin{itemize}

 \item $F(X,Y) \equiv X+Y \mod (X,Y)^{q}$;

 \item $\vphi_{\varpi_K}(X) \equiv \varpi_K u_0 X + T_1 u_1 X^q+\cdots + T_{d-1} u_{d-1}
X^{q^{d-1}} +u_d X^{q^d}$,
où
$u_1,\cdots,u_d$ sont des unités de $\OC_{\hat
K^{nr}}[[T_1,\cdots,T_{d-1}]][[X]]^\times$, cf. \cite{Lubin1}
p106;

 \item $\vphi_{\{ a \} } (X) \equiv aX \mod X^2$ pour tout $a \in \k$.
\end{itemize}
\end{prop}

\begin{defi} Une structure de niveau $n$ sur un $\OC_K$-module de Barsotti-Tate $H/S$ de
hauteur constante
$h$ sur un schéma $S$, est un morphisme de $\OC_K$-modules:
$$\iota_n:(\PC_K^{-n}/\OC_K)^d \longto H[\PC_K^n](S)$$
tel que $\{ \iota_n(x)~:~ x \in (\PC_K^{-n}/\OC_K)^d \}$ soit un système complet de
section de la
$\PC_K^n$-torsion de $H/S$.
\end{defi}

\begin{prop} \label{defi-M-strate} (cf. \cite{h-t})
Soit $H/S$ un $\OC_K$-module de Barsotti-Tate de hauteur constante $h$ sur $S$ connexe
et tel
que l'on ait une suite exacte courte de $\OC_K$-module de Barsotti-Tate:
$$0 \to H^0 \longto H \longto H^{et} \to 0$$
avec $H^0$ formel et $H^{et}$ ind-étale. Un morphisme
$\iota_m:(\PC_K^{-m}/\OC_K)^h \longto H[\PC_K^m](S)$
est alors une structure de niveau de Drinfeld si et seulement s'il existe un facteur
direct $M$ de
$(\PC_K^{-m}/\OC_K)^h$ tel que:

\begin{itemize}
\item $\iota_{m,|M}:M \longto H^0[\PC_K^m](S)$ est une structure de niveau de Drinfeld;

\item $\iota_m$ induit un isomorphisme de $S$-schémas en groupes
$\iota_m^{et}:((\PC_K^{-m}/\OC_K)^h~/~M)_S \longto H^{et}[\PC_K^m]$.
\end{itemize}
\end{prop}

\begin{prop} (cf. \cite{dr1}) Le foncteur qui à un objet $R$ de $C$ associe l'ensemble
des classes d'isomorphismes
des
 déformations sur $R$ de $\Sigma_{K,d}$ munies d'une structure de niveau $n$ est
pro-représentable par un anneau
local
 complet, noethérien, régulier $R_{K,d,n}$ de corps résiduel $\bar \k$ tel que, en
notant
 $\widetilde{\iota_n}$ la structure de niveau universelle sur $R_{K,d,n}$,
$e_1,\cdots,e_d$ la base
 canonique du $\OC_K/\PC_K^n$-module $(\PC_K^{-n}/\OC_K)^d$, et $X$ un paramètre de
 $\widetilde{\Sigma_{K,d}}/R_{K,d}$, alors les $X(\widetilde{\iota_n}(e_i))$ forment un
système de paramètres
 locaux de $R_{K,d,n}$.
 \end{prop}

\begin{defi} Soit $\Psi_{K,l,d,n}^{i}$ le $\bar \Qm_l$-espace vectoriel de dimension
finie associé,
via la théorie des cycles évanescents de Berkovich, au morphisme structural
$\spf R_{K,d,n} \longto \spf \hat \OC_K^{nr}.$
\end{defi}

\marque Cet espace vectoriel est muni entre autre d'une action de $GL_d(\OC_K)$ qui se
factorise par le
morphisme surjectif naturel $GL_d(\OC_K) \longto GL_d(\OC_K/\MC_K^n)$ et on pose
$\Psi_{K,l,d}^{i}= {\DS \lim_{\atop{\longto}{n}}} \Psi_{K,l,d,n}^{i}$
de sorte que pour $\KF_{n}:=\ker (GL_d(\OC_K) \longto GL_d(\OC_K/\PC_K^n))$,
$\Psi_{K,l,d,n}^{i}=(\Psi_{K,l,d}^{i})^{\KF_{n}}$. On introduit le groupe $\NC_K$ (resp.
$\NC_K'$) défini
comme le noyau de
$$(g,\d,c) \in GL_d(K) \times D_{K,d}^\times \times W_K \mapsto
\val(\det(g^{-1})\rn(\d)\art_K^{-1}(c))
\in \Zm$$ (resp. composé avec la projection canonique $\Zm \longto \Zm/d\Zm$), où
$D_{K,d}$ est l'algèbre à
division centrale sur $K$ d'invariant $1/d$ et d'ordre maximal $\DC_{K,d}$.

\marque Pour $\xi$ un caractère d'ordre fini de $K^\times$, on note
$\Psi_{K,l,d,\xi}^{i}$ la
$\xi'$-composante isotypique où $\xi'$ est la restriction de $\xi$ à $\OC_K^\times$.
Ainsi
$\Psi_{K,l,d,\xi}^{i}$ (resp. $\Psi_{K,l,d}^{i}$) est muni d'une action de $\NC_K'$
(resp. de $\NC_K$).

\marque Dans la définition de $R_{K,d,n}$, il est agréable de considérer plutôt les
déformations par
quasi-isogénies ce qui donne un schéma formel $\spf R_{K,d,n,\Zm} \simeq \coprod_\Zm
\spf R_{K,d,n}$ de sorte
que la construction précédente fourni des $\bar \Qm_l$-espaces vectoriels
$\UC_{K,l,d,n}^{i} \simeq
(\UC_{K,l,d}^{i})^{\KF_{n}}$ où
$$\UC_{K,l,d,\xi}^{i} \simeq \ind_{\NC_K}^{GL_d(K) \times D_{K,d}^\times \times
W_{K}}\Psi_{K,l,d,\xi}^{i}$$
est une représentation de $GL_d(K) \times D_{K,d}^\times \times W_{K}$.

Pour toute représentation admissible irréductible $\t$ de $D_{K,d}^\times$ de caractère
central $\xi$, la
réciprocité de Frobenius donne que la composante isotypique
$\UC_{K,l,d,\xi}^{i}(\tau_v)$ est isomorphe à
$\hom_{\DC_{K,d}^\times}(\res^{D_{K,d}^\times}_{\DC_{K,d}^\times}
\t,\Psi_{K,l,d,\xi^{-1}}^{i})$ où l'action
de $(g^c,\s)$ est donnée par celle de $(g,\d,\s) \in \NC_K'$ pour $\d \in
D_{K,d}^\times$ quelconque.

\subsection{Variétés de Shimura simples et systèmes locaux}
\label{rappels-globaux}

Commençons par introduire les variétés de Shimura simples associées à des groupes
unitaires telles qu'elles
sont définies dans \cite{h-t}.

\marque Soit tout d'abord une donnée de type PEL $(F,B,*,V,<,>)$ définie comme suit:

\begin{itemize}

\item soit $E$ une extension quadratique imaginaire pure de $\Qm$ dans laquelle $p$ se
décompose en $u$ et
$u^c$ où $c \in \gal(E/\Qm)$ désigne la conjugaison complexe. Soit alors $F^+$ une
extension totalement
réelle de $\Qm$ de degré $e$ dont on fixe un plongement $\tau:F^+ \hookrightarrow \Rm$.
On pose alors $F=F^+
E$ et on note $v=v_1,\cdots,v_r$ les places de $F$ au dessus de $u$; $F_v$ désigne alors
le complété du
localisé en $v$ de $F$, d'anneau des entiers $\OC_{v}$ et de corps résiduel $\k(v)
\simeq \k_q$ avec
$q=p^{f_1}$.

\item $B$ est une algèbre à division centrale sur $F$ de dimension $d^2$ telle que:

\begin{itemize}
\item $\dim_F B=d^2$;

\item son algèbre opposée $B^{op}$ est isomorphe à $B \otimes_{E,c} E$; soit alors une
involution positive
de seconde espèce $*$ sur $B$, i.e. $*:B \longto B^{op}$ est telle
que $*_{|E}=c$ et $\tr(xx^*) >0$ pour tout $x \in B$.

\item à toute place $x$ de $F$, $B_x$ est soit décomposée soit une algèbre à division;

\item $B$ est décomposée en $v$;

\item si $d$ est pair alors le nombre de places de $F^+$ au dessus desquelles $B$ est
ramifiée  est congru à $1+\frac{d}{2} [F^+/\Qm] \mod 2$.

\end{itemize}

\item $V$ est le $B \otimes B^{op}$-module sous-jacent à $B$. Un accouplement alterné
$(,):V \times V \longto \Qm$
(resp. $(V \otimes \Am^\oo) \times (V \otimes \Am^\oo) \longto \Am^\oo$) $*$-hermitien
pour l'action de $B$ sur $V$
est de la forme
$$(x_1,x_2)_\beta = \tr_{B/\Qm}(x_1 \beta x_2^*)$$
pour $\beta \in B^{*=-1}$ (resp. $B^{*=-1} \otimes \Am^\oo$). On définit alors une
involution de
seconde espèce $\sharp_\beta$ sur $B$ (resp. $B \otimes \Am^\oo$) définie par
$$((b_1 \otimes b_2)x,y)=(x,(b_1^* \otimes b_2^\sharp)y)$$
et on note $G_\beta /\Qm$ (resp. $G_\beta /\Am^\oo$) le groupe algébrique dont les
$R$-points, pour toute
$\Qm$-algèbre (resp. $\Am^\oo$-algèbre) $R$ est l'ensemble des paires
$$(\lambda,g) \in R^\times \times (B^{op} \otimes R)^\times$$
telles que $g g^{\sharp_\beta}=\l$. On note alors $G_{\beta,1}$ le noyau du morphisme
$\nu:G_\beta \longto \Gm_m$
qui envoie $(\lambda,g)$ sur $\lambda$.

\end{itemize}

\marque Pour $i=1,\cdots,r$, on pose $\Lambda_i=\OC_{B_{v_i}} \subset V_{v_i}$ et on
note $\Lambda_i^\vee$
son dual dans $V_{v_i^c}$ de sorte que
$\Lambda=\bigoplus_{i=1}^r \Lambda_i \oplus \bigoplus_{i=1}^r \Lambda_i^\vee \subset V
\otimes_\Qm \Qm_p$
est un $\Zm_p$-réseau de $V \otimes_\Qm \Qm_p$ sur lequel l'accouplement de $V$ se
restreint en un
$\Zm_p$-accouplement. En outre comme $B_v \simeq \Mm_d(F_v)$, en notant $\e_1$
l'idempotent de $B_v$ associé
au premier vecteur de la base canonique, on a $\Lambda_{1,1}:=\e_1 \Lambda_1 \simeq
(\OC_v^\vee)^d$.

\marque \textit{Hypothèse supplémentaire:} soit $\s:F^+ \hookrightarrow \Rm$, on choisit
alors, ce qui est possible,
$\b$ tel que
$$G_{1,\s}(\Rm) = \left \{ \begin{array}{ll} U(1,d-1) & \hbox{si } \s =\t \\ U(0,d) &
\hbox{si } \s \neq \t
\end{array} \right.$$

\marque A $\tau$ est associé une extension
$$(,)_\tau:(V \otimes \Am) \times (V \otimes \Am) \longto \Am$$
ayant pour invariants $(1,d-1)$ en $\tau$ et $(0,d)$ à toutes les autres places infinies. On note $G_\tau$ le
groupe des similitudes associé qui ne dépend que de $\tau$ et pas du choix de $\beta$ bien qu'il y ait
$\sharp \ker^1(\Qm,G_\tau)$ choix possibles où $\ker^1(\Qm,G_\tau)$ désigne le sous-ensemble de
$H^1(\Qm,G_\tau)$ constitué des éléments qui deviennent triviaux dans $H^1(\Qm_{p'},G_\tau)$ pour toute place
$p'$ de $\Qm$.

\begin{rema} \label{rema-ker1}
- On notera que si $d$ est pair alors $\ker^1(\Qm,G_\tau)$ est nul tandis que pour $d$
impair on a
$$\ker^1(\Qm,G_\tau)= \ker \Bigl ( (F^+)^\times/\Qm^\times N_{F/F^+}(F^\times) \longto
\Am_{F^+}^\times/\Am^\times N_{F^+/F}(\Am_F^\times) \Bigr ) $$ Par ailleurs, cf. \cite{Ko1} \S 7, pour $d$
impair l'application $\ker^1(\Qm,Z_{G_\tau}) \longto \ker^1(\Qm,G_\tau) $ est bijective de sorte que toutes
les formes intérieures de $G_\tau$ isomorphes à $G_\tau$ partout localement, sont en fait isomorphes à $G$.

- On a
$G_\tau(\Qm_p)=G(\Qm_p) \simeq (\Qm_p)^\times \times \prod_{i=1}^r
(B_{v_i}^{op})^\times$
et le stabilisateur de $\Lambda \subset V \otimes_\Qm \Qm_p$ est alors sous cet
isomorphisme, identifié à
$\Zm_p^\times \times \prod_{i=1}^r \OC_{B_{v_i}}^\times$.
\end{rema}

\marque Soit $U^p$ un sous-groupe compact de $G(\Am^{\oo,p})$ et $m=(m_1,\cdots,m_r) \in
\Zm_{\geq 0}^r$. On
pose
$$U^p(m)=U^p \times \Zm_p^\times \times \prod_{i=1}^r \ker ( \OC_{B_{v_i}}^\times
\longto
(\OC_{B_{v_i}}/v_i^{m_i})^\times )$$
On considère alors $\XC_{U^p(m)}$ le foncteur de la catégorie des paires
$(S,s)$ où $S$ est un $\OC_o$-schéma connexe localement noethérien et $s$ est un point
géométrique de $S$,
tel que $X_{U^p(m)}(S,s)$ est l'ensemble des classes d'équivalence des $(r+4)$-uplets
$(A,\lambda,i,\bar \mu^p,\iota_{1,m_1},\cdots,\iota_{r,m_r})$ où:

\begin{itemize}

\item $A$ est une variété abélienne sur $S$ de dimension $ed^2$;

\item $\lambda:A \longto A^\vee$ est une polarisation;

\item $i:B \hookrightarrow \End(A) \otimes_\Zm \Qm$ est telle que $(A,i)$ est compatible
au sens de
\cite{h-t} lemme IV.1.2 p94 et $\lambda \circ i(b^*)=i(b)^\vee \circ \lambda$ pour tout
$b \in B$;

\item $\bar \mu^p$ est une $U^p$-orbite $\pi_1(S,s)$-invariante d'isomorphismes de $B
\otimes \Am^{\oo,p}$-modules
$\mu^p: V \otimes_\Qm \Am^{\oo,p} \longto V^p A_s$ qui envoie l'accouplement $<,>$ de $V
\otimes_\Qm \Am^{\oo,p}$ sur
un
$(\Am^{\oo,p})^\times$-multiple de l'accouplement de Weil associé à $\lambda$;

\item $\iota_{1,m_1}: v_1^{-m_1} \Lambda_{1,1}/\Lambda_{1,1} \longto \epsilon_1
A[v_1^{m_1}](S)$ est une structure de
niveau de Drinfeld;

\item pour tout $i=2,\cdots,r$, $\iota_{i,m_i}:(v_i^{-m_i} \Lambda_i/\Lambda_i)_S
\longto A[v_i^{m_i}]$ est
un isomorphisme de $S$-schémas munis d'une action de $\OC_B$.

\end{itemize}

Deux $(r+4)$-uplets $(A,\lambda,i,\bar \mu,\iota_{i,m_i})$ et $(A',\lambda',i',\bar
\mu',\iota_{i,m_i}')$
sont équivalents s'il existe une isogénie $\beta:A \longto A'$ qui envoie $\lambda$ sur
un
$\Qm^\times$-multiple de $\lambda'$, $i$ sur $i'$, $\bar \mu$ sur $\bar \mu'$ et
$\iota_{i,m_i}$ sur
$\iota_{i,m_i}'$.

\marque \rem Si $s'$ est un autre point géométrique de $S$ alors $\XC_{U^p(m)}(S,s)$ est
canoniquement en
bijection avec $X_U(S,s')$ de sorte que l'on peut considérer $\XC_{U^p(m)}$ comme un
foncteur sur les
$\OC_v$-schémas localement noethérien en posant $\XC_{U^p(m)}(\coprod_i S_i)=\prod_i
\XC_{U^p(m)}(S_i)$ où
les $S_i$ sont connexes.

\begin{prop} Si $U^p$ est suffisamment petit\footnote{s'il existe un premier $x$ de
$\Qm$ tel que l'image de
$U^p$ dans $G(\Qm_x)$ ne contient aucun élément non trivial d'ordre fini}, alors
$\XC_{U^p(m)}$ est
représenté par un schéma projectif sur $\spec \OC_v$ que l'on notera $X_{U^p(m)}$ et que
l'on appelle la
variété de Shimura de niveau $U$.
\end{prop}

\rem Soit $h:\Cm \longto \End_B(V)_\Rm$ le morphisme d'algèbres à involutions où
$\End_B(V)$ est munie
de l'adjonction par rapport à la forme symplectique $<,>$ et $\Cm$ de la conjugaison
complexe, définit par
$h(z):=(z,\bar z,\cdots,\bar z) \times (\bar z,\cdots,\bar z) \times \cdots \times (\bar
z,\cdots,\bar z) \in
G_{1,\Rm}$.
On note $X$ la classe de $G_\tau(\Rm)$-conjugaison de $h$ de sorte que $X$ s'identifie à
l'espace symétrique
hermitien
$G_\tau(\Rm)/K_\oo$ où $K_\oo$, un sous-groupe compact maximal modulo le centre de
$G_\tau(\Rm)$, est le
centralisateur
dans le groupe de Lie $G_\tau(\Rm)$. La variété de Shimura $X_{U^p(m)}$ n'est pas en
général associée à la donnée de
Shimura $(G,X)$, on
a en fait la décomposition suivante définie sur $F$, cf. \cite{Ko1} \S 8:

$$X_{U^p(m)}= \coprod_{\ker^1(\Qm,G)} \sh_{U^p(m)} (G',X)$$
où $G'$ parcourt les formes intérieures de $G_\tau$ isomorphes à $G_\tau$ partout
localement et où
$\sh_{U^p(m)}(G',X)$ désigne le modèle
canonique sur $F$ de la variété de Shimura associée à la donnée de Shimura $(G',X)$.
Ainsi d'après le remarque
(\ref{rema-ker1}),
on a
$X_{U^p(m)}= \sh_{U^p(m)}(G_\tau,X)^{\ker^1(\Qm,G_\tau)}$.

\noindent \textbf{Propriété 1}: \textit{quand $U$ varie, les $X_{U^p(m)}$ forment un
système projectif de
schémas projectif sur $\OC_v$ dont les morphismes de transition sont finis et plats:
quand $m_1=m_1'$ ils
sont en plus étales.}

\noindent \textbf{Propriété 2}: \textit{le système projectif $(X_{U^p(m)})_{U^p,m}$ est
naturellement muni
d'une action de $G(\Am^\oo)$.}

\marque Soit $U^p$ un sous-groupe compact suffisamment petit de $G(\Am^{\oo,p})$. On
note $\bar X_{U^p,m}$ la
fibre spéciale de $X_{U^p(m)}$ et soit $\AC$ la variété abélienne universelle sur $\bar
X_U$ naturellement
munie d'une action de $\OC_B$ ce qui donne donc une action de $\OC_{B,p}$ sur
$\AC[p^\oo]$ qui se décompose
alors sous la forme
$$\AC[p^\oo]=\prod_{i=1}^r \AC[v_i^\oo] \times \prod_{i=1}^r \AC[v_i^{c\oo}]$$
Ainsi $\GC:=\e \AC[v^\oo]$ est un $\OC_v$-module de Barsotti-Tate compatible de
dimension $1$ et de hauteur
$d$ sur $\bar X_{U^p,m}$ et les groupes de Barsotti-Tate $\AC[v_i^\oo]$ sont ind-étales
pour $i>1$.

\begin{defi} Pour tout $0 \leq h \leq d-1$, on note $\bar X_{U^p,m}^{[h]}$ et $\bar
X_{U^p,m}^{(h)}$ les strates
respectivement fermées et ouvertes de $\bar X_{U^p,m}$ associées au $\OC_v$-module de
Barsotti-Tate $\GC/\bar
X_ {U^p,m}$, comme dans la proposition (\ref{defi-p-strate}).
\end{defi}

\noindent \textbf{Propriété 3}: \textit{pour tout $0 \leq h \leq d-1$, la strate $\bar
X_{U^p,m}^{[h]}$ est
de pure dimension $h$ et munie d'une action de $G(\Am^\oo)$. Dans le cas de bonne
réduction, i.e. $m_1=0$,
elle est en outre lisse.}

\marque Pour tout point géométrique $x$ de $\bar X_{U^p,m}^{(h)}$, soit $M_x$ le noyau
de l'application
$$(\iota_{1,m_1})_x:v^{-m_1}\Lambda_{1,1}/\Lambda_{1,1} \longto \GC_x[v^{m_1}](\k(x))
\simeq
\GC_x^{et}[v^{m_1}](\k(x))$$
Pour tout $x$, $M_x$ est un facteur direct de rang $d-h$ de
$v^{-m_1}\Lambda_{1,1}/\Lambda_{1,1}$ sur
$\OC_v/v^{m_1}$. De plus quand $x$ varie, les $M_x$ se recollent en un faisceau
localement constant
$\underline M$ sur $\bar X_{U^p,m}^{(h)}$. On obtient ainsi une décomposition de $\bar
X_{U^p,m}^{(h)}$
indexée par les facteurs directs $M$ de $v^{-m_1}\Lambda_{1,1}/\Lambda_{1,1}$, libres de
rang $d-h$ en tant
que $\OC_v/v^{m_1}$-modules:
$\bar X_{U^p,m}^{(h)}=\coprod_M \bar X_{U^p,m,M}^{(h)}$
où pour tout point fermé $x$ de $\bar X_{U^p,m,M}^{(h)}$ on ait $M_x=M$.

\begin{lemm} Supposons $m_1=0$ et $m'_i=m_i$ pour $i>1$ alors $\bar X_{U^p,m',M}^{(h)}
\longto \bar X_{U^p,m}^{(h)}$
est fini et plat.
\end{lemm}

\begin{defi} Soit $M \subset \Lambda_{1,1}$ un $\OC_v$-sous-module libre de rand $d-h$
et facteur direct. Soit $P_M
\subset \aut(\Lambda_{1,1})$ le sous-groupe parabolique maximal qui stabilise $M$. On
pose alors
$$G_M(\Am^\oo)=G(\Am^{\oo,p}) \times \Qm_p^\times \times P_M(F_v) \times \prod_{i=2}^r
(B_{v_i}^{op})^\times$$
et on note $\bar X_{U^p,m,M}^{(h)}:=\bar X_{U^p,m,v^{-m_1}M/M}^{(h)}$.
Par ailleurs on associe à $M$ un élément $\bar g_M \in GL_d(F_v)/P_{d-h,d}(F_v)$ tel que
pour tout relèvement
$g_M$ de ce dernier, $M$ soit l'image par $g_M$ de l'espace engendré par les $d-h$
premiers vecteurs de la
base canonique.
\end{defi}

\noindent \textbf{Propriété 4}: \textit{pour $M$ fixé, le système projectif des $\bar
X_{U^p,m,M}^{(h)}$
admet une action de $G_M(\Am^\oo)$ de sorte que $\bar X_{U^p,m}^{(h)}$ muni de son
action de $G(\Am^\oo)$ est
géométriquement induit à partir de $\bar X_{U^p,m,M}^{(h)}$:}
$$\bar X_{U^p,m}^{(h)} \simeq \bar X_{U^p,m,M}^{(h)} \times_{P_M(\OC_v/v^{m_1})}
GL_d(\OC_v/v^{m_1})$$

\begin{defis} - Pour tout $h$, on notera $M_h$ le facteur direct de $\Lambda_{1,1}$
libre de rang $d-h$ associé
au $d-h$-premiers vecteurs de la base canonique et $G_{(h)}(\Am^\oo)$ le groupe
$G_{M_h}(\Am^\oo)$.

- On notera $\bar X_{U^p,m,M}^{[h]}$ l'adhérence de $\bar X_{U^p,m,M}^{(h)}$ dans $\bar
X_{U^p,m}^{[h]}$.
\end{defis}

\rem On peut montrer que $\bar X_{U^p,m,M}^{[h]}$ est lisse.

\noindent \textbf{Propriété 5}: \textit{pour $\xi$ une représentation de dimension finie
de $G$ sur un $\bar
\Qm_l$-espace vectoriel, il existe pour tout $U^p$ suffisamment petit, un système local
$\LC_\xi$ sur les
$X_{U^p,m}$.}

\marque Dans \cite{h-t}, les auteurs définissent les variétés d'Igusa de première espèce ainsi qu'un
isomorphisme de celles-ci vers $\bar X_{U^p,m,M}^{(h)}$ pour tout $M$. Ces variétés sont munies d'une action
de $G_M^+(\Am^\oo)$ qui agit via le quotient $G_M^+(\Am^\oo) \longto G(\Am^{\oo,p}) \times \Qm_p^\times
\times (\Zm \times GL_h(F_v))^+ \times \prod_{i=2}^r (B_{v_i}^{op})^\times$ où $(\Zm \times GL_h(F_v))^+$
désigne le sous-semi-groupe de $\Zm \times GL_h(F_v)$ formé des éléments $(c,g)$ tels que $\varpi_o^{\lfloor
- c/(n-h) \rfloor} g \in GL_h(\OC_v)$ qui est induit par la surjection $P_M(F_v) \longto \Zm \times
GL_h(F_v)$ qui à $g_v=g_M(g_v^0,g_v^{et})g_M^{-1}$ associe $(v(\det g_v^0),g_v^{et})$.

\noindent \textbf{Propriété 6}: \textit{via la notion de variété d'Igusa de seconde espèce, on peut définir
des systèmes locaux $\FC_{\rho_v,M}$ sur les $\bar X_{U^p,m,M}^{(h)}$ associés à une représentation
irréductible admissible $\rho_v$ des inversibles $D_{v,d-h}^\times$ de l'algèbre à division $D_{v,d-h}$ de
centre $F_v$ et d'invariant $1/(d-h)$. Tout élément $(g^p,g_{p,0},c,g_o^{et},g_{o_i},\d)$ de $G(\Am^{\oo,p})
\times \Qm_p^\times \times (\Zm \times GL_h(F_v))^+ \times \prod_{i=2}^r (B_{v_i}^{op})^\times \times
(D_{v,d-h}^\times/\DC_{v,d-h}^\times)$ définit naturellement un morphisme
$$(g^p,g_{p,0},c,g_v^{et},g_{v_i},\d):(g^p,g_{p,0},c,g_v^{et},g_{v_i},\d)^*
(\FC_{\rho,M} \otimes \LC_\xi) \longto \FC_{\rho,M} \otimes \LC_\xi$$}

\begin{defi} Soit  $H_0/\Qm$ le groupe algébrique forme intérieure de $G_\tau$ telle
que
$H_0(\Rm)$ est compact, $H_{0,v} \simeq D_{v,d}^\times$ et $H_0(\Am^\oo) \simeq
G^{(0)}(\Am^\oo)$.
\end{defi}

\noindent \textbf{Propriété 7}: \textit{l'ensemble des points supersinguliers $\bar X_{U^p,m}^{(0)}$ de $\bar
X_{U^p,m}$ est un sous-schéma de dimension nulle, constitué de $\sharp \ker^1(\Qm,H_0)$ classes d'isogénies
tel que chaque classe d'isogénie $\bar X_{U^p,m}^{(0)}(\bar \Fm_q)_i$ pour $i \in \ker^1(\Qm,H_0)$ est
isomorphe, en tant que $G^{(0)}(\Am^\oo):=G(\Am^{\oo,p}) \times D_{v,d}^\times \times \prod_{i=2}^r
(B_{v_i}^{op})^\times$-ensemble, à (cf. \cite{h-t} lemme V.1.2 p.153):}
$$H_0(\Qm) \backslash G^{(0)}(\Am^\oo)/ (\DC_{v,d}^\times \times U^p(m))$$
\textit{En outre l'action de $c_v \in W_v$ sur cet ensemble est donnée par la translation de $\deg(c_v)$ sur
la composante $\Zm \simeq D_{v,d}^\times / \DC_{v,d}^\times$ qui envoie $1$ sur une uniformisante $\Pi_{v,d}$
de $D_{v,d}$. L'action d'un élément $g^{\oo} \in G^{(0)}(\Am^\oo)$ est donnée par la correspondance}

$$\diagram & H_0(\Qm) \backslash [ \frac{G^{(0)}(\Am^{\oo,v})}{U^{p,v}_1(m')} \times
\Zm]
\dlto^{c_1} \drto^{c_2} \\
H_0(\Qm) \backslash [ \frac{G^{(0)}(\Am^{\oo,v})}{U^{p,v}_2(m)} \times \Zm] & & H_0(\Qm) \backslash [
\frac{G^{(0)}(\Am^{\oo,v})}{U^{p,v}_2(m)} \times \Zm]
\enddiagram$$
\textit{où $U^p_1(m')$ est tel que $U^{p,v}_1(m') \subset U^{p,v}_2(m) \cap (g^{\oo,v})^{-1} U^{p,v}_2(m)
g^{\oo,v}$, avec $c_1$ (resp. $c_2$) induit par l'inclusion $U^{p,v}_1(m') \subset U^{p,v}_2(m)$ (resp. par
la translation à gauche de $g^{\oo,v}$ sur $G^{(0)}(\Am^{\oo,v})$) et la translation de $\val(\det(g_v))$ sur
la composante $\Zm$.}

\rem D'après \cite{Ko7}, on a une bijection canonique $f:\ker^1(\Qm,G_\tau) \longto \ker^1(\Qm,H_0)$ de sorte
que dans la décomposition sur $\OC_v$ $X_{U^p(m)}= \coprod_{\ker^1(\Qm,G_\tau)} \sh_{U^p(m)}(G,X)$ l'ensemble
des points supersinguliers de la composante $\sh_{U^p(m)}(G,X)$ indexée par $i \in \ker^1(\Qm,G_\tau)$
correspond à la classe d'isogénie des points supersinguliers de $X_{U^p(m)}$ associée à $f(i) \in
\ker^1(\Qm,H_0)$.

\begin{defi} \label{defi-fg-s}
Pour tout diviseur $g$ de $d=sg$ et toute représentation irréductible cuspidale $\pi_v$ de $GL_g(F_v)$, on
notera $\FC(g,s,\pi_v)$ le faisceau concentré aux points supersinguliers dont la restriction à toute classe
d'isogénie $i \in \ker^1(\Qm,H_0)$ est
$$\diagram H_0(\Qm) \backslash [(G^{(0)}(\Am^{\oo,v}) / U^{p,v}(m)) \times (\Zm \times
\JL^{-1}([\overleftarrow{s-1}]_{\pi_v}))] \dto \\
\bar X_{U^p,m}^{(0)}(\bar \Fm_q)_i=H_0(\Qm) \backslash
[(G^{(0)}(\Am^{\oo,v})/U^{p,v}(m)) \times \Zm]
\enddiagram$$
où l'action diagonale de $H_0$ est donnée par translation à droite sur
$G^{(0)}(\Am^{\oo,v})=H_0(\Am^{\oo,v})$, par translation de valeur $\val \rn(.)$ sur $\Zm$ et par l'action
naturelle sur $\JL^{-1}([\overleftarrow{s-1}]_{\pi_v})$.
\end{defi}

\begin{defi} Pour tout $0 \leq h < d$, on notera
$$G^{(h)}(\Am^\oo):=G(\Am^{\oo,p}) \times \Qm_p^\times \times D_{v,d-h}^\times \times
GL_h(F_v) \times
\prod_{i=2}^r (B_{v_i}^{op})^\times$$
\end{defi}

\section{Rappels des propriétés cohomologiques}
\label{rappels-coho}
\setcounter{smfthm}{0}

\begin{defi} \label{defi-inert}
Soit $\tau$ une représentation irréductible de $D_{K,h}^\times$, sa restriction à
$\DC_{K,h}^\times$ est une
somme de représentations irréductibles $\rho_{1} \oplus \cdots \oplus \rho_{e_{\tau}}$
et on notera $e_{\tau}$ le nombre de celles ci. Étant donnée une représentation
irréductible $\rho$ de
$\DC_{K,h}^\times$, soient alors $\tau$ et $\tau'$ des sous-représentations
irréductibles de l'induite de
$\DC_{K,h}^\times$ à $D_{K,h}^\times$ de $\rho$: d'après la réciprocité de Frobenius, ce
sont exactement
celles telles que leur restriction à $\DC_{K,h}^\times$ contienne $\rho$. On en déduit
alors que $\tau$ et
$\tau'$ sont inertiellement équivalentes, i.e. $\tau' \simeq \tau \otimes \chi$ avec
$\chi:\d \mapsto
x^{\val(\det \d)}$ pour $x \in \bar \Qm_l^\times$. On note  $\CF_h$ l'ensemble des
classes d'équivalences
inertielles des représentations admissibles et irréductibles du groupe $D_{K,h}^\times$.
De la même façon,
deux représentations $\pi$ et $\pi'$ de $GL_g(K)$ seront dites inertiellement
équivalentes, s'il existe un
caractère $\xi:\Zm \longto \Qm_l^\times$ tel que $\pi \simeq \pi' \circ \xi \circ \val
\circ \det$ et on
notera $\cusp_g$ l'ensemble des classes d'équivalence inertielle des représentations
cuspidales de $GL_g(K)$
ainsi que $e_{\pi}$ le cardinal de la classe d'équivalence inertielle de $\pi$.
\end{defi}

\begin{defi} \label{wo-isotypique}
Soit $\sigma$ une représentation irréductible de $W_K$. Une représenta\-tion de $I_K$
sera dite
$\sigma$-isotypique, si ses sous-représentations irréductibles sont aussi des
sous-représentations
irréductibles de la restriction à $I_K$ de $\sigma$.
\end{defi}

\rem Si une représentation de $I_K$ est $\sigma$-isotypique et $\sigma'$-isotypique
alors $\sigma$ et $\sigma'$ sont
inertiellement équivalentes.

\begin{lemm} \label{epio}
Pour tout entier $g$ et toute représentation irréductible cuspidale $\pi$ de $GL_g(K)$,
$e_{\JL^{-1}([\overleftarrow{s-1}]_{\pi})}=e_{\pi}$.
\end{lemm}

\begin{proof} On rappelle que l'entier en question correspond au nombre de caractères
$\chi$ de $\Zm$ tels que
$\tau \simeq \tau \otimes \chi \circ \val \circ \det$ en notant
$\tau=\JL^{-1}([\overleftarrow{s-1}]_{\pi})$; par Jacquet-Langlands c'est aussi le
nombre de caractères
$\chi$ de $\Zm$ tels que
$[\overleftarrow{s-1}]_{\pi} \simeq [\overleftarrow{s-1}]_{\pi} \otimes (\chi \circ \val
\circ \det) \simeq
[\overleftarrow{s-1}]_{\pi \otimes (\chi \circ \val \circ \det)}$ et donc au nombre de
caractères $\xi$
tels que $\pi \simeq \pi \otimes \xi \circ \val \circ \det$, d'où le résultat.

\end{proof}

\subsection{Sur le modèle local}
\label{rap-loc}

Soit $\tau$ une représentation admissible irréductible de $D_{K,h}^\times$. L'espace que
l'on souhaite
étudier est le $(GL_h(K) \times W_{K})$-module
$$\hom_{D_{K,h}^\times}(\tau,\UC_{K,l,h}^{i})=\UC_{K,l,h}^{i}(\tau) \simeq
\Psi_{K,l,h}^{i}(\tau):=
\hom_{\DC_{K,h}^\times}(\tau,\Psi_{K,l,h}^{i})$$

\marque Pour tout $\tau$, on a un morphisme naturel de $\NC_K$-modules:
$\UC_{K,l,h}^{i}(\tau) \otimes \tau \longto \Psi_{K,l,h}^{i}$
qui envoie $f \otimes v$ sur $f(v)$. On note $\Psi_{K,l,h}^{i}[\tau]$ l'image de ce
morphisme et soit
$\Psi_{K,l,h,m}^{i}[\tau]$ la préimage de $\Psi_{K,l,h}^{i}[\tau]$ dans
$\Psi_{K,l,h,m}^{i}$. Le
sous-module $\Psi_{K,l,h}^{i}[\tau]$ ne dépend que de la classe d'équivalence inertielle
de $\tau$. Le
groupe $\DC_{K,h}^\times$ étant compact, on a
$\Psi_{K,l,h}^{i} = \bigoplus_{\tau \in \CF_h} \Psi_{K,l,h}^{i}[\tau]$.
Soit $\Delta_{\tau}$ un ensemble d'\eles de $D_{K,h}^\times$ tel que les congruences des
$\val (\det \d)$
pour $\d \in \Delta_{\tau}$ forment un système de représentants de $\Zm/e_{\tau} \Zm$.
L'application
$$\begin{array}{rl} \label{action-tordue}
\UC_{K,l,h}^{i}(\tau) \otimes \tau & \longto \bigoplus_{\d \in \Delta_{\tau}}
\Psi_{K,l,h}^{i}[\tau]^\d \\
f \otimes v & \mapsto (f(\d^{-1}v))_\d
\end{array}$$
où $\Psi_{K,l,h}^{i}[\tau]^\d$ est l'espace $\Psi_{K,l,h}^{i}[\tau]$ muni de la
structure de
$\NC_K$-module où $(g,\nu,c)$ agit via $(g,\d^{-1}\nu \d ,c)$, est un isomorphisme de
$\NC_K$-modules.

\begin{theo} (cf. \cite{h-t})
Pour toute représentation irréductible cuspidale $\pi$ de $GL_d(K)$,
$\UC_{K,l,d}^{i}(\JL^{-1}(\pi)^\vee)$ est nul pour $i \neq  d-1$ et
$$\UC_{K,l,d}^{d-1}(\JL^{-1}(\pi)^\vee) \simeq \pi \otimes \rec_{K}^\vee(\pi)
(-\frac{d-1}{2})$$
\end{theo}

\begin{theo} \label{inutile} (cf. \cite{h-t} théorème VII.1.5)
Pour tout diviseur $g$ de $d=sg$ et toute
représentation irréductible cuspidale $\pi$ de $GL_g(K)$, on a
\begin{multline*}
\sum_{i=0}^{d-1} (-1)^i [\UC_{K,l,d}^{d-1-i}(\JL^{-1}(\st_s(\pi)^\vee))]=
\sum_{i=1}^s (-1)^i [\overleftarrow{s-1-i},\overrightarrow{i}]_{\pi} \otimes
\rec_{K}^\vee(\pi)(-\frac{d+s-2-2i}{2})
\end{multline*}
\end{theo}

\subsection{Sur les systèmes locaux d'Harris-Taylor}
\label{ht-prop}

Dans la suite nous ne considérerons plus les schémas sur $\Fm_q$, $\bar X_{U^p,m}$, $\bar X_{U,m}^{(h)}$...
mais plutôt les $\bar \Fm_q$-schémas $\bar X_{U^p,m} \times_{\Fm_q} \bar \Fm_q$, $\bar
X_{U,m}^{(h)}\times_{\Fm_q} \bar \Fm_q$... Afin de ne pas alourdir encore plus les énoncés, nous garderons
les mêmes notations, par exemple $\bar X_{U^p,m}$ désignera ce que l'on devrait noter $\bar X_{U^p,m}
\times_{\Fm_q} \bar \Fm_q$.

\begin{defi}
Pour $\tau_v$ une représentation irréductible admissible de $D_{v,h}^\times$, on note
$\FC_{\tau_v,U^p(m),1}$
le faisceau sur $\bar X_{U^p,m,M_{d-h}}^{(d-h)}$ noté précédemment
$\FC_{\tau_v,U^p(m),M_{d-h}}$ et on définit:
$$\FC_{\tau_v,U^p(m)}:= \FC_{\tau_v,U^p(m),1} \times_{P_{h,d}(\OC_v/v^{m_1})}
GL_d(\OC_v/v^{m_1})$$
le faisceau induit sur toute la strate $\bar X_{U^p,m}^{(d-h)}$.
\end{defi}

\marque Les systèmes locaux $\FC_{\tau_v}=(\FC_{\tau_v,U^p(m)})_{U^p(m)}$ d'Harris-Taylor sont tels que la
restriction à $\bar X_{U^p,m,M_{d-h}}^{(d-h)}$ du $i$-ème faisceau des cycles évanescents $R^i\Psi_v$ vérifie
\begin{equation} \label{iso1}
(R^i\Psi_v)_{|\bar X_{U^p,m,M_{d-h}}^{(d-h)}}^{h} \simeq \bigoplus_{\tau_v \in \CF_h} (\FC_{\tau_v,U^p(m),1}
\otimes \UC_{F_v,l,h,m_1}^{i}(\tau_v))^{h/e_{\tau_v}},
\end{equation}
où l'action d'un élément $(g^{\oo,p}, g_{p,0},g_v^0,g_v^{et},g_{v_i},c_v)$ de $G(\Am^{\oo,p}) \times
\Qm_p^\times \times (GL_{n-h}(F_v) \times GL_h(F_v))^+ \times \prod_{i=2}^r (B_{v_i}^{op})^\times \times W_v$
est donnée par l'action naturelle de
\begin{multline*}
(g^{\oo,p}, g_{p,0}p^{-f_1 v(c_v)},v(\det g_v^0)-v(c_v),g_v^{et},g_{v_i}) \in G(\Am^{\oo,p}) \times
\Qm_p^\times \times (\Zm \times GL_h(F_v))^+ \times \prod_{i=2}^r (B_{v_i}^{op})^\times
\end{multline*}
sur le système $\FC_{\tau_v}$ au dessus de $\bar X_{U^p,m,M_{d-h}}^{(d-h)}$, et par celle de $(g_v^0,c_v)$
sur la tour des $\UC_{F_v,l,h,m_1}^{i}(\tau_v)$. L'isomorphisme (\ref{iso1}) implique alors la proposition
suivante.

\begin{prop} \label{prop-hic}
Pour tout $i,j$, on a un isomorphisme canonique
\begin{multline*}
H^j_c(\bar X_{U^p,m,M_{d-h}}^{(d-h)}, R^i\Psi_v \otimes \LC_\xi)^h \simeq \bigoplus_{\tau_v \in \CF_h}
(H^j_c(\bar X_{U^p,m,M_{d-h}}^{(d-h)},\FC_{\tau_v,U^p(m),1} \otimes \LC_{\xi}) \otimes
\UC_{F_v,l,h,m_1}^{i}(\tau_v))^{h/e_{\tau_v}}
\end{multline*}
tel que l'action de $G(\Am^{\oo,p}) \times \Qm_p^\times \times (GL_{d-h}(F_v) \times GL_h(F_v))^+ \times
\prod_{i=2}^r (B_{v_i}^{op})^\times \times W_v$ se prolonge à $G(\Am^{\oo,p}) \times \Qm_p^\times \times
GL_{d-h}(F_v) \times GL_h(F_v) \times \prod_{i=2}^r (B_{v_i}^{op})^\times \times W_v$ de sorte que l'action
de $(g^{\oo,p}, g_{p,0},g_v^0,g_v^{et},g_{v_i},c_v)$ sur la limite inductive, indexée par $U^p(m)$, du membre
de gauche, induit l'action de $(g^{\oo,p},g_{p,0}p^{-f_1 v(c_v)},v(\det(g_v^0)-v(c_v),g_v^{et},g_{v_i})
\times (g_v^0,c_v)$ sur la limite inductive du membre de droite.
\end{prop}

Dans la suite on couplera la proposition précédente avec la suivante.

\begin{prop} \label{prop-poids}
Soient $1 \leq tg <d$, $\pi_v$ une représentation irréductible cuspidale de $GL_g(F_v)$ et $\Pi_t$ une
représentation de $GL_{tg}(F_v)$. Pour $1 \leq i \leq t$ et pour tout $j$,
$$\lim_{\atop{\to}{U^p(m)}} H^j_c(\bar X_{U^p,m,M_{d-tg}}^{(d-tg)},
\FC_{\JL^{-1}(\st_t(\pi_v)^\vee),U^p(m),1} \otimes \LC_{\xi}) \otimes \Pi_t \otimes \rec_{F_v}^\vee(\pi_v)$$
est naturellement muni d'une action de $G(\Am^{\oo,p}) \times \Qm_p^\times \times GL_{d-tg}(F_v) \times
GL_{tg}(F_v) \times \prod_{i=2}^r (B_{v_i}^{op})^\times \times W_v$ telle que $(g^{\oo,p},
g_{p,0},g_v^0,g_v^{et},g_{v_i},c_v)$ y agisse via l'action naturelle de
$$(g^{\oo,p},g_{p,0}p^{-f_1 v(c_v)},v(\det(g_v^0)-v(c_v),g_v^{et},g_{v_i})$$
sur ${\DS \lim_{\atop{\to}{U^p(m)}}} H^j_c(\bar X_{U^p,m}^{(d-tg)}, \FC_{\JL^{-1}(\st_t(\pi_v)^\vee),U^p(m)}
\otimes \LC_{\xi})$ et celle de $(g_v^0,c_v)$ sur $\Pi_t \otimes \rec_{F_v}^\vee(\pi_v)$.

Ainsi l'espace induit associé
$$\lim_{\atop{\to}{U^p(m)}} H^j_c(\bar X_{U^p,m}^{(d-tg)}, \FC_{\JL^{-1}(\st_t(\pi_v)^\vee),U^p(m)} \otimes \LC_{\xi}) \otimes \Pi_t
\otimes \rec_{F_v}^\vee(\pi_v)$$ en tant que représentation de $GL_d(F_v) \times W_v$ est de la forme
$$\bigoplus_\chi(\ind_{P_{tg,d}(F_v)}^{GL_d(F_v)} \Pi_t \circ \chi(\val(\det))\otimes \pi_\chi)\otimes
\rec_{F_v}^\vee(\pi_v \circ \chi(\val(\det)))$$ où $\chi$ décrit les caractères $\Zm \longto \bar
\Qm_l^\times$ et $\pi_\chi$ est une représentation de $GL_{d-tg}(F_v)$.
\end{prop}

\begin{proof} Le faisceau $\FC_{\JL^{-1}(\st_t(\pi_v)^\vee)}$ étant induit à partir de
$\FC_{\JL^{-1}(\st_t(\pi_v)^\vee),1}$
on a:
\begin{multline} \label{iso-induite1}
\lim_{\atop{\to}{U^p(m)}} H^j_c(\bar X_{U^p,m}^{(d-tg)},\FC_{\JL^{-1}(\st_t(\pi_v)^\vee),U^p(m)} \otimes
\LC_{\xi}) \otimes \Pi_t \otimes \rec_{F_v}^\vee(\pi_v)) \simeq \\
\ind_{P_{tg,d}(F_v)}^{GL_d(F_v)} \lim_{\atop{\to}{U^p(m)}} H^j_c(\bar
X_{U^p,m,M_{d-tg}}^{(d-tg)},\FC_{\JL^{-1}(\st_t(\pi_v)^\vee),U^p(m),1} \otimes \LC_{\xi}) \otimes \Pi_t
\otimes \rec_{F_v}^\vee(\pi_v)
\end{multline}
tel que l'action de $((g_v^0,g_v^{et}),c_v) \in P_{tg,_d}(F_v) \times W_v$ sur le membre de droite de
(\ref{iso-induite1}), soit donnée par l'action de $(g_v^{et},\val(\det g_v^0)-\deg(c_v)) \times (g_v^0,c_v)
\in GL_{d-tg}(F_v) \times \Zm \times GL_{tg}(F_v) \times W_v$ où ${\DS \lim_{\atop{\to}{U^p(m)}}} H^j_c(\bar
X_{U^p,m,M_{d-tg}}^{(d-tg)},\FC_{\JL^{-1}(\st_t(\pi_v)^\vee),U^p(m),1} \otimes \LC_{\xi})$ (resp. $\Pi_t
\otimes \rec_{F_v}^\vee(\pi_v)$) est vue comme un $GL_{d-tg}(F_v) \times \Zm$-module (resp. $GL_{tg}(F_v)
\times W_v$-module). Le résultat découle alors de l'écriture
$$\lim_{\atop{\to}{U^p(m)}} H^j_c(\bar
X_{U^p,m,M_{d-tg}}^{(d-tg)},\FC_{\JL^{-1}(\st_t(\pi_v)^\vee),U^p(m),1} \otimes \LC_{\xi}) \simeq
\bigoplus_{\chi} \chi \otimes \pi_\chi$$ où $\chi$ décrit les caractères $\Zm \longto \bar \Qm_l^\times$ et
où $\pi_\chi$ est une représentation de $GL_{d-tg}(F_v)$, à priori réductible.

\end{proof}

\marque \rem Les systèmes locaux $\FC_{\tau_v}$ ne sont pas irréductibles mais plutôt une somme directe de
$e_{\tau_v}$ systèmes locaux irréductibles. La complexité de l'écriture de (\ref{iso1}), est la contre-partie
de la simplicité de la description de l'action de $GL_d(F_v) \times W_v$ qui tient au fait que l'on a fait
apparaître $\UC_{F_v,l,h}^{i}(\tau_v)$ plutôt que $\hom_{\DC_{v,h}^\times}(\r_v,\Psi_{F_v,l,h}^{i})$ avec
$\rho_v$ une représentation irréductible de $I_v$ de sorte que ce dernier est seulement muni d'une action de
$\NC_v \cap (GL_h(F_v) \times D_{v,h}^{\times})$.

\begin{defi} \label{defi-hyp}
On considère pour une représentation automorphe $\Pi$ de $G_\tau(\Am)$, l'hypothèse
$\hyp(\oo)$ suivante:
$\Pi_\oo$ est cohomologique pour une certaine représentation algébrique $\xi$ sur $\Cm$
de la
restriction des scalaires de $F$ à $\Qm$ de $GL_g$, i.e. il existe $i$ tel que
$$H^i((\lie G_\tau(\Rm)) \otimes_\Rm \Cm,U_\tau,\Pi_\oo \otimes (\xi')^\vee) \neq (0)$$
où $U_\tau$ est un sous-groupe compact modulo le centre de $G_\tau(\Rm)$, maximal, cf.
\cite{h-t} p.92, et où $\xi'$
est le caractère sur $\Cm$ associé à $\xi$ via un isomorphisme $\bar \Qm_l \simeq \Cm$
fixé.
On dira que $\Pi$ vérifie $\hyp(\xi)$ si $\Pi$ vérifie $\hyp(\oo)$ pour $\xi$.
\end{defi}

\marque On note ${\DS [H^*_{h,\xi,\tau_v}]:= \sum_i (-1)^i \lim_{\atop{\longto}{U^p,m}} [H^i_c(\bar
X_{U^p,m,M_{d-h}}^{(d-h)},\FC_{\tau_v,U^p(m),1} \otimes \LC_{\xi})]}$ dans le groupe de Grothendieck des
représentations admissibles de $(G(\Am^{\oo,v}) \times GL_{h}(F_v) \times \Zm$.

\begin{defi} \label{defi-carac}
Pour tout entier $h$, on considère l'identification canonique définie par la valuation de la norme réduite:
$D_{v,h}^\times / \DC_{v,h}^\times \longto \Zm$. On notera $\Xi:\Zm \longto \bar \Qm_l^\times$ le caractère
défini par $\Xi(1)=\frac{1}{q}$
\end{defi}

\begin{prop} \label{prop-somme-alternee}
Pour $\Pi$ une représentation de $G_\tau(\Am)$, on a alors
\begin{equation} \label{somme-alternee}
[H^*_{h,\xi,\tau_v}(\Pi^{\oo,v})]= \left \{
  \begin{array}{ll}
\sharp \ker^1(\Qm,G_\tau)  \e(\Pi)  m(\Pi) \Red_{\tau_v}^h (\Pi_v) & \hbox{si } \Pi_\oo
\hbox{ vérifie } \hyp(\oo) \\
0 & \hbox{sinon}
  \end{array}
\right . \end{equation}
où $\e(\Pi)$ est un signe qui dépend de $\Pi$, $m(\Pi)$ est la multiplicité de $\Pi$
dans l'espace des formes automorphes et $\Red_{\tau_v}^h: \groth(GL_d(F_v)) \longto
\groth(D_{v,h}^\times/\DC_{v,h}^\times \times GL_{d-h}(F_v))$ est défini comme la
composition des deux
homomorphismes suivant:

\begin{itemize}
\item en premier lieu, on a un homomorphisme
$$\begin{array}{l}
\groth(GL_d(F_v)) \longto \groth(GL_h(F_v) \times GL_{d-h}(F_v)) \\
~ [ \Pi_v ] \mapsto [ J_{P_{h,d}}(\Pi_v) \otimes \d_{P_{h,d}}^{1/2} ]
\end{array}$$

\item ensuite on a un homomorphisme
$$\begin{array}{l}
\groth(GL_h(F_v) \times GL_{d-h}(F_v)) \longto \groth(D_{v,h}^\times/\DC_{v,h}^\times
\times GL_{d-h}(F_v)) \\
~ [\a \otimes \b] \mapsto \sum_{\psi} \vol(D_{v,h}^\times/F_v^\times,d \bar h_v)^{-1}
\tr \a
(\phi_{\JL(\tau_v^\vee \otimes \psi)})[\psi \otimes \b],
\end{array}$$
où $\psi$ décrit les caractères de $\Zm \simeq D_{v,h}^\times/\DC_{v,h}^\times$ tels que
$\a$ et $\tau_v^\vee
\otimes \psi$ ont le même caractère central et où l'on considère des mesures de Haar
associées sur
$GL_h(F_v)$ et $D_{v,h}^\times$.
\end{itemize}
\end{prop}

\rem On a $\e(\Pi)=1$ (resp. $\e(\Pi)=(-1)^{s-1}$) s'il existe une place $x$ tel que
$\Pi_x$ est de carré
intégrable (resp. $\Pi_x \simeq \speh_s(\pi_x)$ pour $\pi_x$ irréductible cuspidale de
$GL_g(F_x)$ avec
$d=sg$). \footnote{L'auteur n'est pas certain que le cas $\Pi_x \simeq \speh_s(\pi_x)$
soit connu, de sorte
qu'on en donnera une preuve.}

\marque Pour toute représentation irréductible $\xi$ de $G(\Qm)$, on note
$\CC^\oo_{H_0,\xi}$ l'ensemble des
formes automorphes sur $H_0$ de type $(\xi')^\vee$ à l'infini. De la propriété 7, on en
déduit le résultat
suivant.

\begin{prop} \label{prop-h0-ss}
Pour tout $i$, on a un isomorphisme $G(\Am^\oo) \times W_v$-équivariant
\begin{equation} \label{h0-ss}
\lim_{\atop{\longto}{U^p,m}} H^0(\bar X_{U^p,m}^{(0)},R^i\Psi_v \otimes \LC_{\xi})
\simeq \sharp
\ker^1(\Qm,G_\tau) \hom_{D_{v,d}^\times} (\CC^\oo_{H_0,\xi}, \UC_{F_v,l,d}^{i})
\end{equation}
\end{prop}

\rem Le groupe $H_0(\Rm)$ étant compact, l'espace $\CC^\oo_{H_0,\xi}$ s'écrit comme une
somme directe
$\bigoplus_{\bar \Pi} m(\bar \Pi) \bar \Pi$ où $\bar \Pi$
décrit l'ensemble des classes d'isomorphismes des représentations automorphes de
$H_0(\Am)$ cohomologiques pour
$(\xi')^\vee$.

\part{Filtrations de monodromie: énoncés}

\section{Énonces des théorèmes locaux}

\subsection{Groupes de cohomologie des modèles de Deligne-Carayol}

~ \\

\marque Considérons le complexe à flèches nulles
\begin{multline*}
ML^\bullet (s):=([\overrightarrow{s-1}]_{1} \otimes |\art_{K}^{-1}|^{(s-1)/2},
[\overleftarrow{1},\overrightarrow{s-2}]_{1} \otimes |\art_{K}^{-1}|^{(s-3)/2}
,\cdots,\\
[\overleftarrow{s-1}]_{1} \otimes |\art_{K}^{-1}|^{(1-s)/2}) \otimes \JL^{-1}(\st_s(1)^\vee) \otimes
\rec_{K}^\vee(1)((1-s)/2)
\end{multline*}
où $[\overleftarrow{s-1}]_{1}$ est placé en degré $0$. Pour $\pi$ une représentation cuspidale de $GL_g(K)$,
on pose
\begin{multline*}
\pi \diamond ML^\bullet(s):=([\overrightarrow{s-1}]_{\pi} \otimes
|\art_{K}^{-1}|^{(s-1)/2},
[\overleftarrow{1}, \overrightarrow{s-2}]_{\pi} \otimes
|\art_{K}^{-1}|^{(s-3)/2},\cdots,\\
[\overleftarrow{s-1}]_{\pi} \otimes |\art_{K}^{-1}|^{(1-s)/2}) \otimes \JL^{-1}(\st_s(\pi))^\vee \otimes
\rec_{K}^\vee(\pi)((1-sg)/2)
\end{multline*}

\begin{theo} \label{theo-ripsi-local}
Pour tout $d$, on a $\UC_{K,l,d}^{d-1+\bullet}=\bigoplus_{\atop{g|d}{d=sg}}
\bigoplus_{\pi \in \cusp_g} \pi
\diamond ML^\bullet(s)$ où $\cusp_g$ désigne l'ensemble des classes d'équivalence des
représentations
irréductibles cuspidales de $GL_g(K)$.
\end{theo}

\marque Autrement dit, on a
$$\UC_{K,l,d}^{d-s+i}(\JL^{-1}(\st_{s}(\pi)^\vee)) \simeq
\left \{ \begin{array}{cl} \rec_{K}^\vee(\pi) (-\frac{d-s+2i}{2}) \otimes
[\overleftarrow{i},\overrightarrow{s-i-1}]_{\pi} & 0 \leq i < s \\ 0 & i < 0 \end{array}
\right. $$ où
$[\overleftarrow{i},\overrightarrow{s-i-1}]_{\pi}$ est l'unique quotient irréductible de
l'induite
$$\ind_{P_{ig,sg}(K)}^{GL_{sg}(K)} \st_i(\pi(\frac{(s-i)(g-1)}{2})) \otimes
\speh_{s-i}(\pi(\frac{-i(g-1)}{2})).$$

\rem Le cas Iwahori, pourrait se déduire du résultat principal de \cite{s-s} en
utilisant le théorème de
comparaison de Faltings, cf. \cite{falt}.

\subsection{Filtration de monodromie-locale}

On rappelle l'amélioration du théorème de comparaison de Berkovich donnée par Fargues en
appendice.

\begin{theo-defi} \label{theo-defi-monodromie}
Soit $X \to \spec \OC_K$ un morphisme propre d'un schéma $X$ de type fini de dimension
$d$ sur le trait
$\spec \OC_K$. Le complexe des cycles évanescents $R\Psi_{\eta}(\Qm_l)[d-1]$ est un
faisceau pervers munie
d'une action de $W_K$ qui fournit une filtration de monodromie dont on notera $gr_{k}$
les gradués et on
considère la suite spectrale $E_1^{i,j}=h^{i+j} gr_{-i} \Rightarrow R^{i+j+d-1}
\Psi_{\eta}(\Qm_l).$

Pour tout point géométrique $x$ de la fibre spéciale $X_{s}$, en considérant les germes
en $x$ des $h^i
gr_{k}$, on obtient une suite spectrale $E_{1,x}^{i,j}=(h^{i+j} gr_{-i})_x \Rightarrow
R^{i+j+d-1}
\Psi_{\eta}(\Qm_l)_x$ dont la nature est purement locale de sorte que si $Y \to \spec
\OC_K$ est un autre
schéma avec un point $y$ tel que le complété formel de l'anneau local de $Y$ en $y$ est
isomorphe, en tant
que $\OC_K$-schéma formel, au complété formel de l'anneau local de $X$ en $x$, alors
pour tout $r \geq 1$, on
a $E_{r,x}^{i,j}=E_{r,y}^{i,j}$.

La filtration ainsi obtenue sera dite de monodromie-locale.
\end{theo-defi}

\marque Pour tout $s \geq 1$, on introduit \footnote{Par exemple pour $s=4$ on a
représenté à la figure
(\ref{figure7}) les $MLE_1^{i,j}(4)$ et $MLE_\oo^{i,j}(4)$.} une suite de bicomplexe
$MLE_r^{i,j}(s)$ définit
comme suit:

- pour $r=1$, $MLE_1^{i,j}(s)$ est nul pour $|j|\geq s$ ou $j \equiv s \mod 2$ ou $i+j > 0$ ou $i<(1-s-j)/2$.
Sinon pour $j=1-s+2r$ avec $0 \leq r < s$ et $i=s-1-2r-k$ avec $0 \leq k \leq(s-1-j)/2$, on a
$$MLE_1^{s-1-2r-k,1-s+2r}(s)=[\overleftarrow{s-1-k}]_{1} \overrightarrow{\times}
[\overrightarrow{k-1}]_{1} \otimes |\art_{K}^{-1}|^{(s-1-2r)/2}$$

- les flèches $d_1^{i,j}$ se déduisent des suites exactes courtes
(\ref{suites-exactes});

- pour $r \geq 2$, $MLE_r^{i,j}(s)=MLE_2^{i,j}(s)$ est nul pour $|j| \geq s$ ou $j \equiv s \mod 2$ ou $2i
\neq 1-s-j$, et
$$MLE_2^{1-s+r,s-1-2r}=[\overleftarrow{s-1-r},\overrightarrow{r}]_{1} \otimes
|\art_{K}^{-1}|^{-(s-1-2r)/2}$$

Nous allons en fait prouver l'énoncé suivant dont découle directement le théorème
(\ref{theo-ripsi-local}).

\begin{theo} \label{theo-local-fil}
La filtration de monodromie-locale du complexe $\UC_{K,l,d,n}^\bullet$ est équivariante pour l'action de
l'algèbre de Hecke $\HC_n=\CC^\oo(\KF_{n} \backslash GL_d(K)/ \KF_{n})$ et celle de $I_K$; en outre pour $n$
variant, les filtrations de monodromie locale forment des systèmes inductifs compatibles à l'action des
correspondances de Hecke. On note $gr_{d,n,loc,k}$ les gradués et $E_{r,gr,loc}^{i,j}$ les termes de la suite
spectrale associée avec $E_{1,gr,loc}^{i,j}= h^{i+j} gr_{d,n,loc,-i}$. On a alors:
$$E_{\bullet,gr,loc}^{\bullet,d-1+\bullet}(\frac{d-1}{2})=\bigoplus_{\atop{g|d}{d=sg}} \bigoplus_{\pi \in
\cusp_g} JL^{-1}(\st_s(\pi)^\vee) \otimes \bigl ( \pi \diamond MLE_\bullet^{\bullet,\bullet}(s) \bigr ) $$ où
comme précédemment $\pi \diamond \bigl ( [\overleftrightarrow{s-1}]_{1} \otimes \rec_K^\vee(1) \bigr ) :=
[\overleftrightarrow{s-1}]_{\pi} \otimes \rec_K^\vee(\pi)$ et où $(\frac{d-1}{2})$ désigne la torsion sur la
partie galoisienne.
\end{theo}

\marque Autrement dit on a les propriétés suivantes:
\begin{itemize}
\item[(1)] pour tout $k$, $gr_{d,n,loc,k}$ est une somme directe
$\bigoplus_{\atop{g ~|~ d=sg}{\pi_v \in \cusp_g}} gr_{d,n,loc,k,\pi_v},$
où en tant que $I_v$-module, $gr_{d,n,loc,k,\pi_v}$ est
$\rec_{F_v}^\vee(\pi_v)$-isotypique au sens de la
définition (\ref{wo-isotypique}) \footnote{Sans utiliser le fait que les
$\rec_{F_v}^\vee(\pi_v)$ décrivent
l'ensemble des représentations irréductibles de $W_v$, l'énoncé dit que pour $\sigma_v$
qui ne serait pas de la
forme $\rec_{F_v}^\vee(\pi_v)$, $gr_{d,loc,k,\sigma_v}$ est nul.}

\item[(2)] les $gr_{d,n,loc,k,\pi_v}$ sont nuls pour $|k| \geq s$ et pour $|k| < s$, ses
groupes
de cohomologie $h^i gr_{d,n,loc,k,\pi_v}$ sont nuls pour $i<d-s+|k|$ et pour $i \not
\equiv d-s+k \mod 2$;
pour $d-s+|k| \leq i=d-s+|k|+2r \leq d-1$, $h^{d-s+|k|+2r} gr_{d,loc,k,\pi_v}$ est
l'espace des vecteurs
invariants sous $K_{v,n} \subset GL_d(F_v)$ de l'espace suivant:
$$\JL^{-1}([\overleftarrow{s-1}]_{\pi_v^\vee}) \otimes [\overleftarrow{|k|+2r}]_{\pi_v}
\overrightarrow{\times}
[\overrightarrow{s-2-2r-|k|}]_{\pi_v} \otimes \rec_{F_v}^\vee(\pi_v)(-\frac{i+k}{2});$$

\item[(3)] la suite spectrale $E_{1,gr,loc}^{i,j}=h^{i+j} gr_{d,loc,-i,\pi_v}
\Rightarrow \UC_{F_v,l,d}^{i+j}$ dégénère en $E_2$ et ses $d_1^{i,j}$ se déduisent des suites exactes courtes
\begin{multline} \label{suites-exactes}
0 \to [\overleftarrow{t-1},\overrightarrow{s-t}]_{\pi_v} \to
[\overleftarrow{t-1}]_{\pi_v}
\overrightarrow{\times} [\overrightarrow{s-t-1}]_{\pi_v} \to
[\overleftarrow{t},\overrightarrow{s-t-1}]_{\pi_v} \to 0
\end{multline}
\end{itemize}

\section{Énonces des théorèmes globaux}

\subsection{Rappels sur les faisceaux pervers}
\label{intro}

Pour $X$ un schéma séparé de type fini sur un corps, on considère la $t$-structure
perverse autoduale sur
$D_c^b(X,\bar \Qm_l)$. On reprend les notations et les résultats de \cite{ast}: en
particulier quand il n'y a
pas de confusion relativement à $X$, on note pour $a \in \Zm$, $\lexp p \t^{\leq a}$,
$\lexp p \t^{\geq a}$
les foncteurs de troncations de $D^b_c:=D_c^b(X,\bar \Qm_l)$ dans $D^{\leq
a}:=D_c^b(X,\bar \Qm_l)^{\leq a}$
et $D^{\geq a}:=D_c^b(X,\bar \Qm_l)^{\geq a}$, $\perv:=D^{\geq 0} \cap D^{\leq 0}$,
$\lexp p H^0:= \lexp p
\t^{\leq 0} \lexp p \t^{\geq 0}:D^b_c \longto \perv$.

\begin{defi} Soit $\sigma_v$ une représentation irréductible de $W_v$. Un faisceau
pervers
$P$ irréductible, de la forme $j_{!*} \LC$ pour $j:U \hookrightarrow X$ localement fermé
et $\LC$ un système
local sur $U$, muni d'une action de $I_v$, sera dit $\sigma_v$-isotypique si $\LC$ l'est
au sens de la définition
(\ref{wo-isotypique}).
\end{defi}

On a alors le lemme suivant dont la preuve est immédiate.

\begin{lemm} Soient $\sigma_{v,1}$ et $\sigma_{v,2}$ des représentations irréductibles
non inertiellement équivalentes
de $W_v$
et soient $P_1$, $P_2$ des faisceaux pervers muni d'une action de $I_v$, respectivement
$\sigma_{v,1}$ et
$\sigma_{v,2}$ isotypique. Si $\vphi:P_1 \longto P_2$ est un morphisme de faisceau
pervers $I_v$-équivariant,
pour tout $i$, on a $h^i \vphi=0$.
\end{lemm}

Le lemme (5.3.6) de \cite{ast} donne alors la décomposition en somme directe suivante.

\begin{prop} \label{prop-so}
Tout faisceau pervers $P$ sur $X$ muni d'une action de $I_v$ s'écrit comme une somme
directe
$\bigoplus_{\sigma_v} P_{\sigma_v}$
où la somme porte sur les classes d'équivalence inertielle des représentations
irréductibles de $W_v$ et où
$P_{\sigma_v}$ est $\sigma_v$-isotypique au sens de la définition (\ref{wo-isotypique}).
En particulier pour tout
point géométrique $x:\spec K \longto X$, $(P_{\sigma_v})_x$ est une représentation
$\sigma_v$-isotypique.
\end{prop}

\begin{coro} Soit $P$ un faisceau pervers sur $X$ muni d'une action compatible de $W_v$.
On considère sa filtration
de
monodromie et on note $gr_k P$ son gradué de poids $k$. La suite spectrale
$E_{1,gr,glob}^{i,j}=h^{i+j} gr_{-i} P \Longrightarrow h^{i+j} P$
est la somme directe sur les classes d'équivalence inertielle des représentations
irréductibles $\sigma_v$ de
$W_v$ des suites spectrales
\begin{equation} \label{ss-gr}
E_{1,gr,glob,\sigma_v}^{i,j} =h^{i+j} gr_{-i} P_{\sigma_v} \Longrightarrow h^{i+j}
P_{\sigma_v}
\end{equation}
\end{coro}

\marque Le complexe des cycles évanescents $R\Psi_v[d-1]$ sur $\bar X_{U^p,m}$ est un
faisceau pervers muni
d'une action de $W_{v}$; on considère alors sa filtration de monodromie et on note
$gr_{k,U^p(m)}$ son gradué
de niveau $k$. On note $N$ le logarithme de la partie unipotente de la monodromie de
sorte qu'en particulier
on a
$N^k:gr_{k,U^p(m)} \simeq gr_{-k,U^p(m)}.$
On pose $gr_k=(gr_{k,U^p(m)})_{U^p,m}$, et on dispose alors de la suite spectrale
$E_{1,U^p(m),gr,glob}^{i,j}=h^{i+j} (gr_{-i,U^p(m)}) \Rightarrow R^{i+j+d-1}\Psi_v$
équivariante sous l'action de $G(\Am^{\oo}) \times W_v$, qui, d'après le corollaire
précédent, se scinde en
une somme directe portant sur les classes d'équivalence inertielle des représentations
irréductibles $\sigma_v$
de $W_v$:
\begin{equation} \label{suite-spectrale}
E_{1,U^p(m),gr,glob,\sigma_v}^{i,j}=h^{i+j} gr_{-i,U^p(m),\sigma_v} R^{i+j+d-1} \Psi_v
\end{equation}

\begin{defi} Pour $\xi$ une représentation algébrique irréductible de $G$ sur $\bar
\Qm_l$, on notera
$gr_{k,U^p(m),\xi}$ le
gradué de $(R\Psi_v \otimes \LC_{\xi})[d-1]$ soit $gr_{k,U^p(m),\xi} \simeq
gr_{k,U^p(m)} \otimes \LC_{\xi}$.
Pour $\sigma_v$ une représentation irréductible de $W_v$, on notera
$gr_{k,U^p(m),\xi,\sigma_v}:=gr_{k,U^p(m),\sigma_v} \otimes \LC_{\xi}$. Pour $\pi_v$ une
représentation
irréductible cuspidale de $GL_g(F_v)$, on notera
$gr_{k,U^p(m),\pi_v}$ (resp. $gr_{k,U^p(m),\xi,\pi_v}$)
pour $gr_{k,U^p(m),\rec_{F_v}^\vee(\pi_v)}$ (resp.
$gr_{k,U^p(m),\xi,\rec_{F_v}^\vee(\pi_v)}$) et de manière
générale on mettra un indice $\pi_v$ en lieu et place de $\rec_{F_v}^\vee(\pi_v)$.
\end{defi}

\subsection{Définition de la catégorie des faisceaux pervers de Hecke}
\label{defi-fph}

Il s'agit ici de donner le cadre catégoriel pour les énoncés et les preuves des
résultats des paragraphes
suivants; ce qui suit m'a été suggéré par Jean-François Dat.

\begin{defi}
On considère un groupe $\tilde G$ et on considère un ensemble $\IC$ tel qu'à tout $I \in
\IC$ soit associé un
sous-groupe compact $\KC_I$ de $\tilde G$. On met sur $\IC$ une relation d'ordre partiel
$J \subset I$ si et
seulement si $\KC_J$ est un sous-groupe distingué de $\KC_I$. Un schéma de Hecke pour
$(\tilde G,\IC)$, est
un système projectif de schémas $X_\IC=(X_I)_{I \in \IC}$
relativement à des morphismes finis $r_{J,I}:X_J \longto X_I$ de restriction du niveau,
tel que pour tout $g
\in \tilde G$ et tout $J \subset I$ tels que $g^{-1} \KC_J g \subset \KC_I$, on dispose
d'un morphisme fini
de schémas $[g]_{J,I}:X_J \longto X_I$
vérifiant les propriétés suivantes:

- pour $g \in \KC_I$ et $J \subset I$, $[g]_{J,I}: X_J \longto X_I$ est le morphisme de
restriction du niveau
$r_{J,I}$;

- pour tout $g,g' \in \tilde G$, et tout $K \subset J \subset I$ tels que $g^{-1} \KC_J
g \subset \KC_I$ et
$(g')^{-1} \KC_K g' \subset \KC_J$, on a
$$[gg']_{K,I}= [g]_{J,I} \circ [g']_{K,J}: X_K \longto X_J \longto X_I.$$

- pour tout $J \subset I$, $X_J/\KC_I \simeq X_I$, où $g \in \KC_I$ agit sur $X_J$ via
$[g]_{J,J}$; autrement
dit $r_{J,I}:X_J \longto X_I$ est un bon quotient de $X_J$ par $\KC_I/\KC_J$.
\end{defi}

\marque \textit{Exemples}: avec les notations du \S \ref{rappels-globaux}, $\IC$ est
l'ensemble des
sous-groupes compacts de $\tilde G:=G(\Am^{\oo})$ de la forme $U^p(m)$ avec $X_{U^p,m}$,
ou $\bar
X_{U^p,m,M_{d-h}}^{(d-h)}$ avec $G:=G(\Am^{\oo,v}) \times P_{h,d}(F_v)$.

\begin{defi} Soit $X_\IC=(X_I)_I$ un schéma de Hecke pour $(\tilde G,\IC)$, on définit
alors la catégorie $\FPH_{\tilde G}(X_\IC)$ (resp. $\FH_{\tilde G}(X_\IC)$) des
faisceaux pervers (resp. des
faisceaux) de Hecke sur $X_\IC$ comme la catégorie dont les objets sont les systèmes
$(\FC_I)_{I \in \IC}$ où
$\FC_I$ est un faisceau pervers (resp. faisceau) sur $X_I$ tels que

\begin{itemize}
\item pour tout $g \in \tilde G$ et $J \subset I$ tel que $g^{-1} \KC_J g \subset
\KC_I$, on dispose d'un morphisme
de faisceau sur $X_I$, $u_{J,I}(g):\FC_I \longto [g]_{J,I,*} \FC_J$ soumise à la
condition de cocycle habituelle
$$u_{K,I}(g'g)=[g]_{J,I,*} (u_{K,J}(g')) \circ u_{J,I}(g)$$

\item pour tout $g \in \KC_I$, $u_{J,I}(g)$ se factorise par les $\KC_I$ invariants,
i.e. induit une flèche
$\FC_I \longto r_{J,I,*} \FC_J^{\KC_I}$ où $\FC_J^{\KC_I}$ désigne le sous-faisceau de
$\FC_J$ invariant par
tous les $u_{J,J}(g)$ où $g$ décrit $\KC_I$.
\end{itemize}

Les flèches sont alors les systèmes $(f_I)_{I \in \IC}$ avec $f_I: \FC_I \longto \FC'_I$
tel que le diagramme
suivant soit commutatif:
$$\diagram
\FC_I \rrto^{u_{J,I}(g)} \dto^{f_I} & & [g]_{J,I,*} (\FC_J) \dto^{[g]_{J,I,*} (f_J)} \\
\FC_I' \rrto^{u_{J,I}(g)} & & [g]_{J,I,*} (\FC_J')
\enddiagram$$
\end{defi}

\begin{prop} Pour tout schéma de Hecke $X_\IC$, la catégorie $\FH_{\tilde G}(X_\IC)$
(resp. $\FPH_{\tilde G}(X_\IC)$)
est abélienne (resp. abélienne et artinienne).
\end{prop}

\begin{proof} Elle ne présente aucune difficulté. Montrons par exemple l'existence d'un
noyau pour
$f=(f_I:\FC_I \to \FC_I')_I$. Les catégories à niveau fini étant abélienne, soit pour
tout $I$, $\GC_I$ un
noyau de $f_I$. On a alors le diagramme commutatif suivant:

$$\xymatrix{
\GC_I \ar[r] \ar@{-->}[d] & \FC_I \ar[r]^{f_I} \ar[d]^{u_{J,I}(g)} & \FC_I'
\ar[d]^{u_{J,I}(g)} \\
[g]_{J,I,*} (\GC_J) \ar[r] \ar@{-->}[d] & [g]_{J,I,*} (\FC_J) \ar[r]^{[g]_{J,I,*} (f_J)}
\ar[d]^{[g]_{J,I,*}
(u_{K,J}(g'))} & [g]_{J,I,*} (\FC_J') \ar[d]^{[g]_{J,I,*}( u_{K,J}(g'))} \\
[g'g]_{K,J,*} (\GC_K) \ar[r] & [g'g]_{K,J,*} (\FC_K) \ar[r]^{[g'g]_{K,J,*}( f_K)} &
[g'g]_{K,J,*} (\FC_K')
}$$ De la propriété universelle du noyau $\GC_J$ (resp. $\GC_K$), on en déduit une
flèche
$$u_{J,I}(g):\GC_I \longto [g]_{J,I,*} (\GC_J)  \qquad (\hbox{resp.
}[g]_{K,J,*}(u_{K,J}(g')): [g]_{J,I,*} (\GC_J)
\longto [g'g]_{K,J,*} (\GC_K))$$
la propriété de cocycle pour $\GC$ découle alors de celles pour $\FC$ et $\FC'$ via la
commutativité du
diagramme précédent. De la même façon, en utilisant la propriété universelle du noyau et
la commutativité du
diagramme suivant
$$\xymatrix{
\GC_I \ar[r] \ar@{-->}[d] & \FC_I \ar[r]^{f_I} \ar[d]^{u_{J,I}(g)} & \FC_I'
\ar[d]^{u_{J,I}(g)} \\
r_{J,I,*} (\GC_J)^{\KC_I} \ar[r] & r_{J,I,*} (\FC_J)^{\KC_I} \ar[r]^{r_{J,I,*} (f_J)} &
r_{J,I,*}
(\FC_J')^{\KC_I} }$$ on obtient que pour tout $g \in \KC_I$, $u_{J,I}(g)$ induit une
flèche $\GC_I \longto
r_{J,I,*} \GC_J^{\KC_I}$.

\end{proof}

\marque \rem Dans la suite du texte, les faisceaux de Hecke que l'on considérera seront
tels que pour $g \in
\KC_{I}$, $u_{J,I}(g)$ induit un isomorphisme $\FC_I \longto r_{J,I,*} \FC_J^{\KC_I}$,
de sorte que par
exemple, les invariants sous $\KC_I$ de la cohomologie en niveau $J$, est isomorphe à la
cohomologie en
niveau $I$. Cette propriété est clairement conservée par passage au noyau mais ne l'est
pas pour les
conoyaux. Cela étant, on remarquera que les constructions des paragraphes suivants se
font toujours par des
noyaux à partir d'une flèche entre deux faisceaux vérifiant cette propriété
additionnelle de sorte que tous
les faisceaux que l'on construit la vérifient aussi.

\begin{nota} \label{nota-coho}
Pour $\FC=(\FC_I)_{I \in \IC}$ un objet de $\FPH_{\tilde G}(X_\IC)$ ou de $\FH_{\tilde
G}(X_\IC)$, on notera
$H^i(\FC)$ pour la limite inductive des $H^i(X_I,\FC_I)$.
\end{nota}

\begin{prop} \label{prop-def-fph}
Soient $\bar X_\IC$, $X_\IC$ et $Y_\IC$ des schémas de Hecke pour $(\tilde G,\IC)$ tels
que $j_\IC:X_\IC
\hookrightarrow \bar X_\IC$ soit un système projectif d'inclusions affines
$G$-équivariantes et $Y_\IC = \bar
X_\IC \backslash X_\IC$. On dispose alors des foncteurs
$j_!,Rj_*,j_{!*}:\FPH_G(X_\IC) \longto \FPH_G(\bar X_\IC)$, $\lexp p j^*=\lexp p
j^!:\FPH_G(\bar X_\IC) \longto
\FPH_G(X_\IC)$ et pour $i:Y_\IC:=\bar X_\IC - X_\IC \hookrightarrow \bar X_\IC$ des
foncteurs
$i_*=i_!:\FPH_G(Y_\IC) \longto \FPH_G(\bar X_\IC)$,  $\lexp p i^*, ~\lexp p
Ri^!:\FPH_G(\bar X_\IC) \longto
\FPH_G(Y_\IC)$.
\end{prop}

\begin{proof} On commence par rappeler la proposition suivante tirée de \cite{ast}.

\begin{prop} \label{prop-ast0} (cf. proposition 1.3.17 de \cite{ast})
Pour tout foncteur exact $f$, on note $\lexp p f$ pour $\lexp p H^0 f$ comme foncteur
sur la catégorie des
faisceaux pervers.

\begin{itemize}
\item[(i)] Si $f$ est $t$-exact à gauche (resp. à droite), $\lexp p f$ est alors exact à
gauche (resp. à droite). Par
ailleurs pour $\FC$
dans $D^{\geq 0}$ (resp. $D^{\leq 0}$), on a
$$\lexp p f (\lexp p H^0 \FC) \longmapright{\sim} \lexp p H^0 f(\FC) \quad (\hbox{resp.
} \lexp p H^0 f(\FC)
\longmapright{\sim} \lexp p f (\lexp p H^0 \FC)).$$

\item[(ii)] Soit $(f^+,f_+)$ une paire de foncteurs exacts adjoints; pour que $f^+$ soit
$t$-exact à droite, il faut
et il suffit que $f_+$
soit $t$-exact à gauche et dans ce cas $(\lexp p f^+,\lexp p f_+)$ forment une paire de
foncteurs adjoints.

\item[(iii)] Si $f_1:D^b_{c,1} \to D^b_{c,_2}$ et $f_2:D^b_{c,2} \to D^b_{c,3}$ sont
$t$-exacts à gauche
(resp. à droite), alors $f_2 \circ f_1$ l'est aussi et $\lexp p (f_2 \circ f_1)= \lexp p
f_2 \circ f_1$.

\end{itemize}
\end{prop}

On rappelle que $j_!$, $Rj_*$, $j^*$, $i_*$ sont $t$-exacts et donc égaux à leur version
perverse alors que
$i^*$ est $t$-exact à droite. De même les morphismes $r_{J,I}$ et $[g]_{J,I}$ étant
finis, $r_{J,I}^*$ (resp.
$[g]_{J,I,*}$) est $t$-exact à droite (resp. $t$-exact) de sorte que $(\lexp p
r_{J,I}^*,\lexp p r_{J,I,*})$
et $(\lexp p [g]_{J,I,*}=[g]_{J,I,*},\lexp p [g]_{J,I}^!)$ forment des paires de
foncteurs adjoints avec
$\lexp p r_{J,I}^*$ et $\lexp p [g]_{J,I}^!$ (resp. $\lexp p r_{J,I,*}$, resp. $\lexp p
[g]_{J,I,*}$) exacts
à droite (resp. à gauche, resp. exact). Pour $J \subset I$, considérons le diagramme
suivant:
$$\xymatrix{ X_J \ar@{^{(}->}[rr]^{j_J} \ar[d]^{[g]_{J,I}} & & \bar X_J
\ar[d]_{\overline{[g]}_{J,I}}
& & Y_J \ar@{_{(}->}[ll]_{i_J} \ar[d]^{\widetilde{[g]}_{J,I}} \\
X_I \ar@{^{(}->}[rr]^{j_I} & & \bar X_I & & Y_I \ar@{_{(}->}[ll]_{i_I} }$$ Soient alors
$(\FC_I)_{I \in \IC}$
un objet de $\FPH_{\tilde G}(X_\IC)$ et $g,J,I$ avec $u_{J,I}(g): \FC_I \longto
[g]_{J,I,*} (\FC_I)$. Par
application de $Rj_{\IC,*}$ (resp. $j_{\IC,!}$), on obtient
$Rj_{I,*} \FC_I \longmapright{Rj_{I,*}(u_{J,I}(g))} Rj_{I,*} [g]_{J,I,*}
(\FC_J)=\overline{[g]}_{J,I,*}
(Rj_{J,*} \FC_J)$
(resp. $j_{I,!} \FC_I \longmapright{j_{I,!}(u_{J,I}(g))} j_{I,!}
[g]_{J,I,!}(\FC_J)=\overline{[g]}_{J,I,*} (j_{J,!} \FC_J)$).

En ce qui concerne $j_{\IC,!*}$, on considère le diagramme commutatif de la figure
(\ref{fig-diag}) où la
flèche notée (1) s'obtient à partir du diagramme commutatif suivant
$$\diagram Y_J \rto^{\widetilde{[g]}_{J,I}} \dto_{i_J} & Y_I \dto^{i_I} \\ X_J
\rto^{[g]_{J,I}} & X_I \enddiagram$$
comme le morphisme de changement de base $\bar i_I^* \circ [g]_{J,I,*} \longto
[g]_{J,I,*} \circ i_J^* $.

\begin{figure} \label{fig-diag}
$$\xymatrix{
0 \to j_{I,!*} (\FC_I) \ar[r] \ar[ddddddr] \ar@{-->}[dddddd] & Rj_{I,*} (\FC_I)
\ar[d]_{Rj_{I,*}(u_{J,I}(g))}
\ar[r] & i_{I,*} \circ \lexp p i_I^* \circ Rj_{I,*} (\FC_I) \ar[d]_{i_{I,*}
\circ \lexp p i_I^* \circ Rj_{I,*}(u_{J,I}(g))} \to 0 \\
 & Rj_{I,*} \circ [g]_{J,I,*} (\FC_I) \ar[ddddd]^{=} \ar[r] &
 i_{I,*} \circ \lexp p i_I^* \circ Rj_{I,*} \circ [g]_{J,I,*} (\FC_I) \ar[d]^{=} \\
 & & i_{I,*} \circ \lexp p i_I^* \circ \overline{[g]}_{J,I,*} \circ Rj_{J,*} (\FC_J)
\ar[d]^{=} \\
 & & i_{I,*} \circ \lexp p (i_I^* \circ \overline{[g]}_{J,I,*}) \circ Rj_{J,*} (\FC_J)
\ar[d]^{(1)} \\
 & & i_{I,*} \circ \lexp p (\widetilde{[g]}_{J,I,*} \circ i_J^*) \circ Rj_{J,*} (\FC_J)
\ar[d]^{=} \\
 &  &
i_{I,*} \circ \widetilde{[g]}_{J,I,*} \circ \lexp p i_J^* \circ Rj_{J,*} (\FC_J)
\ar[d]^{=} \\
0 \to \overline{[g]}_{J,I,*} \circ j_{J,!*} (\FC_J) \ar[r] & \overline{[g]}_{J,I,*}
\circ Rj_{J,*} (\FC_J)
\ar[r] & \overline{[g]}_{J,I,*} \circ i_{J,*} \circ \lexp p i_J^* \circ Rj_{J,*} (\FC_J)
\to 0 }$$
\end{figure}

Les cas de $i_*$, $\lexp p i^*$ sont traités ci dessus et ceux de $\lexp p j^*$ et
$\lexp p Ri^!$ se traitent
de manière strictement similaire.

\end{proof}

\marque \textit{Exemples}: étant donné un système local $\LC_\emptyset$ sur
$X_\emptyset$, soit pour $I \in
\IC$, $\LC_I=r_{I,\emptyset}^* \LC_\emptyset$; $(\LC_I[d])_{I \in \IC}$ est alors un
objet de $\FPH_{\tilde
G}(X_\IC)$, où $d$ est la dimension de $X_\IC$. Si en outre les $X_I$ sont la fibre
générique d'un $S$-schéma
$\XC_I$ de dimension relative $d-1$, pour $S$ un trait, alors
$(R\Psi_{\eta}(\LC_I)[d-1])_{I \in \IC}$ est un
objet de $\FPH_{\tilde G}(\bar X_\IC)$ où $\bar X_I$ désigne la fibre spéciale de
$X_I$.

\subsection{Notations}

\begin{defis}\label{defi-type}
\begin{itemize}
\item Pour $1 \leq g \leq d$, on notera $s_g:=\lfloor \frac{d}{g} \rfloor$, la partie
entière
de $d/g$.

\item On introduit les injections
$$i^{h}:\bar X_{U^p,m}^{[d-h]} \hookrightarrow \bar X_{U^p,m}, \qquad j^{ \geq h}: \bar
X_{U^p,m}^{(d-h)}
\hookrightarrow \bar X_{U^p,m}^{[d-h]} \qquad k^{=h}_1:\bar X_{U^p,m,M_{d-h}}^{(d-h)}
\hookrightarrow \bar
X_{U^p,m}^{(d-h)}$$

\item On notera $\FPH(\bar X)$ (resp. $\FPH(X)$, resp. $\FPH(\bar X^{(h)})$,
resp. $\FPH(\bar X^{[h]})$) la catégorie des faisceaux pervers de Hecke sur la tour des
$\bar X_{U^p,m}$
(resp. $X_{U^p,m}$, resp. $\bar X_{U^p,m}^{(h)}$, resp. $\bar X_{U^p,m}^{[h]}$) avec
$\tilde G=G(\Am^\oo)$.

\item On notera aussi $\FPH(\bar X_{M_{d-h}}^{(d-h)})$ (resp. $\FPH(\bar
X_{M_{d-h}}^{[d-h]})$) la catégorie
des faisceaux pervers de Hecke sur la tour des $\bar X_{U^p,m,M_{d-h}}^{(d-h)}$ (resp.
$\bar
X_{U^p,m,M_{d-h}}^{[d-h]}$) avec $\tilde G=G(\Am^{\oo,v}) \times P_{h,d}(F_v)$. On
utilisera des notations
similaires pour les faisceaux de Hecke.

\item Un système projectif de faisceaux pervers ou pas, non nul $\FC=(\FC_I)_{I \in \IC}
\in \FPH(\bar X^{(d-h)})$
(ou dans $\FH(\bar X)$) sera dit \textrm{induit} s'il est de la forme
$$(k^{=h,*}_1 \FC_I)_{I \in \IC} \times_{P_{h,d}^{op}(F_v)} GL_d(F_v)$$
où $\times_{P_{h,d}^{op}(F_v)} GL_d(F_v):\FPH(\bar X_{M_{d-h}}^{(d-h)}) \longto
\FPH(\bar X^{(d-h)})$ est le
foncteur qui à un faisceau $(\FC_{0,I})_I$ associe $(\FC_{0,I}
\times_{P_{h,d}^{op}(\OC_v)} GL_d(\OC_v))_I$
de sorte que pour tout $g \in GL_d(F_v)$, la correspondance associée induit un
isomorphisme de $g^* \FC_{\bar
g} \longmapright{\sim} \FC_0$ où $\FC_{\bar g}$ est la restriction de $\FC$ à la
composante $(\bar
X_{U^p,m,M_g}^{(d-h)})_I$ où $M_g$ est le sous-module associé à l'image de $g$ dans le
quotient
$GL_d(F_v)/P_{h,d}(F_v)$.

\item Pour $tg<d$, on notera $\FC(g,t,\pi_v)_1=(\FC(g,t,\pi_v,I)_1)_{I \in \IC}$ le
système local
$$\FC_{\JL^{-1}([\overleftarrow{t-1}]_{\pi_v^\vee}),1}=(\FC_{\JL^{-1}([\overleftarrow{t-1}]_{\pi_v^\vee}),I,1})_{I
\in \IC}$$ sur la tour des $(\bar X_{U^p,m,M_{d-tg}}^{(d-tg)})_{I=U^p(m)}$ et soit
$\FC(g,t,\pi_v)=(\FC(g,t,\pi_v,I))_I$ le faisceau induit associé, défini donc sur la tour des $(\bar
X_{U^p,m}^{(d-tg)})_{I=U^p(m)}$: $\FC(g,t,\pi_v)[d-tg] \in \FPH(\bar X^{(d-tg)})$.

\item Pour $1 \leq tg \leq d$, un $GL_d(F_v) \times W_v$-faisceau sur la tour $(\bar
X_{U^p,m}^{(d-tg)})_{I=U^p(m)}$ sera dit de type $HT(g,t)$ s'il est induit, et si sa restriction à $(\bar
X_{U^p,m,M_{d-tg}}^{(d-tg)})_{I=U^p(m)}$ est de la forme
$$(\FC(g,t,\pi_v,I)_1 \otimes \Pi_v^{\KF_{v,m_1}} \otimes \chi)_I$$
pour un certain triplet $(\Pi_v,\pi_v,\chi)$ où:

- $\pi_v$ est une représentation irréductible cuspidale de $GL_g(F_v)$,

- $\Pi_v$ est une représentation de $GL_{tg}(F_v)$

- $\chi$ un caractère de $\Zm$,

- et $\KF_{v,m_1}:=\ker (GL_{tg}(\OC_v) \longto GL_{tg}(\OC_v/v^{m_1}))$ \footnote{On ne met pas en indice de
référence à la dimension $tg$ car en général le contexte est clair.}

\noindent telle que l'action d'un élément $(g^{\oo,p},g_{p,0},g_v^0,g_v^{et},g_{v_i},c_v) \in G(\Am^{\oo,p})
\times \Qm_p^\times \times GL_{tg}(F_v) \times GL_{d-tg}(F_v) \times \prod_{i=2}^r (B_{v_i}^{op})^\times
\times W_v$ soit donnée par l'action naturelle de
$$(g^{\oo,p},g_{p,0}p^{-f_1v(c_v)},g_v^{et},g_{v_i},v(\det g_v^0) -v(c_v)) \in G(\Am^{\oo,p}) \times \Qm_p^\times
\times GL_{d-tg}(F_v) \times \Zm$$ sur $\FC_{\JL^{-1}([\overleftarrow{t-1}]_{\pi_v}),1}$ au dessus de $(\bar
X_{U^p,m,M_{d-tg}}^{(d-tg)})_{I=U^p(m)}$, et par celle de $(g_v^0,c_v)$ sur $\Pi_v \otimes \chi$. On le
notera alors
$$HT(g,t,\pi_v,\Pi_v,\chi)=(HT(g,t,\pi_v,\Pi_v,\chi,I))_I, \quad
HT(g,t,\pi_v,\Pi_v,\chi)[d-tg] \in \FPH(\bar
X^{(d-tg)})$$
et on omet $\chi$ si celui-ci est trivial. De même le faisceau $HT(g,t,\pi_v,\Pi_v)
\otimes \LC_{\xi}$ sera
noté
$$HT_{\xi}(g,t,\pi_v,\Pi_v)=(HT_{\xi}(g,t,\pi_v,\Pi_v,I))_I, \quad
HT_{\xi}(g,t,\pi_v,\Pi_v)[d-tg] \in \FPH(\bar
X^{(d-tg)})$$

\item Pour $\pi_v$ une représentation irréductible cuspidale de $GL_g(F_v)$ et $1 \leq t
\leq s_g$, on note
$$\PC(g,t,\pi_v)=(\PC(g,t,\pi_v,I))_I \in \FPH(\bar X)$$
le faisceau pervers défini comme l'extension intermédiaire $j^{\geq tg}_{!*}$ du
faisceau pervers sur $\bar
X_{U^p,m}^{(d-tg)}$, $HT(g,t,\pi_v,[\overleftarrow{t-1}]_{\pi_v})[d-tg] \otimes
\rec_{F_v}^\vee(\pi_v)$, qui
est donc pur de poids $d-tg$ pour $\pi_v$ unitaire \footnote{sinon on rajoute le poids
de
$\rec_{F_v}^\vee(\pi_v)$} et qui ne dépend que de la classe d'équivalence inertielle de
$\pi_v$.
\end{itemize}
\end{defis}

\begin{rema} \label{rema-I} $\FC(g,t,\pi_v)$ et $HT(g,t,\pi_v,\Pi_v)$ (resp.
$\PC(g,t,\pi_v)$) se présentent
sous la forme d'une somme directe de $e_{\pi_v}$ systèmes locaux (resp. faisceaux
pervers) irréductibles:
$$\FC(g,t,\pi_v)=\bigoplus_{i=1}^{e_{\pi_v}} \FC(g,t,\rho_{v,i}) \qquad (\hbox{resp. }
\PC(g,t,\pi_v)=\bigoplus_{i=1}^{e_{\pi_v}} \PC(g,t,\rho_{v,i}))$$ donnés par la
restriction de
$\JL^{-1}([\overleftarrow{t-1}]_{\pi_v^\vee})$ à $\DC_{v,tg}^\times$ qui s'écrit comme
une somme directe
$\bigoplus_{i=1}^{e_{\pi_v}} \rho_{v,i}$ de représentations irréductibles, de sorte que
la différence entre
les faisceaux $\FC(g,t,\rho_{v,i})$ (resp. $\PC(g,t,\rho_{v,i})$) provient de l'action
de $\DC_{v,tg}^\times
\subset \NC_v$ donnée sur chacun comme dans la formule (\ref{action-tordue}).
\end{rema}

\begin{defi} \label{defi-gf}
Soit $\GF=(\GF_I)_I$ le système projectif de groupes de Grothendieck associés à
$\FPH(\bar X_I)$.
\end{defi}

\subsection{Filtration de monodromie}
\label{enonces-glob}

\begin{theo} \label{theo-global0}
Soit $1 \leq g \leq d$ et $s_g$ la partie entière de $d/g$. Pour $\pi_v$ une
représentation irréductible
cuspidale de $GL_g(F_v)$, on a l'égalité suivante dans $\GF$:
$$e_{\pi_v} [R\Psi_{\pi_v}[d-1]]= \sum_{k=1-s_g}^{s_g-1} \sum_{\atop{|k| < t \leq s_g}{t
\equiv k-1 \mod 2}}
\PC(g,t,\pi_v)(-\frac{tg+k-1}{2}),$$ où la torsion concerne l'action de $W_v$.
\end{theo}

\marque Pour tout $s$, soit $MGr_k(s)$ le faisceau pervers nul pour $|k| \geq s$ et égal
à
$$\bigoplus_{\atop{|k| <t \leq s}{t \equiv k-1 \mod 2}} j^{\geq t}_{!*}
\FC(1,t,1_v)[d-t] \otimes
[\overleftarrow{t-1}]_{1_v} \otimes \rec_{F_v}^\vee(1_v) (-(t-1+k)/2)$$
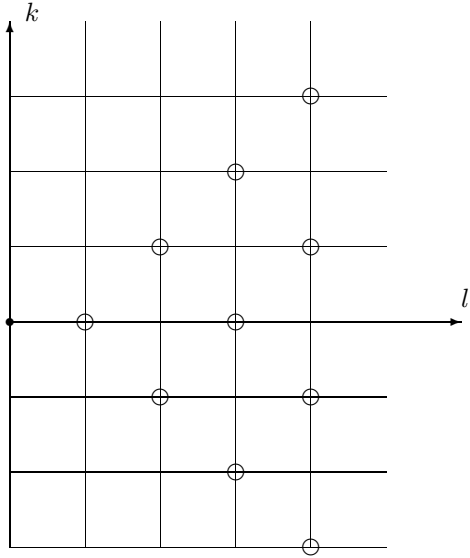
\begin{figure}[!ht]
\caption{\label{figgg0} Pour $s=4$, dans le repère $(t,k)$, on marque sur la ligne
$k=k_0$ les $t$ tels que
$j^{\geq t}_{!*} \FC(1,t,1_v)[d-t] \otimes \st_t(1_v) \otimes \rec_{F_v}^\vee(1_v)
(-(t-1+k)/2)$ soit un
constituant de $MGr_{k_0}(4)$} \setlength{\unitlength}{1cm} \centering
\begin{picture}(6,8.5)(0,-4)
\linethickness{.1pt}

\multiput(0,-3)(0,1){7}{\line(1,0){5}} \multiput(0,-3)(1,0){5}{\line(0,1){7}}

\put(5,0){\vector(1,0){1}} \put(6,.2){$l$}

\put(0,3){\vector(0,1){1}} \put(0.2,4){$k$}

\multiput(1,0)(1,1){4}{\circle{.2}}

\multiput(2,-1)(1,1){3}{\circle{.2}}

\multiput(3,-2)(1,1){2}{\circle{.2}}

\put(4,-3){\circle{.2}}

\put(0,0){\circle*{.1}}

\end{picture}
\end{figure}

\marque Pour tout représentation cuspidale $\pi_v$ de $GL_g(F_v)$, on pose
\begin{multline*}
\pi_v \diamond \bigl ( j^{\geq t}_{!*} \FC(1,t,1_v)[d-l] \otimes
[\overleftarrow{t-1}]_{1_v} \otimes
\rec_{F_v}^\vee(1_v) (-(t-1+k)/2) \bigr ) := \\
j^{\geq tg}_{!*} \FC(g,t,\pi_v)[d-tg] \otimes [\overleftarrow{t-1}]_{\pi_v} \otimes
\rec_{F_v}^\vee(\pi_v)
(-(tg-1+k)/2)
\end{multline*}

\begin{theo} \label{theo-global1}
Soit $1 \leq g \leq d$ et $s_g$ la partie entière de $d/g$. Pour $\pi_v$ une
représentation irréductible
cuspidale de $GL_g(F_v)$, on a $e_{\pi_v} gr_{\bullet,\pi_v} = \pi_v \diamond
MGr_\bullet(s_g)$.
\end{theo}

\marque Autrement dit pour $1 \leq g \leq d$, $\pi_v$ une représentation irréductible
cuspidale de
$GL_g(F_v)$ et pour tout $|k| < s_g$, $e_{\pi_v} gr_{k,\pi_v}$ est égale à
$$\bigoplus_{\atop{|k| < t \leq s_g}{t \equiv k-1 \mod 2}}
\PC(g,t,\pi_v)(-\frac{tg+k-1}{2}),$$
dans $\FPH(\bar X)$, où la torsion concerne l'action de $W_v$.

\begin{defi} \label{defi-pil}
Pour $1 \leq g \leq d$, et on fixe une représentation $\pi_v$ irréductible cuspidale de
$GL_g(F_v)$, que l'on
notera $\chi_v$ pour $g=1$. Pour tout $1 \leq t \leq s_g$, $\Pi_t$ désigne une
représentation quelconque de
$GL_{tg}(F_v)$, qui dans les applications sera elliptique de type $\pi_v$.
\end{defi}

\begin{theo} \label{theo-global2}
Pour $g$, $\pi_v$ et $\Pi_t$ comme ci-dessus, on considère la restriction à la tour des
$(\bar
X_{U^p,m}^{(d-h)})_{I=U^p(m)}$ de $h^i j^{\geq tg}_{!*} HT(g,t,\pi_v,\Pi_t)[d-tg]$;

\begin{itemize}
\item elle est nulle pour $h$ ne s'écrivant pas sous la forme $(t+a)g$ avec $0 \leq a
\leq s_g-t$;

\item pour $h=(t+a)g$ avec $0 \leq a \leq s_g-t$, elle est nulle pour $i \neq
tg-d+a(g-1)$ et sinon elle est
isomorphe dans $\FH(\bar X)$ à \footnote{cf. la définition (\ref{defi-carac})}
$HT(g,t+a,\pi_v,\Pi_t
\overrightarrow{\times} [\overrightarrow{a-1}]_{\pi_v}) \otimes
\Xi^{\frac{a(g-1)}{2}}$.
\end{itemize}
\end{theo}

\begin{coro} Pour $g>1$, $h^{i} j^{\geq tg}_{!*} HT(g,t,\pi_v,\Pi_t)[d-tg]$
est nul pour $i$ ne s'écrivant pas sous la forme $tg-d+a(g-1)$ avec $0 \leq a \leq
s_g-t$ et pour un tel
$i=tg-d+a(g-1)$ il est isomorphe dans $\FH(\bar X)$ à $j^{\geq
(t+a)g}_!HT(g,t+a,\pi_v,\Pi_t
\overrightarrow{\times} [\overrightarrow{a-1}]_{\pi_v}) \otimes
\Xi^{\frac{a(g-1)}{2}}.$
\end{coro}

\begin{coro} \label{coro-higr}
Pour $\pi_v$ une représentation irréductible cuspidale de $GL_g(F_v)$ pour $1 \leq g
\leq d$, on a
\begin{itemize}
\item[(i)] pour $g>1$, $(h^i gr_{k,\pi_v})^{e_{\pi_v}}$ est une somme directe sur tous
les couples $(t,a)$ tels que
$-d+tg+a(g-1)=i$ et $t \equiv k-1 \mod 2$ des faisceaux induits dans $\FH(\bar X)$:
$$j^{\geq (t+a)g}_!HT(g,t+a,\pi_v,[\overleftarrow{t-1}]_{\pi_v} \overrightarrow{\times}
[\overrightarrow{a-1}]_{\pi_v})
\otimes \rec_{F_v}^\vee(\pi_v)(-\frac{tg-1+k+a(g-1)}{2})$$

\item[(ii)] pour $g=1$, la restriction de $h^i gr_{k,\chi_v}$ pour $i \equiv k-1 \mod
2$, à la tour des $\bar
X_{U^p,m}^{(d-t)}$ pour $t \geq d+i$ est égale au faisceau induit dans $\FH(\bar X)$
$$HT(1,t,\chi_v,[\overleftarrow{i-1}]_{\chi_v} \overrightarrow{\times}
[\overrightarrow{t-i-1}]_{\chi_v}) \otimes
\chi_v (-\frac{d+i-1+k}{2})$$ et nulle dans tous les autres cas.

\end{itemize}
\end{coro}

\marque On introduit les bicomplexes
$$MGE_\bullet^{\bullet,\bullet}(s)=\bigoplus_{t=1}^s j^{\geq l}_! \FC(1,t,1_v)[d-t]
\otimes MLE_\bullet^{\bullet,\bullet}(t)$$ de sorte que le point (2) du théorème (\ref{theo-local-fil}) et
donc (\ref{theo-ripsi-local}), découle du résultat suivant.

\begin{theo}  \label{theo-ss}
Soit $1 \leq g \leq d$ et $s_g$ la partie entière de $d/g$. Pour $\pi_v$ une représentation irréductible
cuspidale de $GL_g(F_v)$, on a $e_{\pi_v} E_{\bullet,gr,glob,\pi_v}^{\bullet,\bullet}(\frac{d-1}{2})=\pi_v
\diamond MGE_\bullet^{\bullet,\bullet}(s_g)$.
\end{theo}

\marque Autrement dit, la suite spectrale (\ref{suite-spectrale}) dégénère en $E_2$ et
les applications
$d_1^{i,j}$ sont induites par les suites exactes (\ref{suites-exactes}).

\rem A la figure (\ref{fig0}) on a représenté la suite spectrale (\ref{suite-spectrale})
dans le cas $g=6$ et
$s=4$.

\section{Schéma de la preuve: par récurrence} \label{schema}

\marque Remarquons tout d'abord que le théorème (\ref{theo-global1}) découlerait de
(\ref{theo-global0}) si
on disposait de la conjecture de monodromie-poids (MP) version faisceautique, i.e. la
pureté des gradués de
la filtration de monodromie. D'après le théorème de comparaison de Berkovich-Fargues
(BF), le théorème
(\ref{theo-local-fil}) se déduit de la détermination des germes en un point
supersingulier des gradués de la
filtration de monodromie du complexe des cycles évanescents de la variété de Shimura et
de la suite spectrale
associée, soit donc des théorèmes (\ref{theo-global1}), (\ref{theo-global2}) et
(\ref{theo-ss}). Par
ailleurs, comme on l'a déjà noté, le théorème (\ref{theo-ripsi-local}) découle
directement de
(\ref{theo-local-fil}) puisqu'il s'agit de l'aboutissement de la suite spectrale
associée à la filtration de
monodromie-locale. On résume ces implications dans le tableau suivant.
\begin{figure}[!h]
\begin{center}
\begin{tabular}{|ll|}
\hline (\ref{theo-global0}) & \\ $\Downarrow$ MP? & \\ (\ref{theo-global1}) + &
(\ref{theo-global2}) +
(\ref{theo-ss}) \\ \hline & $\Downarrow$ BF \\ \hfill (\ref{theo-local-fil})& en hauteur
$\leq d$ \\ \hline &
$\Downarrow$ \\ \hfill (\ref{theo-ripsi-local}) & en hauteur $\leq d$ \\ \hline
\end{tabular}
\end{center}
\end{figure}
Le principe est alors de raisonner par récurrence en supposant connu
(\ref{theo-local-fil}) pour tout $d'<d$.

\begin{rema} \label{rema-lim}
A priori, il nous faut travailler à niveau fini et montrer que les résultats obtenus
sont compatibles aux
correspondances de Hecke; cependant devant l'immédiateté de ces vérifications et afin de
ne pas alourdir
encore les preuves, on raisonne souvent sur le système inductif total. En particulier
nous ne marquerons plus
la référence au niveau dans les notations; par exemple on écrira $\FC_{\tau_v}$ au lieu
de $\FC_{\tau_v,U^p(m)}$.
\end{rema}

\marque Dans un premier temps le théorème de comparaison de Berkovich-Fargues (BF) et la
description (HT) des
faisceaux des cycles évanescents en fonction des systèmes locaux d'Harris-Taylor (cf.
(\ref{iso1})), nous
permet, \S \ref{glob0}, de montrer le théorème (\ref{theo-global0}). On montre ensuite,
proposition
(\ref{prop-hij}), le théorème (\ref{theo-global2}) hors des points supersinguliers et on
obtient un contrôle
(C) sur ce qui peut se passer au niveau de ces derniers. Selon le même principe, si on
disposait de
(\ref{theo-local-fil}) pour la hauteur $d$, on en déduirait (\ref{theo-global2}) et
(\ref{theo-ss}). On
démontre, \S \ref{globss}, comment la connaissance des parties de poids minimal de
(\ref{theo-ripsi-local})
permet d'obtenir (\ref{theo-global2}), proposition (\ref{prop-hij-ss-g}) \footnote{de
sorte que
(\ref{theo-global2})=(\ref{prop-hij})+(\ref{prop-hij-ss-g})}, (\ref{theo-ss}), \S
\ref{globss2} et enfin \S
\ref{globalss3},
la conjecture de monodromie-poids version faisceautique i.e. le théorème
(\ref{theo-global1}). Reste alors, \S
\ref{ssce}, à
calculer les parties de poids minimal de (\ref{theo-ripsi-local}) ce qui se fait via
l'étude de la suite
spectrale des cycles évanescents (SSCE) à travers les suites spectrales associées à la
stratification (SSS)
en utilisant le contrôle (C) sur ce qui peut se passer au niveau des points
supersinguliers. Evidemment tout
ceci passe par le calcul, ou au moins le contrôle, des groupes de cohomologie (CHT) des
systèmes locaux
d'Harris-Taylor, qui sont distillés au fur et à mesure et dont nous détaillons
l'enchainement logique après
le tableau récapitulatif suivant.

\begin{figure}[!h]
\begin{center}
\begin{tabular}{|ccl|}
\hline
\begin{tabular}{l} \S \ref{glob0}: (\ref{theo-local-fil}) \\ en hauteur $d'<d$
\end{tabular}
& \begin{tabular}{c} $\Longrightarrow$ \\ BF+HT \end{tabular} & (\ref{theo-global0})  \\
\hline
\begin{tabular}{l} \S \ref{glob2}: (\ref{theo-local-fil}) \\ en hauteur $d'<d$
\end{tabular}
& \begin{tabular}{c} BF+HT \\ $\Longrightarrow$ \\ (\ref{prop-hij}) \end{tabular} &
\begin{tabular}{l} (\ref{theo-global2}) hors des points supersinguliers \\ + C=contrôle
aux points supersinguliers
\end{tabular} \\ \hline
\begin{tabular}{l} (\ref{theo-local-fil}) \\ en hauteur $d' \leq d$ \end{tabular}
& \begin{tabular}{c} $\Longrightarrow$ \\ (\ref{prop-hij-ss-g}) \end{tabular} &
(\ref{theo-global2})+ (\ref{theo-ss})
\\ \hline
\begin{tabular}{l} \S \ref{globss}: (\ref{theo-ripsi-local}) \\ en hauteur $d' < d$ \\
et en poids minimal \\
pour la hauteur $d$ \end{tabular} &
\begin{tabular}{c} $\Longrightarrow$ \\ C+N \end{tabular} &
$\left \{ \begin{array}{l} \hbox{(\ref{theo-global2}):\S \ref{globss1} proposition
(\ref{prop-hij-ss-g})} \\
\hbox{(\ref{theo-ss}): \S \ref{globss2} } \\ \hbox{(\ref{theo-global1}): \S
\ref{globalss3} $\Rightarrow$ MP
faisceautique}
\end{array} \right.$ \\ \hline
\S \ref{ssce0} & \begin{tabular}{c} $\Longrightarrow$ \\ SSCE+SSS \\ +C+CHT
\end{tabular} & \begin{tabular}{l} poids
minimal de (\ref{theo-ripsi-local})\\ en hauteur $d$
\end{tabular} \\ \hline
\end{tabular}
\end{center}
\end{figure}

\marque En ce qui concerne le rôle des divers calculs de groupes de cohomologie, voici
leur enchaînement
logique. Pour commencer on calcule exclusivement des composantes
$\Pi^{\oo,v}$-isotypiques pour $\Pi$ une
représentation irréductible automorphe et cohomologique de $G_\tau(\Am)$, car dans tous
les autres cas tous
les calculs donnent des résultats nuls. On considère plus particulièrement deux
situations selon que $\Pi_v$
est de la forme $[\overleftarrow{s-1}]_{\pi_v}$ ou $[\overrightarrow{s-1}]_{\pi_v}$ pour
$\pi_v$ une
représentation cuspidale de $GL_g(F_v)$ avec $d=sg$. Pour alléger les notations, on
notera $H^i(j^{\geq
tg}_!)$ pour la composante $\Pi^{\oo,v}$-isotypique du $i$-ème groupe de cohomologie du
complexe $j^{\geq
tg}_! \FC(g,t,\pi_v)[d-tg]$ et on utilisera une notation similaire pour $H^i(j^{\geq
tg}_{!*})$.

\marque La situation est largement plus simple dans le cas où
$\Pi_v=[\overleftarrow{s-1}]_{\pi_v}$. On
commence, proposition (\ref{prop-coho1}), par calculer les $H^i(j^{\geq tg}_{!*})$ et on
trouve que ceux-ci
sont tous nuls si $tg \neq d$. En utilisant (\ref{somme-alternee}), on en déduit alors,
corollaire
(\ref{coro-hij-nul}), le calcul des $H^i(j^{\geq tg}_!)$ qui sont nuls pour $i \neq 0$.
Ces calculs
permettent de compléter, corollaire (\ref{coro-pnul}), aux points supersinguliers, la
détermination des
constituants simples du faisceau pervers $j^{\geq tg}_! \FC(g,t,\pi_v)[d-tg]$ obtenu,
corollaire
(\ref{coro-sec-fp}), sur les autres strates grâce à la description des faisceaux des
cycles évanescents en
fonction des systèmes locaux d'Harris-Taylor ainsi qu'au théorème de comparaison de
Berkovich-Fargues. On en
déduit alors, corollaire (\ref{coro-hic}), le calcul des $\Pi^{\oo,v}$-parties des
groupes de cohomologie
$H^j_c(\bar X_{U^p,m}^{(d-h)},R^i\Psi_v \otimes \LC_{\xi})$. A ce moment ci, on dispose
du théorème
(\ref{theo-global1}) de sorte que l'on peut calculer via la suite spectrale des poids
SSP, corollaire
(\ref{coho-global1}), les $\Pi^{\oo,v}$-parties des groupes de cohomologie
$H^i_{\eta_v,\xi}$ de la fibre
générique. On en déduit alors, proposition (\ref{prop-hipsi}), les $\Pi^{\oo,v}$-parties
des termes
$E_2^{p,q}$ de la suite spectrale des cycles évanescents qui nous permet d'obtenir,
corollaire
(\ref{coro-1}), une première liste de possibilités pour les
$\UC_{F_v,l,d}^{i}(\JL^{-1}(\st_s(\pi_v)^\vee))$.

\marque On considère alors $\Pi_v=[\overrightarrow{s-1}]_{\pi_v}$ et on calcule,
proposition
(\ref{prop-not}), les $H^i(j^{\geq tg}_{!*})$ qui sont non nuls pour $i \equiv s-t \mod
2$ et $|i| \leq s-t$.
De ce calcul et de la connaissance des constituants simples des $j^{\geq tg}_!
\FC(g,t,\pi_v)[d-tg]$, on en
déduit, lemme (\ref{lem-rj-combi}), un contrôle sur les $H^i(j^{\geq tg}_!)$. On étudie
ensuite la suite
spectrale des cycles évanescents. On commence par en déterminer l'aboutissement,
corollaire
(\ref{coro-ssce-min}) et proposition (\ref{prop-lrs2}), en utilisant le calcul des
$H^i(j^{\geq tg}_{!*})$
ainsi que le théorème (\ref{theo-global1}). En utilisant le théorème de Lefschetz
difficile ainsi que le
corollaire (\ref{coro-1}) qui contrôle ce qui se passe au niveau des points
supersinguliers, on en déduit,
proposition (\ref{prop-hic-poids}) et corollaire (\ref{strate-alt-p}), le calcul des
$H^i(j^{\geq tg}_!)$. Le
résultat est qu'en ce qui concerne les parties de poids minimal, $s(g-1)$, tout ce passe
au niveau des points
supersinguliers. L'étude de la suite spectrale des cycles évanescents fournit alors la
partie de poids
$s(g-1)$ du théorème (\ref{theo-ripsi-local}) qui d'après ce que l'on a vu, implique les
théorèmes globaux
(\ref{theo-global2}) et (\ref{theo-ss}) qui d'après le théorème de comparaison de
Berkovich-Fargues,
impliquent le théorème local (\ref{theo-local-fil}), ce qui complète la démonstration.
On résume la
discussion précédente dans le tableau suivant.

\begin{figure}[!h]
\begin{center}
\begin{tabular}{|c|r|c|}
\hline $\Pi_v =\st_s(\pi_v)$ & & $\Pi_v=\speh_s(\pi_v)$ \\ \hline

$H^i(j^{\geq tg}_{!*})$ (\ref{prop-coho1}) & & \\
\begin{tabular}{c|c|} SSP $\Downarrow$ & $\Downarrow$ (\ref{somme-alternee}) \\
$H^i_{\eta_v}$ (\ref{coho-global1}) & $H^i(j^{\geq tg}_!)$ (\ref{coro-hij-nul}) \\
\hline  \end{tabular}
$\Rightarrow$ &
\begin{tabular}{c} $j^{\geq tg}_!=\sum j^{\geq (t+r)g}_{!*}$ \\ (\ref{coro-sec-fp}) \\
+ \\  (\ref{coro-pnul})
\end{tabular} $\Rightarrow$ & $H^i(j^{\geq tg}_{!*})$ (\ref{prop-not}) \\

\begin{tabular}{lr} & $\Downarrow$ \end{tabular} & & \begin{tabular}{lr} $\Downarrow$
(\ref{theo-global1}) &
$\Downarrow {\atop{SSP}{LD}}$ \end{tabular} \\

$H^j(\bar X^{(d-h)},R^i\Psi)~h \neq d$ (\ref{coro-hic}) & & \begin{tabular}{l|l}
contrôle & $H^i_{\eta_v}$ \\
$H^i(j^{\geq tg}_!)$ & (\ref{coro-ssce-min}) \\ (\ref{lem-rj-combi}) & (\ref{prop-lrs2})
\\ \hline \end{tabular} \\

$(\ref{coho-global1}) \Downarrow SSS$ & & \begin{tabular}{c} LD+\textbf{(\ref{coro-1})}
\\ $\Downarrow$ \end{tabular}
\\

$H^j(R^i\Psi)~j \neq 0$ (\ref{prop-hipsi}) & & \begin{tabular}{c} \hline $H^i(j^{\geq
tg}_!)$ \\
(\ref{prop-hic-poids}) (\ref{strate-alt-p}) \\ \hline \end{tabular} \\

$\Downarrow$ SSCE & & \begin{tabular}{c} SSCE \\ $\Downarrow$ \end{tabular}  \\

contrôle $\Psi_{F_v}^{d,i}$ \textbf{(\ref{coro-1})} & \begin{tabular}{c|}
(\ref{theo-global2}) \\ (\ref{theo-ss}) \\
\hline \end{tabular} & \begin{tabular}{ll} $\Longleftarrow$ & $\Psi^{d,i}$
(\ref{theo-ripsi-local}) \end{tabular} \\

& \begin{tabular}{c} BF \\ $\Downarrow$ \end{tabular} & \\

& (\ref{theo-local-fil}) & \\ \hline
\end{tabular}
\end{center}
\end{figure}

\begin{rema} \label{rema-zw}
- Afin de ne pas multiplier les situations, nous n'avons considéré dans nos énoncés que
les faisceaux induits
notés $HT(g,t,\pi_v,\Pi_t)$ à partir des systèmes locaux $\FC(g,t,\pi_v)_1 \otimes
\Pi_t$. Ceux-ci sont munis
d'une action par correspondances de $G(\Am^\oo) \times \Zm$; cependant quand on les
considère comme des
constituants de $R\Psi_v$, on les voit comme munis d'une action par correspondances de
$G(\Am^\oo) \times
W_v$. Le lien entre les deux situations est donnée par $c_v \in W_v \mapsto -v(c_v) \in
\Zm$. Afin de
distinguer les deux situations, l'action de $W_v$ sera toujours accompagnée d'un
$\rec_{F_v}^\vee(\pi_v)$
tandis que pour $\Zm$ il sera toujours question du caractère $\Xi$ (cf. plus loin).

- Par ailleurs nous ne considérons pas par la suite l'action totale de $W_v$ mais plutôt
ses frobenius
semi-simplifiés.
\end{rema}


\part{Preuve des théorèmes}

\section{Étude dans le groupe de Grothendieck $\GF$}
\label{glob0}

\subsection{Groupe de Grothendieck des faisceaux pervers: généralités}

~\\

\marque Soit $X$ un $\Fm_q$-schéma $X$. On rappelle que la catégorie $\perv(X)$ des faisceaux pervers sur $X$
est noethérienne et artinienne, ses objets simples étant de la forme $j_{!*} \LC$ où $\LC$ est un système
local irréductible sur $j: U \hookrightarrow X$, où $U$ est un ouvert d'un fermé de $X$.

Pour $P$ un objet de $\perv(X)$, sa dimension $n$ est par définition la plus grande dimension des supports
$U$ des faisceaux pervers simples constituant $P$ de sorte que $-n= \min \{ i~/~ h^i P \neq 0 \}$ et $h^{-n}
P$ a un support de dimension $n$.

\begin{defi} Étant donnée une catégorie triangulée $A$ localement petite, on considère son groupe de Grothendieck
$K(A)$ défini
comme le groupe libre engendré par les classes d'isomorphismes d'objets de $A$ quotienté par les relations:

- $A[1]=-A$;

- $A=B+C$ pour tout triangle distingué $B \longto A \longto C \longmapright{+1}$.
\end{defi}

\begin{prop} Soit $(\DC^{\leq 0},\DC^{\geq 0})$ une catégorie dérivée munie d'une
$t$-structure non dégénérée: on note $\CC$ son coeur qui
est alors une catégorie abélienne de groupe de Grothendieck $\groth(\CC)$. L'application qui à un objet $\FC$ de
$\DC$
associe
$$\sum_i (-1)^i [\lexp p h^i \FC ] \in \groth(\CC)$$
induit un isomorphisme du groupe de Grothendieck $K(\DC)$ de la catégorie triangulée $\DC$, sur $\groth(\CC)$.
\end{prop}

\begin{proof} Pour tout objet $\GC$ de $\DC$ et pour tout $n \in \Zm$, on a un triangle distingué
$$\tau_{\leq n} \GC \longto \GC \longto \tau_{\geq n+1} \GC \longmapright{+1}$$
de sorte que si $\GC$ est un objet de $\DC^{\geq n}$, on a $[\GC]=(-1)^n [\lexp p h^n \GC]+[\tau_{n+1} \GC]$
car $\tau_{\leq n} \GC=\tau_{\leq n} \tau_{\geq n} \GC=(\lexp p h^n \GC)[n]$. Soit alors $a$ et $b$ tel que
$\FC$ soit un objet de $\DC^{[a,b]}$. En appliquant ce qui précède à $\tau_{\geq n} \FC$ pour $n$ variant de
$a$ à $b$, on obtient l'égalité $[\FC]=\sum_i (-1)^i [\lexp p h^i \FC],$ de sorte que l'application $\sum_i
\alpha_i [P_i] \in \groth(\CC) \mapsto \sum_i \alpha_i [P_i] \in K(\DC)$ est inverse de celle de l'énoncé.
\end{proof}

\begin{coro}
Soit $D_b^c(X,\bar \Qm_l)$ la catégorie dérivée des $\bar \Qm_l$-complexes constructibles sur un $\Fm_q$-schéma $X$
que l'on muni de la $t$-structure
de perversité autoduale (resp. de la $t$-structure triviale) de coeur la catégorie $\perv(X)$ (resp. $\const(X)$) des
faisceaux pervers (resp. des faisceaux
constructibles) sur $X$. Pour tout objet $\FC$ de $D_b^c(X,\bar \Qm_l)$, son image dans $\groth(\perv(X))$ est
déterminée par son image
dans $\groth(\const(X))$.
\end{coro}

\begin{lemm} Soit $j:U \hookrightarrow X$ l'inclusion d'un ouvert $U$ d'un $\Fm_q$-schéma $X$; on note $i:Z
\hookrightarrow X$ le fermé complémentaire.
L'image d'un objet $\FC$ de $D_b^c(X,\bar \Qm_l)$ dans $\groth(\const(X))$ est déterminée par l'image de $j^*
\FC$ dans $\groth(\const(U))$ et celle de $i^* \FC$ dans $\groth(\const(Z))$.
\end{lemm}

\begin{proof} Le résultat découle directement de l'existence du triangle distingué
$j_! j^* \FC \longto \FC \longto i_*i^* \FC \longmapright{+1}$.
\end{proof}

On en déduit alors la proposition suivante.

\begin{prop} \label{prop-perv1}
Soit $X$ un schéma muni d'une stratification
$$X^{[0]} \subset X^{[1]} \subset \cdots \subset X^{[d-1]}=X$$
et soit $\PC$ un faisceau pervers sur $X$. Alors l'image de $\PC$ dans le groupe de Grothendieck des faisceaux
pervers
sur $X$ est
déterminée par l'image de $\sum_i (-1)^i [h^i(\PC)_{|X^{(h)}}$ dans le groupe de Grothendieck des faisceaux
localement
constant sur
$X^{(h)}:=X^{[h]}-X^{[h-1]}$ pour tout $0 \leq h < d$.
\end{prop}

\subsection{Image dans $\GF$ de $R\Psi_v[d-1]$}

\begin{prop}
Les complexes $Rj^{\geq tg}_* \FC(g,t,\pi_v)_0[d-tg]$ et $j^{\geq tg}_! \FC(g,t,\pi_v)_0[d-tg]$ sont des
faisceaux pervers.
\end{prop}

\begin{proof} Le résultat découle du fait que $j^{\geq tg}$ est une inclusion affine.

\end{proof}

\begin{prop} \label{prop-libre}
Pour $\pi_v$ une représentation irréductible cuspidale de $GL_g(F_v)$, on a l'égalité suivante dans $\GF$:
\begin{multline} \label{egalite-groth}
e_{\pi_v}[R\Psi_{\pi_v}[d-1]]= \sum_{i=1}^{s_g} \sum_{t=i}^{s_g} (-1)^{t-i} \\
[j^{\geq tg}_! HT(g,t,\pi_v, [\overleftarrow{i-1},\overrightarrow{t-i}]_{\pi_v}) \otimes
\rec_{F_v}^\vee(\pi_v)(-\frac{tg-2+2i-t}{2})]
\end{multline}
\end{prop}

\begin{proof}
L'isomorphisme (\ref{iso1}) et l'égalité (\ref{inutile}) déterminent pour tout $0 \leq h < d$, la somme
$\sum_i (-1)^i [R^i\Psi_{v,|\bar X_{U^p,m}^{(d-h)}}]$
dans le groupe de Grothendieck de la catégorie abélienne des faisceaux constructibles sur $\bar X_{U^p,m}^{(d-h)}$.
Le
résultat découle alors de la proposition (\ref{prop-perv1}).

\end{proof}

\subsection{Décomposition dans $\GF$ des $j_!^{\geq tg} \FC(g,t,\pi_v)[d-tg]$}
\label{bf1}

Rappelons le théorème de comparaison de Berkovich-Fargues

\begin{theo} Soit $z$ un point géométrique de $\bar X_{U^p,m}^{(d-h)}$. Le germe en $z$ de $R\Psi_v$ est égal
à $\Psi_{F_v,l,h,m_1}^\bullet$ en tant qu'objet de la catégorie dérivée filtrée $\Dm^b F(\bar \Qm_l)$.
\end{theo}

Ainsi l'hypothèse de récurrence sur la filtration de monodromie-locale des modèles de Deligne-Carayol de
hauteur strictement inférieure à $d$, donne le corollaire suivant.

\begin{coro} \label{coro-comparaison-grk}
Pour tout point géométrique $z_{I=U^p(m)}$ de $\bar X_{U^p,m}^{(d-h)}$, le germe en $z_I$ de $h^i
gr_{k,\pi_v,I}$ vérifie les propriétés suivantes:

\begin{itemize}
\item[(1)] il est nul si $h$ n'est pas divisible par $g$;

\item[(2)] pour $h=tg$, ils sont nuls pour $|k| \geq t$ quelque soit $i$;

\item[(3)] pour $h=tg$ et $|k| < t$:

(i) ils sont nuls pour $i < tg-d-t+1+|k|$ ou pour $i \not \equiv tg-d-t+1+k \mod 2$ ou pour $i>tg-d$;

(ii) pour $tg-d-t+1+|k| \leq i=tg-d-t+1+|k|+2r \leq tg-d$, la fibre en $z_I$ de $h^i gr_{k,\pi_v,I}$ est
naturellement munie d'une action de $\NC_v \cap \HC_{v,m_1}$ où $\NC_v \subset GL_h(F_v) \times
D_{v,h}^\times \times W_v$ est défini comme au \S \ref{rapel-DC}, où $GL_h(F_v)$ est vu naturellement comme
un sous-groupe du Levi $GL_h(F_v) \times GL_{d-h}(F_v)$ de $GL_d(F_v)$ et où $\HC_{v,m_1}$ est l'algèbre de
Hecke associée à $\KF_{v,m_1} \subset GL_h(\OC_v)$. On obtient alors que la fibre en $z_I$ de $(h^i
gr_{k,\pi_v,I})^{e_{\pi_v}}$ est isomorphe, en tant que $\NC_v \cap \HC_{v,m_1}$ module, aux invariants sous
$\KF_{v,m_1}$ de
$$\JL^{-1}([\overleftarrow{t-1}]_{\pi_v^\vee}) \otimes [\overleftarrow{|k|+2r}]_{\pi_v} \overrightarrow{\times}
[\overrightarrow{t-2r-|k|-2}]_{\pi_v} \otimes \rec_{F_v}^\vee(\pi_v)(-\frac{t(g-1)+|k|+k+2r}{2}).$$

\end{itemize}
\end{coro}

\rem En ce qui concerne le dernier point du corollaire précédent, et donc le fait que la filtration locale
soit équivariante, on peut la déduire du théorème de comparaison de Berkovich-Fargues avec les strates des
variétés de Shimura en rang $d'<d$ en supposant connus, par récurrence, les théorèmes globaux pour celles-ci.

\begin{prop} \label{prop-p}
Pour tout $1 \leq t \leq s_g$, on a l'égalité dans $\GF$
\begin{multline*}
j^{\geq tg}_! HT(g,t,\pi_v,\Pi_t)[d-tg] =  \sum_{i=t}^{s_g} j^{\geq ig}_{!*} HT(g,i,\pi_v,\Pi_t
\overrightarrow{\times}
[\overleftarrow{i-t-1}]_{\pi_v})[d-ig] \otimes \Xi^{\frac{(t-i)(g-1)}{2}} \oplus P_{!,t}
\end{multline*}
où $P_{!,t}$ est une somme de faisceaux pervers concentrés aux points supersinguliers.
\end{prop}

\begin{proof} Dans $\GF$, on écrit $j^{\geq tg}_! HT(g,t,\pi_v,\Pi_t)[d-tg]$ sous la forme
$\sum_i A_{\pi_v,t,i}(\Pi_t)$ où les $A_{\pi_v,t,i}(\Pi_t)$ sont une somme à coefficients positifs de
faisceaux pervers simples de poids $d-tg-i$, que l'on appellera les constituants de $j^{\geq tg}_!
HT(g,t,\pi_v,\Pi_t)[d-tg]$. On raisonne par récurrence descendante sur la dimension des faisceaux pervers qui
interviennent, en traitant tous les $t$ par récurrence de $s_g$ à $1$ et en travaillant sur les systèmes
inductifs comme indiqué dans la remarque (\ref{rema-lim}).

\begin{lemm} \label{lem-fait0}
Pour tout $1 \leq t \leq s_g$, $j^{\geq tg}_{!*} HT(g,t,\pi_v,\Pi_t)[d-tg]$ est le seul constituant de
dimension $d-tg$ de $j^{\geq tg}_! HT(g,t,\pi_v,\Pi_t)[d-tg]$. Par ailleurs tous les autres constituants
de $j^{\geq tg}_! HT(g,t,\pi_v,\Pi_t)[d-tg]$ sont de poids strictement inférieur à $d-tg$.
\end{lemm}

\begin{proof} On a la suite exacte courte de faisceaux pervers
$$0 \to P_{\pi_v,t,0}(\Pi_t) \longto j^{\geq tg}_! HT(g,t,\pi_v,\Pi_t)[d-tg] \longto
j^{\geq tg}_{!*} HT(g,t,\pi_v,\Pi_t)[d-tg] \to 0$$ où $P_{\pi_v,t,0}$ est égale à $i^{\geq tg,*}
j^{\geq tg}_{!*} HT(g,t,\pi_v,\Pi_t)[d-tg-1]$ qui est donc de poids inférieur ou égal à $d-tg-1$ et
de dimension strictement inférieure à $d-tg$, d'où le résultat.

\end{proof}

\marque D'après le lemme précédent, tous les faisceaux pervers qui interviennent dans l'écriture de $j^{\geq
tg}_! HT(g,t,\pi_v,\Pi_t)[d-tg]$ sont de dimension inférieure ou égale à $d-tg$. Ainsi pour $h>d-g$, il n'y a
aucun constituant de dimension $h$ dans les $j^{\geq tg}_! HT(g,t,\pi_v,\Pi_t)[d-tg]$ pour $1 \leq t \leq
s_g$ et pour $h=d-g$, $j^{\geq g}_{!*} HT(g,1,\pi_v,\Pi_1)[d-g]$ est le seul constituant de dimension $d-g$
pour $t=1$. Supposons donc que pour tout $h > h_0$, les constituants des $j^{\geq tg}_!
HT(g,t,\pi_v,\Pi_t)[d-tg]$ pour $1 \leq t \leq s_g$, sont ceux prédits par l'énoncé et traitons la dimension
$h_0$.

\begin{lemm} \label{lem-hok}
Supposons la proposition (\ref{prop-p}) vérifiée pour les faisceaux pervers de dimension strictement
supérieure à $0<h_0 < d-g$, i.e. pour tout $1 \leq t \leq s_g$, les constituants de dimension strictement
plus grande que $h_0$ des $j^{\geq tg}_! HT(g,t,\pi_v,\Pi_t)[d-tg]$ sont les
$$j^{\geq (t+r)g}_{!*} HT(g,t+r,\pi_v,\Pi_t \overrightarrow{\times} [\overleftarrow{r-1}]_{\pi_v})
[d-tg] \otimes \Xi^{\frac{r(g-1)}{2}}$$ avec $d-(t+r)g>h_0$. Pour tout $k$, l'image dans $\GF$ de $e_{\pi_v}
gr_{k,\pi_v}$ est alors égale à
$$\sum_{\atop{|k| < t < (d-h_0)/g}{t \equiv k-1 \mod 2}} \PC(g,t,\pi_v) (-\frac{tg+k-1}{2}) +P_{k,h_0}$$
où $P_{k,h_0}=(P_{k,h_0,I})_I$ est une somme de faisceaux pervers simples de dimension inférieure ou égale à
$h_0$ à support dans la tour des $(\bar X_{U^p,m}^{[h_0]})_{I=U^p(m)}$, i.e. la fibre de $h^iP_{k,h_0,I}$ en
tout point géométrique de $\bar X_{U^p,m}^{(d-h)}$ est nulle pour $h<d-h_0$ et pour tout $i$.
\end{lemm}

\begin{proof} D'après (\ref{prop-libre}), chacun\footnote{On rappelle, cf. la remarque (\ref{rema-I}), que les
$\PC(g,t,\pi_v)$ ne sont pas irréductibles; cependant les arguments fonctionnent de manière strictement
identique pour tous ses constituants simples car la seule différence entre ceux-ci provient de l'action de
$\DC_{v,tg}^\times$.} des $\PC(g,t,\pi_v) (-\frac{tg+k-1}{2})$, pour $d-tg>h_0$ est un constituant d'un
$e_{\pi_v}gr_{r,\pi_v}$; il s'agit alors de montrer que $r=k$. Le résultat découle directement via le
théorème de comparaison de Berkovich-Fargues, de l'hypothèse de récurrence sur le théorème local
(\ref{theo-local-fil}). En effet, pour $tg<d$ les $\pi_v$-parties des $h^0 gr_{tg,loc,k}$ sont pures de poids
$tg-1+k$. L'équivalent du théorème de Serre-Tate et le théorème de Berkovich-Fargues impliquent alors que la
fibre en tout point géométrique $z$ de la strate $tg$ de $h^{tg-d} gr_{k,\pi_v}$ est pure de poids $tg-1+k$.
Or la fibre en un point géométrique de la strate $tg$ de $h^{tg-d} \PC(g,t,\pi_v)(-\frac{tg+k-1}{2})$ est de
poids $tg-1+k$. Ainsi pour $t=1$, $\PC(g,1,\pi_v)((1-g)/2)$ est un constituant de $gr_{0,\pi_v}$ car tous les
autres faisceaux pervers qui interviennent sont de dimension strictement plus petite.

On raisonne alors par récurrence sur $t$; on suppose que pour tout $t <t_0<(d-h_0)/g$, et tout $|k|<t-1$,
$\PC(g,t,\pi_v)(-\frac{tg+k-1}{2})$ est un constituant de $gr_{k,\pi_v}$ et traitons le cas de $t_0$. D'après
ce qui précède on en déduit que pour tout $|k|<t_0$, $\PC(g,t_0,\pi_v)(-\frac{t_0g-1+k}{2})$ est un
constituant de $gr_{r,\pi_v}$ avec $r \geq k$; en effet dans le cas contraire la filtration par le poids de
$gr_r$ donnerait une suite exacte courte
$$0 \to P \longto gr_{r,\pi_v} \longto \PC(g,t_0,\pi_v)(-\frac{t_0g-1+k}{2}) \to 0$$
modulo des faisceaux pervers de dimension strictement plus petite, ce qui donnerait une suite exacte longue
$$ \cdots h^{t_0g-d} gr_{r,\pi_v} \longto h^{t_0g-d} \PC(g,t_0,\pi_v)(-\frac{t_0g-1+k}{2}) \longto h^{t_0g-d+1} P
\cdots$$ or $h^{tg-d+1} P$ a un support de dimension strictement plus petite que $d-t_0g$ de sorte qu'il
existerait un point géométrique de $\bar X_{U^p,m}^{(d-t_0g)}$ tel que la fibre de $h^{t_0g-d} gr_{r,\pi_v}$
serait de poids $t_0g-1+k$ ce qui n'est pas. Ainsi en particulier pour $k=t_0-1$, il existe $r \geq t_0-1$
tel que $gr_{r,\pi_v}$ soit de dimension $d-t_0g$. Pour la même raison que ci-dessus, on en déduit qu'il
s'agit de $r=t_0-1$. On utilise alors l'opérateur $N$ qui finit de placer les
$\PC(g,t_0,\pi_v)(-\frac{t_0g+k-1}{2})$ sur les $gr_{k,\pi_v}$, d'où la récurrence.

Par ailleurs on remarque qu'en tout point géométrique $z$ de $\bar X_{U^p,m}^{(h)}$, pour $h \geq h_0$, la fibre en
$z$ des
faisceaux de cohomologie de $gr_{k,\pi_v,U^p(m)}-P_{k,h_0,U^p(m)}$ est égale à celle du modèle local i.e., d'après le
théorème de
comparaison de Berkovich-Fargues à celle de $gr_{k,\pi_v,U^p(m)}$ de sorte que celle de $P_{k,h_0,U^p(m)}$
est nulle.

\end{proof}

\marque \textit{Retour à la preuve de la proposition (\ref{prop-p})}: le principe est d'étudier l'égalité
(\ref{egalite-groth}) en utilisant le lemme précédent ainsi que (\ref{lem-fait0}) pour $j^{\geq t_0g}_{!}
HT(g,t_0,\pi_v,\Pi_{t_0})[d-t_0g]$ dont le seul constituant de dimension supérieur ou égal à $d-t_0g$ est
$$j^{\geq t_0g}_{!*} HT(g,t_0,\pi_v,\Pi_{t_0})[d-t_0g].$$

\noindent \textbf{Note:} \textit{on suppose, afin de simplifier la rédaction, que $\Pi_t$ est elliptique de
type $\pi_v$.}

\begin{lemm} \label{lem-fp-rpsi}
Pour $1 \leq t \leq s_g$, soit $A_{\pi_v,t,r}(\Pi_t)$ un constituant de $j^{\geq tg}_!
HT(g,t,\pi_v,\Pi_t)[d-tg]$ de poids $d-tg-r$ et de dimension $h_0$. Il existe alors des entiers $k$ et $a
\geq 0$ ainsi que $\Pi_t'$ une représentation elliptique de type $\pi_v$, et donc de même support cuspidal
que $\Pi_t$, telle que $A_{\pi_v,t,r}(\Pi_t') \otimes \rec_{F_v}^\vee(\pi_v)(-\frac{t(g+1)-2(a+1)}{2})$ soit
un constituant de $e_{\pi_v} gr_{k,\pi_v}$.
\end{lemm}

\begin{proof} Soit $A_{\pi_v,t,r}(\Pi_t)$ un constituant de $j^{\geq tg}_! HT(g,t,\pi_v,\Pi_t)[d-tg]$ de sorte que
$A_{\pi_v,t,r}([\overleftarrow{t-1}]_{\pi_v}) \otimes \rec_{F_v}^\vee(\pi_v)(-\frac{t(g+1)-2}{2})$ apparaît
dans le membre de droite de l'égalité de la proposition (\ref{prop-libre}). On considère alors $\d$ maximal
tel qu'il existe une représentation elliptique $\Pi_t'$ de type $\pi_v$ de $GL_{tg}(F_v)$ tel que
$A_{\pi_v,t,r}(\Pi_t') \otimes \rec_{F_v}^\vee(\pi_v)(-\frac{\d}{2})$ apparaisse dans le membre de droite de
l'égalité de la proposition (\ref{prop-libre}). On remarque alors qu'il y apparaît avec un coefficient
positif et aucun négatif de sorte qu'il reste dans la somme ce qui implique le résultat. En effet s'il
apparaissait avec un coefficient négatif, ce serait comme constituant d'un $j^{\geq t'g}_!
HT(g,t',\pi_v,[\overleftarrow{t'-r-1},\overrightarrow{r}]_{\pi_v}) [d-t'g] \otimes \rec_{F_v}^\vee(\pi_v)
(-\frac{t'(g+1)-2(r+1)}{2})$ pour $r$ impair positif et on note que $j^{\geq t'g}_!
HT(g,t',\pi_v,[\overleftarrow{t'-1}]_{\pi_v})[d-t'g] \otimes \rec_{F_v}^\vee(\pi_v)(-\frac{t'(g+1)-2}{2})$
contiendrait un $A_{\pi_v,t,r}(\Pi_t'') \otimes \rec_{F_v}^\vee(\pi_v)(-\frac{\d+2r}{2})$ pour $\Pi_t''$ une
représentation de même support cuspidal que $\Pi_t'$, contredisant la maximalité de $\d$.

\end{proof}

\marque Supposons dans un premier temps que $h_0$ n'est pas de la forme $d-tg$ et qu'il existe $t,r$ tel que
$A_{\pi_v,t,r}(\Pi_t)$ contienne un faisceau pervers simple de dimension $h_0$. Le lemme précédent implique
alors qu'il existe $k \geq 0$ tel que, avec les notations du lemme (\ref{lem-hok}), $P_{k,h_0}$ ait un
constituant simple, disons $B_{\pi_v,k}$, de dimension $h_0$ et de poids $d-1+\delta$ pour un certain entier
$\delta$. On raisonne alors comme dans la preuve du lemme (\ref{lem-hok}). Étant de dimension inférieure ou
égale à $h_0$, $P_{k,h_0}$ est alors de dimension $h_0$ de sorte que $h^{-h_0} P_{k,h_0}$ est de dimension
$h_0$. En outre comme $h^{-h_0} P_{k,h_0,I=U^p(m)}$, d'après le lemme (\ref{lem-hok}), est à support dans
$\bar X_{U^p,m}^{[h_0]}$, on en déduit que pour tout point générique de $\bar X_{U^p,m}^{(h_0)}$, la fibre en
ce point de $h^{-h_0} P_{k,h_0,I}$ est non nulle de sorte qu'il existe un point géométrique $z_I$ de $\bar
X_{U^p,m}^{(h_0)}$ tel que la fibre en $z_I$ de $h^{-h_0} P_{k,h_0,I}$ soit non nulle. Par ailleurs, les
$\pi_v$-parties des $h^{-h_0}gr_{k,loc}$ du modèle local de hauteur $d-h_0$ sont toutes nulles, on en déduit,
d'après le théorème de comparaison de Berkovich-Fargues, que la fibre en $z$ de $h^{-h_0} gr_{k,\pi_v}$ est
nulle ce qui implique que $\delta \neq k$: en effet sinon $e_{\pi_v} gr_{k,\pi_v}$ est une extension de
$$(\bigoplus_{\atop{|k| < t < (d-h_0)/g}{t \equiv k-1 \mod 2}} \PC(g,t,\pi_v) (-\frac{tg+k-1}{2}))
 \oplus B_{k,\pi_v}$$
par des faisceaux pervers de dimension inférieure ou égale à $h_0$ ce qui implique que $(h^{-h_0}
gr_{k,\pi_v})^{e_{\pi_v}}$ admettrait $h^{-h_0} B_{k,h_0}$ comme facteur direct. Selon le même principe, on
doit même avoir $\delta <k$. Comme on ne dispose pas de la pureté des $gr_k$, pour montrer que $\delta=k$ on
raisonne alors comme suit.

\begin{lemm} \label{lem-dualite}
Soit $\PC$ un constituant de $R\Psi_{\pi_v}[d-1]$. On en déduit alors que $D \PC(1-d)$ est un
constituant de $R\Psi_{\pi_v^\vee}[d-1]$.
\end{lemm}

\begin{proof} Le résultat découle simplement de la compatibilité de $R\Psi_v$ avec la dualité de Verdier,
soit $D R\Psi_{\eta_v}(\bar \Qm_l[d-1]) \simeq R\Psi_{\eta_v}(D \bar \Qm_l[d-1])$ et du fait que la fibre
générique de $X_{U^p,m}$ est lisse de sorte que $D \bar \Qm_l[d-1] \simeq \bar \Qm_l[d-1](1-d)$.

\end{proof}

On considère alors un constituant simple de dimension $h_0$ de poids maximal, disons $d-1+r$ de
$R\Psi_{\pi_v}[d-1]$. C'est alors un constituant de $gr_{k,\pi_v}$ pour $k > r$ d'après ce qui précède. Par
ailleurs en utilisant le lemme (\ref{lem-dualite}), on obtient un constituant simple de dimension $h_0$ de
$gr_{-k,\pi_v^\vee}$ et de poids $d-1-r$, de sorte que l'opérateur de monodromie $N^k$ fournit un constituant
simple de dimension $h_0$ de $gr_{k,\pi_v}$ et de poids $d-1-r+2k > d-1+r$, d'où la contradiction par
maximalité de $r$.

\marque Supposons désormais que $h_0=d-t_0g$ et supposons avoir montré par récurrence que pour tout $t_1 \leq
t \leq s_g$, le seul constituant de dimension $d-t_0g$ de $j^{\geq tg}_! HT(g,t,\pi_v,\Pi_t)[d-tg]$, est
$j^{\geq t_0g}_{!*} HT(g,t_0g,\pi_v,\Pi_t \overrightarrow{\times} [\overleftarrow{t_0-t-1}]_{\pi_v})[d-t_0g]
\otimes \Xi^{\frac{(t_0-t)(g-1)}{2}}$. Le résultat est vérifié pour $t_1=t_0$, supposons le donc vérifié
jusqu'au rang $t_1$ et traitons le cas de $t_1-1$. On étudie alors les faisceaux pervers de dimension
$d-t_0g$ dans le membre de droite de (\ref{egalite-groth}), en particulier ceux de poids $d-t_0$, ce qui
donne $j^{\geq t_0g}_{!*} HT(g,t_0,\pi_v,\Pi)[d-t_0g] \otimes \rec_{F_v}^\vee(\pi_v) (-\frac{t_0(g-1)}{2})$
avec
$$(-1)^{t_0-1}\Pi=[\overrightarrow{t_0-1}]_{\pi_v} - \sum_{i=1}^{t_0-t_1} (-1)^i [\overrightarrow{t_0-i-1}]_{\pi_v}
\overrightarrow{\times} [\overleftarrow{i-1}]_{\pi_v}=(-1)^{t_1-t_0}
[\overrightarrow{t_1-1},\overleftarrow{t_0-t_1}]_{\pi_v}.$$ Or ce dernier ne peut pas apparaître dans le
membre de gauche de (\ref{egalite-groth}). En effet si $t_1$ est pair c'est évident car à gauche il ne peut y
avoir que des coefficients positifs. Sinon de manière générale, on raisonne comme suit. Le lemme
(\ref{lem-dualite}) donnerait que $j^{\geq t_0g}_{!*} HT(g,t_0,\pi_v^\vee,
[\overleftarrow{t_0-t_1},\overrightarrow{t_1-1}]_{\pi_v^\vee})[d-t_0g] \otimes
\rec_{F_v}(\pi_v)(-\frac{t_0(g+1)-2}{2})$ serait un constituant de $e_{\pi_v}[R\Psi_{\pi_v}[d-1]]$ qui, vu le
poids, ne pourrait provenir que de $j^{\geq t_0g}_!
HT(g,t_0,\pi_v^\vee,[\overleftarrow{t_0-1}]_{\pi_v^\vee})[d-t_0g] \otimes \rec_{F_v}(\pi_v)
(-\frac{t_0(g+1)-2}{2}),$ ce qui n'est pas d'après le lemme (\ref{lem-fait0}).

Ainsi $j^{\geq t_0g}_{!*} HT(g,t_0,\pi_v,[\overrightarrow{t_1-1},\overleftarrow{t_0-t_1}]_{\pi_v})[d-t_0g]
\otimes \rec_{F_v}^\vee(\pi_v) (-\frac{t_0(g-1)}{2})$ doit être un constituant d'un $j^{\geq tg}_!
HT(g,t,\pi_v,[\overleftarrow{i},\overrightarrow{t-i-1}]_{\pi_v})[d-tg] \otimes
\rec_{F_v}^\vee(\pi_v)(-\frac{t(g-1)+2i}{2})$ pour $0 \leq i < t$ et $t>t_0$ et $t \equiv t_1-1 \mod 2$. Le
résultat, i.e. $t=t_1-1$ et $r=0$, découle alors des trois lemmes suivants.

\begin{lemm} \label{lem-induite}
Pour un point géométrique $z$ de $\bar X_{U^p,m,M_{d-t_0g}}^{(d-t_0g)}$, les fibres en $z$ des faisceaux de
cohomologies des $j^{\geq tg}_{!*} HT(g,t,\pi_v,\Pi_t)[d-tg] \otimes \rec_{F_v}^\vee(\pi_v)$, pour $t \geq
t_0$, sont, en tant que représentation de $GL_{t_0g}(F_v) \times GL_{d-t_0g}(F_v) \times W_v$ de la forme
$\bigoplus_\chi (\Pi_t(\chi) \times \pi_\chi) \otimes \pi_\chi' \otimes \rec_{F_v}^\vee(\pi_v)(\chi)$ où
$\chi$ décrit les caractères de $\Zm$ et $\pi_\chi$ (resp. $\pi_\chi'$) est une représentation de
$GL_{(t_0-t)g}(F_v)$ (resp. $GL_{d-t_0g}(F_v)$).
\end{lemm}

\begin{proof} C'est évident en utilisant que les strates sont induites, i.e.
$$j^{\geq tg}_{!*} HT(g,t,\pi_v,\Pi_t)= j^{\geq tg}_{1,!*} \FC(g,t,\pi_v)_1 \otimes \Pi_t
\times_{P_{tg,t_0g,d}(F_v)} P_{t_0g,d}(F_v)$$ en tant que $P_{t_0g,d}(F_v) \times \Zm$-module. Les torsions
découlent alors de l'action telle qu'elle est décrite au \S \ref{ht-prop}.

\end{proof}

\begin{lemm} Supposons que pour tout $1 \leq t \leq t_0$, les constituants simples de dimension strictement supérieure
à $d-t_0g$ des $j^{\geq tg}_{!*} HT(g,t,\pi_v,\Pi_t)[d-tg]$ sont ceux prévus par la proposition
(\ref{prop-p}) et supposons qu'il existe $t<t_0$ tel que
$$j^{\geq tg}_! HT(g,t,\pi_v,[\overleftarrow{r},\overrightarrow{t-r-1}]_{\pi_v})[d-tg] \otimes
\rec_{F_v}^\vee(\pi_v)(-\frac{t(g-1)+2r}{2})$$ admette $j^{\geq t_0g}_{!*}
HT(g,t_0,\pi_v,[\overrightarrow{t_1-1},\overleftarrow{1},\overleftrightarrow{t_0-t_1}]_{\pi_v})[d-t_0g]
\otimes \rec_{F_v}^\vee(\pi_v) (-\frac{t_0(g-1)}{2})$ comme constituant. On a alors $t \leq t_1$ et $r=0$.
\end{lemm}

\begin{proof} On considère la filtration par le poids de
$j^{\geq tg}_! HT(g,t,\pi_v,[\overleftarrow{r},\overrightarrow{t-r-1}]_{\pi_v})[d-tg] \otimes
\rec_{F_v}^\vee(\pi_v)(-\frac{t(g-1)+2r}{2})$ et la suite spectrale qui s'en déduit
\begin{multline*}
E_1^{i,j}=h^{i+j} gr_{-i}(t,r) \Rightarrow  h^{i+j} j^{\geq tg}_! HT(g,t,\pi_v,[\overleftarrow{r},
\overrightarrow{t-r-1}]_{\pi_v})[d-tg] \otimes \rec_{F_v}^\vee(\pi_v)(-\frac{t(g-1)+2r}{2})
\end{multline*}
où $gr_k(t,r)$ désigne le gradué de poids $k$ de $j^{\geq tg}_!
HT(g,t,\pi_v,[\overleftarrow{r},\overrightarrow{t-r-1}]_{\pi_v})[d-tg] \otimes
\rec_{F_v}^\vee(\pi_v)(-\frac{t(g-1)+2r}{2})$. D'après le lemme précédent toutes les fibres aux points
géométriques de $\bar X_{U^p,m,M_{d-t_0g}}^{(d-t_0g)}$ des $E_1^{i,j}$ avec $i+j<t_0g-d$ sont, en tant que
$GL_{t_0g}(F_v) \times W_v$-module, de la forme $\bigoplus_\chi [\overleftarrow{r},
\overrightarrow{t-r-1}]_{\pi_v(\chi)} \times \pi_\chi \otimes
\rec_{F_v}^\vee(\pi_v(\chi))(-\frac{t(g-1)+2r}{2})$ alors que celles des faisceaux de cohomologies de
$j^{\geq t_0g}_{!*}
HT(g,t_0,\pi_v,[\overrightarrow{t_1-1},\overleftarrow{1},\overleftrightarrow{t_0-t_1}]_{\pi_v})[d-t_0g]
\otimes \rec_{F_v}^\vee(\pi_v) (-\frac{t_0(g-1)}{2})$ sont de la forme $\bigoplus_\chi
[\overrightarrow{t_1-1},\overleftarrow{1},\overleftrightarrow{t_0-t_1}]_{\pi_v(\chi)} \times \pi'_\chi
\otimes \rec_{F_v}^\vee(\pi_v(\chi))(-\frac{t(g-1)}{2})$. En remarquant que les $E_1^{i,t_0g-d+1-i}$ ont un
support de dimension strictement inférieur à $d-t_0g$ et que les $E_\oo^n$ sont nuls pour $n \neq tg-d$, on
en déduit alors que $r=0$ et $t \leq t_1$.
\end{proof}

\begin{lemm} \label{lem-combi}
Pour tout $1 \leq t < t_0$ et toute représentation $\Pi$ elliptique de type $\pi_v$ de $GL_{(t_0-t)g}(F_v)$
distincte de $[\overleftarrow{t_0-t-1}]_{\pi_v}$, $j^{\geq tg}_! HT(g,t,\pi_v,\Pi_t)[d-tg]$ ne contient pas
$$j^{\geq t_0g}_{!*} HT(g,t_0,\pi_v,\Pi_t \overrightarrow{\times} \Pi)[d-tg] \otimes \Xi^{(t_0-t)(g-1)/2}.$$
\end{lemm}

\begin{proof} Dans le cas contraire considérons $t$ minimal pour cette propriété. En l'appliquant à
$\Pi_t=[\overrightarrow{t-1}]_{\pi_v}$, on en déduit que $j^{\geq t_0g}_{!*}
HT(g,t_0,\pi_v,[\overrightarrow{t-1},\overleftarrow{1},\overleftrightarrow{t_0-t-1})[d-tg] \otimes
\rec_{F_v}^\vee(\pi_v)(-\frac{t_0(g-1)}{2})$ reste dans le membre de droite de (\ref{egalite-groth}), où
$[\overleftrightarrow{t_0-t-1}]_{\pi_v}$ désigne $\Pi$. En effet il y apparaît via $j^{\geq tg}_{!*}
HT(g,t,\pi_v,[\overrightarrow{t-1}]_{\pi_v})[d-tg] \otimes \rec_{F_v}^\vee(\pi_v)(-\frac{t(g-1)}{2})$ et
n'est pas compensé car d'après le lemme précédent ce ne pourrait qu'être pour un $t'<t$ ce qui contredirait
la minimalité de $t$. Si le signe est négatif on obtient directement la contradiction, sinon on argument
comme suit: par application de la dualité de Verdier, on en déduit que $e_{\pi_v} [gr_{k,\pi_v^\vee}]$ pour
un certain $k$, devrait contenir $j^{\geq t_0g}_{!*} HT(g,t_0,\pi_v^\vee,\lceil
[\overrightarrow{t-1}]_{\pi_v} \overrightarrow{\times} \Pi_t \rceil^\vee)[d-t_0g] \otimes \rec_{F_v}(\pi_v)
(-\frac{t_0(g+1)-2}{2})$ qui est de poids $d+t_0-2$. Or tous les constituants de dimension strictement
supérieur à $d-t_0g$ de $R\Psi_{\pi_v^\vee}[d-1]$ sont de poids strictement inférieur à $d+t_0-2$ de sorte
qu'il existerait un point géométrique de $\bar X_{U^p,m}^{(d-t_0g)}$ tel que la fibre de
$h^{t_0g-d}gr_{k,\pi_v^\vee}$ en ce point admettrait un facteur direct de poids $d+t_0-2$ de la forme $\lceil
[\overrightarrow{t-1}]_{\pi_v} \overrightarrow{\times} \Pi_t \rceil^\vee$ ce qui n'est pas d'après le
corollaire (\ref{coro-comparaison-grk}).

\end{proof}

\marque \textit{Fin de la preuve de la proposition (\ref{prop-p})}: pour tout $k$, on pose dans $\GF$:
$$Q_{k,\pi_v,t_0}:=e_{\pi_v} gr_{k,\pi_v}-\sum_{\atop{|k| < t \leq t_0}{t \equiv k-1 \mod 2}} \PC(g,t,\pi_v)
(-\frac{tg+k-1}{2})$$
Le corollaire (\ref{coro-comparaison-grk}) et ce qui précède, prouvent alors que pour
tout point générique $z$ de dimension supérieure ou égale à $d-t_0g$, $\sum_i (-1)^i (h^i Q_{k,\pi_v,t_0})_z$
est nulle. On en déduit donc d'après la proposition (\ref{prop-perv1}), que $Q_{k,\pi_v,t_0}$ est de dimension
strictement inférieure à $d-t_0g$, ce qui conclut la récurrence.

\end{proof}

\begin{defi}
Pour $t \leq t'$ et $\Pi_t$ une représentation de $GL_{tg}(F_v)$, soit $\PC_-(g,t',\pi_v,t,\Pi_t) \in
\FPH(\bar X)$ \footnote{Comme dans la remarque (\ref{rema-I}), ceux-ci ne sont pas simples mais plutôt une
somme directe de $e_{\pi_v}$ faisceaux pervers simples où la différence entre eux provient de l'action de
$\DC_{v,tg}^\times$. Par ailleurs pour ce qui est des poids cf. la remarque (\ref{rema-zw}).} le faisceau
$j^{\geq t'g}_{!*} HT(g,t',\pi_v,\Pi_t \overrightarrow{\times} [\overleftarrow{t'-t-1}]_{\pi_v}) [d-t'g]
\otimes \Xi^{\frac{(t'-t)(g - 1)}{2}}$ pur de poids $d-t'g -2(t'-t)$ et $\PC_+(g,t',\pi_v,t,\Pi_t)=j^{\geq
t'g}_{!*} HT(g,t',\pi_v,\Pi_t \overleftarrow{\times} [\overleftarrow{t'-t-1}]_{\pi_v}) [d-t'g] \otimes
\Xi^{\frac{(t'-t)(g + 1)}{2}}$, pur de poids $d-t'g + 2(t'-t)$.
\end{defi}

La filtration par le poids donne alors le corollaire suivant.

\begin{coro} \label{coro-sec-fp}
\begin{itemize}

\item[(i)] Il existe, pour $0 \leq i \leq s_g-t$, des faisceaux pervers $P_{\pi_v,t,i}(\Pi_t)\in \FPH(\bar X)$
tels que l'on ait, dans $\FPH(\bar X)$, les suites exactes suivantes:
$$0 \to P_{\pi_v,t,0}(\Pi_t) \longto j^{\geq tg}_! HT(g,t,\pi_v,\Pi_t)[d-tg] \longto j^{\geq tg}_{!*}
HT(g,t,\pi_v,\Pi_t)[d-tg] \to 0$$
$$0 \to P_{\pi_v,t,1}(\Pi_t) \longto P_{\pi_v,t,0}(\Pi_t) \longto \PC_-(g,t+1,\pi_v,t,\Pi_t) \oplus
A_{\pi_v,t,0}(\Pi_t)
\to 0$$
$$\cdots$$
$$0 \to P_{\pi_v,t,r}(\Pi_t) \longto P_{\pi_v,t,r-1}(\Pi_t) \longto \PC_-(g,t+r,\pi_v,t,\Pi_t) \oplus
A_{\pi_v,t,r-1}(\Pi_t)
 \to 0$$
$$\cdots$$
$$0 \to P_{\pi_v,t,s_g-t}(\Pi_t) \longto P_{\pi_v,t,s_g-t-1}(\Pi_t) \longto \PC_-(g,s_g,\pi_v,t,\Pi_t) \oplus
A_{\pi_v,t,s_g-t-1}(\Pi_t) \longto 0$$ avec $A_{\pi_v,t,i}(\Pi_t)$, pour $1 \leq i \leq s_g-t$, des faisceaux
pervers concentrés aux points supersinguliers purs de poids $d-tg-i-1$ et $P_{\pi_v,t,s_g-l}(\Pi_t)$ un
faisceau pervers concentré aux points supersinguliers de poids inférieur ou égal à $d-tg-s_g+t$.

\item[(ii)] Dualement pour la dualité de Verdier, le complexe $Rj^{\geq tg}_* HT(g,t,\pi_v,\Pi_t)[d-tg]$
est un objet de $\FPH(\bar X)$ qui s'insère dans les suites exactes courtes
\begin{multline*}
0 \to j^{\geq tg}_{!*} HT(g,t,\pi_v,\Pi_t)[d-tg] \longto Rj^{\geq tg}_* HT(g,t,\pi_v,\Pi_t)[d-tg] \longto
DP_{\pi_v,t,0}(\Pi_t)(tg-d) \to 0
\end{multline*}
\begin{multline*} 0 \to DA_{\pi_v,t,0}(\Pi_t)(tg-d) \oplus \PC_+(g,t+1,\pi_v,t,\Pi_t) \longto \\
DP_{\pi_v,t,0}(\Pi_t)(tg-d) \longto DP_{\pi_v,t,1}(\Pi_t)(tg-d) \to 0
\end{multline*}
$$\cdots$$
\begin{multline*} 0 \to DA_{\pi_v,t,r-1}(\Pi_t)(tg-d) \oplus \PC_+(g,t+r,\pi_v,t,\Pi_t) \longto \\
DP_{\pi_v,t,r-1}(\Pi_t)(tg-d) \longto DP_{\pi_v,t,r}(\Pi_t)(tg-d) \to 0
\end{multline*}
$$\cdots$$
\begin{multline*} 0 \to DA_{\pi_v,t,s_g-t-1}(\Pi_t)(tg-d) \oplus \PC_+(g,s_g,\pi_v,t,\Pi_t) \longto \\
DP_{\pi_v,t,s_g-t-1}(\Pi_t)(tg-d) \longto DP_{\pi_v,t,s_g-t}(\Pi_t)(tg-d) \to 0
\end{multline*}

\end{itemize}
\end{coro}

\rem Hors des points supersinguliers, nous avons montré que les $gr_k$ étaient purs de poids $d-1+k$, ce qui
n'a rien d'impressionnant puisque finalement on l'a déduite de la pureté locale qui nous est donnée d'après
l'hypothèse de récurrence. Nous verrons à la proposition (\ref{prop-mono}) comment la démontrer pour les
faisceaux pervers supportés par les points supersinguliers.

\begin{coro} \label{coro-grk}
Pour tout $|k| < s_g$, l'image de $e_{\pi_v} gr_{k,\pi_v}$ dans $\GF$ est égale à
$$(\sum_{\atop{|k| < t \leq s_g}{t \equiv k-1 \mod 2}} \PC(g,t,\pi_v)(-\frac{tg+k-1}{2})) + P_k$$
où comme ci-dessus, $P_k$ est une somme de faisceaux pervers concentrés aux points supersinguliers. Dans tous
les autres cas $gr_{k,\pi_v}$ est de dimension nulle, concentré aux points supersinguliers.
\end{coro}

\subsection{Étude aux points supersinguliers des $j_!^{\geq tg} \FC(g,t,\pi_v)[d-tg]$}
\label{fp-ponctuel}

Le but de ce paragraphe est de déterminer les faisceaux pervers ponctuels non précisés dans la proposition
(\ref{prop-p}). On rappelle que le raisonnement du paragraphe précédent ne s'appliquait pas au niveau des
points supersinguliers car nous n'y connaissons pas l'aboutissement de (\ref{suite-spectrale}). Une idée
naive est que pour connaître un faisceau ponctuel, on peut commencer par calculer son groupe de cohomologie
$H^0$.

\begin{prop} \label{prop-coho1}
Soit $1 \leq g < d$ ne divisant pas $d$ (resp. $g$ divisant $d=sg$), et soit $\Pi$ une représentation globale
de $G_\tau(\Am)$ telle que $\Pi$ vérifie $\hyp(\xi)$ avec $\Pi_v \simeq \st_s(\pi_v)$ pour $\pi_v$ une
représentation irréductible cuspidale de $GL_g(F_v)$. Pour tout $i$ et $1 \leq t \leq s_g$ (resp. $1 \leq t <
s$), la composante $\Pi^{\oo,v}$-isotypique des groupes de cohomologie des faisceaux pervers $j^{\geq
tg}_{!*} HT_{\xi}(g,t,\pi_v,\Pi_t)[d-tg]$ est nulle, soit avec les notations de (\ref{nota-coho}):
$$H^i(j^{\geq tg}_{!*} HT_{\xi}(g,t,\pi_v,\Pi_t)[d-tg])[\Pi^{\oo,v}]=0.$$
\end{prop}

\begin{proof} On raisonne par récurrence pour $t$ variant de $s-1$ à $1$; l'initialisation de la récurrence se fait
d'elle même dans la preuve de l'induction qui suit. On reprend les suites exactes courtes de faisceaux
pervers du corollaire (\ref{coro-sec-fp}). On montre tout d'abord, par récurrence sur $i$ de $s-t-1$ à $0$,
que pour tout $j \neq 0$, les groupes de cohomologie $H^j(P_{\pi_v,t,i}(\Pi_t)\otimes
\LC_{\xi})[\Pi^{\oo,v}]$ sont nuls et que
\begin{multline*}
H^0(P_{\pi_v,t,i}(\Pi_t)\otimes \LC_{\xi})[\Pi^{\oo,v}] =
\sum_{k=i}^{s_g-t-1} H^0(A_{\pi_v,t,k}(\Pi_t) \otimes \LC_{\xi})[\Pi^{\oo,v}] + H^0(P_{\pi_v,t,s_g-t}(\Pi_t)
\otimes \LC_{\xi})[\Pi^{\oo,v}]
\end{multline*}
Le résultat est clairement vrai pour $i=s-t-1$ car $P_{\pi_v,t,s-t-1}(\Pi_t)$ est un faisceau pervers
ponctuel. Supposons donc le résultat acquis jusqu'au rang $i+1$ et montrons le au rang $i$. La suite exacte
longue de cohomologie associée à
$$0 \to P_{\pi_v,t,i}(\Pi_t) \longto P_{\pi_v,t,i-1}(\Pi_t) \longto \PC_-(g,t+i,\pi_v,t,\Pi_t) \oplus
A_{\pi_v,t,i-1}(\Pi_t) \to 0$$ fournit, pour $j \neq 0$, les isomorphismes
$H^j(P_{\pi_v,t,i}(\Pi_t)\otimes \LC_{\xi})[\Pi^{\oo,v}]
\simeq H^j(P_{\pi_v,t,i-1}(\Pi_t)\otimes \LC_{\xi})[\Pi^{\oo,v}]$ car d'après l'hypothèse de récurrence
portant sur les $t$,
$H^j(\PC(g,t+i,\pi_v,t,\Pi_t)\otimes \LC_{\xi})[\Pi^{\oo,v}]$
est nul ce qui implique la nullité de $H^j(\PC_-(g,t+i,\pi_v,t,\Pi_t)\otimes \LC_{\xi})[\Pi^{\oo,v}]$.
Par ailleurs pour $j=0$, on a la suite exacte courte
\begin{multline*}
0 \to H^0(P_{\pi_v,t,i}(\Pi_t) \otimes \LC_{\xi})[\Pi^{\oo,v}] \longto
H^0(P_{\pi_v,t,i-1}(\Pi_t) \otimes \LC_{\xi})
\longto H^0(A_{\pi_v,t,i-1}(\Pi_t) \otimes \LC_{\xi}) \to 0
\end{multline*}
d'où le résultat.

On considère alors la suite exacte longue de cohomologie associée à
$$0 \to P_{\pi_v,t,0}(\Pi_t) \longto j^{\geq tg}_! HT(g,t,\pi_v,\Pi_t)[d-tg] \longto
j^{\geq tg}_{!*} HT(g,t,\pi_v,\Pi_t)[d-tg] \to 0$$ qui s'écrit
$$0 \to H^{-1}(j^{\geq tg}_!) \longto H^{-1}(j^{\geq tg}_{!*}) \longto H^0(P_{\pi_v,t,s-t-1})
\longto H^0(j^{\geq tg}_!) \longto H^0(j^{\geq tg}_{!*}) \to 0$$ et pour tout $i \neq -1,0$, $H^i(j^{\geq
tg}_!) \simeq H^i(j^{\geq tg}_{!*})$ où pour alléger les notations, on a omis d'écrire
$HT_{\xi}(g,t,\pi_v,\Pi_t)[d-tg]$ ainsi que $[\Pi^{\oo,v}]$. On en déduit alors que pour $i\neq 0$,
$H^i(j^{\geq tg}_!)$ est pur de poids $d-1+i$ alors que $H^0(j^{\geq tg}_!)$ est mixte de poids inférieur ou
égal à $d-1$. Le calcul de la somme alternée $\sum_i (-1)^i H^i(j^{\geq tg}_!)$, laquelle est nulle pour $g
\nmid d$ et pour $g | d$ est constituée d'un seul terme de poids $d-1-(s-t)$, implique alors la nullité des
$H^i(j^{\geq tg}_!)$ pour $i >0$ et celle des $H^i(j^{\geq tg}_{!*}
HT_{\xi}(g,t,\pi_v,\Pi_t)[d-tg])[\Pi^{\oo,v}]$ pour $i \geq 0$ et pour tout $\pi_v$. La dualité de Poincaré
donne alors la nullité des $H^i(j^{\geq tg}_{!*}
HT_{\xi^\vee}(g,t,\pi_v^\vee,\Pi_t^\vee)[d-tg])[\Pi^{\oo,v}]$ pour tout $\pi_v$ et tout $i \leq 0$ et donc
finalement la nullité des $H^i(j^{\geq tg}_{!*}HT_{\xi}(g,t,\pi_v,\Pi_t)[d-tg])[\Pi^{\oo,v}]$ pour tout $i$
et tout $\pi_v$.
\end{proof}

\marque \rem Il est possible de faire des calculs strictement similaires pour $\pi_v$ quelconque. Par
exemple, au lemme (\ref{prop-mono}) on traite le cas de
$\pi_v \simeq \st_{n_1}(\xi_1) \boxplus \cdots \boxplus \st_{n_r}(\xi_r),$
le cas général étant traité dans la preuve du théorème (\ref{theo-mono-glob}).

\marque D'après le corollaire (\ref{coro-sec-fp}) et en remarquant qu'un faisceau pervers ponctuel n'a de la
cohomologie qu'en degré zéro, les corollaires suivants découlent directement de (\ref{somme-alternee}).

\begin{coro} \label{coro-hij-nul}
Pour $g$ divisant $d=sg$ (resp. $1 \leq g < d$ ne divisant pas $d$), soit $\Pi$ une représentation
irréductible de $G(\Am)$ vérifiant $\hyp(\xi)$ telle que $\Pi_v \simeq \st_s(\pi_v)$. On a alors:

(i) pour tout $i \neq d-tg$ (resp. pour tout $i$), $H^i_c(\bar X_{U^p,m,M_{d-tg}}^{(d-tg)},\FC(g,t,\pi_v)_1
\otimes \LC_{\xi})[\Pi^{\oo,v}]$ et donc $H^i_c(\bar X_{U^p,m}^{(d-tg)}, \FC(g,t,\pi_v,I)\otimes
\LC_{\xi})[\Pi^{\oo,v}]$ sont nuls et pour $i=d-tg$
\begin{multline*}
\lim_{\atop{\to}{U^p(m)}} H^{d-tg}_c(\bar X_{U^p,m,M_{d-tg}}^{(d-tg)},\FC(g,t,\pi_v)_1 \otimes
\LC_{\xi})[\Pi^{\oo,v}]= \\
\sharp \ker^1(\Qm,G_\tau) m(\Pi) [\overleftarrow{s-t-1}]_{\pi_v(-t(g-1)/2)} \otimes
(\Xi^{\frac{(s-t)(g-1)}{2}} \bigoplus_{\chi \in \AF(\pi_v)} \chi^{-1})
\end{multline*}
\begin{multline*} \lim_{\atop{\longto}{U^p(m)}} H^{d-tg}_c(\bar X_{U^p,m}^{(d-tg)},HT_{\xi}(g,t,\pi_v,\Pi_t))
[\Pi^{\oo,v}] \\ =\sharp \ker^1(\Qm,G_\tau) m(\Pi) (\Pi_t \overrightarrow{\times} [\overleftarrow{s-t-1}]_{\pi_v})
\otimes
(\Xi^{(s-t)(g-1)/2} \bigoplus_{\chi \in \AF(\pi_v)} \chi^{-1})
\end{multline*}
en tant que représentation de $GL_{(s-t)g}(F_v) \times \Zm$ et de $GL_d(F_v) \times \Zm$, où
$\AF(\pi_v)$ est l'ensemble des caractères $\chi:\Zm \longto \Qm_l^\times$, tels que $\pi_v \otimes \chi \circ
\val(\det) \simeq \pi_v$ et $m(\Pi)$ est la multiplicité de $\Pi$ dans l'espace des formes automorphes.

(ii) pour tout $1 \leq t \leq s$, on a
\begin{multline*}
H^0(j^{\geq tg}_! HT(g,t,\pi_v,\Pi_t)[d-tg] \otimes \rec_{F_v}^\vee(\pi_v)) \\ =\sharp \ker^1(\Qm,G_\tau) m(\Pi)
e_{\pi_v}(\Pi_t
\overrightarrow{\times} [\overleftarrow{s-t-1}]_{\pi_v}) \otimes \rec_{F_v}^\vee(\pi_v) |-|^{-(s-t)(g-1)/2}
\end{multline*}

\end{coro}

\marque On rappelle que $H_0$ est une forme intérieure de $G_\tau$ telle que $H_0(\Am^{\oo,v}) \simeq
G_\tau(\Am^{\oo,v})$,
$H_0(\Rm)$ est compact et $H_{0,v} \simeq D_{v,d}^\times$.

\noindent \textbf{Hypothèse}: pour la preuve de l'énoncé suivant on considère une représentation automorphe
cohomologique et irréductible $\Pi$ de $G_\tau(\Am)$ telle que:
\begin{itemize}
\item $\Pi_v \simeq \st_s(\pi_v)$ pour $\pi_v$ une représentation irréductible cuspidale de $GL_g(F_v)$ avec $d=sg$;

\item l'ensemble $\AF_{H_0,\xi}(\Pi)$ des représentations irréductibles automorphes $\bar \Pi$ de $H_0(\Am)$ telle
que
$(\bar \Pi^{\oo,v})^\vee \simeq \Pi^{\oo,v}$, est réduit à un élément avec $m(\bar \Pi)=m(\Pi)$, $\bar \Pi_v^\vee
\simeq
\JL^{-1}(\Pi_v)$.
\end{itemize}

\rem L'existence d'une telle représentation $\Pi$ est assurée par Harris-Labesse (cf. \cite{H-L} théorème
3.2.1).

\begin{coro} \label{coro-pnul}
Pour tout $1 \leq t \leq s_g$, les faisceaux pervers $A_{\pi_v,t,i}(\Pi_t)$ du corollaire (\ref{coro-sec-fp})
sont tous nuls et pour $g$ ne divisant pas $d$ (resp. $g|d=sg$) $P_{\pi_v,t,s_g-t}(\Pi_t)$ est nul (resp.
$P_{\pi_v,t,s-t-1}(\Pi_t)$ est le faisceau pervers ponctuel de support l'ensemble des points supersinguliers
tel que $P_{\pi_v,t,s-t-1}(\Pi_t)$ est isomorphe à $\FC(g,s,\pi_v)(-\frac{(s-t)(g-1)}{2}) \otimes (\Pi_t
\overrightarrow{\times} [\overleftarrow{s-t-1}]_{\pi_v})$).
\end{coro}

\begin{proof} D'après la proposition précédente on a
\begin{multline*}
\sum_{i=0}^{s_g-t-1} \lim_{\atop{\to}{U^p(m)}} H^0(\bar X_{U^p,m},A_{\pi_v,t,i}(\Pi_t) \otimes
\LC_{\xi})[\Pi^{\oo,v}]
+
 \lim_{\atop{\to}{U^p(m)}} H^0(\bar X_{U^p,m},P_{\pi_v,t,s_g-t}(\Pi_t)\otimes \LC_{\xi})[\Pi^{\oo,v}] \\
=\lim_{\atop{\to}{U^p(m)}} H^{d-tg}(\bar X_{U^p,m},j^{\geq tg}_! HT_{\xi}(g,t,\pi_v,\Pi_t,I))[\Pi^{\oo,v}].
\end{multline*}
Le membre de droite est d'après (\ref{somme-alternee}) soit nul, pour $g$ ne divisant pas $d$, soit de poids
$d-tg-s_g+t$ pour $g$ divisant $d$.

- Les faisceaux pervers $A_{\pi_v,t,i}(\Pi_t)$, pour $0 \leq i < s_g-t$ étant de dimension zéro et purs de
poids $d-tg-i$, la nullité de $H^0(\bar X_{U^p,m,M_{d-tg}}^{[d-tg]},A_{\pi_v,t,i}(\Pi_t)\otimes
\LC_{\xi})[\Pi^{\oo,v}]$ implique celle des $A_{\pi_v,t,i}(\Pi_t)$.

- De même pour $g$ ne divisant pas $d$, $P_{\pi_v,t,s-t}$, faisceau pervers concentré aux points
supersinguliers, est nul car son groupe de cohomologie en degré $0$ l'est.

- Pour $g$ divisant $d$, on rappelle, cf. le corollaire (\ref{coro-hij-nul}), que
$H^i_c(\bar X_{U^p,m,M_{d-tg}}^{(d-tg)},\FC(g,t,\pi_v)_1 \otimes \LC_{\xi})[\Pi^{\oo,v}]$
est nul pour $i \neq d-tg$, et
\begin{multline*}
\lim_{\atop{\to}{U^p(m)}} H^{d-tg}_c(\bar X_{U^p,m,M_{d-tg}}^{(d-tg)},\FC(g,t,\pi_v)_1 \otimes
\LC_{\xi}))[\Pi^{\oo,v}]= \\
\sharp \ker^1(\Qm,G_\tau) m(\Pi) \bigoplus_{\chi \in \AF(\pi_v)} [\overleftarrow{s-t-1}]_{\pi_v(-t(g-1)/2)} \otimes
\Xi^{\frac{(s-t)(g-1)}{2}} \chi
\end{multline*}
où $\AF(\pi_v)$ est l'ensemble des caractères $\chi:\Zm \longto \Qm_l^\times$, tels que $\pi_v \otimes \chi
\circ \val(\det) \simeq \pi_v$. Ainsi pour toute représentation elliptique $\Pi_t$ de type $\pi_v$ de
$GL_{tg}(F_v)$,
\begin{multline*}
\lim_{\atop{\to}{U^p(m)}} H^0(\bar X_{U^p,m}^{[d-tg]},j^{\geq tg}_! HT_{\xi}(g,t,\pi_v,\Pi_t)[d-tg])
[\Pi^{\oo,v}]=
\\ \sharp \ker^1(\Qm,G_\tau) m(\Pi) e_{\pi_v} \Pi_t \overrightarrow{\times} [\overleftarrow{s-t-1}]_{\pi_v}
\otimes \Xi^{\frac{(s-t)(g-1)}{2}}
\end{multline*}
en tant que représentation de $GL_d(F_v) \times \Zm$. On étudie alors comme précédemment l'égalité
(\ref{egalite-groth}) et tout particulièrement les faisceaux pervers de dimension nulle. En ce qui concerne
les faisceaux pervers simples de poids $s(g+1)-4$, le membre de droite de (\ref{egalite-groth}) est égal à
\begin{multline*}
(P_{\pi_v,s-1,1}([\overleftarrow{s-2}]_{\pi_v})(-\frac{(s-1)(g+1)-2}{2}) -
\FC(g,s,\pi_v) \otimes [\overleftarrow{s-2},\overrightarrow{1}]_{\pi_v}) \otimes \rec_{F_v}^\vee(\pi_v)
(-\frac{s(g+1)-4}{2})
\end{multline*}
de sorte que $P_{\pi_v,s-1,1}([\overleftarrow{s-2}]_{\pi_v})$ contient $\FC(g,s,\pi_v)(-\frac{g-1}{2})
\otimes [\overleftarrow{s-2},\overrightarrow{1}]_{\pi_v}$ et donc, vu que les strates sont induites,
$\FC(g,s,\pi_v) (-\frac{g-1}{2}) \otimes [\overleftarrow{s-2}]_{\pi_v} \overrightarrow{\times}
[\overleftarrow{0}]_{\pi_v}$. Par ailleurs d'après (\ref{h0-ss}), on a le lemme suivant

\begin{lemm} \label{lem-pts-ss}
Pour toute représentation $\Pi_s$ de $GL_d(F_v)$, on a
$$\lim_{\atop{\to}{U^p(m)}} H^0(\bar X_{U^p,m}^{(0)},\FC(g,s,\pi_v) \otimes \LC_{\xi} \otimes \Pi_{s})=\sharp
\ker^1(\Qm,G_\tau)
e_{\pi_v} \bigl ( \CC^{\oo}_{H_0,\xi} [\JL^{-1}(\overleftarrow{(s)}_{\pi_v})] \bigr )^\vee \otimes \Pi_{s}$$ de sorte
que sa
partie $\bar \Pi^{\oo,v}$-isotypique est $e_{\pi_v} m(\bar \Pi) \Pi_{s}$.
\end{lemm}

\marque \rem Sans utiliser \cite{H-L}, on arriverait au fait suivant: pour toute représentation irréductible
automorphe $\bar \Pi$ de $H_0(\Am)$, cohomologique pour $(\xi')^\vee$ de multiplicité $m(\bar \Pi)$, il
existe une unique représentation irréductible automorphe $\Pi$ de $G_\tau(\Am)$ telle que $\Pi^{\oo,v} \simeq
(\bar \Pi^{\oo,v})^\vee$ et $\Pi_v \simeq \JL(\bar \Pi_v^\vee)$, avec $m(\Pi) \geq m(\bar \Pi)$. On renvoie
au \S\ref{corres-jl} pour des correspondances de Jacquet-Langlands générales.

\marque En appliquant ce lemme à $\Pi_s=\Pi_{s-1} \overrightarrow{\times} [\overleftarrow{0}]_{\pi_v}$ et en
utilisant les propriétés imposées à $\Pi$, on obtient que la partie $\bar \Pi^{\oo,v}$-isotypique de
${\DS \lim_{\atop{\to}{U^p(m)}}} H^0(\bar X_{U^p,m}^{(0)},\FC(g,s,\pi_v) \otimes \Pi_{s-1} \overrightarrow{\times}
[\overleftarrow{0}]_{\pi_v})$
est égale à la partie $\Pi^{\oo,v}$-isotypique de
${\DS \lim_{\atop{\to}{U^p(m)}}} H^0_c(\bar X_{U^p,m}^{(g)},HT(g,s-1,\pi_v,\Pi_{s-1})[g]$
qui d'après ce qui précède est égale à celle de
${\DS \lim_{\atop{\to}{U^p(m)}}} H^0(\bar X_{U^p,m},P_{\pi_v,s-1,1}(\Pi_{s-1})).$
On en déduit donc que $P_{\pi_v,s-1,1}(\Pi_{s-1})$ est égal à $\FC(g,s,\pi_v)(-\frac{g-1}{2}) \otimes \Pi_{s-1}
\overrightarrow{\times} [\overleftarrow{0}]_{\pi_v}$.

On raisonne alors par récurrence sur $t$ de $s-1$ à $1$, en supposant que pour tout $t > t_0$,
$P_{\pi_v,t,s-t-1}(\Pi_t)$ est isomorphe à $\FC(g,s,\pi_v)(-\frac{(s-t)(g-1)}{2}) \otimes \Pi_t
\overrightarrow{\times} [\overleftarrow{s-t-1}]_{\pi_v}$. On regarde alors les faisceaux pervers de dimension
nulle et de poids $s(g+1)-2(s-t_0+1)$ dans l'égalité (\ref{egalite-groth}), soit
\begin{multline*}
P_{\pi_v,t_0,s-t_0-1} ( [\overleftarrow{t_0-1}]_{\pi_v}) \otimes (-\frac{t_0(g+1)-2}{2})+ \FC(g,s,\pi_v) \otimes \\
(\sum_{i=1}^{s-t_0} (-1)^i [\overleftarrow{t_0-1},\overrightarrow{i}]_{\pi_v} \overrightarrow{\times}
[\overleftarrow{s-t_0-i-2}]_{\pi_v}) \otimes \rec_{F_v}^\vee(\pi_v) (-\frac{s(g+1)-2(s-t_0+1)}{2})
\end{multline*}
ce qui donne
\begin{multline*}
P_{\pi_v,t_0,s-t_0-1}([\overleftarrow{t_0-1}]_{\pi_v}) -  \FC(g,s,\pi_v) \otimes
[\overleftarrow{t_0-1},\overrightarrow{1},\overleftarrow{s-t_0-1}]_{\pi_v} \otimes \rec_{F_v}^\vee(\pi_v)
(-\frac{(s-t_0)(g-1)}{2})
\end{multline*}
positif. On en déduit alors que $\FC(g,s,\pi_v) (-\frac{(s-t_0)(g-1)}{2}) \otimes
[\overleftarrow{t_0-1},\overrightarrow{1},\overleftarrow{s-t_0-1}]_{\pi_v}$
est contenu dans $P_{\pi_v,t_0,s-t_0-1}([\overleftarrow{t_0-1}]_{\pi_v})$
et donc vu le caractère
induit des strates, que $P_{\pi_v,t_0,s-t_0-1}(\Pi_{t_0})$ contient $\FC(g,s,\pi_v)(-\frac{(s-t_0)(g-1)}{2})
\otimes \Pi_{t_0} \overrightarrow{\times} [\overleftarrow{s-t_0-1}]_{\pi_v}$. L'égalité provient alors de
l'égalité des $H^0$ comme dans le cas $t_0=s-1$, traité ci-avant, d'où le résultat.

\end{proof}

\begin{coro} \label{global0-ok}
Le théorème (\ref{theo-global0}) est vrai.
\end{coro}

\begin{proof} En effet le résultat découle directement de la proposition (\ref{prop-libre}) en réinjectant les
égalités de la proposition (\ref{prop-p}) où les $P_{!,t}$ sont nuls d'après le corollaire (\ref{coro-pnul}).
On obtient ainsi pour tout $1 \leq i \leq s_g$:

\begin{multline*}
\sum_{t=i}^{s_g} (-1)^{t-i} j^{\geq tg}_!
HT(g,t,\pi_v,[\overleftarrow{i-1},\overrightarrow{t-i}]_{\pi_v})[d-tg] \otimes \rec_{F_v}^\vee(\pi_v)
(-\frac{tg-2+2i-t}{2}) = \sum_{l=i}^{s_g} (-1)^{t-i} \sum_{r=0}^{s_g-t} \\ j^{\geq (t+r)g}_{!*}
HT(g,t+r,\pi_v,[\overleftarrow{i-1},\overrightarrow{t-i}]_{\pi_v} \overrightarrow{\times}
[\overleftarrow{r-1}]_{\pi_v})[d-(t+r)g]  \otimes \rec_{F_v}^\vee(\pi_v) (-\frac{(t+r)(g-1)+2(i-1)}{2})
\end{multline*}
laquelle somme est alors égale à
$$\sum_{t=i}^{s_g} j^{\geq tg}_{!*} HT(g,t,\pi_v,\Pi_{t,i})[d-tg] \otimes \rec_{F_v}^\vee(\pi_v)
(-\frac{t(g-1)+2(i-1)}{2})$$ où $\Pi_{t,i}:=\sum_{r=i}^t (-1)^{r-i}
[\overleftarrow{i-1},\overrightarrow{r-i}]_{\pi_v} \overrightarrow{\times} [\overleftarrow{t-r-1}]_{\pi_v}=
[\overleftarrow{t-1}]_{\pi_v}$ avec
$$\begin{array}{rl} e_{\pi_v}[R\Psi_v[d-1]] = & \sum_{i=1}^{s_g} \sum_{t=i}^{s_g}
j^{\geq tg}_{!*}
HT(g,t,\pi_v,\Pi_{t,i})[d-tg] \\ & \hfill \otimes \rec_{F_v}^\vee(\pi_v) (-\frac{t(g-1)+2(i-1)}{2}) \\
= & \sum_{t=1}^{s_g} j^{\geq tg}_{!*} HT(g,t,\pi_v,[\overleftarrow{t-1}]_{\pi_v}) \\ & \hfill \otimes
\rec_{F_v}^\vee(\pi_v)
(-\frac{tg-1}{2}) (\sum_{\atop{|k| < t}{k \equiv t-1 \mod 2}} (-\frac{k}{2}) \\
 = & \sum_{k=1-s}^{s-1} \sum_{\atop{|k| < t \leq s_g}{t \equiv k-1 \mod 2}} j^{\geq tg}_{!*}
HT(g,t,\pi_v,[\overleftarrow{t-1}]_{\pi_v})[d-tg] \\
& \hfill \otimes \rec_{F_v}^\vee(\pi_v) (-\frac{tg+k-1}{2})
\end{array}$$
d'où le résultat.

\end{proof}

\begin{coro} \label{coro-hic}
Pour tout $g|d=sg$, et $\pi_v$ une représentation cuspidale de $GL_g(F_v)$, on considère les groupes de
cohomologie $H^j_c(\bar X_{U^p,m}^{(d-h)},R^i\Psi_{\pi_v}\otimes \LC_{\xi})$ ainsi que leur limite inductive
sur tous les idéaux $I$ de $A$, limite que l'on notera $H^j_{=h,i,\pi_v,\xi}$. Pour $\Pi$ une représentation
automorphe de $G(\Am)$ vérifiant $\hyp(\xi)$, on a les résultats suivant:

\begin{itemize}
    \item[(1)] pour $g$ ne divisant pas $d$, les $H^j_{=h,i,\pi_v,\xi}[\Pi^{\oo,v}]$ sont nuls pour tous $j,h,i$;

    \item[(2)] pour $g$ divisant $d=sg$ et $\Pi_v \simeq [\overleftarrow{s-1}]_{\pi_v}$, les
    $H^j_{=h,i,\pi_v,\xi}[\Pi^{\oo,v}]$ sont nuls pour $h$ qui n'est pas de la forme
    $tg$ avec $1 \leq t \leq s$;

    \item[(3)] pour $g$ divisant $d=sg$ et $\Pi_v \simeq [\overleftarrow{s-1}]_{\pi_v}$, les
    $H^j_{=tg,i,\pi_v,\xi}[\Pi^{\oo,v}]$ sont nuls pour $j \neq d-tg$;

    \item[(4)] pour $g$ divisant $d=sg$ et $\Pi_v \simeq [\overleftarrow{s-1}]_{\pi_v}$, les
    $H^{d-tg}_{=tg,i,\pi_v,\xi}[\Pi^{\oo,v}]$ sont nuls pour $i$ ne vérifiant pas
    $t(g-1) \leq i \leq tg-1$. Si $1 \leq t < s$ et $i=tg-r$ avec $1 \leq r \leq t$,
    $H^{d-tg}_{=tg,tg-r,\pi_v,\xi}[\Pi^{\oo,v}]$ est isomorphe à
    $$\sharp \ker^1(\Qm,G_\tau) m(\Pi)
    ([\overleftarrow{t-r},\overrightarrow{r-1}]_{\pi_v} \overrightarrow{\times} [\overleftarrow{s-t-1}]_{\pi_v})
    \otimes \rec_{F_v}^\vee(\pi_v) (-\frac{s(g-1)-2(r-t)}{2})$$ en tant que représentation de $GL_d(F_v) \times W_v$,
où $m(\Pi)$
    est la multiplicité de $\Pi$ dans l'espace des formes automorphes.

\end{itemize}
\end{coro}

\begin{proof} On rappelle (cf. la proposition (\ref{prop-hic})) que pour tout $i,j$, on a un isomorphisme canonique
$$H^j_c(\bar X_{U^p,m,M_{d-h}}^{(d-h)}, R^i\Psi_v \otimes \LC_{\xi})^h \simeq \bigoplus_{\tau_v \in \CF_h}
(H^j_c(\bar X_{U^p,m,M_{d-h}}^{(d-h)},\FC_{\tau_v} \otimes \LC_{\xi}) \otimes
\UC_{F_v,l,h,m}^{i}(\tau_v))^{h/e_{\tau_v}}$$ de sorte que les résultats découle directement du corollaire
(\ref{coro-hij-nul}).

\end{proof}

\section{Faisceaux de cohomologie des $j^{\geq tg}_{!*} \FC(g,t,\pi_v)$}
\label{glob2}

Rappelons, cf \S \ref{schema}, que d'après le théorème de comparaison de Berkovich-Fargues, la filtration de
monodromie-locale du complexe $\Psi_{F_v,l,d}^\bullet$ sera donnée par les germes en un point supersinguliers
des faisceaux de cohomologie des $gr_k$. On propose alors de calculer tous les faisceaux de cohomologie des
$gr_k$. D'après la proposition (\ref{prop-p}) précisée par le corollaire (\ref{coro-pnul}), il nous suffit de
déterminer les faisceaux de cohomologie des $j^{\geq tg}_{!*} \FC(g,t,\pi_v)[d-tg]$. Nous verrons que cela ne
pose aucun problème en dehors des points supersinguliers car on peut utiliser l'hypothèse de récurrence dans
le théorème (\ref{theo-local-fil}). Au niveau des points supersinguliers on ne dispose d'aucun renseignement
local, cependant la proposition (\ref{prop-p}) nous permet de restreindre efficacement les possibilités pour
les germes en un point supersingulier de ces faisceaux de cohomologie.

\subsection{Une écriture dans $\GF$ de $j^{\geq tg}_{!*} \FC(g,t,\pi_v)[d-tg]$}

\begin{lemm} \label{lem-hij0}
Pour tout $1 \leq t \leq s_g$, on a l'égalité dans $\GF$:
\begin{multline*}
\PC(g,t,\pi_v)= \sum_{r=0}^{s_g-t} (-1)^r j^{\geq (t+r)g}_! HT(g,t+r,\pi_v,[\overleftarrow{t-1}]_{\pi_v}
\overrightarrow{\times} [\overrightarrow{r-1}]_{\pi_v})[d-(t+r)g] \\
\otimes \rec_{F_v}^\vee(\pi_v) (-\frac{r(g-1)}{2})
\end{multline*}
\end{lemm}

\begin{proof} On démontre le résultat par récurrence sur $t$ de $s_g$ à $1$, en utilisant les corollaires
(\ref{coro-sec-fp}) et (\ref{coro-pnul}). Le cas $t=s_g$ est directement donné par loc. cit. Supposons le
résultat acquis jusqu'au rang $t+1$ et traitons le cas de $t$. D'après loc. cit. on a
\begin{multline} \label{egalite-hij}
j^{\geq tg}_! HT(g,t,\pi_v,[\overleftarrow{t-1}]_{\pi_v})[d-tg] \otimes \rec_{F_v}^\vee(\pi_v) =
\PC(g,t,\pi_v)
\\ + \sum_{i=1}^{s-t} j^{\geq (t+i)g}_{!*} HT(g,t+i,\pi_v,[\overleftarrow{t-1}]_{\pi_v}
\overrightarrow{\times} [\overleftarrow{i-1}]_{\pi_v})[d-(t+i)g] \otimes \rec_{F_v}^\vee(\pi_v)
(-\frac{ig-i}{2})
\end{multline}

D'après l'hypothèse de récurrence on a
\begin{multline*}
j^{\geq (t+i)g}_{!*} HT(g,t+i,\pi_v,[\overleftarrow{t-1}]_{\pi_v} \overrightarrow{\times}
[\overleftarrow{i-1}]_{\pi_v})[d-(t+i)g]
\otimes \rec_{F_v}^\vee(\pi_v) (-\frac{ig-i}{2}) = \\
\sum_{r=0}^{s-t-i} (-1)^r j^{\geq (t+i+r)g}_! HT(g,t+i+r,\pi_v,[\overleftarrow{t-1}]_{\pi_v}
\overrightarrow{\times} [\overleftarrow{i-1}]_{\pi_v} \overrightarrow{\times}
[\overrightarrow{r-1}]_{\pi_v})[d-(t+i+r)g] \\ \otimes \rec_{F_v}^\vee(\pi_v) (-\frac{(i+r)g-i-r}{2})
\end{multline*}
de sorte que (\ref{egalite-hij}) s'écrit
\begin{multline*}
\PC(g,t,\pi_v)-j^{\geq tg}_! HT(g,t,\pi_v,[\overleftarrow{t-1}]_{\pi_v})[d-tg] \otimes \rec_{F_v}^\vee(\pi_v)= \\
\sum_{r=1}^{s-t} (-1)^{r} j^{\geq (t+r)g}_! HT(g,t+r,\pi_v,[\overleftarrow{t-1}]_{\pi_v}
\overrightarrow{\times} \Pi_r) [d-(t+r)g] \otimes \rec_{F_v}^\vee(\pi_v) (-\frac{(rg-r}{2})
\end{multline*}
où $\Pi_r=(-1)^{r-1} [\overleftarrow{r-1}]_{\pi_v} + \sum_{i=1}^{r-1} (-1)^{i-1}
[\overleftarrow{i-1}]_{\pi_v} \overrightarrow{\times}
[\overrightarrow{r-i-1}]_{\pi_v}=[\overrightarrow{r-1}]_{\pi_v}$ d'où le résultat.
\end{proof}

\subsection{Filtration de monodromie-locale en hauteur non maximale}
\label{bf2}

On rappelle que d'après l'hypothèse de récurrence, le théorème (\ref{theo-local-fil}) est connu pour les
modèles locaux de Deligne-Carayol de hauteur strictement inférieure à $d$. Ainsi l'aboutissement de la suite
spectrale (\ref{suite-spectrale}) est connue en dehors des points supersinguliers tandis qu'on ne connaît que
les germes en des points non supersinguliers des termes $E_1$ (cf. le théorème-définition
(\ref{theo-defi-monodromie})).

\begin{prop} \label{prop-hij}
Pour $g \neq 1$ et $1 \leq t \leq s_g$, les faisceaux de cohomologie $h^i \PC(g,t,\pi_v)$ sont nuls pour tout
$i < t-s_g$ et $i$ qui n'est pas de la forme $tg-d+a(g-1)$. Pour $i=tg-d+a(g-1)$ avec $0 \leq a < s_g-t$, ils
sont égaux dans $\FH(\bar X)$ à
$$j^{\geq (t+a)g}_! HT(g,t+a,\pi_v, [\overleftarrow{t-1}]_{\pi_v} \overrightarrow{\times}
[\overrightarrow{a-1}]_{\pi_v}) \otimes \rec_{F_v}^\vee(\pi_v) (-\frac{a(g-1)}{2}).$$
\end{prop}

\begin{proof} Pour $i < tg-d$, les $h^i \PC(g,t,\pi_v)$ sont tous nuls car $\PC(g,t,\pi_v)$ est de dimension
$d-tg$. On a la suite exacte courte de faisceaux pervers
$$0 \to P_{\pi_v,1,0}(\pi_v) \otimes \rec_{F_v}^\vee(\pi_v) \longto j^{\geq g}_! \FC(g,t,\pi_v,\pi_v)[d-tg] \otimes
\rec_{F_v}^\vee(\pi_v)
\longto \PC(g,t,\pi_v) \to 0$$
où $P_{\pi_v,1,0}(\pi_v)$ est un faisceau pervers de dimension $d-(t+1)g$ de sorte que
$h^{tg-d} j^{\geq g}_{!*} \FC(g,1,\pi_v)[d-g]$ est isomorphe à $j^{\geq g}_! \FC(g,1,\pi_v)$.
Le principe est alors d'utiliser le théorème de
comparaison de Berkovich-Fargues couplé avec le lemme (\ref{lem-hij0}).
Remarquons tout d'abord que d'après le corollaire (\ref{coro-grk}), pour tout $k$, $gr_{k,\pi_v}$ est pur
hors des points supersinguliers, de sorte qu'en ce qui concerne les strates non supersingulières, on a
$$e_{\pi_v} h^i gr_{k,\pi_v}=\bigoplus_{\atop{|k| < t \leq s_g}{t \equiv k-1 \mod 2}} h^i
\PC(g,t,\pi_v)(-\frac{tg-1+k}{2})$$
Ainsi d'après le corollaire (\ref{coro-comparaison-grk}), on en déduit le lemme suivant.

\begin{lemm} \label{lem-hij1}
Les $h^i \PC(g,t,\pi_v)$ vérifient les propriétés suivantes:

\begin{itemize}
\item hormis les points supersinguliers, ils sont à support dans les strates $\bar X_{U^p,m}^{(d-t'g)}$ pour
$t \leq t' \leq s_g$ avec $t'g-d-t'+t-1' \leq i \leq t'g-d$ et $i \equiv t'g-d-t'+t-1 \mod 2$;

\item pour $i=t'g-d-t'+t-1+2r$, la fibre en un point géométrique de $\bar X_{U^p,m}^{(d-t'g)}$ de $h^i
\PC(g,t,\pi_v)$ est un facteur direct de
$$[\overleftarrow{t+2r-1}]_{\pi_v} \overrightarrow{\times}
[\overrightarrow{t'-t-2r-1}]_{\pi_v} \otimes \rec_{F_v}^\vee(\pi_v) (-\frac{t'(g-1)+2(t-1)+2r}{2})$$
\end{itemize}
\end{lemm}

\marque Soit $z$ un point géométrique de $\bar X_{U^p,m}^{(d-t'g)}$. Les strates étant induites, on en déduit
que la fibre en $z$ de $h^i \PC(g,t,\pi_v)(-\frac{t(g-1)}{2})$ est de la forme:
$$\bigoplus_{\chi} ([\overleftarrow{t-1}]_{\pi_v \circ \chi((t'-t)(g-1)/2)} \overrightarrow{\times} \pi_\chi)
\otimes \rec_{F_v}^\vee(\pi_v \circ \chi) (-\frac{t(g-1)}{2})$$ où $\pi_\chi$ est une représentation de
$GL_{(t'-t)g}(F_v)$. Pour $\chi=\Xi^{(-(t'-t)(g-1)-2r)/2}$ avec $r>0$, on remarque que si
$[\overleftarrow{t-1}]_{\pi_v(r)} \overrightarrow{\times} \pi_\chi$ contient un des deux constituants de
$[\overleftarrow{t+2r-1}]_{\pi_v} \overrightarrow{\times} [\overrightarrow{t'-t-2r-1}]_{\pi_v}$
alors il contient aussi tous les constituants de
$([\overleftarrow{t-1}]_{\pi_v} \overrightarrow{\times} [\overrightarrow{t'-t-2r-1}]_{\pi_v})
\overleftarrow{\times} [\overleftarrow{2r-1}]_{\pi_v}$ alors que, par exemple
$[\overleftarrow{2r-1},\overrightarrow{1},\overleftarrow{t-1},\overrightarrow{t'-t-2r}]_{\pi_v} \otimes
\rec_{F_v}^\vee(\pi_v) (-\frac{t'(g-1)+2r}{2})$ n'est pas un constituant de $e_{\pi_v} h^{t'g-d-t'+t-1+2r}
gr_{1-t,\pi_v}$ d'où la contradiction.

On en déduit ainsi que $h^{t'g-d-t'+t-1} \PC(g,t,\pi_v)$ est le seul faisceau de cohomologie ayant un support
d'intersection non vide avec la strate $t'g$. Le raisonnement étant valide pour tout $t'$, on en déduit aussi
que le support de $h^{t'g-d-t'+t-1} \PC(g,t,\pi_v)$ est contenu dans la strate $t'g$. Le lemme
(\ref{lem-hij0}) donne alors que la restriction à la strate $t'g$ de $h^{t'g-d-t'+t-1} \PC(g,t,\pi_v)$ est
isomorphe à $HT(g,t',\pi_v, [\overleftarrow{t-1}]_{\pi_v} \overrightarrow{\times}
[\overrightarrow{t'-t-1}]_{\pi_v}) \otimes \rec_{F_v}^\vee(\pi_v) (-\frac{(t'-t)(g-1)}{2}),$ de sorte que
$$h^{t'g-d-t'+t-1} \PC(g,t,\pi_v) \simeq j^{\geq t'g}_! HT(g,t',\pi_v, [\overleftarrow{t-1}]_{\pi_v}
\overrightarrow{\times} [\overrightarrow{t'-t-1}]_{\pi_v}) \otimes \rec_{F_v}^\vee(\pi_v)
(-\frac{(t'-t)(g-1)}{2}).$$

\end{proof}

\begin{prop} \label{prop-hij1}
Pour $g=1$, les $h^i \PC(1,t,\chi_v)$, sont nuls pour $i \neq t-d$ et pour $i=t-d$ leur restriction à $\bar
X_{U^p,m}^{(d-t-r)}$, pour $0 \leq r < d-t$ est égale dans $\FH(\bar X^{d-t-r})$ à
$HT(1,t+r,\chi_v,[\overleftarrow{t-1}]_{\chi_v} \overrightarrow{\times} [\overrightarrow{r-1}]_{\chi_v})$.
\end{prop}

\begin{proof} Elle est strictement identique à celle pour $g \neq 1$, en considérant à chaque étape les restrictions
aux strates $\bar X_{U^p,m}^{(d-t-r)}$.

\end{proof}

\begin{prop} \label{prop-hij-ss}
Pour $g$ ne divisant pas $d$, $h^{s_g g-d-(s_g-l)} \PC(g,t,\pi_v)$ est, dans $\FH(\bar X)$, isomorphe à
$$j^{\geq s_g g}_! HT(g,s_g,\pi_v,([\overleftarrow{t-1}]_{\pi_v} \overrightarrow{\times}
[\overrightarrow{s_g-t-1}]_{\pi_v})) \otimes \rec_{F_v}^\vee(\pi_v)(-\frac{(s_g-t)(g-1)}{2}).$$ Par ailleurs
les $h^i \PC(g,t,\pi_v)$ sont nuls pour $s_g g -d -(s_g-t) < i \leq 0$.
\end{prop}

\begin{proof} Tant qu'on est en dehors des points supersinguliers, les arguments précédents s'adaptent, en utilisant
le théorème de comparaison de Berkovich-Fargues avec la connaissance de la filtration de monodromie-locale
des modèles de Deligne-Carayol en hauteur strictement inférieure à $d$. En ce qui concerne les points
supersinguliers, on raisonne par récurrence sur $t$ de $s_g$ à $1$. Pour $t=s_g$, le corollaire
(\ref{coro-sec-fp}), joint au corollaire (\ref{coro-pnul}), donne
$$\PC(g,s_g,\pi_v)=j^{\geq s_g g}_! HT(g,s_g,\pi_v,[\overleftarrow{s_g-1}]_{\pi_v})[d-s_gg] \otimes
\rec_{F_v}^\vee(\pi_v)$$
d'où le résultat. Supposons donc le résultat acquis jusqu'au rang $t+1$ et traitons le cas de $t$. On
considère les suites exactes courtes du corollaire (\ref{coro-sec-fp}), où l'on rappelle que d'après
(\ref{coro-pnul}), les $A_{\pi_v,t,g}$ sont nuls:
\begin{multline*}
0 \to P_{\pi_v,t,0}([\overleftarrow{t-1}]_{\pi_v}) \longto j^{\geq tg}_! HT(g,t,\pi_v,
[\overleftarrow{t-1}]_{\pi_v})[d-tg] \longto \\ j^{\geq tg}_{!*}
HT(g,t,\pi_v,[\overleftarrow{t-1}]_{\pi_v})[d-tg] \to 0
\end{multline*}
$$\cdots$$
$$0 \to P_{\pi_v,t,s_g-l}([\overleftarrow{t-1}]_{\pi_v}) \longto
P_{\pi_v,t,s_g-t-1}([\overleftarrow{t-1}]_{\pi_v}) \longto \PC_-(g,s_g,\pi_v,t,[\overleftarrow{t-1}]_{\pi_v})
\to 0$$ où $P_{\pi_v,t,s_g-t}([\overleftarrow{t-1}]_{\pi_v})$ est le faisceau pervers nul. On démontre alors
par récurrence sur $r$, de $s_g-t$ à $0$, que le germe en un point supersingulier de $h^i
P_{\pi_v,t,r}([\overleftarrow{t-1}]_{\pi_v})$ est nul pour tout $i$, d'où le résultat.

\end{proof}

En direction du théorème (\ref{theo-global2}), on a le résultat suivant.

\begin{lemm} \label{lem-hij3}
Pour $g \neq 1$ divisant $d=sg$ (resp. $g=1$), la fibre en un point supersingulier de $h^i \PC(g,t,\pi_v)$
(resp. de $h^i \PC(1,t,\chi_v)$) est nulle pour $i<t-s$. Par ailleurs les $h^{t-s+i} \PC(g,t,\pi_v)$ sont
concentrés aux points supersinguliers pour $0 \leq i < s-t$, de fibre (resp. la fibre en un point
supersingulier de $h^{t-s+i} \PC(1,t,\chi_v)$ est), en tant que $GL_d(F_v) \times W_v$-module, un sous-espace,
éventuellement nul, de
$$[\overleftarrow{t-1}]_{\pi_v} \overrightarrow{\times} [\overrightarrow{s-t-2-i}]_{\pi_v}
\overrightarrow{\times} [\overleftarrow{i}]_{\pi_v} \otimes \rec_{F_v}^\vee(\pi_v) (-\frac{(s-t)(g-1)}{2}).$$
\end{lemm}

\begin{proof} Nous ne traitons que le cas $g \neq 1$, le raisonnement pour $g=1$ étant strictement identique. On
raisonne par récurrence sur $t$ de $s$ à $1$, le cas $t=s$ étant trivial. Supposons donc le résultat acquis
jusqu'au rang $t+1$ et traitons le cas de $t$. On considère alors les suites exactes courtes de faisceaux
pervers de (\ref{coro-sec-fp}). On remarque tout d'abord que toutes les fibres en un point supersingulier des
faisceaux de cohomologie des $\PC_-(g,t+r,\pi_v,t,\Pi_t)$ sont de poids $(s-t)(g-1)$ de sorte qu'il en est de
même pour les $P_{\pi_v,t,r}(\Pi_t)$ et donc pour $\PC(g,t,\pi_v)$. Le résultat découle alors de l'étude des
suites exactes longues associées. Une façon plus visuelle et plus directe d'obtenir le résultat est de
considérer l'isomorphisme $Dj^{\geq tg}_{!*} \FC(g,t,\pi_v)[d-tg] \simeq j^{\geq tg}_{!*}
\FC(g,t,\pi_v^\vee)[d-tg](d-tg)$ et de regarder la suite spectrale
\begin{equation} \label{ss-dualite}
E_2^{p,q}=R^p Hom(h^{-q} \PC(g,t,\pi_v),K_{s_o}) \Rightarrow h^{p+q} \PC(g,t,\pi_v^\vee)(d-1)
\end{equation}
où $K_{s_o}$ désigne le complexe dualisant sur $\bar X_{U^p,m}$ \footnote{On pourra voir les termes
$E_2^{p,q}$ à la figure (\ref{figure1}) (resp. (\ref{figure2})), pour $g=7$, $s=5$ et $t=4$ (resp. $t=3$,
$t=2$ et $t=1$).}. On rappelle que pour un faisceau $\LC$ sur $\bar X_{U^p,m}^{(d-(l+r)g)}$, par adjonction
on a
$$RHom(i^{(t+r)g}_* j^{\geq (t+r)g}_! \LC,f^! \bar \Qm_l) \simeq  i^{(t+r)g}_* Rj^{\geq (t+r)g}_*
RHom(\LC,j^{\geq (t+r)g,!} i^{(t+r)g,!} f^! \bar \Qm_l)$$ et comme $\bar X_{U^p,m}^{(d-(t+r)g)}$ est lisse,
on a
$j^{\geq (t+r)g,!} i^{(t+r)g,!} f^! \bar \Qm_l \simeq \bar \Qm_l[2(d-(t+r)g)](d-(t+r)g)$
soit pour $p \geq -2(d-(t+r)g)$,
\begin{multline} \label{etoile}
R^p Hom (j^{\geq (t+r)g}_! HT(g,t+a,\pi_v,  [\overleftarrow{t-1}]_{\pi_v} \overrightarrow{\times}
[\overleftarrow{a-1}]_{\pi_v}) \otimes |\art_{F_v}^{-1}|^{-a(g-1)/2}, K_{s_o})  \simeq \\ (i^{(t+a)g}_*
R^{p+2(d-(t+a)g)} j^{\geq (l+a)g}_* HT(g,t+a,\pi_v^\vee,  [\overleftarrow{t-1}]_{\pi_v^\vee}
\overleftarrow{\times} [\overleftarrow{a-1}]_{\pi_v^\vee}) \otimes |\art_{F_v}^{-1}|^{-a(g+1)/2})(tg-d)
\end{multline}

Ainsi pour $z$ un point supersingulier, $(E_2^{p,q})_z$, pour $q$ de la forme $d-tg-a(g-1)$ (resp. $q=a+1$)
avec $0 \leq a < s-t$, et $p=-(s-t-a)g-i$ avec $0 \leq i \leq s-t-a$ (resp. $p=0$), s'il est non nul, est
mixte de poids $(s-t-a-r)(g-1)/2+(a+r)(g+1)/2$ avec $0 \leq r <s-t-a$ si $p<-(s-t-a)g$ et pur de poids
$(s-t)(g+1)/2$ si $p=(s-t-a)g$. Dans ce dernier cas on obtient alors
$$[\overleftarrow{t-1}]_{\pi_v} \overrightarrow{\times} [\overrightarrow{s-t-a-1}]_{\pi_v} \overleftarrow{\times}
[\overleftarrow{a-1}]_{\pi_v} \otimes |\art_{F_v}^{-1}|^{-(s-t)(g+1)/2}$$ (resp. $(h^{-q} j^{\geq tg}_{!*}
\FC(g,t,\pi_v)_0)_z^\vee (-tg) \otimes |\art_{F_v}^{-1}|^{-(s-t)(g+1)/2}$ où le dual est pris en tant que
représentation de $GL_{(s-t)g}(F_v)$). Dans la figure (\ref{fig-dualite2}) on représente ces $(E_2^{p,q})_z$
de poids $(s-t)(g+1)$. Le résultat découle alors du fait que les $E_\oo^i$ sont tous nuls pour $i \geq 0$.

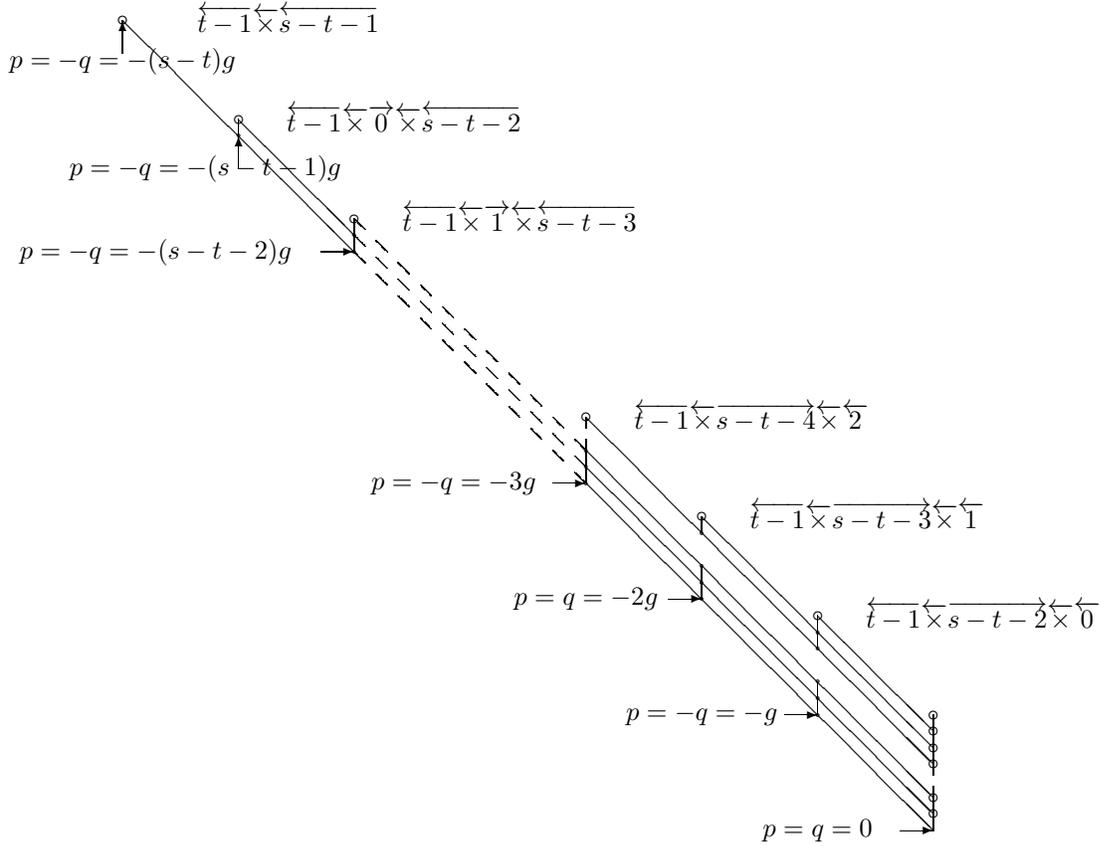
\begin{figure}[!ht]
\setlength{\unitlength}{.22cm}
\begin{picture}(50,53)(-5,-50)
\linethickness{.1pt}

\put(0,0){\circle{.4}} \put(10,0){\makebox(0,0){$\overleftarrow{t-1} \overleftarrow{\times}
\overleftarrow{s-t-1}$}}

\put(0,-2){\vector(0,1){2}} \put(0,-2.5){\makebox(0,0){$p=-q=-(s-t)g$}} \put(0,0){\line(1,-1){14}}

\put(7,-7){\circle{.1}} \put(7,-7){\line(0,1){1}} \put(7,-6){\circle{.4}} \put(7,-6){\line(1,-1){7}}
\put(7,-9){\vector(0,1){2}} \put(5,-9){\makebox(0,0){$p=-q=-(s-t-1)g$}}

\put(17,-6){\makebox(0,0){$\overleftarrow{t-1} \overleftarrow{\times} \overrightarrow{0}
\overleftarrow{\times} \overleftarrow{s-t-2}$}} \put(14,-14){\circle{.1}} \put(14,-14){\line(0,1){2}}
\put(12,-14){\vector(1,0){2}} \put(2,-14){\makebox(0,0){$p=-q=-(s-t-2)g$}} \dashline{1}(14,-14)(28,-28)
\put(14,-13){\circle{.1}} \dashline{1}(14,-13)(28,-27) \put(14,-12){\circle{.4}}

\put(24,-12){\makebox(0,0){$\overleftarrow{t-1} \overleftarrow{\times} \overrightarrow{1}
\overleftarrow{\times} \overleftarrow{s-t-3}$}} \dashline{1}(14,-12)(28,-26) \put(28,-28){\line(0,1){2}}
\put(26,-28){\vector(1,0){2}} \put(20,-28){\makebox(0,0){$p=-q=-3g$}} \dashline{1}(28,-26)(28,-24)
\multiput(28,-28)(0,1){3}{\circle{.1}} \multiput(28,-28)(0,1){3}{\line(1,-1){21}} \put(28,-24){\circle{.4}}
\put(28,-24){\line(1,-1){21}}

\put(38,-24){\makebox(0,0){$\overleftarrow{t-1} \overleftarrow{\times} \overrightarrow{s-t-4}
\overleftarrow{\times} \overleftarrow{2}$}} \put(35,-30){\circle{.4}} \put(35,-30){\line(1,-1){14}}
\put(35,-35){\line(0,1){2}} \put(33,-35){\vector(1,0){2}} \put(28,-35){\makebox(0,0){$p=q=-2g$}}
\multiput(35,-35)(0,1){3}{\circle{.1}} \put(35,-31){\circle{.1}} \put(35,-31){\line(0,1){1}}

\put(45,-30){\makebox(0,0){$\overleftarrow{t-1} \overleftarrow{\times} \overrightarrow{s-t-3}
\overleftarrow{\times} \overleftarrow{1}$}} \put(42,-36){\circle{.4}} \put(42,-36){\line(1,-1){7}}
\put(42,-42){\line(0,1){2}} \multiput(42,-42)(0,1){3}{\circle{.1}} \put(40,-42){\vector(1,0){2}}
\put(35,-42){\makebox(0,0){$p=-q=-g$}} \put(42,-38){\line(0,1){2}} \multiput(42,-38)(0,1){2}{\circle{.1}}

\put(52,-36){\makebox(0,0){$\overleftarrow{t-1} \overleftarrow{\times} \overrightarrow{s-t-2}
\overleftarrow{\times} \overleftarrow{0}$}} \put(49,-49){\line(0,1){2}} \put(47,-49){\vector(1,0){2}}
\put(42,-49){\makebox(0,0){$p=q=0$}} \dashline{1}(49,-47)(49,-45) \put(49,-45){\line(0,1){3}}
\multiput(49,-48)(0,1){2}{\circle{.4}} \multiput(49,-45)(0,1){4}{\circle{.4}}

\end{picture}
\caption{\label{fig-dualite2} parties de poids $(s-t)(g+1)$ de la fibre aux points supersinguliers des
$E_2^{p,q}$ de la suite spectrale (\ref{ss-dualite})}
\end{figure}

\end{proof}

\section{Preuve des théorèmes globaux sous (9.1.1)}
\label{globss}

On rappelle, cf. \S \ref{schema}, que le théorème (\ref{theo-local-fil}) découle, d'après le théorème de
comparaison de Berkovich-Fargues, des théorèmes globaux (\ref{theo-global1}), (\ref{theo-global2}) et
(\ref{theo-ss}). En outre le théorème (\ref{theo-local-fil}) implique, grâce à (\ref{theo-ss}) via
Berkovich-Fargues, le théorème (\ref{theo-ripsi-local}). Par ailleurs le principe de la preuve de la
proposition (\ref{prop-hij}) montre que (\ref{theo-local-fil}) implique les théorèmes globaux. Le but de ce
paragraphe est de montrer que le théorème (\ref{theo-ripsi-local}) implique les théorèmes globaux
(\ref{theo-global1}), (\ref{theo-global2}) et (\ref{theo-ss}) et donc le théorème local
(\ref{theo-local-fil}). Par souci d'efficacité, on montrera, en utilisant l'opérateur $N$, qu'il nous suffit
en fait de connaître les parties de poids $s(g-1)$ de (\ref{theo-ripsi-local}).

\subsection{Preuve de (\ref{theo-global2}) sous (\ref{cas-m})}
\label{globss1}

\begin{prop} \label{cas-m}
Pour tout diviseur $g$ de $d=sg$ et $\pi_v$ une représentation cuspidale irréductible de $GL_g(F_v)$, la
partie de poids $d-s$ de $\UC_{F_v,l,d,\xi_v}^{d-s+i}(\JL^{-1}([\overleftarrow{s-1}]_{\pi_v^\vee})$ est nulle
pour $0 < i <s$ et égale à $[\overrightarrow{s-1}]_{\pi_v} \otimes \rec_{F_v}^\vee(\pi_v) (-\frac{d-s}{2})$
pour $i=0$.
\end{prop}

La preuve de cette proposition sera donnée au paragraphe suivant; montrons comment le théorème
(\ref{theo-local-fil}) en découle. On commence par prouver le théorème (\ref{theo-global2}) qui d'après la
proposition (\ref{prop-hij}) découle alors de la proposition suivante.

\begin{prop} \label{prop-hij-ss-g}
Sous le résultat de la proposition (\ref{cas-m}), pour $g \neq 1$ divisant $d$, $h^{t-s+i} \PC(g,t,\pi_v)$
est nul pour $i \neq 0$ et sinon est concentré aux points supersinguliers de fibre isomorphe à
$[\overleftarrow{t-1}]_{\pi_v} \overrightarrow{\times} [\overrightarrow{s-t-1}]_{\pi_v} \otimes
\rec_{F_v}^\vee(\pi_v) (-\frac{(s-t)(g-1)}{2})$ pour $i=0$.
\end{prop}

\begin{proof} D'après le lemme (\ref{lem-hij3}) tous les germes en un point supersingulier des $\PC(g,t,\pi_v)$,
sont de poids $(s-t)(g-1)$. Le raisonnement de la preuve de la proposition (\ref{prop-hij}) s'applique alors
tel quel en étudiant les $\PC(g,t,\pi_v)(-\frac{t(g-1)}{2})$ pour $1 \leq t \leq s$ et en utilisant la
connaissance des parties de poids $s(g-1)$ de (\ref{theo-local-fil}).

\end{proof}

\subsection{Preuve de (\ref{theo-ss}) sous (\ref{cas-m})}
\label{globss2}

\begin{proof} On considère la suite spectrale
(\ref{suite-spectrale}) associée à la filtration de monodromie. D'après la proposition précédente et le point
(ii) de la proposition (\ref{prop-p}), le germe en un point supersingulier $z$ de $E_1^{i,j}$ vérifie les
propriétés suivantes (cf. la figure (\ref{figure7}) dans le cas $s=4$ et $g=3$):

\begin{itemize}
\item il est nul pour $i,j$ ne vérifiant pas la condition suivante: il existe $0 \leq k \leq s-1$ tel
que $j=1-s+2k$ et $-k \leq i \leq s-1-2k$;

\item pour $j=1-s+2k$ avec $0 \leq k \leq s-1$, le germe en $z$ de $E_1^{s-1-2k-r,1-s+2k}$ pour $0 \leq r \leq s-1-k$
est isomorphe à
$[\overleftarrow{s-r-1}]_{\pi_v} \overrightarrow{\times} [\overrightarrow{r-1}]_{\pi_v} \otimes
\rec_{F_v}^\vee(\pi_v)
(-\frac{s(g-1)+2k}{2}).$
\end{itemize}

En outre d'après (\ref{cas-m}), pour $j=1-s$ et pour tout $0 \leq i < s-1$, l'application
$d_1^{i,1-s}:E_1^{i,1-s} \longto E_1^{i+1,1-s}$ induit en $z$ la flèche non triviale
$$[\overleftarrow{i}]_{\pi_v} \overrightarrow{\times} [\overrightarrow{s-i-2}]_{\pi_v} \longto
[\overleftarrow{i+1}]_{\pi_v} \overrightarrow{\times} [\overrightarrow{s-i-3}]_{\pi_v}$$ dont le noyau est
$[\overleftarrow{i},\overrightarrow{s-i-1}]_{\pi_v}$ et le conoyau est
$[\overleftarrow{i+2},\overrightarrow{s-i-3}]_{\pi_v}$. Le logarithme $N^k$ de la partie unipotente de la
monodromie induit un diagramme commutatif
$$\xymatrix{ E_1^{i,1-s} \rrto^{d_1^{i,1-s}} & & E_1^{i+1,1-s} \\
E_1^{i-2k,1-s+2k} \uto^{N^k} \rrto^{d_1^{i-2k,1-s+2k}} & & E_{i-2k+1,1-s+2k} \uto^{N^k}}$$ En remarquant que
pour $0 \leq k \leq -i$, $N^k$ induit des isomorphismes $E_1^{i-2k,1-s+2k} \simeq E_1^{i,1-s}$ et
$E_1^{i+1-2k,1-s+2k} \simeq E_1^{i+1,16s}$, on en déduit que la flèche $d_1^{i-2k,1-s+2k}$ est aussi induite
par la flèche non triviale
$[\overleftarrow{i}]_{\pi_v} \overrightarrow{\times} [\overrightarrow{s-i-2}]_{\pi_v} \longto
[\overleftarrow{i+1}]_{\pi_v} \overrightarrow{\times} [\overrightarrow{s-i-1}]_{\pi_v}$.

Finalement pour tout $1 \leq t' < t \leq s$ et tout $k \equiv t'-1 \mod 2$ avec $|k| \leq t'-1$, les flèches
$$h^{-d+tg-(t-t')} gr_{k,\pi_v} \longto h^{-d+tg-(t-t'-1)} gr_{k-1,\pi_v}$$
de la suite spectrale (\ref{suite-spectrale}), sont non nulles, et se déduisent pour $t' < t-1$, des suites
exactes
\begin{multline*} \label{suites-exactes}
0 \to [\overleftarrow{t'-1},\overrightarrow{t-t'}]_{\pi_v} \to [\overleftarrow{t'-1}]_{\pi_v}
\overrightarrow{\times} [\overrightarrow{t-t'-1}]_{\pi_v} \to \\ \to [\overleftarrow{t'}]_{\pi_v}
\overrightarrow{\times} [\overrightarrow{t-t'-2}]_{\pi_v} \to
[\overleftarrow{t'+1},\overrightarrow{t-t'-2}]_{\pi_v} \to 0
\end{multline*}
et pour $t'=t-1$ de la suite exacte $0 \to [\overleftarrow{t-2},\overrightarrow{1}]_{\pi_v} \to
[\overleftarrow{t-2}]_{\pi_v} \overrightarrow{\times} [\overrightarrow{0}]_{\pi_v} \to
[\overleftarrow{t-1}]_{\pi_v} \to 0$.
\end{proof}

\rem On verra plus loin (cf. les remarques (\ref{rema-N}) et (\ref{rema-N2})) qu'on pourrait en fait montrer
(\ref{theo-local-fil}) directement, sans utiliser l'opérateur de monodromie $N$.

\begin{coro}
Pour $0 \leq i \leq d-tg$ et $g >1$, $R^ij^{\geq tg}_* HT(g,t,\pi_v,\Pi_t)[d-tg]$ est, dans $\FH(\bar X)$,
une somme directe sur tous les couples $(n,r)$ tels que $ng+r(g-1)=i$ des faisceaux de type $HT(g,t+n+r)$

$$j^{\geq (t+n+r)g}_!HT(g,t+n+r,\pi_v,(\Pi_t \overrightarrow{\times} [\overleftarrow{n-1}]_{\pi_v})
\overleftarrow{\times} [\overrightarrow{r-1}]_{\pi_v}) \otimes \Xi^{\frac{n(g+1)+r(g-1)}{2}}$$ Pour $g=1$,
$R^ij^{\geq h}_* HT(1,h,\chi_v,\Pi_h)[d-h]$ est à support dans la tour des $\bar X_{U^p,m}^{(d-h-i)}$ et sa
restriction à $\bar X_{U^p,m}^{(d-h-i-r)}$ est le faisceau: $HT(1,h+i+r,\chi_v, (\Pi_h
\overrightarrow{\times} [\overleftarrow{i-1}]_{\chi_v}) \overleftarrow{\times}
[\overleftarrow{r-1}]_{\chi_v}) \otimes \Xi^{i}$.
\end{coro}

\subsection{Pureté de la filtration de monodromie de $R\Psi_v[d-1]$}
\label{globalss3}

La proposition suivante correspond au théorème (\ref{theo-global1}) qui rappelons le découle, en utilisant la
pureté de la filtration de monodromie, du théorème (\ref{theo-global0}) qui est prouvé au corollaire
(\ref{global0-ok}).

\begin{prop} \label{prop-grk-final}
Pour tout $|k| < s_g$, $e_{\pi_v} gr_{k,\pi_v}$ est, dans $\FPH(\bar X)$, égal à
$$\sum_{\atop{|k| < t \leq s_g}{\equiv k-1 \mod 2}} \PC(g,t,\pi_v) (-\frac{tg+k-1}{2})$$
Dans tous les autres cas $gr_{k,\pi_v}$ est nul.
\end{prop}

\begin{proof} Par rapport au corollaire (\ref{coro-grk}), il s'agit de traiter les points supersinguliers et donc,
d'après le corollaire (\ref{coro-pnul}) de montrer que pour tout $r \equiv s-1 \mod 2$ et $|r| \leq s-1$,
$\PC(g,s,\pi_v)(-\frac{sg+r-1}{2})$ est un constituant de $e_{\pi_v} gr_{r,\pi_v}$. \footnote{On remarquera
que le raisonnement suivant est valable pour toutes les strates de sorte qu'en raisonnant par récurrence
l'argument du corollaire (\ref{coro-grk}) en utilisant le théorème de comparaison de Berkovich-Fargues
n'était pas strictement nécessaire.}

On commence par remarquer, en utilisant la partie unipotente $N$ de la monodromie, qu'il suffit en fait de
montrer que $\PC(g,s,\pi_v) (-\frac{s(g-1)}{2})$ est un constituant de $e_{\pi_v} gr_{1-s,\pi_v}$.
Considérons alors $\PC(g,s,\pi_v) (-\frac{s(g-1)}{2})$ et soit $r$ tel qu'il soit un constituant de
$e_{\pi_v} gr_{r,\pi_v}$. Soit $z_I$ un point géométrique de $\bar X_{U^p,m}^{(0)}$. D'après la proposition
(\ref{prop-hij-ss-g}), les germes en $z_I$ des $h^{-1} gr_{k,\pi_v}$ sont tous de poids strictement plus
grand que $s(g-1)$ sauf pour $k=2-s$. Par ailleurs le germe en $z_I$ de $R^{d-1}\Psi_{\pi_v}$ n'a aucun
constituant de poids $s(g-1)$ de sorte que $r \leq 2-s$.

\noindent - Si on avait $r \leq -s$, la partie unipotente $N$ de la monodromie donnerait que
$\PC(g,s,\pi_v)(-\frac{s(g-1)-2r}{2})$
serait un constituant de $e_{\pi_v} R\Psi_{\pi_v}[d-1]$ ce qui n'est pas d'après (\ref{prop-libre})
(\ref{prop-p}) et (\ref{coro-pnul}), car $1-s-2r > s-1$.

\noindent - Supposons $r=2-s$: pour $s>2$, par application de la dualité de Verdier et de l'opérateur de monodromie
$N$, on en
déduit que $\PC(g,s,\pi_v) \otimes \rec_{F_v}^\vee(\pi_v)(-\frac{s(g+1)-2}{2})$ et. $\PC(g,s,\pi_v)
(-\frac{s(g+1)-4}{2})$ sont des constituants de $e_{\pi_v}gr_{s-2,\pi_v}$. Ainsi par une nouvelle application de $N$,
$\PC(g,s,\pi_v) (-\frac{s(g+1)-4}{2})$ est aussi un constituant de $e_{\pi_v} gr_{s-5,\pi_v}$ et donc devrait
apparaître deux fois ce qui n'est pas. Dans le cas $s=2$, la situation défavorable correspondrait à
$\PC(g,2,\pi_v)(-g)$ et $\PC(g,2,\pi_v)(1-g)$ constituants de $e_{\pi_v} gr_{0,\pi_v}$ de sorte que dans la
cohomologie globale, la monodromie serait nulle ce qui est en contradiction avec la proposition suivante pour
$s=2$.

\end{proof}

\begin{prop} \label{prop-mono}
Pour tout diviseur $g$ de $d=sg$ et toute représentation irréductible cuspidale de $GL_g(F_v)$, il existe une
représentation irréductible automorphe $\Pi$ de $D_\Am^\times$ vérifiant $\hyp(\oo)$ telle que $\Pi_v \simeq
[\overleftarrow{s-1}]_{\pi_v}$ telle que l'opérateur de monodromie $N$ sur $H^{d-1}_{\eta_v}[\Pi^\oo]$ soit
non nul où $H^i_{\eta_v}={\DS \lim_{\atop{\to}{U^p(m)}}} H^i(X_{U^p,m} \otimes_{\OC_v} \bar F_v,\LC_\xi)$.
\end{prop}

\begin{proof} Le principe, qui nous a été suggéré par M. Harris, est de se ramener au cas Iwahori par changement de
base résoluble. Commençons donc par traiter le cas Iwahori. Notons que récemment Yoshida et Taylor, cf.
\cite{y-t}, ont rédigé ce résultat qui par ailleurs m'avait été expliqué par A. Genestier. Cependant dans
notre cas, vu que l'on ne s'intéresse qu'au cas $s=2$ et donc $g=1$, le résultat a déjà été prouvé par
Carayol, cf. \cite{ca}.

Considérons désormais le cas général de $\pi_v$ une représentation irréductible cuspidale de $GL_g(F_v)$ pour
$g>1$. On reprend les notations de \cite{h-t}. La somme alternée de la cohomologie de la variété globale à
valeur dans $\LC_\xi$, dans le groupe de Grothendieck correspondant, y est écrite sous la forme
$$\sum_{\pi^\oo} \pi^\oo \otimes [R_\xi(\pi^\oo)] \qquad [R_\xi(\pi^\oo)]=\sum_i (-1)^i R^i_\xi(\pi^\oo)$$
où $\pi^\oo$ décrit l'ensemble des représentations irréductibles de $G(\Am^\oo)$. Par ailleurs, corollaire
VI.2.7 de loc. cit., s'il existe une représentation irréductible $\pi_\oo$ de $G_\tau({\mathbb R})$ telle que
$\pi:= \pi^\oo \otimes \pi_\oo$ vérifie $\hyp(\oo)$ avec $BC(\pi)=(\psi,\Pi)$ où $\JL(\Pi)$ est cuspidale,
avec les notations de loc. cit., alors $R^i_\xi(\pi^\oo)$ est nulle pour $i \neq d-1$.

\noindent Soit alors $\sigma_v=\rec_{F_v}(\pi_v)$ et soit $L_v$ une extension de $F_v$ telle que la restriction
de $\sigma_v$ à $L_v$ soit non ramifiée. On rappelle que l'extension $L_v/F_v$ est résoluble; cf.\cite{se} IV-2.
On globalise alors la situation comme dans \cite{h-t}: soit $(F')^+/F^+$ une extension résoluble de corps
totalement réels telle que:

-la place $v$ de $F^+$ soit inerte dans $(F')^+$;

- l'extension $(F')^+_v / F^+_v$ est isomorphe à l'extension $L_v / F_v$.

\noindent On pose alors $F'=E(F')^+$. D'après le corollaire VI.2.6 de loc. cit., on peut choisir une
représentation automorphe cuspidale $\Pi$ de $GL_d(\Am_F)$ telle que $\Pi$ vérifie:

- $\Pi^c \simeq \Pi^\vee$;

- $\Pi_\oo$ a le même caractère central qu'une représentation algébrique de $Res^G_\Qm (GL_d)$;

- $\Pi_v \simeq \st_s(\pi_v \otimes \psi)$ pour un certain caractère $\psi$ de $F'_\omega$.

\noindent Soit alors le changement de base $\Pi':=BC_{F'/F}(\Pi)$; d'après loc. cit. (théorèmes VI.1.1 et
VI.2.9) on associe à $\Pi$ et $\Pi'$ des représentations $\pi$ et $\pi'$ de respectivement $G_\tau(\Am_F)$ et
$G_\tau(\Am_{F'})$ qui vérifient $\hyp(\oo)$ avec $\pi'_v \simeq [\overleftarrow{s-1}]_{\zeta} \boxplus
\cdots \boxplus [\overleftarrow{s-1}]_{\zeta}$, pour $\zeta$ un caractère de $F'_o$. Par ailleurs d'après le
théorème VII.1.9 de loc. cit. appliqué aux bonnes places de $F$, et en utilisant le théorème de densité de
Cebotarev, on en déduit que la représentation galoisienne $R_{\xi'}(\pi^{',\oo})$ est isomorphe à
$R_\xi(\pi^\oo)_{|\gal(\bar F'/F')}$ avec
$R_{\xi'}(\pi^{',\oo})_\omega \simeq (\sp_s \otimes \zeta^{-1})^g$
de sorte que $R_\xi(\pi^\oo) \simeq \sp_s \otimes \rec_{F_v}^\vee(\pi_v)$, d'où le résultat.

\end{proof}

\begin{prop} \label{prop-mono-1}
Soit $\Pi$ une représentation irréductible automorphe de $G_\tau(\Am)$ obtenue via le théorème VI.2.9 de
\cite{h-t} et telle que $\Pi_v$ est l'induite irréductible $[\overleftarrow{s_1-1}]_{\xi_1} \boxplus \cdots
\boxplus [\overleftarrow{s_r-1}]_{\xi_r}$ où les $\xi_i$ sont des caractères de $F_v^\times$. En tant que
représentation de $W_v$, $H^{d-1}_{\eta_v}[\Pi^\oo]$ est la somme directe des $\sharp \ker^1(\Qm,G_\tau)
m(\Pi) \bigoplus_{i=1}^r \sp_{s_i} \otimes \xi_i$ où $\sp_s$ est la représentation de dimension $s$,
$|-|^{(1-s)/2} \oplus \cdots \oplus |-|^{(s-1)/2}$ où l'indice de nilpotence de l'opérateur de monodromie $N$
est égal à $s$.
\end{prop}

\begin{proof} Le cas $s=2$ et $g=1$ de la proposition (\ref{prop-mono}) est prouvé par Carayol dans \cite{ca} de
sorte que la proposition précédente est vraie pour $g=1$ et $s$ quelconque. On étudie ensuite la suite
spectrale
\begin{equation} \label{ss-poids-hi1}
E_1^{i,j}[\Pi^\oo]=H^{i+j}(gr_{-i})[\Pi^\oo] \Rightarrow H^{i+j}_{\eta_v}[\Pi^\oo]
\end{equation}
où pour tout $i$, $H^i_{\eta_v}$ (resp. $H^i_{\eta_v,\xi}$) désigne la limite projective sur $U^p(m)$ des
groupes de cohomologie de la fibre générique du faisceau constant $\bar \Qm_l$ (resp. $\LC_\xi$) de
$X_{U^p,m}$. Cette étude est faite dans le cas général au théorème (\ref{theo-mono-glob}); le point est qu'en
ce qui concerne les $\Pi^\oo$-parties, nous ne devons considérer que les représentations $\pi_v$ de
$GL_1(F_v)$, i.e. le cas $g=1$, qui rappelons le est connu grâce à Carayol, d'où le résultat.

\end{proof}

\rem Pour des résultats généraux sur les composantes locales $\Pi_v$ ainsi que sur la partie
$\Pi^\oo$-isotypique des groupes de cohomologie de la fibre générique, on renvoie respectivement aux
théorèmes (\ref{prop-compo-locale}) et (\ref{theo-mono-glob}).

\section{Étude de la suite spectrale des cycles évanescents}
\label{ssce0}

Le but de ce paragraphe est de prouver la proposition (\ref{cas-m}). Le principe de la preuve est d'étudier
la suite spectrale des cycles évanescents en y intégrant, via les suites spectrales associées à la
stratification, la connaissance des $H^i_c(\bar X_{U^p,m}^{(d-tg)}, \FC(g,t,\pi_v)\otimes \LC_{\xi})$.

\subsection{Cas où $\Pi_v \simeq \st_s(\pi_v)$}

Soit $1 \leq g \leq d$ et $\pi_v$ une représentation irréductible cuspidale de $GL_g(F_v)$. On fixe dans la
suite une représentation automorphe irréductible $\Pi$ de $G_\tau(\Am)$ vérifiant $\hyp(\xi)$ et telle que
$\Pi_v \simeq [\overleftarrow{s-1}]_{\pi_v}$ pour $\pi_v$ une représentation irréductible cuspidale de
$GL_g(F_v)$ et on note $m(\Pi)$ la multiplicité de $\Pi$ dans l'espace des formes automorphes. On commence
par un résultat déjà présent dans \cite{lrs}.

\begin{coro} \label{coho-global1}
Les groupes de cohomologie de la fibre générique
$H^i_{\eta_v,\xi}[\Pi^{\oo,v}]$ sont nuls pour $i \neq d-1$ et
$$H^{d-1}_{\eta_v,\xi}[\Pi^{\oo,v}]=\sharp \ker^1(\Qm,G_\tau) m(\Pi) [\overleftarrow{s-1}]_{\pi_v} \otimes
\rec_{F_v}^\vee([\overleftarrow{s-1}]_{\pi_v})
(-\frac{sg-1}{2})$$
\end{coro}

\begin{proof} On considère la suite spectrale
$$E_{1,\xi}^{i,j}[\Pi^{\oo,v}]:=H^{i+j}(gr_{-i,\xi})[\Pi^{\oo,v}] \Rightarrow H^{i+j+d-1}_{\eta_v,\xi}[\Pi^{\oo,v}]$$
Le théorème (\ref{theo-global1}) découle de la proposition (\ref{prop-mono}). La proposition
(\ref{prop-coho1}) donne alors la nullité de $E_{1,\xi}^{i,j}[\Pi^{\oo,v}]$ pour: $i+j \neq 0$, ou $|i| \geq s$ ou
$i \equiv s \mod 2$. Pour $|i|<s$ et $i=s-1-2r$, on a
\begin{multline*}
E_{\oo,\xi}^{s-1-2r,2r+1-s}[\Pi^{\oo,v}]=E_{1,\xi}^{s-1-2r,2r+1-s}[\Pi^{\oo,v}]=\\
\sharp \ker^1(\Qm,G_\tau) m(\Pi) [\overleftarrow{s-1}]_{\pi_v} \otimes \rec_{F_v}^\vee(\pi_v) (-\frac{s(g-1)+2r}{2})
\end{multline*}
d'où le résultat.

\end{proof}

\begin{prop} \label{prop-hipsi}
Pour tout $0 \leq i \leq d-1$, les $GL_d(F_v) \otimes W_v$-modules
$H^j(R^i\Psi_{\pi_v,\xi})[\Pi^{\oo,v}]$
vérifient les propriétés suivantes:
\begin{itemize}
\item[(1)] ils sont nuls si $g$ n'est pas un diviseur de $d$;

\item[(2)] pour $g$ un diviseur de $d=sg$, ils sont nuls si $j$ n'est pas de la forme $d-tg$ pour $1 \leq t
\leq s$;

\item[(3)] pour $g$ un diviseur de $d=sg$ et $j=d-tg$ avec $1 \leq t \leq s$, ils sont nuls si $i$ n'est pas
de la forme $tg-r$ avec $1 \leq r \leq t$;

\item[(4)] pour $g$ un diviseur de $d=sg$ et $1 \leq t < s$,
$H^{d-tg}(R^{tg-r}\Psi_{\pi_v})[\Pi^{\oo,v}]$ est isomorphe à
$$\sharp \ker^1(\Qm,G_\tau) m(\Pi) [\overleftarrow{t-r},
\overrightarrow{r-1}]_{\pi_v} \overrightarrow{\times} [\overleftarrow{s-t-1}]_{\pi_v} \otimes
\rec_{F_v}^\vee(\pi_v) (-\frac{s(g-1)+2(t-r)}{2})$$

\end{itemize}
\end{prop}

\begin{proof} On utilise la suite spectrale associée à la stratification
\begin{equation} \label{sss}
E_{1,,\xi,U^p(m),\pi_v}^{p,q;i}=H_c^{p+q}(\bar X_{U^p,m}^{(d-p+1)},R^i\Psi_{\pi_v,\xi}) \Rightarrow H_c^{p+q}(\bar
X_{U^p,m},R^i\Psi_{\pi_v,\xi})
\end{equation}
On rappelle que d'après le corollaire (\ref{coro-hic}), $E_{1,\xi,U^p(m),\pi_v}^{p,q;i}[\Pi^{\oo,v}]$ est non nul
si et seulement si:

\begin{itemize}
\item $g$ est un diviseur de $d=sg$,

\item $p-1=tg$ pour $1 \leq t \leq s$,

\item $p+q=d-tg$,

\item $i=tg-r$ avec $1 \leq r \leq t$.
\end{itemize}
Les points (1), (2) et (3) en découlent alors directement. Par ailleurs on a
\begin{multline*}
\lim_{\atop{\longto}{U^p(m)}} E_{1,\xi,U^p(m),\pi_v}^{tg+1,d-2tg-1;tg-r}[\Pi^{\oo,v}]\simeq \sharp
\ker^1(\Qm,G_\tau) m(\Pi) [\overleftarrow{t-r},
\overrightarrow{r-1}]_{\pi_v} \overrightarrow{\times} [\overleftarrow{s-t-1}]_{\pi_v} \\
\otimes \rec_{F_v}^\vee(\pi_v) (-\frac{s(g-1)-2(r-t)}{2})
\end{multline*}
de sorte que pour tout $k\geq 1$, les flèches $d_k^{p,q;i}:E_{k,\xi,\pi_v}^{p,q;i}[\Pi^{\oo,v}] \longto
E_{k,\xi,\pi_v}^{p+k,q+k-1}[\Pi^{\oo,v}]$ de (\ref{sss}) sont toutes nulles. En effet pour que
$E_{k,\xi,\pi_v}^{p,q;i}[\Pi^{\oo,v}]$ (resp. $E_{k,\xi,\pi_v}^{p+k,q+k-1;i}[\Pi^{\oo,v}]$) soit non nul, il faut
qu'il existe $1 \leq t_1 \leq s$ et $1 \leq r_1 \leq t_1$ (resp. $1 \leq t_2 \leq s$ et $1 \leq r_2 \leq
t_2$) avec
$$(p,q,i)=(t_1g+1,d-2t_1g-1,t_1g-r_1) ~(\hbox{resp. } (p+k,q+k-1,i)=(t_2g+1,d-2t_2g-1,t_2g-r_2)).$$
Ce qui donne $1=3g(t_2-t_1)$; on voit alors que pour tout $k\geq 1$, $E_{k,\xi,\pi_v}^{p,q;i}[\Pi^{\oo,v}]$ et
$E_{k,\xi,\pi_v}^{p+k,q+k-1}[\Pi^{\oo,v}]$ ne peuvent pas être tous deux non nuls de sorte que
$E_{\oo,\xi,U^p(m),\pi_v}^{p,q;i}[\Pi^{\oo,v}]=E_{1,\xi,U^p(m),\pi_v}^{p,q;i}[\Pi^{\oo,v}]$ d'où le résultat d'après
le
corollaire (\ref{coro-hic}).

\end{proof}

\begin{coro} \label{coro-1}
La partie de poids $d-s$ de $\UC_{F_v,l,d}^{d-i}(\JL^{-1}([\overleftarrow{s-1}]_{\pi_v^\vee}))$ est un
constituant de $\Pi_i \otimes \rec_{F_v}^\vee(\pi_v) (-\frac{d-s}{2})$ où
$$\Pi_i= \left \{ \begin{array}{cl} [\overleftarrow{s-1}]_{\pi_v} & \hbox{pour } i=1 \\
{[}\overrightarrow{0}]_{\pi_v} \overrightarrow{\times} [\overleftarrow{s-2}]_{\pi_v} & \hbox{pour } i=2 \\
\cdots \\
{[}\overrightarrow{i-2}]_{\pi_v} \overrightarrow{\times} [\overleftarrow{s-i}]_{\pi_v} & \hbox{pour } i \\
\cdots \\
{[}\overrightarrow{s-2}]_{\pi_v} \overrightarrow{\times} [\overleftarrow{0}]_{\pi_v} & \hbox{pour } i=s
\end{array} \right.$$
En outre $\UC_{F_v,l,d}^{d-s}(\JL^{-1}([\overleftarrow{s-1}]_{\pi_v^\vee}))$ contient
$[\overrightarrow{s-1}]_{\pi_v} \otimes \rec_{F_v}^\vee(\pi_v) (-\frac{s(g-1)}{2})$ et la partie de poids
$s(g-1)$ de $\sum_{i=1}^s (-1)^i \UC_{F_v,l,d}^{d-i}(\JL^{-1}([\overleftarrow{s-1}]_{\pi_v^\vee}))$ est égale
à
$(-1)^s [\overrightarrow{s-1}]_{\pi_v} \otimes \rec_{F_v}^\vee(\pi_v) (-\frac{s(g-1)}{2}).$

\end{coro}

\begin{proof} On étudie la $\Pi^{\oo,v}$-partie de la suite spectrale des cycles évanescents pour $\Pi$ vérifiant
les propriétés du début de ce paragraphe:
\begin{equation} \label{ssce}
E_{2,\xi,\pi_v}^{p,q}[\Pi^{\oo,v}]=\lim_{\atop{\longto}{U^p(m)}} H^p(\bar
X_{U^p,m},R^a\Psi_{\pi_v,\xi})[\Pi^{\oo,v}] \Rightarrow H^{p+q}_{\eta_v,\pi_v,\xi}[\Pi^{\oo,v}]
\end{equation}

\begin{lemm} \label{lemm-prec}
Pour tout $1 < r < s$, la partie de poids $s(g-1)$ de $E_{2,\xi,\pi_v}^{0,d-r}[\Pi^{\oo,v}]$
est un constituant de $E_{2,\xi,\pi_v}^{d-(r-1)g,(r-1)(g-1)}[\Pi^{\oo,v}]$, i.e. de
$\sharp \ker^1(\Qm,G_\tau) m(\Pi)[\overrightarrow{r-2}]_{\pi_v} \overrightarrow{\times} [\overleftarrow{s-r}]_{\pi_v}
\otimes \rec_{F_v}^\vee(\pi_v)
(-\frac{s(g-1)}{2});$ celle de $E_{2,\xi,\pi_v}^{0,d-1}[\Pi^{\oo,v}]$ est nulle.
\end{lemm}

\begin{proof} D'après la proposition précédente, les $E_{2,\xi,\pi_v}^{p,q}[\Pi^{\oo,v}]$ de poids $s(g-1)$ pour $p
\neq
0$, non nuls, sont
$$E_{2,\xi,\pi_v}^{d-rg,r(g-1)}[\Pi^{\oo,v}] \simeq \sharp \ker^1(\Qm,G_\tau) m(\Pi) [\overrightarrow{r-1}]_{\pi_v}
\overrightarrow{\times}
[\overleftarrow{s-r-1}]_{\pi_v} \otimes \rec_{F_v}^\vee(\pi_v) (-\frac{s(g-1)}{2})$$ Or, d'après le
corollaire (\ref{coho-global1}), pour $1 \leq r < s$, la partie de poids $s(g-1)$ de
$E_{\oo,\xi,\pi_v}^{d-r}[\Pi^{\oo,v}]$ est nulle d'où le résultat.

\end{proof}

\begin{lemm}
Pour $g >1$, le terme $E_{2,\xi,\pi_v}^{0,d-r}[\Pi^{\oo,v}]$ de la suite spectrale des cycles évanescents est
égal à
$\sharp \ker^1(\Qm,H_0) \hom_{D_{v,d}^\times}(\CC^\oo_{H_0,\xi},\UC^{d-r}_{F_v,l,d})[\Pi^{\oo,v}]$.
\end{lemm}

\begin{proof} On étudie la suite spectrale associée à la stratification pour le faisceau $R^{d-r}\Psi_{\pi_v}$.
D'après la proposition précédente, pour tout $1 \leq t < s$, $H^0_c(\bar
X_{U^p,m}^{(d-tg)},R^{d-r}\Psi_{\pi_v,\xi})[\Pi^{\oo,v}]$ est nul; pour $g \neq 1$, il en est de même de
$H^1_c(\bar X_{U^p,m}^{(d-tg)},R^{d-r}\Psi_{\pi_v,\xi})[\Pi^{\oo,v}]$, d'où le résultat.

\end{proof}

\begin{lemm} Pour $g=1$, $\hom_{D_{v,d}^\times}(\CC^\oo_{H_0,\xi},\UC^{d-r}_{F_v,l,d})[\Pi^{\oo,v}]$
est un constituant de
$$m(\Pi)[\overrightarrow{r-2}]_{\chi_v} \overrightarrow{\times} [\overleftarrow{d-r}]_{\chi_v} \otimes \chi_v^\vee.$$
\end{lemm}

\begin{proof} Pour $g=1$, on remarque que le seul $H^1_c(\bar X_{U^p,m}^{(d-h)},R^{d-r}\Psi_{\chi_v,\xi})$
ayant une partie de poids $0$ non nulle, est pour $h=d-1$ et $r=d$, de $\Pi^{\oo,v}$-composante isotypique
égale à $\sharp \ker^1(\Qm,G_\tau) m(\Pi)[\overrightarrow{d-2}]_{\chi_v} \overrightarrow{\times}
[\overleftarrow{0}]_{\chi_v} \otimes \chi_v^\vee$ ce
qui correspond à la contribution de $E_{2,\xi,\chi_v}^{1,0}[\Pi^{\oo,v}]$ dans le lemme (\ref{lemm-prec}), d'où le
résultat.

\end{proof}

Ainsi pour tout $1 \leq i \leq s$, la partie de poids $d-s$ de
$\UC_{F_v,l,d}^{d-i}(\JL^{-1}([\overleftarrow{s-1}]_{\pi_v^\vee}))$ est un constituant de $\Pi_i \otimes
\rec_{F_v}^\vee(\pi_v) (-\frac{d-s}{2})$ où $\Pi_i$ est comme dans l'énoncé. Par ailleurs on remarque que
$E_{2,\xi,\pi_v}^{d-s+1,(s-1)(g-1)}[\Pi^{\oo,v}]$ contient $[\overrightarrow{s-1}]_{\pi_v} \otimes
\rec_{F_v}^\vee(\pi_v) (-\frac{s(g-1)}{2})$ alors que $E_{2,\pi_v}^{d-s+2,(s-2)(g-1)}[\Pi^{\oo,v}]$ non; on
en déduit alors que $\UC_{F_v,l,d}^{d-i}(\JL^{-1}(\pi_v)^\vee)$ contient $[\overrightarrow{s-1}]_{\pi_v}
\otimes \rec_{F_v}^\vee(\pi_v) (-\frac{s(g-1)}{2})$. Le calcul sur la somme alternée correspond à
(\ref{inutile}).

\end{proof}

\subsection{Involution de Zelevinski et première preuve de (\ref{cas-m})}

Il est possible de prouver la proposition (\ref{cas-m}) sans plus d'étude cohomologique en utilisant le résultat
suivant.

\begin{theo} \label{involution}
Pour toute représentation irréductible cuspidale $\pi_v$ de $GL_g(F_v)$, on a, pour tout $s \geq 1$ et pour tout $i$,
un isomorphisme canonique
$$\left ( \UC_{F_v,l,sg}^{sg-1+\bullet}(\JL^{-1}(\st_s(\pi_v))) \right )^{\vee,\iota}(d-1) \simeq
\UC_{F_v,l,sg}^{sg+s-2-2 \bullet}(\JL^{-1}(\st_s(\pi_v^\vee)))$$
où dans le membre de droite l'exposant $\vee,\iota$ désigne le dual composé avec l'involution
de Zelevinski, $\iota$, sur $GL_{sg}(F_v)$.
\end{theo}

\rem Laurent Fargues a une preuve de ce résultat qui utilise, tout ou partie, l'isomorphisme de Faltings, à partir
d'un résultat similaire du coté de
l'espace de Drinfeld.

Ainsi avec les notations du corollaire (\ref{coro-1}), s'il existait $1 \leq i \leq s$ tel que $\Pi_i$ admette un
constituant autre que
$[\overrightarrow{s-1}]_{\pi_v}$, on en déduirait que
$\UC_{F_v,l,sg}^{sg-s-1+i}(\JL^{-1}([\overleftarrow{s-1}]_{\pi_v}))$ aurait un constituant de poids
$sg+s-2$ autre que $[\overleftarrow{s-1}]_{\pi_v^\vee}$. Or, d'après les corollaires (\ref{coro-sec-fp}) et
(\ref{coro-pnul}), pour tout $k$, les germes
aux points supersinguliers des $h^i gr_{k,\pi_v^\vee}$ sont tous de poids
strictement plus petit que $sg+s-2$ ou alors égales à $[\overleftarrow{s-1}]_{\pi_v^\vee} \otimes \rec_{F_v}(\pi_v)
(-\frac{sg+s-2}{2})$, d'où le résultat.

\subsection{Cas  $\Pi_v \simeq \speh_s(\pi_v)$}

Nous allons prouver la proposition (\ref{cas-m}) sans utiliser le théorème (\ref{involution}). Pour cela nous
revenons
à l'étude de la suite spectrale
des cycles évanescents.

\begin{figure}[!ht]

\setlength{\unitlength}{.5cm}
\begin{picture}(19,12)(-5,-2)
\linethickness{.1pt}

\put(0,2){\vector(0,1){5}}

\dashline{.2}(0,0)(0,2)

\put(0,0){\vector(1,0){6}}

\put(0,2){\line(1,0){5}}

\put(7,2){\makebox(0,0){$j=g-1$}}

\put(0,5){\circle{.3}}

\put(-2,5){\makebox(0,0){$\genfrac{}{}{0pt}{}{[\overleftarrow{0}]_{\pi_v} \overrightarrow{\times}
[\overrightarrow{0}]_{\pi_v}}{\otimes \rec_{F_v}^\vee(\pi_v)|-|^{-2(g-1)/2}}$}}

\put(0,6){\circle{.3}}

\put(2,6){\vector(-1,0){2}}

\put(2,6){$i+j=d-1$}

\put(2,8){\vector(-1,-1){2}}

\put(2,9){\makebox(0,0){$\genfrac{}{}{0pt}{}{[\overleftarrow{1}]_{\pi_v} \otimes \rec_{F_v}^\vee(\pi_v)}{ (
|-|^{-2(g+1)/2} \oplus |-|^{-2(g-1)/2})}$}}

\put(0,6){\line(1,-1){4}}

\put(4,2){\circle{.3}}

\put(5,3){\vector(-1,-1){1}}

\put(5,4){\makebox(0,0){$\genfrac{}{}{0pt}{}{[\overleftarrow{0}]_{\pi_v} \overrightarrow{\times}
[\overrightarrow{0}]_{\pi_v}}{\otimes \rec_{F_v}^\vee(\pi_v)|-|^{-2(g-1)/2}}$}}

\put(0,5){\vector(4,-3){4}} \dashline{.2}(4,0)(4,2)

\put(4,-1){\vector(0,1){1}}

\put(4,-1.5){\makebox(0,0){$i=d-g$}}

\put(0,-.5){\makebox(0,0){$(0,0)$}}

\multiput(0,2)(1,0){5}{\circle*{.1}}

\multiput(0,2)(0,1){5}{\circle*{.1}}

\put(12,2){\vector(0,1){5}}
\dashline{.2}(12,0)(12,2)

\put(12,0){\vector(1,0){6}}

\put(12,2){\line(1,0){6}}

\put(20,2){\makebox(0,0){$j=g-1$}}

\put(12,5){\circle{.3}}

\put(10,5){\makebox(0,0){$\genfrac{}{}{0pt}{}{[\overleftarrow{0}]_{\pi_v} \overrightarrow{\times}
[\overrightarrow{0}]_{\pi_v}}{\otimes \rec_{F_v}^\vee(\pi_v) |-|^{-2(g-1)/2}}$}}

\put(12,6){\circle{.3}}

\put(14,6){\vector(-1,0){2}}

\put(14,6){$i+j=d-1$}

\put(14,8){\vector(-1,-1){2}}

\put(14,9){\makebox(0,0){$\genfrac{}{}{0pt}{}{[\overleftarrow{1}]_{\pi_v} \otimes \rec_{F_v}^\vee(\pi_v)}
{(|-|^{-2(g-1)/2} \oplus |-|^{-2(g+1)/2})}$}}

\put(12,6){\line(1,-1){4}}

\put(17,2){\circle{.3}}

\put(18,3){\vector(-1,-1){1}}

\put(18,4){\makebox(0,0){$\genfrac{}{}{0pt}{}{[\overleftarrow{0}]_{\pi_v} \overleftarrow{\times}
[\overrightarrow{0}]_{\pi_v}}{\otimes \rec_{F_v}^\vee(\pi_v) |-|^{-2(g+1)/2}}$}}

\put(12,6){\line(5,-4){5}}

\put(17,2){\vector(1,-1){0}}

\dashline{.2}(17,0)(17,2)

\put(17,-1){\vector(0,1){1}}

\put(17,-1.5){\makebox(0,0){$i=d-g+1$}}

\put(12,-.5){\makebox(0,0){$(0,0)$}}

\multiput(12,2)(1,0){5}{\circle*{.1}}

\multiput(12,2)(0,1){5}{\circle*{.1}}

\end{picture}
\caption{\label{figu1} $E_2^{p,q}[\Pi^{\oo,v}]$ de (\ref{ssce}) pour $\Pi_v \simeq \st_2(\pi_v)$ et $\Pi_v \simeq
\speh(\pi_v)$}
\end{figure}
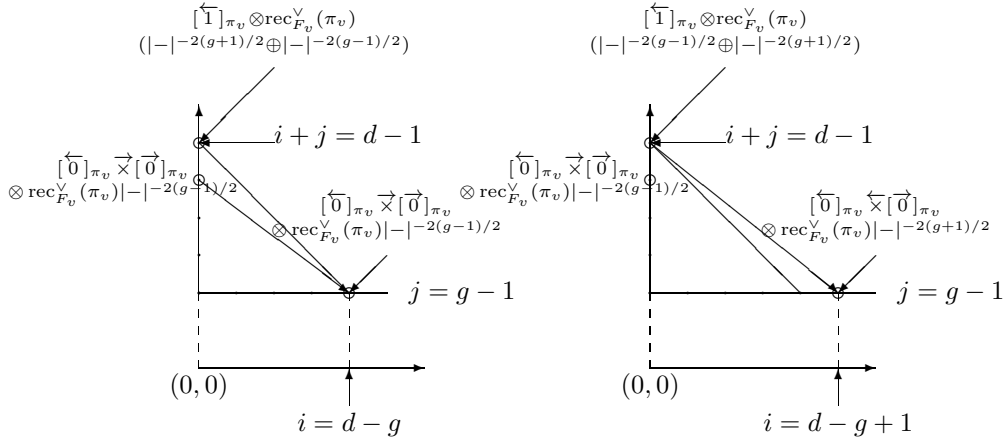

\marque On peut remarquer que la connaissance des $E_{2,\chi,\pi_v}^{p,q}[\Pi^{\oo,v}]$ de la suite spectrale
(\ref{ssce}) ne nous fournit pas le théorème local. La partie gauche de la figure (\ref{figu1}) illustre ce fait dans
le cas
$s=2$, où pour $\pi_v$ une représentation irréductible cuspidale de $GL_{d/2}(F_v)$, l'on n'arrive pas à
exclure le cas:
\begin{itemize}

\item $\UC_{F_v,l,d}^{d-1}(\JL^{-1}([\overleftarrow{1}]_{\pi_v^\vee}))=[\overleftarrow{1}]_{\pi_v} \otimes
(\rec_{F_v}^\vee(\pi_v)(-\frac{d-1}{2}) \otimes \sp_2)$;

\item $\UC_{F_v,l,d}^{d-2}(\JL^{-1}([\overleftarrow{1}]_{\pi_v^\vee}))= [\overleftarrow{0}]_{\pi_v}
\overrightarrow{\times} [\overleftarrow{0}]_{\pi_v} \otimes \rec_{F_v}^\vee(\pi_v) (-\frac{d-2}{2})$.
\end{itemize}

\noindent Si la flèche indiquée est un isomorphisme, on obtient bien le bon aboutissement.

\marque Dans le cas où l'on considère une représentation automorphe $\Pi$ de $G_\tau(\Am)$ vérifiant
$\hyp(\xi)$ et telle que $\Pi_v \simeq [\overrightarrow{s-1}]_{\pi_v}$, on calculera les termes
$E_{2,\xi}^{p,q}[\Pi^{\oo,v}]$ de la suite spectrale des cycles évanescents et on montrera alors que le
théorème local en découle. Dans la partie de droite de la figure (\ref{figu1}), on illustre comment on exclut
le cas défavorable ci-avant, en remarquant que $[\overleftarrow{0}]_{\pi_v} \overrightarrow{\times}
[\overleftarrow{0}]_{\pi_v}$ n'est pas isomorphe à $[\overleftarrow{0}]_{\pi_v} \overleftarrow{\times}
[\overleftarrow{0}]_{\pi_v}$, contredisant le théorème de Lefschetz difficile. On pourra aussi se référer aux
figures (\ref{figure5}) et (\ref{figure6}) qui détaillent pour $s=4$ et $g=2$, les
$E_{2,\pi_v}^{p,q}[\Pi^{\oo,v}]$ respectivement dans les cas $\Pi_v=[\overleftarrow{s-1}]_{\pi_v}$ et
$\Pi_v=[\overrightarrow{s-1}]_{\pi_v}$.

\marque On considère \footnote{cf. aussi les conditions qui précèdent le corollaire (\ref{coro-pnul})} dans
la suite $\Pi$ une représentation irréductible automorphe de $G_\tau(\Am)$ vérifiant $\hyp(\xi)$ et telle que
$\Pi_v \simeq [\overrightarrow{s-1}]_{\pi_v}$, pour $\pi_v$ une représentation irréductible cuspidale de
$GL_g(F_v)$. On suppose par ailleurs que l'ensemble des représentations irréductibles automorphes $\bar \Pi$
de $H_0(\Am)$ telles que $\bar \Pi^{\oo,v} \simeq \Pi^{\oo,v}$ est réduite à une unique représentation, en
particulier $m(\Pi)=m(\bar \Pi)$, et $\bar \Pi_v \simeq \JL^{-1}([\overleftarrow{s-1}]_{\pi_v})$. L'objectif
est de calculer les termes $E_2^{p,q}[\Pi^{\oo,v}]$ de la suite spectrale des cycles évanescents.

\begin{prop} \label{prop-not}
Pour tout $1 \leq t \leq s$, $H^i(j^{\geq tg}_{!*}
HT_\xi(g,t,\pi_v,[\overleftarrow{t-1}]_{\pi_v})[d-tg])[\Pi^{\oo,v}]$ est nul pour $|i|
> s-t$ ou $i \not \equiv s-t \mod 2$ et sinon, il est isomorphe à
\begin{multline*}
\sharp \ker^1(\Qm,G_\tau) m(\Pi)([\overleftarrow{t-1}]_{\pi_v} \overrightarrow{\times}
[\overrightarrow{\frac{s-t-i}{2}-1}]_{\pi_v}) \overleftarrow{\times}
[\overrightarrow{\frac{s-t+i}{2}-1}]_{\pi_v}  \otimes (\Xi^{\frac{(s-t)g+i}{2}} \otimes \bigoplus_{\chi \in
\AF(\pi_v)} \chi^{-1})
\end{multline*}
en tant que représentation de $GL_d(F_v) \times \Zm$, où
$\AF(\pi_v)$ est l'ensemble des caractères $\chi:\Zm \longto \Qm_l^\times$, tels que $\pi_v \otimes \chi \circ
\val(\det) \simeq \pi_v$.
\end{prop}

\begin{proof} On raisonne par récurrence pour $t$ variant de $s$ à $1$, le cas $t=s$ est donné par le lemme
(\ref{lem-pts-ss}) et les conditions imposées ci-dessus à $\Pi$. Les suites exactes courtes
\begin{multline*}
0 \to P_{\pi_v,t,0}([\overleftarrow{t-1}]_{\pi_v}) \longto j^{\geq tg}_! HT(g,t,\pi_v,
[\overleftarrow{t-1}]_{\pi_v})[d-tg]  \longto \\ j^{\geq tg}_{!*}
HT(g,t,\pi_v,[\overleftarrow{t-1}]_{\pi_v})[d-tg] \to 0
\end{multline*}
$$0 \to P_{\pi_v,t,i}([\overleftarrow{t-1}]_{\pi_v}) \longto P_{\pi_v,t,i-1}([\overleftarrow{t-1}]_{\pi_v})
\longto \PC_-(g,t+i,\pi_v,t,[\overleftarrow{t-1}]_{\pi_v}) \to 0$$ pour $1 \leq i \leq s-t-1$, fournissent
\begin{multline} \label{egalite-purete}
\sum_i (-1)^i H^i(j_{!*}^{\geq tg} HT_{\xi}(g,t,\pi_v,[\overleftarrow{t-1}]_{\pi_v})[d-tg])[\Pi^{\oo,v}]=
\sum_{r=1}^{s-t} (-1)^r \sum_i (-1)^i  H^i([\overleftarrow{t-1}]_{\pi_v} \\ \overrightarrow{\times}
[\overleftarrow{r-1}]_{\pi_v})[\Pi^{\oo,v}] \otimes \Xi^{r/2}  +
\sum_i (-1)^i H^i(j_!^{\geq tg} HT_{\xi}(g,t,\pi_v,[\overleftarrow{t-1}]_{\pi_v})[d-tg])[\Pi^{\oo,v}]
\end{multline}
où par simplification, on écrit $H^i([\overleftarrow{t-1}]_{\pi_v} \overrightarrow{\times}
[\overleftarrow{r-1}]_{\pi_v})$ pour
$$\lim_{\atop{\to}{U^p(m)}} H^i(\bar X_{U^p,m}^{[d-(t+r)g]}, j^{\geq (t+r)g}_{!*} HT_{\xi}(g,t+r,\pi_v,
[\overleftarrow{t-1}]_{\pi_v} \overrightarrow{\times} [\overleftarrow{r-1}]_{\pi_v})[d-(t+r)g])$$

D'après l'hypothèse de récurrence et (\ref{somme-alternee}), pour tout $i \equiv 1 \mod 2$, la partie de
poids $(s-t)(g-1)+i$ du membre de droite de (\ref{egalite-purete}) est nulle; il en est de même pour
$i>2(s-t)$ ou $i<0$, tandis que pour $0 \leq i=2k < 2(s-t)$ (resp. $i=2(s-t)$) celle-ci est égale à
$$\sum_{r=1}^{s-t-k} H^{t-s+2k+r}([\overleftarrow{t-1}]_{\pi_v} \overrightarrow{\times}
[\overleftarrow{r-1}]_{\pi_v})[\Pi^{\oo,v}] \otimes \Xi^{r/2}$$ (resp. à la partie de poids $(s-t)(g+1)$ de
(\ref{somme-alternee})), ce qui donne
$$\sharp \ker^1(\Qm,G_\tau) m(\Pi) ~\pi \otimes (\Xi^{\frac{s(g-1)+2k}{2}} \otimes \bigoplus_{\chi \in \AF(\pi_v)}
\chi^{-1})$$
$$\begin{array}{rl} \pi = & \bigl ( [\overleftarrow{t-1}]_{\pi_v} \overrightarrow{\times} ((-1)^{s-t-k}
[\overleftarrow{s-t-k-1}]_{\pi_v} \\ & +\sum_{r=1}^{s-t-k-1} (-1)^{r-1} [\overleftarrow{r-1}]_{\pi_v}
\overrightarrow{\times} [\overrightarrow{s-t-r-k-1}]_{\pi_v}) \bigr )
\overleftarrow{\times} [\overrightarrow{k-1}]_{\pi_v} \\
= & ([\overleftarrow{t-1}]_{\pi_v} \overrightarrow{\times} [\overrightarrow{s-t-k-1}]_{\pi_v})
\overleftarrow{\times} [\overrightarrow{k-1}]_{\pi_v} \\
& (\hbox{resp. } [\overleftarrow{t-1}]_{\pi_v} \overleftarrow{\times} [\overrightarrow{s-t-1}]_{\pi_v}).
\end{array}$$
Or $H^i(j_{!*}^{\geq tg} HT_{\xi}(g,t,\pi_v,[\overleftarrow{t-1}]_{\pi_v}) [d-tg])$ est pur de poids
$(s-t)g+i$ de sorte que son semi-simplifié est égal à celui de l'énoncé. On conclut alors à l'égalité des
représentations, et pas seulement de leur semi-simplifiée, en remarquant que les strates étant induites,
l'espace précédent, en tant que représentation de $GL_d(F_v)$ est de la forme
$\ind_{P_{tg,d}(F_v)}^{GL_d(F_v)} [\overleftarrow{t-1}]_{\pi_v((s-t)g-i)/2} \otimes \pi'$
pour une certaine représentation $\pi'$ de $GL_{d-tg}(F_v)$.
\end{proof}

\begin{coro} \label{coro-ssce-min}
Pour tout $i \neq d-s$, la partie de poids $s(g-1)$ de $H^i_{\eta_v,\xi}[\Pi^{\oo,v}]$ est nulle alors que
pour $i=d-s$ elle est égale à $\sharp \ker^1(\Qm,G_\tau) m(\Pi)[\overrightarrow{s-1}]_{\pi_v} \otimes
\rec_{F_v}^\vee(\pi_v)
(-\frac{s(g-1)}{2})$.
\end{coro}

\begin{rema} \label{rema-N} En fait on peut à ce stade déterminer complètement les $H^i_{\eta_v,\xi}$, cependant
comme
on l'a
déjà remarqué seule la connaissance des parties de poids $s(g-1)$ nous est nécessaire. Pour un énoncé
complet, on pourra voir la proposition (\ref{prop-lrs2}).
\end{rema}

\begin{proof} On écrit $i$ sous la forme $d-1-\d$ et on étudie la suite spectrale
\begin{equation} \label{ss-poids}
E_{1,\xi}^{i,j}:=H^{i+j}(gr_{-i,\xi}) \Rightarrow H^{d-1+i+j}_{\eta_v,\xi}
\end{equation}
qui, d'après la pureté, dégénère en $E_2$. On pourra se référer à la figure (\ref{figure8}) où l'on a
représenté les $H^i(gr_{k,\xi})[\Pi^{\oo,v}]$ pour $s=4$. D'après les propositions (\ref{prop-p}) (ii) et
(\ref{prop-grk-final}), on a
$$e_{\pi_v} H^i(gr_{k,\xi})[\Pi^{\oo,v}]=\bigoplus_{\atop{|k|< t \leq s}{t \equiv k+1 \mod 2}} H^i(\PC(g,t,\pi_v)
(-\frac{tg-1+k}{2}))[\Pi^{\oo,v}]$$ Ainsi d'après la proposition (\ref{prop-not}), pour $\d>0$, la partie de
poids $s(g-1)$ de $H^{d-1+\d}_{\eta_v,\xi}[\Pi^{\oo,v}]$ est nulle; précisément pour $\d >0$, les poids de
$H^{d-1+\d}_{\eta_v,\xi}[\Pi^{\oo,v}]$ sont parmi les $k=s(g-1)+2 \d +2r$ avec $0 \leq r < s-\d$ et sa partie
de poids $s(g-1)+2\d$, que l'on notera avec un indice, est un quotient de
\begin{multline*}
E_{1,\xi}^{1-s+\d,s-1}[\Pi^{\oo,v}]_{s(g-1)+2\d}=\frac{1}{e_{\pi_v}} H^{\d}(\PC(g,s-\d,\pi_v)
(-\frac{(s-\d)(g-1)}{2})) \\
=\sharp \ker^1(\Qm,G_\tau) m(\Pi)[\overleftarrow{s-\d-1}]_{\pi_v} \overleftarrow{\times}
[\overrightarrow{\d-1}]_{\pi_v} \otimes
\rec_{F_v}^\vee(\pi_v) (-\frac{s(g-1)+2\d}{2})
\end{multline*}

Par ailleurs, pour $\d>0$, d'après (\ref{prop-not}), la partie de poids $s(g-1)$ de
$H^{d-s-\d}_{\eta_v,\xi}[\Pi^{\oo,v}]$ est un sous-quotient de
\begin{multline*}
E_{1,\xi}^{1-s+\d,s-1-2\d}[\Pi^{\oo,v}]_{s(g-1)}=\frac{1}{e_{\pi_v}} H^{-\d}(\PC(g,s-\d,\pi_v)
(-\frac{(s-\d)(g-1)}{2})) \\
=\sharp \ker^1(\Qm,G_\tau) m(\Pi)[\overleftarrow{s-\d-1}]_{\pi_v} \overrightarrow{\times}
[\overrightarrow{\d-1}]_{\pi_v} \otimes
\rec_{F_v}^\vee(\pi_v) (-\frac{s(g-1)}{2})
\end{multline*}

On utilise alors le théorème de Lefschetz difficile dont on rappelle l'énoncé ci-après.

\begin{theo} Le fibré canonique sur la fibre générique de $X_{U^p,m}$ est ample et équivariant pour l'action
de $G(\Am^\oo)$ et $W_v$; il induit alors une classe $h \in H^2_{\eta_v,\xi} (1)$,
telle que les applications itérées du cup produit
$h^i: H^{d-1-i}_{\eta_v,\xi} \longto H^{d-1+i}_{\eta_v,\xi} (i)$
sont des isomorphismes.
\end{theo}

\marque Ainsi pour $s>2$, on observe que si la partie de poids $s(g-1)$ de
$H^{d-1-\d}_{\eta_v}[\Pi^{\oo,v}]$, pour $1 < \d <s-1$, est non nulle, ses constituants sont de la forme
$$\sharp \ker^1(\Qm,G_\tau) m(\Pi)[\overleftrightarrow{s-2},\overrightarrow{1}]_{\pi_v} \otimes
\rec_{F_v}^\vee(\pi_v)
(-\frac{s(g-1)}{2})$$
alors qu'un constituant non nul de poids $s(g-1)+2\d$ de $H^{d-1+2\d}_{\eta_v,\xi}[\Pi^{\oo,v}]$, s'il
existe, est de la forme
$$\sharp \ker^1(\Qm,G_\tau) m(\Pi)[\overleftrightarrow{s-2},\overleftarrow{1}]_{\pi_v} \otimes \rec_{F_v}^\vee(\pi_v)
(-\frac{s(g-1)+2\d}{2}).$$
La contradiction découle alors du théorème de Lefschetz difficile et de
l'observation de l'orientation de la dernière flèche.

\marque Pour $\d=1$ et $s>2$, la partie de poids $s(g-1)+2$ de $H^{d}_{\eta_v,\xi}[\Pi^{\oo,v}]$ est un
quotient de $\sharp \ker^1(\Qm,G_\tau) m(\Pi)[\overleftarrow{s-2}]_{\pi_v} \overleftarrow{\times}
[\overrightarrow{0}]_{\pi_v} \otimes \rec_{F_v}^\vee(\pi_v) (-\frac{s(g-1)+2}{2})$ alors que la partie de
poids $s(g-1)$ de $H^{d-2}_{\eta_v,\xi}[\Pi^{\oo,v}]$ est un constituant de $\sharp \ker^1(\Qm,G_\tau) m(\Pi)
[\overleftarrow{s-2}]_{\pi_v} \overrightarrow{\times} [\overrightarrow{0}]_{\pi_v} \otimes
\rec_{F_v}^\vee(\pi_v) (-\frac{s(g-1)}{2}).$ La contradiction découle alors du théorème de Lefschetz
difficile et du fait que $[\overleftarrow{s-1}]_{\pi_v}$ n'est pas un quotient de
$[\overleftarrow{s-2}]_{\pi_v} \overleftarrow{\times} [\overrightarrow{0}]_{\pi_v}$.

\marque Pour $\d=s-1$ et $s \geq 2$, on observe que
$E_{2,\xi}^{d-s}[\Pi^{\oo,v}]=H^{d-s}_{\eta_v,\xi}[\Pi^{\oo,v}]$ est un sous-espace de
$E_{1,\xi}^{0,1-s}[\Pi^{\oo,v}]$ qui est égal à
\begin{multline*}
\frac{1}{e_{\pi_v}} H^{1-s}(\PC(g,1,\pi_v)) =\sharp \ker^1(\Qm,G_\tau) m(\Pi)[\overleftarrow{0}]_{\pi_v}
\overrightarrow{\times} [\overrightarrow{s-2}]_{\pi_v} \otimes \rec_{F_v}^\vee(\pi_v) (-\frac{s(g-1)}{2})
\end{multline*}
alors que $E_{2,\xi}^{d-s}[\Pi^{\oo,v}]=H^{d+s-2}_{\eta_v,\xi}[\Pi^{\oo,v}]$ est un quotient de
$E_{1,\xi}^{0,s-1}[\Pi^{\oo,v}]$ qui est égal à
\begin{multline*}
\frac{1}{e_{\pi_v}} H^{s-1}(\PC(g,1,\pi_v)) =\sharp \ker^1(\Qm,G_\tau) m(\Pi)[\overleftarrow{0}]_{\pi_v}
\overleftarrow{\times} [\overrightarrow{s-2}]_{\pi_v} \otimes \rec_{F_v}^\vee(\pi_v) (-\frac{s(g+1)-2}{2})
\end{multline*}
Ainsi d'après Lefschetz difficile, s'ils sont non nuls, ils
doivent être égaux à $[\overrightarrow{s-1}]_{\pi_v}$. Par ailleurs on remarque que ce dernier n'est pas un
constituant de $H^{d-s+1}_{\eta_v,\xi}[\Pi^{\oo,v}]$, ni de $H^{d+s-3}_{\eta_v,\xi}[\Pi^{\oo,v}]$, de sorte
que $\sharp \ker^1(\Qm,G_\tau) m(\Pi)[\overrightarrow{s-1}]_{\pi_v} \otimes \rec_{F_v}^\vee(\pi_v)
(-\frac{s(g-1)}{2})$ est effectivement un
constituant de $H^{d-s}_{\eta_v,\xi}$, d'où le résultat.

\end{proof}

\begin{prop} \label{prop-ssce-poids} Pour tout $p \neq 0$ et tout $q$, sont nulles, les parties de poids $s(g-1)$ des
$E_{2,\xi}^{p,q}[\Pi^{\oo,v}]$ de la suite spectrale des cycles évanescents
\begin{equation} \label{ssce2}
E_{2,\xi}^{p,q}=\lim_{\atop{\longto}{U^p(m)}} H^p(\bar X_{U^p,m},R^q\Psi_v \otimes \LC_{\xi}) \Rightarrow
H^{p+q}_{\eta_v,\xi}
\end{equation}
\end{prop}

\begin{proof} Le principe est d'étudier les $E_{2,\xi}^{p,q}[\Pi^{\oo,v}]$ via les suites spectrales associées à la
stratification
\begin{equation} \label{sss2}
E_{1,U^p(m),\xi}^{p,q;i}=H_c^{p+q}(\bar X_{U^p,m}^{(d-p+1)},R^i\Psi \otimes \LC_{\xi}) \Rightarrow
H_c^{p+q}(\bar X_{U^p,m},R^i\Psi \otimes \LC_{\xi})
\end{equation}
Ainsi le résultat découle simplement de la proposition suivante.

\begin{prop} \label{prop-hic-poids}
Pour tout $1 \leq t < s$ et pour tout $i$, la partie de poids $(s-t)(g-1)$ de $H^i_c(\bar
X_{U^p,m}^{(d-tg)},\FC(g,t,\pi_v) \otimes \LC_{\xi})[\Pi^{\oo,v}]$ est nulle.
\end{prop}

\begin{proof} Il s'agit dans un premier temps d'étudier les parties de poids $(s-t)(g-1)$ des
$H^i_c(\bar X_{U^p,m}^{(d-tg)},HT_{\xi}(g,t,\pi_v,\Pi_t))[\Pi^{\oo,v}]$ pour une représentation $\Pi_t$
quelconque de $GL_{tg}(F_v)$.

\begin{lemm} \label{lem-rj-combi}
Pour tout $1 \leq t < s$, les
${\DS \lim_{\atop{\to}{U^p(m)}}} H^i(\bar X_{U^p,m}^{[d-tg]},(j^{\geq tg}_! HT_{\xi}(g,t,\pi_v,\Pi_t)[(s-t)g]
))[\Pi^{\oo,v}]$ vérifient les propriétés suivantes:

\begin{itemize}
\item[(i)] ils sont de la forme $\bigoplus_\chi(\ind_{P_{tg,d}(F_v)}^{GL_d(F_v)} (\Pi_t \otimes \chi) \otimes
\pi_\chi)
\otimes \chi$, en tant que représentation de $GL_d(F_v) \times \Zm$, où $\chi$ décrit les caractères $\Zm
\longto \bar \Qm_l^\times$ et $\pi_\chi$ est une représentation de $GL_{(s-t)g}(F_v)$;

\item[(ii)] ils sont nuls pour $|i| \geq s-t+1$;

\item[(iii)] ils sont en général mixtes de poids $(s-t)(g+1)-2(k-1)$ pour $1 \leq k \leq s-t+1$ vérifiant
$s-t-2(k-1) \leq i \leq s-t-(k-1)$;

\item[(iv)] soit $t-s \leq i_0 < 0$, le plus petit indice $i$ tel que la partie de poids
$(s-t)(g-1)$ de ${\DS \lim_{\atop{\to}{U^p(m)}}} H^i(\bar X_{U^p,m}^{[d-tg]},j^{\geq tg}_!
HT_{\xi}(g,t,\pi_v,\Pi_t)[(s-t)g]) [(s-t)g])[\Pi^{\oo,v}]$ soit non nulle\footnote{Si un tel $i_0$, n'existe
pas l'énoncé est vide.}. Cette dernière est alors égale à
$$\sharp \ker^1(\Qm,G_\tau) m(\Pi)\Pi_t \overrightarrow{\times}
[\overleftarrow{i_0-t+s-1},\overrightarrow{-i_0}]_{\pi_v} \otimes (\Xi^{\frac{(s-t)(g-1)}{2}} \otimes
\bigoplus_{\chi \in \AF(\pi_v)} \chi^{-1}).$$
\end{itemize}
\end{lemm}

\begin{proof} Le point (i) correspond à la proposition (\ref{prop-poids}) qui découle directement de l'action de
$GL_{tg}(F_v)$ sur la strate via $\val(\det)$ et du fait que les strates non supersingulières sont induites.

Pour les points (ii)-(iii), on considère les suites exactes
\begin{equation} \label{secfpi}
0 \to P_{\pi_v,t,i+1}(\Pi_t) \longto P_{\pi_v,t,i}(\Pi_t) \longto \PC_-(g,t+i+1,\pi_v,t,\Pi_t) \to 0
\end{equation}
et on reprend les notations simplifiées de la preuve de la proposition (\ref{prop-not}) en notant
$P_{\pi_v,t,-1,\xi}(\Pi_t):=j^{\geq tg}_! HT_{\xi}(g,t,\pi_v,\Pi_t)[d-tg]$. Par facilité on notera
$P_{\pi_v,t,i}(\Pi_t)$ pour $P_{\pi_v,t,i}(\Pi_t) \otimes \LC_{\xi}$. Le résultat découle alors directement
du cas $i=-1$ dans le lemme suivant.

\begin{lemm} \label{lem-utile}
Pour tout $-1 \leq i \leq s-t-1$, les
${\DS \lim_{\atop{\to}{U^p(m)}}} H^j(\bar X_{U^p,m}^{[d-(t+i+1)g]},P_{\pi_v,t,i}(\Pi_t) \otimes
\LC_\xi)[\Pi^{\oo,v}]$
vérifient les points suivants:
\begin{itemize}

\item[(1)] ils sont nuls pour $|j| \geq s-t-i$;

\item[(2)] pour $|j| < s-t-i$, ils sont mixtes de poids $(s-t)(g+1)-2(k+i)$ pour $1 \leq k \leq s-t-i$ vérifiant
$s-t-i-1-2(k-1) \leq j \leq s-t-i-1-(k-1)$;

\end{itemize}
\end{lemm}

\begin{proof} (1)-(2) On raisonne par récurrence descendante, le cas $i=s-t-1$ étant évident car
$P_{\pi_v,t,s-t-1}(\Pi_t)$ est le faisceau concentré aux points supersinguliers
$$\FC(g,s,\pi_v) \otimes \Pi_t((s-t)/2) \times [\overleftarrow{s-t-1}]_{\pi_v(-t/2)} \otimes
|\art_{F_v}^{-1}|^{-\frac{(s-t)(g-1)}{2}}$$ Supposons donc le résultat acquis jusqu'au rang $i+1$ et traitons
le cas de $P_{\pi_v,t,i}(\Pi_t)$. On considère la suite exacte longue de cohomologie associée à la suite
exacte courte (\ref{secfpi}). D'après l'hypothèse de récurrence les $H^j(P_{\pi_v,t,i+1}(\Pi_t))$ sont nuls
pour $|j| \geq s-t-i-1$ et d'après la proposition (\ref{prop-not}), les $H^j(\PC_-(g,t+i+1,\pi_v,t,\Pi_t))$
sont nuls pour $|j| \geq s-t-i$, d'où le point (1). En ce qui concerne le point (2), la suite exacte longue
en question s'écrit alors
\begin{multline} \label{sel1}
0 \to H^{-(s-t-i-1)}(P_{\pi_v,t,i}(\Pi_t)) \longto H^{-(s-t-i-1)}(j^{\geq (t+i+1)g}_{!*}) \longto
H^{-(s-t-i-2)}(P_{\pi_v,t,i+1}(\Pi_t))  \longto \\ H^{-(s-t-i-2)} P_{\pi_v,t,i}(\Pi_t)  \longto 0 \cdots
\cdots 0 \to H^{(s-t-i-1)-2r}(P_{\pi_v,t,i+1}(\Pi_t)) \longto \\ H^{(s-t-i-1)-2r}(P_{\pi_v,t,i}(\Pi_t)) \longto
H^{(s-t-i-1)-2r}(j^{\geq (t+i+1)g}_{!*}) \longto
 H^{(s-t-i-1)-2r+1}(P_{\pi_v,t,i+1}(\Pi_t)) \\ \longto H^{(s-t-i-1)-2r+1}(P_{\pi_v,t,i}(\Pi_t)) \to 0
\cdots \\
0 \to H^{(s-t-i-1)-2}(P_{\pi_v,t,i+1}(\Pi_t)) \longto H^{(s-t-i-1)-2}(P_{\pi_v,t,i}(\Pi_t)) \longto \\
H^{(s-t-i-1)-2}(j^{\geq (t+i+1)g}_{!*})  \longto
H^{s-t-i-2}(P_{\pi_v,t,i+1}(\Pi_t))   \longto H^{s-t-i-2}(P_{\pi_v,t,i}(\Pi_t)) \to 0 \\
0 \to H^{s-t-i-1}(P_{\pi_v,t,i}(\Pi_t)) \longto H^{s-t-i-1}(j^{\geq (t+i+1)g}_{!*}) \to 0
\end{multline}
On rappelle que $H^j(j^{\geq (t+i+1)g}_{!*})$ est pur de poids $(s-t)g-(i+1)+j$ de sorte qu'en utilisant
l'hypothèse de récurrence, les $H^j(P_{\pi_v,t,i}(\Pi_t))$ sont de poids $(s-t)(g+1)-2(i+1+k-1)$ avec $1 \leq
k \leq s-t-i$.

Soit alors $j$ de la forme $s-t-1-i-(2r+1)$ avec $0 \leq 2r+1 \leq 2(s-t-1-i)$; la suite exacte longue
ci-dessus montre alors que les poids de $H^j(P_{\pi_v,t,i}(\Pi_t))$ sont ceux de
$H^j(P_{\pi_v,t,i+1}(\Pi_t))$, i.e. $(s-t)(g+1)-2(i+1+k)$ pour $1 \leq k \leq s-t-i-1$ vérifiant $-2(k-1)
\leq j-(s-t-2-i) \leq -(k-1)$. Le changement de variable $k'=k+1$ donne alors le résultat, i.e.
$H^j(P_{\pi_v,t,i}(\Pi_t))$ est de poids $(s-t)(g+1)-2(i+k')$ avec $1 \leq k' \leq s-t-i$ vérifiant
$-2(k'-1)+1 \leq j-(s-t-1-i) \leq -(k'-1)$ soit ce qui est prévu car $j-(s-t-1-i)$ est impair\footnote{Le cas
$j-(s-t-1-i)=-2(k'-1)$ n'est pas à considérer.}.

Pour $j$ de la forme $s-t-1-i-2r$ avec $0 \leq 2r \leq 2(s-t-1-i)$, la suite exacte longue précédente montre
que les poids de $H^j(P_{\pi_v,t,i}(\Pi_t))$ sont, à priori, ceux de $H^j(P_{\pi_v,t,i+1}(\Pi_t))$ ainsi que
celui de $H^j(j^{\geq (t+i+1)g}_{!*})$ soit $(s-t)(g+1)-2(i+1+r)$. D'après l'hypothèse de récurrence,
$H^j(P_{\pi_v,t,i}(\Pi_t))$ est de poids $(s-t)(g+1)-2(i+1+k)$ avec $1 \leq k \leq s-t-i-1$ vérifiant
$-2(k-1) \leq j-(s-t-2-i) \leq -(k-1)$, de sorte que $H^j(P_{\pi_v,t,i}(\Pi_t))$ est de poids
$(s-t)(g+1)-2(i+k')$ avec $1 \leq k' \leq s-t-i$ vérifiant $-2(k'-1)+1 \leq j-(s-t-1-i) \leq -(k'-1)$ soit ce
qui est prévu car le cas $j-(s-t-1-i)=-2(k'-1)$ est justement donné par le poids de $H^j(j^{\geq
(t+i+1)g}_{!*})$ soit $(s-t)(g+1)-2(i+1+r)$.

\end{proof}

\marque \textit{Suite de la preuve du lemme (\ref{lem-rj-combi})}: (iv) Dans la suite on ne considère que les
parties de poids $(s-t)(g-1)$. D'après le lemme (\ref{lem-utile}), $H^{t-s}(P_{\pi_v,t,0}(\Pi_t))$ est nul de
sorte que la suite exacte longue (\ref{sel1}) s'écrit
\begin{multline*}
0 \to H^{t-s}(j^{\geq tg}_!) \longto \Pi_t \overrightarrow{\times} [\overrightarrow{s-t-1}]_{\pi_v} \otimes
(\Xi^{(s-t)(g-1)/2} \otimes \bigoplus_{\chi \in \AF(\pi_v)} \chi^{-1})\\
\longto H^{1+t-s}(P_{\pi_v,t,0}(\Pi_t)) \longto H^{1+t-s}(j^{\geq tg}_!) \to 0 \to \cdots
\end{multline*}

\marque Le cas $i_0=t-s$ découle alors du fait que les strates sont induites, i.e. si $H^{t-s}(j^{\geq
tg}_!)$ est un sous-espace de
$\sharp \ker^1(\Qm,G_\tau) m(\Pi) \Pi_t \overrightarrow{\times} [\overrightarrow{s-t-1}]_{\pi_v} \otimes
(\Xi^{\frac{(s-t)(g-1)}{2}} \otimes \bigoplus_{\chi \in \AF(\pi_v)} \chi^{-1})$
alors il est égal à tout l'espace.

\marque Pour $i_0=1+t-s$, la suite exacte longue (\ref{sel1}) s'écrit
\begin{multline*}
0 \to \sharp \ker^1(\Qm,G_\tau) m(\Pi)\Pi_t \overrightarrow{\times} [\overrightarrow{s-t-1}]_{\pi_v} \otimes
(\Xi^{\frac{(s-t)(g-1)}{2}}
\otimes \bigoplus_{\chi \in \AF(\pi_v)} \chi^{-1}) \longto \\
H^{1+t-s}(P_{\pi_v,t,0}(\Pi_t)) \longto H^{1+t-s}(j^{\geq tg}_!) \to 0
\end{multline*}
tandis que celle associée à $0 \to P_{\pi_v,t,1}(\Pi_t) \longto P_{\pi_v,t,0}(\Pi_t) \longto
\PC_-(g,t+1,\pi_v,t,\Pi_t) \to 0$ s'écrit
\begin{multline*}
0 \to H^{1+t-s}(P_{\pi_v,t,0}(\Pi_t)) \longto \sharp \ker^1(\Qm,G_\tau) m(\Pi) \\
 \Pi_t \overrightarrow{\times} [\overrightarrow{0}]_{\pi_v}
\overrightarrow{\times} [\overrightarrow{s-t-2}]_{\pi_v} \otimes (\Xi^{\frac{(s-t)(g-1)}{2}} \otimes
\bigoplus_{\chi \in \AF(\pi_v)} \chi^{-1}) \to \cdots
\end{multline*}
 Ainsi si $H^{1+t-s}(j^{\geq tg}_!)$ est non nul,
alors $H^{1+t-s}(P_{\pi_v,t,0}(\Pi_t))$ contient strictement
$$\sharp \ker^1(\Qm,G_\tau) m(\Pi)\Pi_t
\overrightarrow{\times} \overrightarrow{(s-t)}_{\pi_v} \otimes (\Xi^{\frac{(s-t)(g-1)}{2}} \otimes
\bigoplus_{\chi \in \AF(\pi_v)} \chi^{-1})$$
et étant de la forme $\sharp \ker^1(\Qm,G_\tau) m(\Pi)\Pi_t
\overrightarrow{\times} \pi \otimes (\Xi^{\frac{(s-t)(g-1)}{2}} \otimes \bigoplus_{\chi \in \AF(\pi_v)}
\chi^{-1})$ ainsi qu'un sous-espace de $\sharp \ker^1(\Qm,G_\tau) m(\Pi)\Pi_t \overrightarrow{\times}
[\overrightarrow{0}]_{\pi_v} \overrightarrow{\times} [\overrightarrow{s-t-2}]_{\pi_v} \otimes
(\Xi^{\frac{(s-t)(g-1)}{2}} \otimes \bigoplus_{\chi \in \AF(\pi_v)} \chi^{-1}),$ on en déduit qu'il est égal
à ce dernier de sorte que $H^{1+t-s}(j^{\geq tg}_!)$ est isomorphe à
$$\sharp \ker^1(\Qm,G_\tau) m(\Pi)\Pi_t \overrightarrow{\times}
[\overleftarrow{1},\overrightarrow{s-t-1}]_{\pi_v} \otimes (\Xi^{\frac{(s-t)(g-1)}{2}} \otimes
\bigoplus_{\chi \in \AF(\pi_v)} \chi^{-1}),$$ d'où le résultat.

\marque Par ailleurs on remarque de la même façon que si $i_0 >1+t-s$, la partie de poids $(s-t)(g-1)$ de
$H^{1+t-s}(P_{\pi_v,t,0}(\Pi_t))$ est alors égale à
$$\sharp \ker^1(\Qm,G_\tau) m(\Pi)\Pi_t \overrightarrow{\times} [\overrightarrow{s-t-1}]_{\pi_v} \otimes
(\Xi^{\frac{(s-t)(g-1)}{2}} \otimes
\bigoplus_{\chi \in \AF(\pi_v)} \chi^{-1}).$$

\marque Supposons alors $i_0 \geq 2+t-s$. La suite exacte longue (\ref{sel1}) donne l'égalité pour tout $i
\geq 2$, des parties de poids $(s-t)(g-1)$ de $H^{t-s+i}(P_{\pi_v,t,0}(\Pi_t))$ et de $H^{t-s+i}(j^{\geq
tg}_!)$. On va montrer que, pour tout $2 \leq i \leq i_0-(t-s)$, la partie de poids $(s-t)(g-1)$ de
$H^{t-s+r}(P_{\pi_v,t,i-2}(\Pi_t))$ est nulle pour $i \leq r < i_0-l+s$ et que celle de
$H^{i+t-s-1}(P_{\pi_v,t,i-2}(\Pi_t))$ est égale à
$$\sharp \ker^1(\Qm,G_\tau) m(\Pi)\Pi_t \overrightarrow{\times} [\overleftarrow{i-1},\overrightarrow{s-t-i}]_{\pi_v}
\otimes (\Xi^{\frac{(s-t)(g-1)}{2}} \otimes
\bigoplus_{\chi \in \AF(\pi_v)} \chi^{-1}).$$ D'après ce que l'on vient de voir, c'est vrai pour $i=2$.
Supposons donc le résultat acquis jusqu'au rang $i$ et traitons le cas de $i+1$. La suite exacte longue de
cohomologie associée à la suite exacte courte $0 \to P_{\pi_v,t,i-1}(\Pi_t) \longto P_{\pi_v,t,i-2}(\Pi_t)
\longto \PC_-(g,t+i-1,\pi_v,t,\Pi_t) \to 0$, s'écrit
\begin{multline*} 
0 \to H^{i-1+t-s}(P_{\pi_v,t,i-2}(\Pi_t)) \longto  \sharp \ker^1(\Qm,G_\tau) m(\Pi)\Pi_t \overrightarrow{\times}
[\overleftarrow{i-2}]_{\pi_v}
\overrightarrow{\times} [\overrightarrow{s-t-i}]_{\pi_v} \\ \otimes (\Xi^{\frac{(s-t)(g-1)}{2}} \otimes
\bigoplus_{\chi \in \AF(\pi_v)} \chi^{-1})
\longto H^{i+t-s}(P_{\pi_v,t,i-1}(\Pi_t)) \longto H^{i+t-s}(P_{\pi_v,t,i-2}(\Pi_t)) \to 0 \cdots
\end{multline*}
ainsi que l'égalité des parties de poids $(s-t)(g-1)$ des espaces $H^{i+t-s+r}(P_{\pi_v,t,i-1}(\Pi_t))$ et
$H^{i+t-s+r}(P_{\pi_v,t,i-2}(\Pi_t))$ pour tout $r>0$. Le cas $i+1$ découle alors de la nullité de la partie
de poids $(s-t)(g-1)$ de $H^{i+t-s}(P_{\pi_v,t,i-2}(\Pi_t))$ et de l'isomorphisme
$$H^{i+t-s-1}(P_{\pi_v,t,i-2}(\Pi_t)) \simeq \sharp \ker^1(\Qm,G_\tau) m(\Pi)\Pi_t \overrightarrow{\times}
[\overleftarrow{i-2}, \overrightarrow{s-t-i+1}]_{\pi_v} \otimes (\Xi^{\frac{(s-t)(g-1)}{2}} \otimes
\bigoplus_{\chi \in \AF(\pi_v)} \chi^{-1}).$$

\marque Considérons
$0 \to P_{\pi_v,t,i_1+1}(\Pi_t) \longto P_{\pi_v,t,i_1}(\Pi_t) \longto \PC_-(g,t+i_1+1,\pi_v,t,\Pi_t) \to 0$
avec $i_1=i_0-t+s-2$, et la suite exacte longue de cohomologie associée qui s'écrit
\begin{multline} \label{sel3}
0 \to H^{i_0-1}(P_{\pi_v,t,i_1}(\Pi_t)) \longto  \sharp \ker^1(\Qm,G_\tau) m(\Pi) \Pi_t \overrightarrow{\times}
[\overleftarrow{i_1}]_{\pi_v}
\overrightarrow{\times} [\overrightarrow{s-t-i_1-2}]_{\pi_v} \\ \otimes (\Xi^{\frac{(s-t)(g-1)}{2}} \otimes
\bigoplus_{\chi \in \AF(\pi_v)} \chi^{-1})  \longto H^{i_0}(P_{\pi_v,t,i_1+1}(\Pi_t)) \longto
H^{i_0}(P_{\pi_v,t,i_1}(\Pi_t)) \to 0 \cdots
\end{multline}
avec $H^{i_0}(P_{\pi_v,t,i_1}(\Pi_t))=H^{i_0}(j^{\geq tg}_!)$ non nul par hypothèse. La suite exacte longue
associée à
$0 \to P_{\pi_v,t,i_1+2}(\Pi_t) \longto P_{\pi_v,t,i_1+1}(\Pi_t) \longto \PC_-(g,i_1+2,\pi_v,t,\Pi_t) \to 0$
s'écrit
\begin{multline*}
0 \to H^{i_0}(P_{\pi_v,t,i_1+1}(\Pi_t)) \longto \sharp \ker^1(\Qm,G_\tau) m(\Pi) \\ \Pi_t \overrightarrow{\times}
[\overleftarrow{i_1+1}]_{\pi_v}
\overrightarrow{\times} [\overrightarrow{s-t-i_1-3}]_{\pi_v}  \otimes (\Xi^{\frac{(s-t)(g-1)}{2}} \otimes
\bigoplus_{\chi \in \AF(\pi_v)} \chi^{-1}) \longto \cdots
\end{multline*}
Ainsi si on veut que $H^{i_0}(P_{\pi_v,t,i_1}(\Pi_t))$ soit non nul, il faut que
$H^{i_0}(P_{\pi_v,t,i_1+1}(\Pi_t))$ soit égal à $\sharp \ker^1(\Qm,G_\tau) m(\Pi)\Pi_t
\overrightarrow{\times} [\overleftarrow{i_1+1}]_{\pi_v} \overrightarrow{\times}
[\overrightarrow{s-t-i_1-3}]_{\pi_v} \otimes (\Xi^{\frac{(s-t)(g-1)}{2}} \otimes \bigoplus_{\chi \in
\AF(\pi_v)} \chi^{-1})$ et donc
$$H^{i_0}(P_{\pi_v,t,i_1+1}(\Pi_t)) \simeq \sharp \ker^1(\Qm,G_\tau) m(\Pi)  \Pi_t \overrightarrow{\times}
[\overleftarrow{i_1+2},\overrightarrow{s-t-i_1-3}]_{\pi_v} \otimes (\Xi^{\frac{(s-t)(g-1)}{2}} \otimes
\bigoplus_{\chi \in \AF(\pi_v)} \chi^{-1})$$ d'où le résultat.
\end{proof}

\marque \textit{Retour sur la preuve de la proposition (\ref{prop-hic-poids}):}

\rem Si on savait que les strates $\bar X_{U^p,m}^{(d-h)}$ étaient affines, la proposition
(\ref{prop-hic-poids}) découlerait directement du point (iii) du lemme (\ref{lem-rj-combi}) pour $k=s-t+1$
car seul $H^{d-tg}_c(\bar X_{U^p,m}^{(d-h)},\FC(g,t,\pi_v) \otimes \LC_{\xi})[\Pi^{\oo,v}]$ pourrait être de
poids $t-s$ ce qui serait contradictoire avec (\ref{somme-alternee}). Ne disposant pas de l'affinité des
strates $\bar X_{U^p,m}^{(d-h)}$, on entre plus précisément dans la combinatoire.

On raisonne par récurrence sur $t$ de $1$ à $s_g$ puis pour $t$ fixé, par récurrence sur $i$ de $t-s$ à $0$.
L'initialisation de la récurrence pour $t=1$ se traite comme le passage de $t$ à $t+1$; on suppose alors le
résultat acquis pour tout $1 \leq t' < t$ \footnote{Pour $t=1$, l'hypothèse de récurrence est vide.}.

\marque Pour $i=t-s$, si l'espace en question était non nul, on aurait d'après le point (iv) du lemme
(\ref{lem-rj-combi}),
\begin{multline*}
\lim_{\atop{\to}{U^p(m)}} H^{(s-t)(g-1)}_c(\bar X_{U^p,m},HT_{\xi}(g,t,\pi_v,\Pi_t))[\Pi^{\oo,v}] =
\\ \sharp \ker^1(\Qm,G_\tau) m(\Pi) \Pi_t \overrightarrow{\times} [\overrightarrow{s-t-1}]_{\pi_v} \otimes
(\Xi^{\frac{(s-t)(g-1)}{2}} \otimes
\bigoplus_{\chi \in \AF(\pi_v)} \chi^{-1}).
\end{multline*}
On considère alors la suite spectrale (\ref{sss}) pour $i=t(g-1)$. On remarque que toutes les parties de
poids $s(g-1)$ de $E_1^{p,q;t(g-1)}[\Pi^{\oo,v}]$ sont nulles pour $p-1 \neq tg$; en effet pour $p-1 < tg$
cela découle de l'hypothèse de récurrence et pour $p-1 > tg$ du fait que $t'(g-1)>t(g-1)$ pour $p-1=t'g$. On
obtient alors
\begin{multline*}
\lim_{\atop{\to}{U^p(m)}} H^{(s-t)(g-1)}(\bar X_{U^p,m},R^{t(g-1)}\Psi_{\pi_v} \otimes
\LC_{\xi})[\Pi^{\oo,v}] \simeq \\ \sharp \ker^1(\Qm,G_\tau) m(\Pi) [\overleftarrow{t-1}]_{\pi_v}
\overrightarrow{\times} [\overrightarrow{s-t-1}]_{\pi_v} \otimes \rec_{F_v}^\vee(\pi_v) (-\frac{s(g-1)}{2})
\end{multline*}
On considère alors la suite spectrale (\ref{ssce}), de sorte que $E_{2,\xi}^{(s-t)(g-1),t(g-1)}[\Pi^{\oo,v}]$
est isomorphe à l'espace ci-dessus. Or comme d'après l'hypothèse de récurrence toutes les parties de poids
$s(g-1)$ des $E_{2,\xi}^{p,q}[\Pi^{\oo,v}]$ pour $q<l(g-1)$ sont nulles et que d'après le lemme
(\ref{lem-rj-combi}) il en est de même pour $p+q < d-s$, on en déduit que $E_{\oo,\xi}^{d-s}[\Pi^{\oo,v}]$
admettrait $\sharp \ker^1(\Qm,G_\tau) m(\Pi)[\overleftarrow{t-1}]_{\pi_v} \overrightarrow{\times}
[\overrightarrow{s-t-1}]_{\pi_v} \otimes
\rec_{F_v}^\vee(\pi_v) (-\frac{s(g-1)}{2})$ comme sous-espace ce qui n'est pas d'après la proposition
(\ref{prop-ssce-poids}).

\marque Supposons donc le résultat vérifié pour tout $t-s \leq i=t-s+\d \leq i_0 < -1$ et traitons le cas de
$i_0$ D'après le point (iv) du lemme (\ref{lem-rj-combi}), on obtiendrait
$$\sharp \ker^1(\Qm,G_\tau) m(\Pi) [\overleftarrow{t-1}]_{\pi_v} \overrightarrow{\times}
[\overleftarrow{\d},\overrightarrow{s-t-1-\d}]_{\pi_v} \otimes (\Xi^{\frac{(s-t)(g-1)}{2}} \otimes
\bigoplus_{\chi \in \AF(\pi_v)} \chi^{-1})$$ Comme précédemment, la suite spectrale (\ref{sss}) pour
$i=t(g-1)$, donne que la partie de poids $s(g-1)$ de $E_{2,\xi}^{(s-t)(g-1)+\d,t(g-1)}[\Pi^{\oo,v}]$ dans la
suite spectrale (\ref{ssce}) est isomorphe à
$\sharp \ker^1(\Qm,G_\tau) m(\Pi) [\overleftarrow{t-1}]_{\pi_v} \overrightarrow{\times}
[\overleftarrow{\d},\overrightarrow{s-t-1-\d}]_{\pi_v} \otimes \rec_{F_v}^\vee(\pi_v) (-\frac{s(g-1)}{2}).$
On remarque alors que $E_{\oo,\xi}^{d-s+\d}[\Pi^{\oo,v}]$ admettrait
$\sharp \ker^1(\Qm,G_\tau) m(\Pi) [\overleftarrow{t+\d},\overrightarrow{s-t-1-\d}]_{\pi_v} \otimes
\rec_{F_v}^\vee(\pi_v) (-\frac{s(g-1)}{2})$
comme sous-espace car ce dernier n'apparaît pas dans les $E_{2,\xi}^{p,q}[\Pi^{\oo,v}]$ pour $q<l(g-1)$
d'après l'hypothèse de récurrence, ni pour $q=t'(g-1)>t(g-1)$ d'après le lemme (\ref{lem-rj-combi}) (i) ainsi
que le corollaire (\ref{coro-1}).

\end{proof}

\subsection{Preuve de (\ref{cas-m})}

Il s'agit donc de prouver la proposition (\ref{cas-m}). En étudiant les parties de poids $s(g-1)$ des
composantes $\Pi^{\oo,v}$ isotypiques pour $\Pi$ irréductible automorphe tel que $\Pi_v \simeq
\speh_s(\pi_v)$, on s'est ramené d'après le paragraphe précédent, à une situation similaire à celle de
\cite{boy}, i.e. dans la suite spectrale des cycles évanescents, seuls les points supersinguliers
contribuent. La fin de la preuve procède alors exactement comme dans loc. cit.

\marque Précisément, la partie de poids $s(g-1)$ de
$\UC^{d-s}_{F_v,l,d}(\JL^{-1}([\overleftarrow{s-1}]_{\pi_v^\vee}))$, d'après le corollaire (\ref{coro-1}),
est un constituant de $[\overrightarrow{s-2}]_{\pi_v} \overrightarrow{\times} [\overrightarrow{0}]_{\pi_v}
\otimes \rec_{F_v}^\vee(\pi_v) (-\frac{s(g-1)}{2})$. On considère alors la suite spectrale (\ref{sss2}) pour
$i=d-s$, où, d'après le corollaire (\ref{coro-hic}), tous les $E_{1,\xi}^{p,q;d-s}$ sont nuls pour $p+q \neq
0$ ou $p-1 \neq d$ de sorte que le terme $E_{2,\xi}^{0,d-s}[\Pi^{\oo,v}]$ de (\ref{ssce2}), et donc
$E_{\oo,\xi}^{0,d-s}[\Pi^{\oo,v}]$ d'après (\ref{prop-hic-poids}), est égal à cet espace qui est donc,
d'après (\ref{prop-ssce-poids}), égal à $[\overrightarrow{s-1}]_{\pi_v} \otimes \rec_{F_v}^\vee(\pi_v)
(-\frac{s(g-1)}{2})$.

\marque On suppose avoir montré par récurrence que pour tout $0 \leq r < r_0< s-1$, les parties de poids
$s(g-1)$ de $\UC_{F_v,l,d}^{d-s+r}(\JL^{-1}([\overleftarrow{s-1}]_{\pi_v^\vee}))$ sont comme prévues, i.e.
nulles pour $r \neq 0$, et égales à $[\overrightarrow{s-1}]_{\pi_v} \otimes \rec_{F_v}^\vee(\pi_v)
(-\frac{s(g-1)}{2})$ pour $r=0$. D'après le corollaire (\ref{coro-1}), la partie de poids $s(g-1)$ de
$\UC_{F_v,l,d}^{d-s+r_0}(\JL^{-1}([\overleftarrow{s-1}]_{\pi_v^\vee}))$ est un constituant de
$[\overrightarrow{s-2-r_0}]_{\pi_v} \overrightarrow{\times} [\overleftarrow{r_0}]_{\pi_v} \otimes
\rec_{F_v}^\vee(\pi_v) (-\frac{s(g-1)}{2})$. Or comme la partie de poids $s(g-1)$ de $\sum_{i=0}^{s-1} (-1)^i
\UC_{F_v,l,d}^{d-s+i}(\JL^{-1}([\overleftarrow{s-1}]_{\pi_v^\vee}))$ est égale à
$[\overrightarrow{s-1}]_{\pi_v} \otimes \rec_{F_v}^\vee(\pi_v) (-\frac{s(g-1)}{2})$, en remarquant, d'après
le corollaire (\ref{coro-1}), que $[\overrightarrow{s-1-r_0},\overleftarrow{r_0+1}]_{\pi_v} \otimes
\rec_{F_v}^\vee(\pi_v) (-\frac{s(g-1)}{2})$ ne peut pas être un constituant de
$\UC_{F_v,l,d}^{d-s+r}(\JL^{-1}([\overleftarrow{s-1}]_{\pi_v^\vee}))$ pour $r>r_0$, on en déduit que si la
partie de poids $s(g-1)$ de $\UC_{F_v,l,d}^{d-s+r_0}(\JL^{-1}([\overleftarrow{s-1}]_{\pi_v^\vee}))$ est non
nulle, elle est alors égale à $[\overrightarrow{s-2-r_0},\overleftarrow{r_0+2}]_{\pi_v} \otimes
\rec_{F_v}^\vee(\pi_v) (-\frac{s(g-1)}{2})$. Comme précédemment ce dernier espace est aussi égal à la partie
de poids $s(g-1)$ de $\frac{1}{\sharp
\ker^1(\Qm,G_\tau)m(\Pi)}E_{\oo,\xi}^{d-s+r_0}[\rec_{F_v}^\vee(\pi_v)][\Pi^{\oo,v}]$ ce qui n'est pas d'après
(\ref{prop-ssce-poids}).

Finalement le cas $r_0=s-1$ est donné en utilisant que la partie de poids $s(g-1)$ de $\sum_{i=0}^{s-1}
(-1)^i \UC_{F_v,l,d}^{d-s+i}(\JL^{-1}([\overleftarrow{s-1}]_{\pi_v^\vee}))$ est égale à
$[\overrightarrow{s-1}]_{\pi_v} \otimes \rec_{F_v}^\vee(\pi_v) (-\frac{s(g-1)}{2})$.
\end{proof}

\begin{rema} \label{rema-N2}
Comme dans la remarque (\ref{rema-N}), on aurait pu traiter tous les poids des $\UC_{F_v,l,d}^{i}$. Il suffit
pour cela de montrer un analogue du point (iv) de la proposition (\ref{prop-hic-poids}) pour tous les poids.
\end{rema}


\part{Compléments sur la cohomologie globale}

\section{Parties sans monodromie de la cohomologie globale}

\begin{propb} \label{prop-lrs2}
Soient $g$ un diviseur de $d=sg$ et $\pi_v$ une représentation irréductible cuspidale unitaire de
$GL_g(F_v)$. On considère une représentation automorphe $\Pi$ de $G_\tau(\Am)$ vérifiant $\hyp(\xi)$ telle
que $\Pi_v \simeq [\overrightarrow{s-1}]_{\pi_v}$. Le $GL_d(F_v) \times W_v$-module
$H^{d+s-2-2i}_{\eta_v,\xi}[\Pi^{\oo,v}]$ est alors isomorphe à
$\sharp \ker^1(\Qm,G_\tau) m(\Pi)[\overrightarrow{s-1}]_{\pi_v} \otimes \rec_{F_v}^\vee(\pi_v) (-\frac{d+s-2-2i}{2})$
pour $0 \leq i <s$.
\end{propb}

\begin{proof} La proposition (\ref{prop-ssce-poids}) joint au théorème de Lefschetz difficile implique que pour tout
$0 \leq r \leq s-1$, $H^{d-s+2r}_{\eta_v,\xi}[\Pi^{\oo,v}]$ admet $\sharp
\ker^1(\Qm,G_\tau)m(\Pi)[\overrightarrow{s-1}]_{\pi_v} \otimes \rec_{F_v}^\vee(\pi_v) (-\frac{s(g-1)+2r}{2})$
comme facteur direct. On reprend la suite spectrale (\ref{ss-poids}) en utilisant la proposition
(\ref{prop-not}) ainsi que le théorème de Lefschetz difficile. On pourra se reporter à la figure
(\ref{figure8}) où l'on a représenté les $H^i(gr_{k,\xi})[\Pi^{\oo,v}]$ dans le cas $s=4$. On rappelle tout
d'abord que pour $1< \d <s$, on a
\begin{multline} \label{higrk}
\frac{H^{-\d}(gr_{k,\xi})[\Pi^{\oo,v}]}{\sharp \ker^1(\Qm,G_\tau)m(\Pi)}= \bigoplus_{\atop{|k|< t \leq s}{t
\equiv k+1 \equiv s+\d \mod 2}} ([\overleftarrow{t-1}]_{\pi_v} \overrightarrow{\times}
[\overrightarrow{\frac{s-t+\d-2}{2}}]_{\pi_v}) \overleftarrow{\times}
[\overrightarrow{\frac{s-t-\d-2}{2}}]_{\pi_v} \\  \otimes \rec_{F_v}^\vee(\pi_v) (-\frac{d-1-\d+k}{2})
\end{multline}

\begin{multline}
\frac{H^{\d}(gr_{k,\xi})[\Pi^{\oo,v}]}{\sharp \ker^1(\Qm,G_\tau)m(\Pi)}= \bigoplus_{\atop{|k|< t \leq s}{t
\equiv k+1 \equiv s + \d \mod 2}} ([\overleftarrow{t-1}]_{\pi_v} \overrightarrow{\times}
[\overrightarrow{\frac{s-t-\d-2}{2}}]_{\pi_v}) \overleftarrow{\times}
[\overrightarrow{\frac{s-t+\d-2}{2}}]_{\pi_v} \\ \otimes \rec_{F_v}^\vee(\pi_v) (-\frac{d-1+\d+k}{2})
\end{multline}
On remarque ainsi que tous les $[\overrightarrow{s-1}]_{\pi_v}$ des $E_{1,\xi}^{p,q}[\Pi^{\oo,v}]$ de
(\ref{ss-poids}) restent dans l'aboutissement; il nous faut alors montrer que tous les autres disparaissent.
Dans la suite on ne considérera plus les $[\overrightarrow{s-1}]_{\pi_v}$.

\marque Prenons dans un premier temps $\d>1$, de sorte qu'outre les $[\overrightarrow{s-1}]_{\pi_v}$, les
éventuels constituants non nuls de $E_{\oo,\xi}^{d-s-\d}[\Pi^{\oo,v}]$ (resp. de
$E_{\oo,\xi}^{d-s+\d}[\Pi^{\oo,v}]$) sont de la forme
$$\sharp \ker^1(\Qm,G_\tau) m(\Pi)[\overrightarrow{r},\overleftarrow{t-1},\overrightarrow{s-t-r}]_{\pi_v} \otimes
\rec_{F_v}^\vee(\pi_v)
(-\frac{d-1-\d-k}{2})$$ avec $s-t-r=r+\d,r+\d+1,r+\d-1$ (resp. $s-t-r=r-\d,r-\d-1,r-\d+1$) et certains
entiers $k$ qu'il n'est pas nécessaire de préciser. Pour $\d \geq 2$, on remarque que les ensembles $\{
\d,\d+1,\d-1 \}$ et $\{ -\d,-\d+1,-\d-1 \}$ sont disjoints de sorte que d'après le théorème de Lefschetz
difficile, l'éventuel constituant est forcément nul. \footnote{On peut aussi argumenter en remarquant que ce
ne sont pas des constituants locaux licites d'une représentation automorphe alors que dans l'aboutissement
seules celles ci doivent apparaître.}

\marque Pour $\d=1$, le raisonnement est plus fin et demande de distinguer les sous-espaces des quotients
dans nos induites. Si on reprend le raisonnement précédent dans le cas $\d=1$, le théorème de Lefschetz
difficile impose que les éventuels constituants de $H^{d-1 \pm 1}_{\eta_v,\xi}[\Pi^{\oo,v}]$ sont de la forme

\begin{equation} \label{forme}\sharp \ker^1(\Qm,G_\tau) m(\Pi)
[\overrightarrow{\frac{s-t}{2}},\overleftarrow{t-1},\overrightarrow{\frac{s-t}{2}}]_{\pi_v}  \otimes
\rec_{F_v}^\vee(\pi_v) (-\frac{d-1-k}{2})
\end{equation}
pour $1 \leq t \leq s$, $0 \leq r \leq (s-t)/2$ et certains entiers $k$ qu'il n'est pas nécessaire de
préciser. On revient sur (\ref{higrk}) et sur la suite spectrale (\ref{ss-poids}). Ainsi on a
\begin{multline*}
\frac{E_{1,\xi}^{-k,-1+k}[\Pi^{\oo,v}]}{\sharp \ker^1(\Qm,G_\tau)m(\Pi)}= \bigoplus_{\atop{|k|< t \leq s}{t
\equiv k+1 \equiv s-1 \mod 2}} ([\overleftarrow{t-1}]_{\pi_v} \overrightarrow{\times}
[\overrightarrow{\frac{s-t-1}{2}}]_{\pi_v}) \overleftarrow{\times} [\overrightarrow{\frac{s-t-3}{2}}]_{\pi_v}
\\ \otimes \rec_{F_v}^\vee(\pi_v) (-\frac{d-2+k}{2})
\end{multline*}
\begin{multline*}
\frac{E_{1,\xi}^{-k-1,-1+k}[\Pi^{\oo,v}]}{\sharp \ker^1(\Qm,G_\tau) m(\Pi)}= \bigoplus_{\atop{|k|< t \leq
s}{t \equiv k+1 \equiv s-2 \mod 2}} ([\overleftarrow{t-1}]_{\pi_v} \overrightarrow{\times}
[\overrightarrow{\frac{s-t}{2}}]_{\pi_v}) \overleftarrow{\times} [\overrightarrow{\frac{s-t-4}{2}}]_{\pi_v}
\\ \otimes \rec_{F_v}^\vee(\pi_v) (-\frac{d-2+k}{2})
\end{multline*}
\begin{multline*}
\frac{E_{1,\xi}^{-k+1,-1+k}[\Pi^{\oo,v}]}{\sharp \ker^1(\Qm,G_\tau) m(\Pi)}= \bigoplus_{\atop{|k|< t \leq
s}{t \equiv k+1 \equiv s \mod 2}} ([\overleftarrow{t-1}]_{\pi_v} \overrightarrow{\times}
[\overrightarrow{\frac{s-t-2}{2}}]_{\pi_v}) \overleftarrow{\times} [\overrightarrow{\frac{s-t-2}{2}}]_{\pi_v}
\\ \otimes \rec_{F_v}^\vee(\pi_v) (-\frac{d-2+k}{2})
\end{multline*}

\marque Pour $k=s-1$, on remarque que la partie de poids $s(g+1)-2$ de $H^{d-2}_{\eta_v,\xi}[\Pi^{\oo,v}]$
est un sous-espace de $\sharp \ker^1(\Qm,G_\tau) m(\Pi)[\overleftarrow{s-2}]_{\pi_v} \overrightarrow{\times}
[\overrightarrow{0}]_{\pi_v} \otimes \rec_{F_v}^\vee(\pi_v) (-\frac{s(g+1)-2}{2})$. Or
$[\overleftarrow{s-1}]_{\pi_v}$ n'étant pas un sous-espace de $[\overleftarrow{s-2}]_{\pi_v}
\overrightarrow{\times} [\overrightarrow{0}]_{\pi_v}$, on en déduit, en utilisant (\ref{forme}), que la
partie de poids $s(g+1)-2$ de $H^{d-2}_{\eta_v,\xi}[\Pi^{\oo,v}]$ est nulle.

\marque Pour $0 \leq k \leq s-3$, on a:
\begin{itemize}
\item $E_{1,\xi}^{-k,-1+k}[\Pi^{\oo,v}]=E_{1,\xi}^{-k-2,k+1}[\Pi^{\oo,v}] \otimes |-|^{2} \oplus \sharp
\ker^1(\Qm,G_\tau) m(\Pi)
V_{k,-1} \otimes \rec_{F_v}^\vee(\pi_v) (-\frac{d-2+k}{2}) $, avec
$V_{k,-1}=([\overleftarrow{k}]_{\pi_v}
\overrightarrow{\times} [\overrightarrow{\frac{s-k-2}{2}}]_{\pi_v}) \overleftarrow{\times}
[\overrightarrow{\frac{s-k-4}{2}}]_{\pi_v};$

\item $E_{1,\xi}^{-k-1,-1+k}[\Pi^{\oo,v}]=E_{1,\xi}^{-k-3,k+1}[\Pi^{\oo,v}] \otimes |-|^{2} \oplus \sharp
\ker^1(\Qm,G_\tau) m(\Pi)
V_{k,-2} \otimes \rec_{F_v}^\vee(\pi_v) (-\frac{d-2+k}{2})$, avec
$V_{k,-2}=([\overleftarrow{k+1}]_{\pi_v}
\overrightarrow{\times} [\overrightarrow{\frac{s-k-2}{2}}]_{\pi_v}) \overleftarrow{\times}
[\overrightarrow{\frac{s-k-6}{2}}]_{\pi_v};$

\item $E_{1,\xi}^{-k+1,-1+k}[\Pi^{\oo,v}]=E_{1,\xi}^{-k-1,k+1}[\Pi^{\oo,v}] \otimes |-|^{2} \oplus \sharp
\ker^1(\Qm,G_\tau) m(\Pi)
V_{k,0} \otimes \rec_{F_v}^\vee(\pi_v) (-\frac{d-2+k}{2})$ avec $V_{0,0}=0$ et
$V_{k,0}=([\overleftarrow{k-1}]_{\pi_v}
\overrightarrow{\times} [\overrightarrow{\frac{s-k-4}{2}}]_{\pi_v}) \overleftarrow{\times}
[\overrightarrow{\frac{s-k-4}{2}}]_{\pi_v} \hbox{ pour }0 < k \leq s-3.$
\end{itemize}

On suppose alors par récurrence que $E_\oo^{-k-2,k+1}[\Pi^{\oo,v}]=(\ker d_1^{-k-2,k+1} / \im
d_1^{-k-3,k+1})[\Pi^{\oo,v}]$ est nulle; l'opérateur de monodromie $N$ implique alors que
$E_\oo^{-k,-1+k}[\Pi^{\oo,v}]$ est égal à $\ker d_+ /\im d_- \otimes \rec_{F_v}^\vee(\pi_v)
(-\frac{d-2+k}{2})$ avec $V_{k,-2} \longmapright{d_-} V_{k,-1} \longmapright{d_+} V $. Ainsi d'après
(\ref{forme}), si $E_\oo^{-k,-1+k}[\Pi^{\oo,v}]$ était non nul, il serait égal à
$$\sharp \ker^1(\Qm,G_\tau)
m(\Pi)[\overrightarrow{\frac{s-k}{2}},\overleftarrow{k-1},\overrightarrow{\frac{s-k}{2}}]_{\pi_v} \otimes
\rec_{F_v}^\vee(\pi_v) (-\frac{d-2+k}{2}).$$ En remarquant qu'aucun des constituants de
$[\overrightarrow{\frac{s-k}{2}},\overleftarrow{k-1}]_{\pi_v} \overrightarrow{\times}
[\overrightarrow{\frac{s-k-2}{2}}]_{\pi_v}$ n'est un constituant de $V_{k,-2}$, on en déduit que quelque soit
$d_-$, on a une surjection
$$V_{k,-1} / \im d_- \twoheadrightarrow \sharp \ker^1(\Qm,G_\tau) m(\Pi)
[\overrightarrow{\frac{s-k}{2}},\overleftarrow{k-1}]_{\pi_v}
\overrightarrow{\times} [\overrightarrow{\frac{s-k-2}{2}}]_{\pi_v} \otimes \rec_{F_v}^\vee(\pi_v)
(-\frac{d-2+k}{2})$$ alors que
$[\overrightarrow{\frac{s-k}{2}},\overleftarrow{k-1},\overrightarrow{\frac{s-k}{2}}]_{\pi_v}$ n'est pas un
sous-espace de $[\overrightarrow{\frac{s-k}{2}},\overleftarrow{k-1}]_{\pi_v} \overrightarrow{\times}
[\overrightarrow{\frac{s-k-2}{2}}]_{\pi_v}$ de sorte que $E_\oo^{-k,-1+k}[\Pi^{\oo,v}]$ est nul.

\marque Pour $1-s \leq k \leq 0$ on raisonne de manière strictement identique en étudiant les suites exactes
$E_{1,\xi}^{-k-1,k+1} \longmapright{d_1^{-k-1,k+1}} E_{1,\xi}^{-k,k+1} \longmapright{d_1^{-k,k+1}}
E_{1,\xi}^{-k+1,k+1}$ et en remplaçant l'étude des sous-espaces par celle des quotients.

\end{proof}

\begin{corob} \label{strate-alt-p}
Soit $\Pi$ une représentation irréductible automorphe de $G_\tau(\Am)$ telle que $\Pi_v \simeq
[\overrightarrow{s-1}]_{\pi_v}$ où $\pi_v$ est une représentation irréductible cuspidale unitaire de
$GL_g(F_v)$. Si $[H^i_{h,\xi,\tau_v}(\Pi^{\oo,v})]$ muni de son action naturelle de $GL_d(F_v)$ est non nul
dans le groupe de Grothendieck $\groth(GL_{d-h}(F_v) \times (D_{v,h}^\times/\DC_{v,h}^\times))$, on a alors
$\tau_v=\JL^{-1}([\overleftarrow{s-1}]_{\pi_v^\vee})$, $h=tg$ avec $1 \leq t \leq s$ et $i=(s-t)(g+1)$. Dans
ce cas, ils sont donnés par
\begin{multline*}
\lim_{\atop{\to}{U^p(m)}} H^{(s-t)(g+1)}_c(\bar X_{U^p,m,M_{d-tg}}^{(d-tg)},\FC_{\tau_v} \otimes
\LC_{\xi})[\Pi^{\oo,v}] \simeq \\ \sharp \ker^1(\Qm,G_\tau) m(\Pi) [\overrightarrow{s-t-1}]_{\pi_v(t(g+1)/2)}
\otimes (\Xi^{\frac{(s-t)(g+1)}{2}} \bigoplus_{\chi \in \DF(\pi_v)} \chi^{-1})
\end{multline*}
en tant que représentation de $GL_{(s-t)g}(F_v) \times (D_{v,tg}^\times/\DC_{v,tg}^\times)$, où $\DF(\pi_v)$
est l'ensemble des caractères $\chi$ de $\Zm \simeq D_{v,tg}^\times /\DC_{v,tg}^\times$ tels que
$\JL^{-1}([\overleftarrow{s-1}]_{\pi_v^\vee}) \otimes \chi^{-1} \simeq
\JL^{-1}([\overleftarrow{s-1}]_{\pi_v^\vee}).$
\end{corob}

\begin{proof} On rappelle que d'après le lemme (\ref{lem-rj-combi}), les composantes $\Pi^{\oo,v}$-isotypiques
des $H^i_c(\bar X_{U^p,m}^{(d-tg)},\FC(g,t,\pi_v) \otimes \LC_{\xi}) \otimes \Pi_t$ sont mixtes de poids
$(s-t)(g-1)+2k$ avec $0 \leq k < s-t$ et de la forme $\sharp \ker^1(\Qm,G_\tau) m(\Pi)(\Pi_t
\overrightarrow{\times} \pi_+) \overleftarrow{\times} \pi_- \otimes \rec_{F_v}^\vee(\pi_v)
(-\frac{(s-t)(g-1)+2k}{2})$ où $\pi_+$ (resp. $\pi_-$) est une représentation elliptique de
$GL_{(s-t-k)g}(F_v)$ (resp. de $GL_{kg}(F_v)$). On raisonne alors par récurrence sur $k$, de $0$ à $s-t-2$,
afin de montrer que pour tout $i$, les parties de poids $(s-t)(g-1)+2k$ des $H^i_c(\bar
X_{U^p,m}^{(d-tg)},\FC(g,t,\pi_v,I) \otimes \LC_{\xi})[\Pi^{\oo,v}]$ sont nulles. Le cas $k=0$ a été traité
dans la proposition (\ref{prop-hic-poids}). Supposons le résultat acquis jusqu'au rang $k<s-t-2$ et montrons
le au rang $k+1$. On raisonne alors par récurrence sur $t$ de $1$ à $s-1$. L'initialisation de la récurrence
se prouve comme le cas général; on suppose donc le résultat acquis jusqu'au rang $t<s-2$ et prouvons le au
rang $t+1$.

\marque Supposons qu'il existe $j$ tel que la partie de poids $(s-t)(g-1)+2(k+1)$ de $H^j_c(\bar
X_{U^p,m}^{(d-tg)},\FC(g,t,\pi_v) \otimes \LC_{\xi})[\Pi^{\oo,v}]$ soit non nulle. On étudie ensuite la suite
spectrale (\ref{sss2}) associée à la stratification pour $i=t(g-1)$. Ainsi la partie de poids
$m_k:=s(g-1)+2(k+1)$ de $E_{1,\xi}^{tg+1,j-tg-1;t(g-1)}[\Pi^{\oo,v}]$ est non nulle de la forme $\sharp
\ker^1(\Qm,G_\tau) m(\Pi)([\overrightarrow{t-1}]_{\pi_v} \overrightarrow{\times} \pi_+)
\overleftarrow{\times} \pi_- \otimes \rec_{F_v}^\vee(\pi_v) (-\frac{m_k}{2})$ où $\pi_+$ (resp. $\pi_-$) est
une représentation elliptique de $GL_{(s-t-k)g}(F_v)$ (resp. de $GL_{kg}(F_v)$). Or d'après l'hypothèse de
récurrence pour tout $1 \leq t' <t$, les parties de poids $m_k=(s-t')(g-1)+t'(g-1)+2(k+1)$ des
$E_{1,\xi}^{t'g+1,j'-t'g-1;t(g-1)}[\Pi^{\oo,v}]$ sont nulles; en effet celles-ci sont données par les parties
de poids $(s-t')(g-1)+2k'$ avec $s(g-1)+2(k+1)=(s-t')(g-1)+t'(g-1)+2\d+2k'$ avec $t(g-1)=t'(g-1)+\d$, $0 < \d
< t'<t$, soit $k'=k+1-\d$. Pour $t'>t$, les parties de poids $m_k$ des
$E_{1,\xi}^{t'g+1,j'-t'g-1;t(g-1)}[\Pi^{\oo,v}]$ sont nulles car $t(g-1)$ ne s'écrit pas sous la forme
$t'(g-1)+\d$ avec $0 \leq \d < t'$.

\marque On étudie la suite spectrale (\ref{ssce2}) des cycles évanescents. D'après ce qui précède, la partie
de poids $m_k$ de $E_{2,\xi}^{j,t(g-1)}[\Pi^{\oo,v}]$ est de la forme $\sharp \ker^1(\Qm,G_\tau)
m(\Pi)([\overrightarrow{t-1}]_{\pi_v} \overrightarrow{\times} \pi_+) \overleftarrow{\times} \pi_- \otimes
\rec_{F_v}^\vee(\pi_v) (-\frac{m_k}{2}).$ Par ailleurs les parties de poids $m_k$ des
$E_{2,\xi}^{j+r+1,t(g-1)-r}[\Pi^{\oo,v}]$, pour $r \geq 1$ sont nulles; en effet celles-ci proviendraient à
travers la suite spectrale (\ref{sss2}), des parties de poids $(s-t')(g-1)+2k'$ de la composante
$\Pi^{\oo,v}$-isotypique de $H^{j+r-1}_c(\bar X_{U^p,m}^{(d-l'g)},\FC(g,t',\pi_v) \otimes \LC_{\xi}) \otimes
[\overleftarrow{k'}, \overrightarrow{t-1-k'}]_{\pi_v},$ pour $1 \leq t' <t$ avec
$s(g-1)+2(k+1)=(s-t')(g-1)+t'(g-1)+2\d+2k'$, soit $k'=k+1-\d$ avec $0 < \d < t'$ qui sont nulles d'après
l'hypothèse de récurrence. En ce qui concerne les parties de poids $m_k$ des
$E_{2,\xi}^{j-r-1,t(g-1)+r}[\Pi^{\oo,v}]$ pour $r >0$, elles proviennent à nouveau des parties de poids
$(s-t')(g-1)+2k'$ des composantes $\Pi^{\oo,v}$-isotypiques des $H^{j-r-1}_c(\bar
X_{U^p,m}^{(d-t'g)},\FC(g,t',\pi_v) \otimes \LC_{\xi}) \otimes
[\overleftarrow{k'},\overrightarrow{t-1-k'}]_{\pi_v}$ pour $1 \leq t' \leq s$ avec $k'=k+1-\d$ et $0 \leq
\d=(t-t')(g-1)+r < t'$, ce qui impose d'après l'hypothèse de récurrence $t'>t$, avec $\d=0$ et $k'=k+1$. On
remarque alors que $\sharp \ker^1(\Qm,G_\tau)
m(\Pi)[\overleftrightarrow{k},\overrightarrow{t-1},\overleftarrow{1},\overleftrightarrow{s-t-k-1}]_{\pi_v}
\otimes \rec_{F_v}^\vee(\pi_v) (-\frac{s(g-1)+2(k+1)}{2})$ serait un constituant de
$E_{\oo,\xi}^j[\Pi^{\oo,v}]$, ce qui n'est pas d'après la proposition (\ref{prop-lrs2}).

\end{proof}

\rem On pourra voir la figure (\ref{figure6}), où l'on a représenté pour $s=4$ et $g=2$, les
$E_2^{p,q}[\Pi^{\oo,v}]$ de la suite spectrale des cycles évanescents.

\section{Correspondances de Jacquet-Langlands globales}
\label{corres-jl}

Le but de se paragraphe est d'apporter l'amélioration suivante à la proposition (3.2.1) de \cite{H-L}.

\begin{propb} \label{theo-jl}
Il existe une bijection dite de Jacquet-Langlands entre:

\noindent - les représentations irréductibles automorphes $\bar \Pi$ de $H_0(\Am)$ cohomologiques pour $\xi'$;

\noindent - les représentations irréductibles automorphes $\Pi$ de $G_\tau(\Am)$ cohomologiques pour $\xi'$ vérifiant
l'une des deux
conditions suivantes:
\begin{itemize}
\item[(a)] $\Pi_v$ est une représentation essentiellement de carré intégrable, i.e. $\Pi_v \simeq
[\overleftarrow{s-1}]_{\pi_v}$
pour $\pi_v$ une représentation irréductible cuspidale de $GL_g(F_v)$ avec $d=sg$;

\item[(b)] $\Pi_v \simeq [\overrightarrow{s-1}]_{\pi_v}$,
\end{itemize}
compatible aux correspondances de Jacquet-Langlands locales, soit $\bar \Pi^{\oo,v} \simeq \Pi^{\oo,v}$ et
$\Pi_v \simeq \JL(\bar \Pi_v)$ dans le cas (a) et $\Pi_v \simeq \iota(\JL(\bar \Pi_v))$ dans le cas (b) où
$\iota$ désigne l'involution de Zelevinski. En outre on a
$m(\Pi)=m(\bar \Pi)$.

Par ailleurs soit $\Pi^v$ une représentation de $G_\tau(\Am^v)$ telle que pour toute représentation
$\Pi_v$ de $GL_d(F_v)$ avec $\Pi:=\Pi^v \Pi_v$ vérifiant $\hyp(\xi)$, $\Pi_v$ n'est pas de la forme
$[\overleftarrow{s-1}]_{\pi_v}$ ou $[\overrightarrow{s-1}]_{\pi_v}$ pour $\pi_v$ une représentation
irréductible cuspidale de $GL_g(F_v)$ avec $d=sg$. Alors il n'existe pas de représentation irréductible
automorphe $\bar \Pi$ de $H_0(\Am)$ telle que $\bar \Pi^{\oo,v} \simeq \Pi^{\oo,v}$.
\end{propb}

\begin{proof} Soit $\Pi$ une représentation automorphe de $G_\tau(\Am)$ cohomologique pour $\xi$ et telle que $\Pi_v
\simeq
[\overleftarrow{s'-1}]_{\pi_v}$ pour une certaine représentation cuspidale unitaire $\pi_v$ de
$GL_{g'}(F_v)$. D'après la proposition (\ref{prop-coho1}), pour tout $1 \leq t < s'$, on a
\begin{multline*}
\lim_{\atop{\to}{U^p(m)}} H^0(\bar X_{U^p,m}^{[d-tg']},j^{\geq tg'}_! HT_{\xi}(g',t,\pi_v,\Pi_t)[d-tg'])
[\Pi^{\oo,v}]= \\
\sharp \ker^1(\Qm,G_\tau) m(\Pi) \Pi_t \overrightarrow{\times} [\overleftarrow{s'-t-1}]_{\pi_v} \otimes
(\Xi^{\frac{(s'-t)(g'-1)}{2}} \otimes \bigoplus_{\chi \in \AF(\pi_v)} \chi^{-1})
\end{multline*}
Par ailleurs ce dernier est aussi égal à
\begin{multline*}
\lim_{\atop{\to}{U^p(m)}} H^0(\bar X_{U^p,m}^{(0)} ,\FC(g',s',\pi_o) \otimes \LC_{\xi} \otimes \Pi_t
\overrightarrow{\times} [\overleftarrow{s'-t-1}]_{\pi_v})[\Pi^{\oo,v}] \otimes (\Xi^{\frac{(s'-t)(g'-1)}{2}}
\otimes \bigoplus_{\chi \in \AF(\pi_v)} \chi^{-1})
\end{multline*}
qui est égal à $\sum_{\bar \Pi \in \UF_{H_0,\xi}(\Pi^{\oo,v})} \sharp \ker^1(\Qm,H_0) m(\bar \Pi) \Pi_t
\overrightarrow{\times} [\overleftarrow{s'-t-1}]_{\pi_v} \otimes (\Xi^{\frac{(s'-t)(g'-1)}{2}} \otimes
\bigoplus_{\chi \in \AF(\pi_v)} \chi^{-1})$ d'après le lemme (\ref{lem-pts-ss}), où
$\UF_{H_0,\xi}(\Pi^{\oo,v})$ désigne l'ensemble des classes d'isomorphismes des représentations irréductibles
$\bar \Pi$ automorphes cohomologiques pour $(\xi')^\vee$ et telles que $(\bar \Pi^{\oo,v})^\vee \simeq
\Pi^{\oo,v}$ avec $\bar \Pi_v^\vee \simeq \JL^{-1}([\overleftarrow{s'-1}]_{\pi_v})$ et où on rappelle que
$\sharp \ker^1(\Qm,G_\tau)=\sharp \ker^1(\Qm,H_0)$.

\marque De la même façon, soit $\Pi$ une représentation irréductible automorphe de $G_\tau(\Am)$
cohomologique pour $\xi$ et telle que $\Pi_v \simeq [\overrightarrow{s'-1}]_{\pi_v}$. On reprend la preuve de
la proposition (\ref{prop-not}). Pour $t=s'$, le lemme (\ref{lem-pts-ss}) donne
\begin{multline*}
\lim_{\atop{\to}{U^p(m)}} H^0(\bar X_{U^p,m}^{(s'g')},\FC(g',s',\pi_v) \otimes
[\overleftarrow{s'-1}]_{\pi_v})[\Pi^{\oo,v}]= \\ \sharp \ker^1(\Qm,H_0) \Bigl ( \sum_{\bar \Pi \in
\UF_{H_0,\xi}(\Pi^{\oo,v})} m(\bar \Pi) \Bigr ) [\overleftarrow{s'-1}]_{\pi_v} \otimes \bigoplus_{\chi \in
\AF(\pi_v)} \chi^{-1}
\end{multline*}
Pour $t=s'-1$, on obtient alors
\begin{multline*}
\frac{1}{\sharp \ker^1(\Qm,G_\tau)} \sum_i (-1)^i H^i(j^{\geq (s'-1)g'}_{!*})=  - \Bigl ( \bigl ( \sum_{\bar
\Pi \in \UF_{H_0,\xi}(\Pi^{\oo,v})} m(\bar \Pi) \bigr ) [\overleftarrow{s'-2}]_{\pi_v}
\overrightarrow{\times} [\overleftarrow{0}]_{\pi_v} \otimes \Xi^{-1/2} +  \\ \bigl ( \sum_{\Pi \in
\UF_{G_\tau,\xi}(\Pi^{\oo,v})} m(\Pi) \bigr )  [\overleftarrow{s'-2}]_{\pi_v} \overleftarrow{\times}
[\overleftarrow{0}]_{\pi_v} \otimes \Xi^{1/2} \Bigr ) \otimes (\Xi^{g'/2} \otimes \bigoplus_{\chi \in
\AF(\pi_v)} \chi^{-1})
\end{multline*}
qui d'après la pureté, doit être égal à $\frac{-1}{\sharp \ker^1(\Qm,G_\tau)} (H^{-1}(j^{\geq
(s'-1)g'}_{!*})+H^1(j^{\geq (s'-1)g'}_{!*}))$ de sorte que la dualité de Verdier donne le résultat.

\marque Supposons que $\Pi_v$ ne soit ni de la forme $[\overleftarrow{s'-1}]_{\pi_v}$ ni
$[\overrightarrow{s'-1}]_{\pi_v}$ pour tout diviseur $s'$ de $d=s'g'$ et toute représentation cuspidale
$\pi_v$ de $GL_{g'}(F_v)$. Supposons qu'il existe une représentation $\bar \Pi$ cohomologique pour $\xi$
telle que $\bar \Pi^{\oo,v} \simeq \Pi^{\oo,v}$ et soit $s'$, $\pi_v$ tel que $\bar \Pi_v \simeq
\JL^{-1}([\overleftarrow{s'-1}]_{\pi_v})$. Le lemme (\ref{lem-pts-ss}) donne comme précédemment
$$\frac{1}{\sharp \ker^1(\Qm,H_0)}H^0(\FC(g,s',\pi_v) \otimes
[\overleftarrow{s'-1}]_{\pi_v}) \simeq \Bigl ( \sum_{\bar \Pi \in \UF_{H_0,\xi}(\Pi^{\oo,v})} m(\bar \Pi)
\Bigr ) [\overleftarrow{s'-1}]_{\pi_v} \otimes \bigoplus_{\chi \in \AF(\pi_v)} \chi^{-1}$$ de sorte que
\begin{multline*}
\frac{1}{\sharp \ker^1(\Qm,G_\tau)}\sum_i (-1)^i H^i(j^{\geq (s'-1)g'}_{!*})[\Pi^{\oo,v}]= \Bigl ( -
\bigl ( \sum_{\Pi \in \UF_{G_\tau,\xi}(\Pi^{\oo,o})} m(\Pi) \bigr ) (\ind_{P_{d-g',d}(F_v)}^{GL_d(F_v)}
[\overleftarrow{s'-2}]_{\pi_v} \\
\otimes \red_{\JL^{-1}([\overleftarrow{s'-2}]_{\pi_v})}^{d-g'}(\Pi_v)) + \bigl ( \sum_{\bar \Pi \in
\UF_{H_0,\xi}(\Pi^{\oo,v})} m(\bar \Pi) \bigr ) [\overleftarrow{s'-2}]_{\pi_v} \overrightarrow{\times}
[\overleftarrow{0}]_{\pi_v} \otimes \Xi^{(g'-1)/2} \Bigr ) \otimes \bigoplus_{\chi \in \AF(\pi_v)} \chi^{-1})
\end{multline*}
La condition de pureté et la compatibilité à la dualité de Verdier, impose alors que $\Pi_v$ est soit
isomorphe à $[\overleftarrow{s'-1}]_{\pi_v}$ soit à $[\overrightarrow{s'-1}]_{\pi_v}$ d'où la contradiction.

\marque Ce dernier raisonnement montre en outre qu'étant donné une représentation $\bar \Pi$ comme dans
l'énoncé, il lui correspond une représentation $\Pi$ vérifiant les conditions de l'énoncé. Dans le cas où
$\Pi_v \simeq [\overrightarrow{s'-1}]_{\pi_v}$, l'égalité $s=s'$ découle du corollaire
(\ref{coro-compo-locale}) ci-après.

\end{proof}

\section[Composantes locales]{Composantes locales des représentations automorphes}

\begin{propb} \label{prop-compo-locale}
Soit $\Pi$ une représentation irréductible automorphe de $G_\tau(\Am)$ vérifiant $\hyp(\oo)$. Pour tout $0
\leq i < d$, il existe alors un réel $\d$ ainsi que des entiers $n_i$ \footnote{$n_i$ représente la dimension
de la partie primitive de la représentation galoisienne de $H^i_{\eta_v}$} tels que:

\begin{itemize}
\item $\sum_{i=0}^{d-1} n_i(d-i)=d$,

\item pour $i$ tel que $n_i \neq 0$, $n_j=0$ pour $j \equiv i+1 \mod 2$, tels que pour toute place $v$ telle que
$D_{v,d}^\times \simeq GL_d(F_v)$, il existe:
\begin{itemize}
\item des entiers $t_j>1$ et $g_j$ pour $1 \leq j \leq u$, ainsi que des entiers
$t'_k$ et $g'_k$ pour $1 \leq k \leq u'$, avec $\sum_{j=1}^u t_jg_j+\sum_{k=1}^{u'} t'_k g'_k=d,$

\item des représentations irréductibles cuspidales $\pi_{v,j}$ et $\pi_{v,k}'$ de respectivement $GL_{g_j}(F_v)$ et
$GL_{g'_k}(F_v)$
telles que $\pi_{v,j} (\delta)$ et $\pi'_{v,k} (\delta)$ soient unitaires
\end{itemize}
vérifiant les conditions suivantes:

\begin{itemize}
\item pour tout $0 \leq i < d-1$, $\sum_{k~/~ t'_k=d-1-i} g'_k=n_i$;

\item $\Pi_v$ est l'induite irréductible $[\overleftarrow{t_1-1}]_{\pi_{v,1}} \times \cdots \times
[\overleftarrow{t_u-1}]_{\pi_{v,u}} \times [\overrightarrow{t_1'-1}]_{\pi'_{v,1}} \times \cdots \times
[\overrightarrow{t'_{u'}-1}]_{\pi'_{v,u'}}$, avec:

\begin{itemize}
\item si $u >0$ alors pour tout $j$, $s_j \equiv 1 \mod 2$;

\item si $u=0$ alors tous les $s_i$ ont la même parité donnée par $(-1)^{s_i-1}=\e(\Pi)$.
\end{itemize}

\end{itemize}
\end{itemize}
\end{propb}

\begin{corob} Toute représentation automorphe de $G_\tau(\Am)$ cohomologique, vérifie la conjecture de
Ramanujan-Peterson.
\end{corob}

\begin{corob} \label{coro-compo-locale}
Soit $\Pi$ une représentation irréductible automorphe de $G_\tau(\Am)$ cohomologique et supposons qu'il
existe une place $v_0$ telle que $G_{v_0}^\times \simeq GL_d(F_{v_0})$ et $\Pi_{v_0}$ tempérée, i.e. avec les
notations du théorème précédent $t'_i=1$ pour tout $1 \leq i \leq u'$ (resp. $\Pi_{v_0} \simeq
[\overrightarrow{s-1}]_{\pi_{v_0}}$ pour $\pi_{v_0}$ une représentation cuspidale de $GL_g(F_{v_0})$ avec
$d=sg$). On en déduit alors, en utilisant les notations du théorème précédent, que pour toute place $v$ non
ramifiée, pour tout $1 \leq k \leq u'$ (resp. pour tout $1 \leq j \leq u$) $t'_k=0$ (resp. $t_j=0$), i.e.
$\Pi_v$ est tempérée et $n_i=0$ pour tout $0 \leq i <d-1$ (resp. $n_i=0$ pour $i \neq s$ et $\Pi_v \simeq
[\overrightarrow{s-1}]_{\pi'_{v,1} \times \cdots \times \pi'_{v,u'}}$).
\end{corob}

\begin{proof} \textit{de la proposition (\ref{prop-compo-locale}):} A torsion par un caractère près,
$\Pi_v$ est unitaire et donc d'après la classification de Tadic de la forme
$[\overrightarrow{s_1-1}]_{[\overleftarrow{t_1-1}]_{\pi_{v,1}(\lambda_1)}} \times \cdots \times
[\overrightarrow{s_u-1}]_{[\overleftarrow{t_u-1}]_{\pi_{v,u}(\lambda_u)}}$ où $\pi_{v,i}$ sont des
représentations irréductibles cuspidales unitaires de $GL_{g_i}(F_v)$ avec $\lambda_i \in ]
\frac{-1}{2},\frac{1}{2}[$ et $\sum_{i=1}^u g_is_it_i=d$. L'induite étant irréductible, l'ordre des facteurs
n'importe pas. Soit alors $r \geq 1$ tel que $\pi_v:=\pi_{v,1}(\lambda_1) \simeq \cdots \simeq
\pi_{v,r}(\lambda_r)$ et $\pi_{v,i}(\lambda_i)$ n'est pas isomorphe à $\pi_v$ pour $i
>r$. On note $t_0:=max_{1 \leq i \leq r} \{ s_i,t_i \}$. On note $g$ l'entier tel que $\pi_v$ est une
représentation cuspidale de $GL_g(F_v)$. Par ailleurs soit $0\leq k_1$ (resp. $k_1 \leq k_2$, resp. $k_2 \leq
k_3 \leq r$) tel que pour tout $1 \leq i \leq k_1$ (resp. $k_1 < i \leq k_2$, resp. $k_2 < i \leq k_3$) on
ait $s_i=t_i=t_0$ (resp. $t_i < s_i=t_0$, resp. $s_i < t_i=t_0$) et $\max_{k_3 < i \leq r} \{ t_i,s_i \} <
t_0$.

\begin{lemb} \label{lem-red1}
Pour tout $\max \{ s ,t \}< r \leq st$, $\red_{\JL^{-1}([\overleftarrow{r-1}]_{\pi_v})}
([\overrightarrow{s-1}]_{[\overleftarrow{t-1}]_{\pi_v}})=(0)$.
\end{lemb}

\begin{proof} Le résultat découle essentiellement de \cite{ze} et en particulier du lemme suivant.

\begin{lemb} \label{lem-ze}
Pour tout $s,t$, les constituants de $J_{N_{g,2g,\cdots,stg}^{op}}
([\overrightarrow{s-1}]_{[\overleftarrow{t-1}]_{\pi_v}})$ sont de la forme $\pi_v(\frac{2-t-s}{2}+\s(1))
\otimes \cdots \otimes \pi_v(\frac{2-t-s}{2}+\s(st))$ où le multi-ensemble $R:=\{ \s(1),\cdots, \s(st) \}$
est tel que $\{ \pi_v(\frac{2-t-s}{2}+r)~/~r \in R \}$ soit égal au support cuspidal de
$[\overrightarrow{s-1}]_{[\overleftarrow{t-1}]_{\pi_v}}$. Pour tout $0 \leq i \leq st-1$, on pose
$\s^{-1}(i)=\{ n_i(1) < n_i(2) < \cdots \}$. On a alors les conditions suivantes:

\noindent (i) pour $0 \leq k <s-1$ et $0 \leq i < t$, on a $n_{t-1+k-i}(\min(k+1,i+1)) >
n_{t+k-i}(\min(k+2,i+1))$;

\noindent (ii) pour $0 \leq k < s$ et $0 \leq i < j <t$, on a $n_{t-1+k-i}(\min(k+1,i+1)) >
n_{t-1+k-j}(\min(k+1,j+1))$.
\end{lemb}

\rem Quand on regarde le support cuspidal de
$\overbrace{[\overleftarrow{t-1}]_{\pi_v(\frac{1-s}{2})} \boxplus \cdots
[\overleftarrow{t-1}]_{\pi_v(\frac{s-1}{2})}}^s$ dans l'ordre de gauche à droite, pour tout $0 \leq k < s-1$
et $0 \leq i < t$, le $\pi_v(\frac{2-t-s}{2}+t-1+k-i)$ du $k$-ème facteur est le $\min(k+1,i+1)$-ème. La
condition (ii) du lemme précédent affirme que dans un même facteur, on doit prendre les $\pi_v(i)$ de gauche
à droite, tandis que la condition (i) précise que le $\pi_v(i)$ du $k$-ème facteur arrive après le
$\pi_v(i+1)$ du $k+1$-ème facteur.

\begin{proof} Comme $([\overrightarrow{s-1}]_{[\overleftarrow{t-1}]_{\pi_v}})$ s'injecte dans les
induites $[\overleftarrow{t-1}]_{\pi_v(\frac{1-s}{2})} \times
[\overrightarrow{s-2}]_{[\overleftarrow{t-1}]_{\pi_v(1/2)}}$ et $[\overleftarrow{t-1}]_{\pi_v(\frac{1-s}{2})}
\boxplus [\overleftarrow{t-1}]_{\pi_v(\frac{3-s}{2})} \times
\overrightarrow{(s-2)}_{[\overleftarrow{t-1}]_{\pi_v(1)}}$, on se ramène aisément, par récurrence sur $s$, au
cas $s=2$ soit à étudier $[\overleftarrow{t-1}]_{\pi_v(-1/2)} \boxplus [\overleftarrow{t-1}]_{\pi_v(1/2)}$.
Or dans le groupe de Grothendieck, cette dernière est égale à $[\overleftarrow{t-1}]_{\pi_v(-1/2)} \times
[\overleftarrow{t-1}]_{\pi_v(1/2)} - [\overleftarrow{t}]_{\pi_v} \times [\overleftarrow{t-2}]_{\pi_v}$. La
condition (ii) constitue l'unique restriction sur les $\s$ pour l'induite totale
$[\overleftarrow{t-1}]_{\pi_v(-1/2)} \times [\overleftarrow{t-1}]_{\pi_v(1/2)}$. Il suffit alors de remarquer
que les $\s$ vérifiant (i) et pas (ii) sont obtenues à partir de $[\overleftarrow{t}]_{\pi_v} \times
[\overleftarrow{t-2}]_{\pi_v}$.

\end{proof}

Ainsi si on obtenait un constituant de la forme $[\overleftarrow{r-1}]_{\pi_v} \otimes \pi$, on en déduirait,
par transitivité du foncteur de Jacquet, avec les notations du lemme précédent, que la bijection $\s_0$ telle
que $\s_0(i)=\s(1)-i+1$ pour $1 \leq i \leq r$ vérifierait les conditions (i) et (ii) ce qui impose $r \leq
t$.

De même si on obtenait $[\overrightarrow{r-1}]_{\pi_v} \otimes \pi$, la bijection $\s_0$ telle que
$\s_0(i)=\s_0(1)+i-1$ pour $1 \leq i \leq r$ vérifierait les conditions (i) et (ii) ce qui impose $r \leq s$,
d'où le résultat.

\end{proof}

On en déduit alors que pour tout $t > t_0$, $\red_{\JL^{-1}([\overleftarrow{t-1}]_{\pi_v^\vee})} (\Pi_v)$ est
nul de sorte que d'après (\ref{somme-alternee}), que pour tout $t_0 < t \leq s$, $\sum_i (-1)^i H^i(j^{\geq
tg}_! HT_\xi(g,t,\pi_v,\Pi_t))[\Pi^\oo]=0$. D'après le lemme (\ref{lem-hij0}), on en déduit alors que pour
tout $i$, $H^i(j^{\geq tg}_{!*} HT_\xi(g,t,\pi_v,\Pi_t))[\Pi^\oo]$ est nul. Ainsi en reprenant l'argument de
la preuve du lemme (\ref{lem-red1}) et en utilisant le lemme (\ref{lem-ze}), on obtient:

\begin{lemb} \label{lem-red2}
Soit $r:=\max \{ s,t \}$. Si $r>t$ (resp. $r>s$), alors $\red_{\JL^{-1}([\overleftarrow{r-1}]_{\pi_v^\vee})}
([\overrightarrow{s-1}]_{[\overleftarrow{t-1}]_{\pi_v}})$ est isomorphe à $(-1)^{r-1}
\Xi^{\frac{(t-1)(sg-1)}{2}} \otimes [\overrightarrow{s-1}]_{[\overleftarrow{t-2}]_{\pi_v(-\frac{sg-1}{2})}}$
(resp. $\Xi^{\frac{(s-1)(tg+1)}{2}} \otimes
[\overrightarrow{s-2}]_{[\overleftarrow{t-1}]_{\pi_v(-\frac{tg+1}{2})}})$) Dans le cas $r=s=t$, on obtient la
somme des deux termes ci-dessus.
\end{lemb}

Ainsi d'après (\ref{somme-alternee}),
$\frac{\e(\Pi)}{\sharp \ker^1(\Qm,G_\tau)} \sum_i (-1)^i H^i (j^{\geq t_0g}_!
HT_\xi(g,t_0,\pi_v,\Pi_l)[d-t_0g])[\Pi^\oo]$ qui d'après ce qui précède est égal à $\frac{\e(\Pi)}{\sharp
\ker^1(\Qm,G_\tau)}\sum_i (-1)^i
H^i(j^{\geq t_0g}_{!*} HT_\xi(g,t,\pi_v,\Pi_l)[d-t_0g])$, est donné par
\begin{multline*}
\Xi^{\frac{d-t_0g}{2}} ((\Xi^{\frac{t_0-1}{2}} (-1)^{t_0-1} \Pi_l (\frac{t_0-1}{2}) \times
([\overrightarrow{t_0-2}]_{[\overleftarrow{t_0-1}]_{\pi_v(-1/2)}})^{k_1}+ \\
\Xi^{-\frac{t_0-1}{2}} (\Pi_t (\frac{1-t_0}{2}) \times
[\overrightarrow{t_0-1}]_{[\overleftarrow{t_0-2}]_{\pi_v(1/2)}})^{k_1}) \times
[\overrightarrow{s_2-1}]_{[\overleftarrow{t_2-1}]_{\pi_{v,2}(\lambda_2)}} \times \cdots \times
[\overrightarrow{s_u-1}]_{[\overleftarrow{t_u-1}]_{\pi_{v,u}(\lambda_u)}} + \\ \sum_{i=k_1+1}^{k_2}
\Xi^{-\frac{t_i-1}{2}} (-1)^{t_0-1} ([\overrightarrow{s_1-1}]_{[\overleftarrow{t_1-1}]_{\pi_{v,1}
(\lambda_1)}} \times \cdots \times
[\overrightarrow{s_{i-1}-1}]_{[\overleftarrow{t_{i-1}-1}]_{\pi_{v,i-1}(\lambda_{i-1})}} \times \\
\Pi_t (\frac{1-t_i}{2}) \times [\overrightarrow{t_0-1}]_{[\overleftarrow{t_i-2}]_{\pi_v(-1/2)}} \times
[\overrightarrow{s_{i+1}-1}]_{[\overleftarrow{t_{i+1}-1}]_{\pi_{v,i+1}(\lambda_{i+1})}} \times \cdots \times
[\overrightarrow{s_u-1}]_{[\overleftarrow{t_u-1}]_{\pi_{v,u}(\lambda_u)}}) +
\end{multline*}
\begin{multline*}
+ \sum_{i=k_2+1}^{k_3} \Xi^{\frac{s_i-1}{2}}
(([\overrightarrow{s_1-1}]_{[\overleftarrow{t_1-1}]_{\pi_{v,1}(\lambda_1)}} \times \cdots \times
[\overrightarrow{s_{i-1}-1}]_{[\overleftarrow{t_{i-1}-1}]_{\pi_{v,i-1}(\lambda_{i-1})}}
\times \\
\Pi_t (\frac{s_i-1}{2}) \times [\overrightarrow{s_i-2}]_{[\overleftarrow{t_0-1}]_{\pi_v(1/2)}}) \times
[\overrightarrow{s_{i+1}-1}]_{[\overleftarrow{t_{i+1}-1}]_{\pi_{v,i+1}(\lambda_{i+1})}} \times \cdots \times
[\overrightarrow{s_u-1}]_{[\overleftarrow{t_u-1}]_{\pi_{o,u}(\lambda_u)}})
\end{multline*}

En utilisant que pour tout $i$, $H^i(j^{\geq t_0g}_{!*} HT_\xi(g,t_0,\pi_v,\Pi_l)[d-t_0g])$ est pur de poids
$d-t_0g+i$, on en déduit alors que ceux-ci sont alors isomorphes à
\begin{multline*} \Pi_t (\frac{t_0-1}{2}) \times
([\overrightarrow{t_0-1}]_{[\overleftarrow{t_0-2}]_{\pi_v(-1/2)}})^{k_1} \times
[\overrightarrow{s_2-1}]_{[\overleftarrow{t_2-1}]_{\pi_{o,2}(\lambda_2)}} \times \cdots \times
[\overrightarrow{s_u-1}]_{[\overleftarrow{t_u-1}]_{\pi_{o,u}(\lambda_u)}} \\ \hbox{ pour } i=1-t_0
\end{multline*}
\begin{multline*}
\Pi_t (\frac{1-t_0}{2}) \times ([\overrightarrow{t_0-2}]_{[\overleftarrow{t_0-1}]_{\pi_v(1/2)}})^{k_1} \times
[\overrightarrow{s_2-1}]_{[\overleftarrow{t_2-1}]_{\pi_{v,2}(\lambda_2)}} \times \cdots \times
[\overrightarrow{s_u-1}]_{[\overleftarrow{t_u-1}]_{\pi_{v,u}(\lambda_u)}} \\ \hbox{ pour } i=t_0-1
\end{multline*}
\begin{multline*}
\bigoplus_{k_1<k \leq k_2~/~1-t_k=i} (-1)^{t_0-1}
([\overrightarrow{s_1-1}]_{[\overleftarrow{t_1-1}]_{\pi_{v,1}(\lambda_1)}} \times \cdots \times
[\overrightarrow{s_{i-1}-1}]_{[\overleftarrow{t_{i-1}-1}]_{\pi_{v,i-1}(\lambda_{i-1})}} \times  \Pi_t
(\frac{1-t_i}{2}) \times \\ [\overrightarrow{t_0-1}]_{[\overleftarrow{t_i-2}]_{\pi_v(-1/2)}} \times
[\overrightarrow{s_{i+1}-1}]_{[\overleftarrow{t_{i+1}-1}]_{\pi_{v,i+1}(\lambda_{i+1})}} \times \cdots \times
[\overrightarrow{s_u-1}]_{[\overleftarrow{t_u-1}]_{\pi_{v,u}(\lambda_u)}}) \\ \hbox{ pour } 0 \leq -i < t_0-1
\end{multline*}
\begin{multline*}
\bigoplus_{k_2<k \leq k_3 ~/~s_k-1=i}
(([\overrightarrow{s_1-1}]_{[\overleftarrow{t_1-1}]_{\pi_{v,1}(\lambda_1)}} \times \cdots \times
[\overrightarrow{s_{i-1}-1}]_{[\overleftarrow{t_{i-1}-1}]_{\pi_{v,i-1}(\lambda_{i-1})}} \times \Pi_t
(\frac{s_i-1}{2}) \\ \times [\overrightarrow{s_i-2}]_{[\overleftarrow{t_0-1}]_{\pi_v(1/2)}}) \times
[\overrightarrow{s_{i+1}-1}]_{[\overleftarrow{t_{i+1}-1}]_{\pi_{v,i+1}(\lambda_{i+1})}} \times \cdots \times
[\overrightarrow{s_u-1}]_{[\overleftarrow{t_u-1}]_{\pi_{v,u}(\lambda_u)}}) \\ \hbox{ pour } 0 \leq i <t_0-1
\end{multline*}

Si on veut que ces résultats soient compatibles à la dualité de Verdier, il faut alors que pour tout $i$,
$\min \{ s_i,t_i \}=1$. En outre s'il existe $i$ tel que $t_i>1$ alors $\e(\Pi)=1$ et pour tout $j$, $s_j
\equiv 1 \mod 2$. Si pour tout $i$, $t_i=1$ alors tous les $s_i$ ont la même parité donnée par
$(-1)^{s_i-1}=\e(\Pi)$.

Par ailleurs les $H^i(j^{\geq tg}_{!*} HT_\xi(g,t,\pi_v,[\overleftarrow{t-1}]_{\pi_v})[d-tg]) \otimes
\rec_{F_v}^\vee(\pi_v) (1-tg)$ servent à calculer la cohomologie de la fibre générique via la suite spectrale
des poids. On en déduit alors que tous les poids $\rec_{F_v}^\vee(\pi_{v,i}(\lambda_i))$ s'obtiennent à
partir de ceux obtenus par la cohomologie de la fibre générique par une translation entière, d'où le
résultat.

\end{proof}

\section{Pureté de la filtration de monodromie des $H^i_{\eta_v}$}

\begin{theob} \label{theo-mono-glob}
Pour toute représentation irréductible $\Pi$ et pour tout $i$, les gradués, $gr_{k,\xi}[\Pi^\oo]$ de la
filtration de monodromie de $H^i_{\eta_v,\xi}[\Pi^\oo]$ sont purs de poids $d-1+k$ \footnote{plus le poids de
$\xi$ que l'on supposera nul}. En fait pour $i \neq d-1$, l'opérateur de monodromie $N$ agit trivialement sur
les $H^i_{\eta_v,\xi}$.
\end{theob}

\begin{proof} Soit $v$ une place non ramifiée. Remarquons tout d'abord que s'il n'existe pas de $\Pi_v'$
\footnote{On peut en outre supposer que ce dernier soit de même support cuspidal que $\Pi_v$} tel que $\Pi_v'$
ne vérifie pas $\hyp(\xi)$, alors tous les $H^i_{\eta_v,\xi}[\Pi^\oo]$ sont nuls. Considérons alors $\Pi$
irréductible vérifiant $\hyp(\xi)$ de sorte que, d'après ce qui précède, $\Pi_v$ est de la forme
$$[\overrightarrow{s_1-1}]_{\pi_{v,1}} \times \cdots \times [\overrightarrow{s_u-1}]_{\pi_{v,u}} \times
[\overleftarrow{t_1-1}]_{\pi_{v,1}'} \times \cdots \times [\overleftarrow{t_{u'}-1}]_{\pi_{v,u'}'}$$ avec
$\pi_{v,i}$ et $\pi_{v,j}'$ irréductibles cuspidales unitaires de respectivement $GL_{g_i}(F_v)$ et
$GL{g_j'}(F_v)$. \footnote{En outre tous les $s_i$ ont la même congruence modulo $2$.} On étudie alors la
suite spectrale
\begin{equation} \label{ss-gr-3}
E_1^{i,j}=H^{i+j}(gr_{-i,\xi}) \Rightarrow H^{i+j+d-1}_{\eta_v,\xi}
\end{equation}
D'après le théorème (\ref{theo-global1}), les $H^r(gr_{k,\xi})[\Pi^{\oo,v}]$ se calculent à partir des
$H^r(\PC(g_i,t,\pi_{v,i}))[\Pi^{\oo,v}]$ pour $1 \leq i \leq u$ et $1 \leq tg_i \leq d$ ainsi que des
$H^r(\PC(g_i',t,\pi_{o,i}'))[\Pi^{\oo,v}]$ pour $1 \leq i \leq u'$ et $1 \leq tg_i' \leq d$. Le calcul de ces
derniers s'obtient par récurrence descendante sur $t$ via le calcul des
$\red_{\JL^{-1}([\overleftarrow{t-1}]_{\pi_{v,i}})}^{tg_i}(\Pi_v)$.

\marque A priori, pour une preuve rigoureuse on devrait introduire des notations comme dans la preuve de la
proposition (\ref{prop-compo-locale}) en regroupant les $\pi_{v,i}$ qui sont isomorphes, puis en les triant
selon leur $s_i$. Cependant comme le foncteur de Jacquet est additif sur les induites et vu que les
$H^i(j^{\geq tg}_{!*} HT_\xi(g,t,\pi_v,\Pi_l))[\Pi^{\oo,v}]$ se calculent à partir de ceux-ci, par
récurrence, on remarque alors que les divers $\pi_{v,i}$ n'interagissent pas entre eux de sorte que chacun
des $H^i(j^{\geq tg}_{!*} HT_\xi(g,t,\pi_v,\Pi_l)[d-tg]) [\Pi^{\oo,v}]$ sera de la forme
$$\bigoplus_{\atop{k/\pi_{v,k} \simeq \pi_v}{s_k \geq t}} H(i,\pi_v,t,s_k,-) \oplus
\bigoplus_{\atop{k/\pi_{v,k}' \simeq \pi_v}{t_k \geq t}} H(i,\pi_v,t,t_k,+)$$ où $H(i,\pi_v,t,s_k,-)$ (resp.
$H(i,\pi_v,t,t_k,+)$) correspond au calcul de $H^i(j^{\geq tg}_{!*}
HT_\xi(g,t,\pi_v,\Pi_t)[d-tg])[\Pi^{\oo,v}]$ effectué comme si $k$ est tel que pour tout $1 \leq i \neq k
\leq u$ et pour tout $1 \leq j \leq u'$ (resp. pour tout $1 \leq i \leq u$ et $1 \leq j \neq k \leq u'$),
$\pi_{v,i}$ et $\pi_{v,j}'$ ne sont pas isomorphes à $\pi_v$. Plaçons nous alors dans une telle situation,
i.e. $\Pi_v \simeq [\overleftarrow{s-1}]_{\pi_v} \times \pi_{v,1}$ (resp. $\Pi_v \simeq
[\overrightarrow{s-1}]_{\pi_v} \times \pi_{v,1}$) où $\pi_{v,1}$ est une représentation irréductible de
$GL_{d-sg}(F_v)$ de support cuspidal disjoint de la droite associée à $\pi_v$ irréductible cuspidale de
$GL_g(F_v)$. Les calculs sont alors exactement similaires à ceux déjà effectués, rappelons rapidement comment
cela s'articule.

- pour $t > s$, les $H^i(j^{\geq tg}_{!*} HT_\xi(g,t,\pi_v,\Pi_t))[\Pi^{\oo,v}]$ sont nuls: on procède par
récurrence descendante sur $t$ de $s_g$ à $s+1$. Le cas $t=s_g$ quand $s_g g =d$ découle du lemme
(\ref{lem-pts-ss}) en utilisant une correspondance de Jacquet-Langlands globale comme au théorème
(\ref{theo-jl}). Dans les autres situations le résultat découle du fait que
$\red_{\JL^{-1}([\overleftarrow{t-1}]_{\pi_v})}^{tg} (\Pi_v)$ est nul;

- pour $t=s$, on a $H^0(j^{\geq sg}_{!*} HT_\xi(g,s,\pi_v,[\overleftarrow{s-1}]_{\pi_v})[d-sg])[\Pi^{\oo,v}]
\simeq \Pi_v$ \footnote{Dans le cas respé on utilise que $(-1)^{s-1}=\e(\Pi)$.}: le résultat découle
directement du point précédent et du fait que
$\red_{\JL^{-1}([\overleftarrow{s-1}]_{\pi_v^\vee})}^{sg}(\Pi_v)=\pi_{v,1}(-sg/2) \otimes
\Xi^{\frac{d-sg}{2}}$;

- pour $1 \leq t < s$ et $\Pi_v \simeq [\overleftarrow{s-1}]_{\pi_v} \times \pi_{v,1}$, tous les $H^i(j^{\geq
tg}_{!*} HT(g,t,\pi_v,[\overleftarrow{t-1}]_{\pi_v}))[\Pi^{\oo,v}]$ sont nuls: la preuve est strictement
identique à celle de la proposition (\ref{prop-coho1});

- pour $1 \leq t < s$ et $\Pi_v \simeq [\overrightarrow{s-1}]_{\pi_v} \times \pi_{v,1}$, les $H^i(j^{\geq
tg}_{!*} HT(g,t,\pi_v,[\overleftarrow{t-1}]_{\pi_v})[d-tg])[\Pi^{\oo,v}]$ sont nuls pour $|i| > s-t$ ou $i
\equiv s-t-1 \mod 2$. Dans les autres cas, il est isomorphe à
\begin{multline*}
\sharp \ker^1(\Qm,G_\tau)m(\Pi)([\overleftarrow{t-1}]_{\pi_v} \overrightarrow{\times}
[\overrightarrow{\frac{s-t+i - 4}{2}}]_{\pi_v}) \overleftarrow{\times}
[\overrightarrow{\frac{s-t-i-4}{2}}]_{\pi_v}  \otimes (\Xi^{\frac{(s-t)g+i}{2}} \otimes \bigoplus_{\chi \in
\AF(\pi_v)} \chi^{-1})
\end{multline*}
La preuve est identique à celle de la proposition
(\ref{prop-not}).

\marque On revient alors à l'étude de la suite spectrale (\ref{ss-gr-3}) dont on vient de déterminer les
termes $E_1^{i,j}[\Pi^{\oo,v}]$. La détermination des flèches $d_1^{i,j}$ se fait alors comme dans la preuve
de la proposition (\ref{prop-lrs2}) car à nouveau les divers $\pi_{v,i}$ n'interagissent pas entre eux dans
les arguments de dualité. Considérons en effet le cas de deux contributions: si celles-ci correspondent à des
Steinberg généralisées alors tout est concentré sur les $H^0$ de sorte que toutes les flèches sont nulles.
Supposons donc avoir à traiter le cas de deux $\speh$ \footnote{Le cas d'une $\speh$ avec une $\st$ se traite
de la même façon.}. Si elles sont de même longueur, on se retrouve alors dans la situation de la preuve de
(\ref{prop-lrs2}) mais avec des multiplicité $2$, ce qui ne change en rien les arguments. Supposons les alors
de longueur distinctes, $t'<t$ et reprenons les arguments de la preuve de (\ref{prop-lrs2}). Tant qu'on
étudie des $H^i(gr_k)$ avec $\min \{ |i|, |k| \}  >t$, il n'y a rien à rajouter car n'intervient que la
composante de plus grande longueur. Dans les autres cas, l'argument consiste à comparer par dualité les
termes $H^i$ et $H^{-i}$. Pour $i>0$, $H^i$ (resp. $H^{-i}$) est une somme directe constituée de
représentations de la forme
\begin{equation} \label{type1}
\sharp \ker^1(\Qm,G_\tau)m(\Pi)[\overrightarrow{t-1}]_{\pi_v} \times ([\overrightarrow{r_1-1}]_{\pi_v}
\overrightarrow{\times}
[\overleftarrow{t'-r_1-r_2-1}]_{\pi_v} \overrightarrow{\otimes} [\overrightarrow{r_2-1}]_{\pi_v})
\end{equation}
avec $r_2 \leq t'/2$ (resp. $r_1 \leq t'/2$) \footnote{On peut être plus précis, cf. la proposition
(\ref{prop-not}).}
\begin{equation} \label{type2}
\sharp \ker^1(\Qm,G_\tau) m(\Pi) ([\overrightarrow{r_1-1}]_{\pi_v} \overrightarrow{\times}
[\overleftarrow{t-r_1-r_2-1}]_{\pi_v}
\overrightarrow{\times} [\overrightarrow{r_2-1}]_{\pi_v} ) \times [\overrightarrow{t'-1}]_{\pi_v}
\end{equation}
avec $r_2 \leq t/2$ (resp. $r_1 \leq t/2$). On remarque alors que les composants de (\ref{type1}) sont
disjoints de ceux de (\ref{type2}), sauf éventuellement $[\overrightarrow{t-1}]_{\pi_v} \times
[\overrightarrow{t'-1}]_{\pi_v}$. Ce sont donc les seuls qui peuvent rester dans l'aboutissement \footnote{On
le savait déjà car seules les représentations automorphes subsistent dans l'aboutissement.}. Pour voir que
ces derniers restent effectivement dans l'aboutissement il suffit de remarquer qu'ils n'apparaissent que dans
les termes $H^i(gr_{0,\xi})[\Pi^{\oo,v}]$ de sorte qu'il ne peut y avoir de simplifications, d'où le résultat.

\end{proof}


\part{Récapitulatifs}

\section{sur les faisceaux}


\marque - $R \Psi_v(\LC_{\xi}) \simeq R \Psi_v(\bar \Qm_l) \otimes \LC_{\xi}$;

- $R\Psi_v[d-1]$ est un objet de $\FPH(\bar X)$;

- $R\Psi_v[d-1]=\bigoplus_{\pi_v \in \cusp(d)} R\Psi_{\pi_v}[d-1]$ avec les notations de la proposition
(\ref{prop-so}), où $\cusp(d)$ désigne l'ensemble des classes d'équivalence inertielle des représentations
irréductibles cuspidales de $GL_g(F_v)$ pour tout $1 \leq g \leq d$;

- $R\Psi_{\pi_v}[d-1]$ est mixte de poids $d-1$; sa filtration de monodromie est égale à celle par le poids
au décalage de $d-1$ près, on note $gr_{k,\pi_v}$ le gradué de poids $k$ de sa filtration de monodromie de
sorte qu'en particulier $N^k: gr_{k,\pi_v} \longmapright{\sim} gr_{-k,\pi_v}$;

- pour $1 \leq g \leq d$ et $\pi_v$ une représentation irréductible cuspidale de $GL_g(F_v)$; pour tout $|k|
< s_g:=\lfloor d/g \rfloor$, dans $\FPH(\bar X)$ on a (cf. \ref{theo-global1}):
$$e_{\pi_v} gr_{k,\pi_v} = \bigoplus_{\atop{|k| < t \leq s_g}{t \equiv k-1 \mod 2}} 
\PC(g,t,\pi_v)(-\frac{tg+k-1}{2})$$
où $\PC(g,t,\pi_v):=j^{\geq tg}_{!*} HT(g,t,\pi_v,[\overleftarrow{t-1}]_{\pi_v})[d-tg] \otimes
\rec_{F_v}^\vee(\pi_v)$;

- la suite spectrale des poids (\ref{suite-spectrale}) dégénère en $E_2$ et les applications $d_1^{i,j}$ sont
induites par les suites exactes (\ref{suites-exactes}).

\marque $j^{\geq tg}_{!*} HT(g,t,\pi_v,\Pi_t)[d-tg] \in \FPH(\bar X)$ est pur de poids $d-tg$ pour $\pi_v$
unitaire; ses faisceaux de cohomologies $h^i j^{\geq tg}_{!*} HT(g,t,\pi_v,\Pi_t)[d-tg]$ sont, pour $g>1$
\footnote{pour $g=1$ cf. (\ref{theo-global2})}, nuls pour $i$ ne s'écrivant pas sous la forme $tg-d+a(g-1)$
avec $0 \leq a \leq s_g-t$ et pour un tel $i=tg-d+a(g-1)$, isomorphe dans $\FH(\bar X)$ à
$$j^{\geq (t+a)g}_!HT(g,t+a,\pi_v,\Pi_t \overrightarrow{\times} [\overrightarrow{a-1}]_{\pi_v}) \otimes
\Xi^{\frac{a(g-1)}{2}}$$

\marque Pour tout $1 \leq t \leq s_g$, on a, cf. (\ref{prop-p}), l'égalité dans $\GF$
$$j^{\geq tg}_! HT(g,t,\pi_v,\Pi_t)[d-tg] = \sum_{i=t}^{s_g} j^{\geq ig}_{!*} HT(g,i,\pi_v,\Pi_t 
\overrightarrow{\times}
[\overleftarrow{i-t-1}]_{\pi_v})[d-ig] \otimes \Xi^{\frac{(t-i)(g-1)}{2}}$$ Dualement on a

$$Rj^{\geq tg}_* HT(g,t,\pi_v,\Pi_t)[d-tg] = \sum_{i=t}^{s_g} j^{\geq ig}_{!*} HT(g,i,\pi_v,\Pi_t 
\overleftarrow{\times}
[\overleftarrow{i-t-1}]_{\pi_v})[d-ig] \otimes \Xi^{\frac{(t-i)(g+1)}{2}} $$ On en déduit alors que, au moins
pour $g>1$,

\begin{multline*}
R^i j^{\geq tg}_* HT(g,t,\pi_v,\Pi_t)[d-tg]=\bigoplus_{(r,a)~/~i=(t+r+a)g-d-a} \\
 j^{\geq (t+r+a)g}_! HT(g,t+r+a,\pi_v,(\Pi_t \overleftarrow{\times} [\overleftarrow{r-1}]_{\pi_v})
\overrightarrow{\times} [\overrightarrow{a-1}]_{\pi_v}) \otimes \Xi^{((r+a)(g-1)+2r)/2}
\end{multline*}

\section{sur les calculs de groupe de cohomologie}


\marque \textit{Cas où $\Pi_v$ est de la forme $\st_s(\pi_v)$}:

\begin{itemize}
\item[(1)] d'après la proposition (\ref{prop-coho1}),
$H^i(j^{\geq tg}_{!*} HT_{\xi}(g,t,\pi_v,\Pi_t))[\Pi^\oo]$ est nul pour tout $i$.

\item[(2)] d'après le corollaire (\ref{coro-hij-nul}),
$H^i(j^{\geq tg}_{!} HT_{\xi}(g,t,\pi_v,\Pi_t))[\Pi^\oo]$ est donné par
$$\left \{ \begin{array}{cl} 0 & \hbox{si } i \neq 0 \\
\sharp \ker^1(\Qm,G_\tau) m(\Pi) \Pi_t \overrightarrow{\times} [\overleftarrow{s-t-1}]_{\pi_v} \otimes ( 
\Xi^{\frac{(s-t)(g-1)}{2}}
\otimes \bigoplus_{\chi \in \AF(\pi_v)} \xi^{-1}) & i=0 \end{array} \right.$$ où $\AF(\pi_v)$ est l'ensemble
des caractères $\chi:\Zm \longto \Qm_l^\times$, tels que $\pi_v \otimes \chi^{-1} \circ \val(\det) \simeq
\pi_v$.

\item[(3)] par application de la dualité de Verdier, $H^i(Rj^{\geq tg}_{*} HT_{\xi}(g,t,\pi_v,\Pi_t))[\Pi^\oo]$
est donné par
$$\left \{ \begin{array}{cl} 0 & \hbox{si } i \neq 0 \\
\sharp \ker^1(\Qm,G_\tau) m(\Pi)\Pi_t \overleftarrow{\times} [\overleftarrow{s-t-1}]_{\pi_v} \otimes 
(\Xi^{(s-t)(g+1)/2} \otimes
\bigoplus_{\chi \in \AF(\pi_v)} \chi^{-1}) & i=0
\end{array} \right.$$
\end{itemize}

\marque \textit{Cas où $\Pi_v$ est de la forme $\speh_s(\pi_v)$}:

\begin{itemize}
\item[(1)] d'après la proposition (\ref{prop-not}), $H^i(j^{\geq tg}_{!*}
HT_{\xi}(g,t,\pi_v,\Pi_t))[\Pi^\oo]$ est donné par

$$\left \{
\begin{array}{l} \sharp \ker^1(\Qm,G_\tau) m(\Pi) \Pi_t \overrightarrow{\times} [\overrightarrow{s-t-1}]_{\pi_v} 
\otimes
(\Xi^{\frac{(s-t)(g-1)}{2}} \otimes \bigoplus_{\chi \in \AF(\pi_v)} \chi^{-1}) \\ \hfill i=l-s \\
\sharp \ker^1(\Qm,G_\tau) m(\Pi) (\Pi_t \overrightarrow{\times} [\overrightarrow{\frac{s-t-i}{2}-1}]_{\pi_v}) 
\overleftarrow{\times}
[\overrightarrow{\frac{s-t+i}{2}-1}]_{\pi_v} \otimes (\Xi^{\frac{(s-t)g+i}{2}} \otimes \bigoplus_{\chi \in
\AF(\pi_v)} \chi^{-1}) \\ \hfill \left \{ \begin{array}{l} l-s < i < s-t \\  i  \equiv l-s \mod 2 \end{array} \right. 
\\
\sharp \ker^1(\Qm,G_\tau) m(\Pi) \Pi_t \overleftarrow{\times}
[\overrightarrow{s-t-1}]_{\pi_v} \otimes (\Xi^{\frac{(s-t)(g+1)}{2}} \otimes \bigoplus_{\chi \in \AF(\pi_v)}
\chi^{-1}) \\ \hfill i=s-t
\end{array} \right. $$

\item[(2)] d'après le corollaire (\ref{strate-alt-p}),
$H^i(j^{\geq tg}_{!} HT_{\xi}(g,t,\pi_v,\Pi_t))[\Pi^\oo]$ est donné par
$$\left \{ \begin{array}{cl} 0 & \hbox{si } i \neq s-t \\
\sharp \ker^1(\Qm,G_\tau) m(\Pi) \Pi_t \overleftarrow{\times} [\overrightarrow{s-t-1}]_{\pi_v} \otimes 
(\Xi^{\frac{(s-t)(g+1)}{2}}
\bigoplus_{\chi \in \AF(\pi_v)} \chi^{-1}) & i=s-t \end{array} \right.$$

\item[(3)] par application de la dualité de Verdier, $H^i(Rj^{\geq tg}_{*} HT_{\xi}(g,t,\pi_v,\Pi_t))[\Pi^\oo]$ est 
donné par
$$\left \{ \begin{array}{cl} 0 & \hbox{si } i \neq 0 \\
\sharp \ker^1(\Qm,G_\tau) m(\Pi) \Pi_t \overleftarrow{\times} [\overleftarrow{s-t-1}]_{\pi_v} \otimes 
(\Xi^{(s-t)(g+1)/2} \otimes
\bigoplus_{\chi \in \AF(\pi_v)} \chi^{-1}) & i=0
\end{array} \right.$$
\end{itemize}

\section{sur la suite spectrale des cycles évanescents}


A partir du calcul des groupes de cohomologie à support compact des $HT(g,t,\pi_v,\Pi_t)$, l'isomorphisme
(cf. la proposition (\ref{prop-hic}))
$$H^j_c(\bar X_{U^p,m,M_{d-h}}^{(d-h)}, R^i\Psi_v \otimes \LC_{\xi}))^h \simeq \bigoplus_{\tau_v \in \CF_h}
(H^j_c(\bar X_{U^p,m,M_{d-h}}^{(d-h)},\FC_{\tau_v} \otimes \LC_{\xi}) \otimes
\widetilde{\UC_{F_v,l,h,m}^{i}(\tau_v)})^{h/e_{\tau_v}}$$ permet de calculer les
$$H^j_{=h,i,\pi_v,\xi}[\Pi^{\oo,v}]:=\lim_{\atop{\longto}{U^p(m)}} H^j_c(\bar X_{U^p,m}^{(d-h)}, R^i\Psi_{\pi_v} 
\otimes
\LC_{\xi}))$$ Les suites spectrales associées à la stratification
\begin{equation}
E_{1,U^p(m),\pi_v}^{p,q;i}=H_c^{p+q}(\bar X_{U^p,m}^{(d-p+1)}, R^i\Psi_{\pi_v}) \Rightarrow H_c^{p+q}(\bar
X_{U^p,m}, R^i\Psi_{\pi_v})
\end{equation}
permettent alors de calculer les termes $E_2^{i,j}[\Pi^{\oo,v}]$ de la suite spectrale des cycles
évanescents. On résume les résultats obtenus dans

\bigskip

\marque \textit{le cas où $\Pi_v \simeq \st_s(\pi_v)$}: le corollaire (\ref{coro-hic}) donne tout d'abord:

\begin{coro} Les $H^j_{=h,i,\pi_v,\xi}[\Pi^{\oo,v}]$ vérifient les propriétés suivantes:
\begin{itemize}
    \item[(1)] pour $g$ divisant $d=sg$ et $\Pi_v \simeq [\overleftarrow{s-1}]_{\pi_v}$, les
    $H^j_{=h,i,\pi_v,\xi}[\Pi^{\oo,v}]$ sont nuls pour $h$ qui n'est pas de la forme
    $tg$ avec $1 \leq t \leq s$;

    \item[(2)] pour $g$ divisant $d=sg$ et $\Pi_v \simeq [\overleftarrow{s-1}]_{\pi_v}$, les
    $H^j_{=tg,i,\pi_v,\xi}[\Pi^{\oo,v}]$ sont nuls pour $j \neq d-tg$;

    \item[(3)] pour $g$ divisant $d=sg$ et $\Pi_v \simeq [\overleftarrow{s-1}]_{\pi_v}$, les
    $H^{d-tg}_{=tg,i,\pi_v,\xi}[\Pi^{\oo,v}]$ sont nuls pour $i$ ne vérifiant pas
    $t(g-1) \leq i \leq tg-1$. Si $1 \leq t < s$ et $i=tg-r$ avec $1 \leq r \leq t$, on a
    \begin{multline*} H^{d-tg}_{=tg,tg-r,\pi_v,\xi}[\Pi^{\oo,v}] \simeq  \\ \sharp \ker^1(\Qm,G_\tau) m(\Pi)
    ([\overleftarrow{t-r},\overrightarrow{r-1}]_{\pi_v} \overrightarrow{\times} [\overleftarrow{s-t-1}]_{\pi_v})
    \otimes \rec_{F_v}^\vee(\pi_v) (-\frac{s(g-1)-2(r-t)}{2})
    \end{multline*}
    en tant que représentation de $GL_d(F_v) \times W_v$, où $m(\Pi)$
    est la multiplicité de $\Pi$ dans l'espace des formes automorphes.

\end{itemize}
\end{coro}

Ainsi, cf. la preuve de la proposition (\ref{prop-hipsi}), les suites spectrales associées à la
stratification sont triviales, i.e. dégénèrent en $E_1$ de sorte que les aboutissements sont donnés par la
proposition (\ref{prop-hipsi}) dont on rappelle l'énoncé ci-dessous:

\begin{prop}
Pour tout $0 \leq i \leq d-1$, les $H^j(R^i\Psi_{\pi_v,\xi})[\Pi^{\oo,v}]$, en tant que $GL_d(F_v) \otimes
W_v$-modules, vérifient les propriétés suivantes:
\begin{itemize}
\item[(1)] ils sont nuls si $g$ n'est pas un diviseur de $d$;

\item[(2)] pour $g$ un diviseur de $d=sg$, ils sont nuls si $j$ n'est pas de la forme $d-tg$ pour $1 \leq t
\leq s$;

\item[(3)] pour $g$ un diviseur de $d=sg$ et $j=d-tg$ avec $1 \leq t \leq s$, ils sont nuls si $i$ n'est pas
de la forme $tg-r$ avec $1 \leq r \leq t$;

\item[(4)] pour $g$ un diviseur de $d=sg$ et $1 \leq t < s$,
$H^{d-tg}(R^{tg-r}\Psi_{\pi_v,\xi})[\Pi^{\oo,v}]$ est isomorphe à
$$\sharp \ker^1(\Qm,G_\tau) m(\Pi) [\overleftarrow{t-r},
\overrightarrow{r-1}]_{\pi_v} \overrightarrow{\times} [\overleftarrow{s-t-1}]_{\pi_v} \otimes
\rec_{F_v}^\vee(\pi_v) (-\frac{s(g-1)-2(r-t)}{2})$$

\end{itemize}
\end{prop}

Enfin la suite spectrale des cycles évanescents dégénère en $E_3$ et l'aboutissement
$H^i_{\eta_v}[\Pi^{\oo,v}]$ est donné par le corollaire (\ref{coho-global1}):

\begin{coro}
Les groupes de cohomologie de la fibre générique $H^i_{\eta_v,\xi}[\Pi^{\oo,v}]$ sont nuls pour $i \neq d-1$ et
$$H^{d-1}_{\eta_v}[\Pi^{\oo,v}]=\sharp \ker^1(\Qm,G_\tau) m(\Pi)[\overleftarrow{s-1}]_{\pi_v} \otimes 
\rec_{F_v}^\vee([\overleftarrow{s-1}]_{\pi_v})
(-\frac{sg-1}{2})$$
\end{coro}

\marque \textit{Cas où $\Pi_v \simeq \speh_s(\pi_v)$}:

\begin{coro} Les $H^j_{=h,i,\pi_v,\xi}[\Pi^{\oo,v}]$ vérifient les propriétés suivantes:
\begin{itemize}
    \item[(1)] pour $g$ divisant $d=sg$ et $\Pi_v \simeq [\overrightarrow{s-1}]_{\pi_v}$, les
    $H^j_{=h,i,\pi_v,\xi}[\Pi^{\oo,v}]$ sont nuls pour $h$ qui n'est pas de la forme
    $tg$ avec $1 \leq t \leq s$;

    \item[(2)] pour $g$ divisant $d=sg$ et $\Pi_v \simeq [\overrightarrow{s-1}]_{\pi_v}$, les
    $H^j_{=tg,i,\pi_v,\xi}[\Pi^{\oo,v}]$ sont nuls pour $j \neq (s-t)(g+1)$;

    \item[(3)] pour $g$ divisant $d=sg$ et $\Pi_v \simeq [\overrightarrow{s-1}]_{\pi_v}$, les
    $H^{(s-t)(g+1)}_{=tg,i,\pi_v,\xi}[\Pi^{\oo,v}]$ sont nuls pour $i$ ne vérifiant pas
    $t(g-1) \leq i \leq tg-1$. Si $1 \leq t < s$ et $i=tg-r$ avec $1 \leq r \leq t$, on a
    \begin{multline*} H^{d-tg}_{=tg,tg-r,\pi_v,\xi}[\Pi^{\oo,v}] \simeq \sharp \ker^1(\Qm,G_\tau) m(\Pi) \\
    ([\overleftarrow{t-r},\overrightarrow{r-1}]_{\pi_v} \overleftarrow{\times} [\overrightarrow{s-t-1}]_{\pi_v})
    \otimes \rec_{F_v}^\vee(\pi_v) (-\frac{s(g+1)-2r}{2})
    \end{multline*}
    en tant que représentation de $GL_d(F_v) \times W_v$, où $m(\Pi)$
    est la multiplicité de $\Pi$ dans l'espace des formes automorphes.

\end{itemize}
\end{coro}

\begin{proof} Le raisonnement est identique à celui de la preuve de la proposition (\ref{prop-hipsi}). On remarque
alors que pour tout $k\geq 1$, les flèches
$$d_k^{p,q;i}:E_{k,\pi_v,\xi}^{p,q;i}[\Pi^{\oo,v}] \longto E_{k,\pi_v,\xi}^{p+k,q+k-1}[\Pi^{\oo,v}]$$
des suites spectrales associées à la stratification sont toutes nulles. En effet, d'après le corollaire
(\ref{strate-alt-p}), pour que $E_{k,\pi_v}^{p,q;i}[\Pi^{\oo,v}]$ (resp.
$E_{k,\pi_v,\xi}^{p+k,q+k-1;i}[\Pi^{\oo,v}]$) soit non nul, il faut qu'il existe $1 \leq t_1 \leq s$ et $1 \leq
r_1 \leq t_1$ (resp. $1 \leq t_2 \leq s$ et $1 \leq r_2 \leq t_2$) avec
$$(p,q,i)=(t_1g+1,d-2t_1g-1+s-t_1,t_1g-r_1)$$
$$(\hbox{resp. } (p+k,q+k-1,i)=(t_2g+1,d-2l_2g-1+s-t_2,t_2g-r_2)).$$
Ce qui donne $1=(3g+1)(t_2-t_1)$; on voit alors que pour tout $k\geq 1$, $E_{k,\pi_v,\xi}^{p,q;i}[\Pi^{\oo,v}]$
et $E_{k,\pi_v,\xi}^{p+k,q+k-1}[\Pi^{\oo,v}]$ ne peuvent pas être tous deux non nuls de sorte que
$E_{\oo,I,\pi_v,\xi}^{p,q;i}[\Pi^{\oo,v}]=E_{1,I,\pi_v,\xi}^{p,q;i}[\Pi^{\oo,v}]$ d'où le résultat d'après le
corollaire précédent.

\end{proof}

Enfin la suite spectrale des cycles évanescents dégénère en $E_3$ et l'aboutissement
$H^i_{\eta_v,\xi}[\Pi^{\oo,v}]$ est donnée par la proposition (\ref{prop-lrs2}) dont on rappelle l'énoncé.

\begin{prop} Les groupes de cohomologie de la fibre générique $H^i_{\eta_v,\xi}[\Pi^{\oo,v}]$
sont nuls pour $|d-1-i| \geq s$ ou pour $i \not \equiv d-s \mod 2$ et sinon
$$H^{d+s-2-2i}_{\eta_v,\xi}[\Pi^{\oo,v}] \simeq \sharp \ker^1(\Qm,G_\tau) m(\Pi) [\overrightarrow{s-1}]_{\pi_v} 
\otimes \rec_{F_v}^\vee(\pi_v)
(-\frac{d+s-2-2i}{2})$$
\end{prop}






\part{Figures}

\noindent \textit{Aide à la lecture des figures (\ref{figure5}) et (\ref{figure6}):}
\begin{itemize}

\item on indique les coordonnées $(p,q)$ telles que $E_2^{p,q} [\Pi^{\oo,v}]$ est non nulle avec $\Pi$ une
représentation irréductible automorphe de $G(\Am)$ vérifiant $\hyp(\oo)$ et telle que $\Pi_v \simeq
\st_4(\pi_v)$ (resp. $\speh_4(\pi_v)$) dans la figure 5 (resp. figure 6) pour $\pi_v$ une représentation
cuspidale de $GL_2(F_v)$. Les représentations obtenues sont des représentations elliptiques de type $\pi_v$.
Pour chacune de ces coordonnées, on indique qu'elle est la strate qui la donne dans la suite spectrale de
stratification correspondante;

\item le nombre en dessous des $\overleftrightarrow{3}$ indique le poids;

\item toutes les flèches indiquées sont les seules non triviales et découlent des suites exactes 
(\ref{suites-exactes});

\item on note avec un cercle plus grand les coordonnées qui contribuent à l'aboutissement de la suite spectrale.
\end{itemize}

\begin{figure}
\centering \setlength{\unitlength}{1cm}
\begin{picture}(9,19)(0,-10)
\linethickness{.1pt}

\put(0,0){\vector(1,0){9}}
\put(0,0){\vector(0,1){9}}
\put(0,0){\circle*{.1}}

\put(0,1){\line(1,-1){1}}
\multiput(0,1)(1,-1){2}{\circle*{.1}}

\put(0,2){\line(1,-1){2}}
\multiput(0,2)(1,-1){3}{\circle*{.1}}

\put(0,3){\line(1,-1){3}}
\multiput(0,3)(1,-1){4}{\circle*{.1}}

\put(0,4){\line(1,-1){4}}
\multiput(0,4)(1,-1){5}{\circle*{.1}}

\put(0,5){\line(1,-1){5}}
\multiput(0,5)(1,-1){6}{\circle*{.1}}

\put(0,6){\line(1,-1){6}}
\multiput(0,6)(1,-1){7}{\circle*{.1}}

\put(0,7){\line(1,-1){7}}
\multiput(0,7)(1,-1){8}{\circle*{.1}}

\put(1,8){\vector(-1,-1){1}}
\put(1,8){\vector(-1,-2){1}}
\put(1,8){\vector(-1,-3){1}}
\put(1,8){\vector(-1,-4){1}}
\put(2,8.5){\makebox(0,0){$h=4 \times 2$}}

\put(0,7){\circle{.3}} \put(-.5,7){\makebox(0,0){$\genfrac{}{}{0pt}{}{\overleftarrow{3}}{10}$}}

\put(0,6){\circle{.2}}
\put(-.5,6){\makebox(0,0){$\genfrac{}{}{0pt}{}{[\overleftarrow{2},\overrightarrow{1}]}{8}$}}
\put(0,6){\vector(2,-1){2}}

\put(0,5){\circle{.2}}
\put(-.5,5){\makebox(0,0){$\genfrac{}{}{0pt}{}{[\overleftarrow{1},\overrightarrow{2}]}{6}$}}
\put(0,5){\vector(2,-1){2}}

\put(0,4){\circle{.2}} \put(-.5,4){\makebox(0,0){$\genfrac{}{}{0pt}{}{\overrightarrow{3}}{4}$}}
\put(0,4){\vector(2,-1){2}}

\put(3,6){\vector(-1,-1){1}}
\put(3,6){\vector(-1,-2){1}}
\put(3,6){\vector(-1,-3){1}}
\put(3,6.5){\makebox(0,0){$h=3 \times 2$}}

\put(2,5){\line(0,-1){2}} \put(2,5){\circle{.3}}
\put(1.5,3){\makebox(0,0){$\genfrac{}{}{0pt}{}{\overrightarrow{2} \overrightarrow{\times}
\overrightarrow{0}]}{4}$}}

\put(2,4){\circle{.2}} \put(1.5,4){\makebox(0,0){$\genfrac{}{}{0pt}{}{[\overleftarrow{1},\overrightarrow{1}]
\overrightarrow{\times} \overrightarrow{0}}{6}$}} \put(2,4){\vector(2,-1){2}}

\put(2,3){\circle{.2}} \put(1.5,5){\makebox(0,0){$\genfrac{}{}{0pt}{}{\overleftarrow{2}
\overrightarrow{\times} \overrightarrow{0}}{8}$}} \put(2,3){\vector(2,-1){2}}

\put(5,4){\vector(-1,-1){1}}
\put(5,4){\vector(-1,-2){1}}
\put(5,4.5){\makebox(0,0){$h=2 \times 2$}}

\put(4,3){\line(0,-1){1}} \put(4,3){\circle{.3}}
\put(4.7,3){\makebox(0,0){$\genfrac{}{}{0pt}{}{\overleftarrow{1} \overrightarrow{\times}
\overleftarrow{1}}{6}$}}

\put(4,2){\circle{.2}} \put(4.7,2){\makebox(0,0){$\genfrac{}{}{0pt}{}{\overrightarrow{1}
\overrightarrow{\times} \overleftarrow{1}}{4}$}} \put(4,2){\vector(2,-1){2}}

\put(7,2){\vector(-1,-1){1}}
\put(7,2.5){\makebox(0,0){$h=1 \times 2$}}

\put(6,1){\circle{.3}} \put(6.5,1){\makebox(0,0){$\genfrac{}{}{0pt}{}{\overleftarrow{0}
\overrightarrow{\times} \overleftarrow{2}}{4}$}}

\put(0,-10){\vector(1,0){9}} \put(0,-10){\vector(0,1){9}} \put(0,-10){\circle*{.1}}

\put(0,-9){\line(1,-1){1}} \multiput(0,-9)(1,-1){2}{\circle*{.1}}

\put(0,-8){\line(1,-1){2}} \multiput(0,-8)(1,-1){3}{\circle*{.1}}

\put(0,-7){\line(1,-1){3}} \multiput(0,-7)(1,-1){4}{\circle*{.1}}

\put(0,-6){\line(1,-1){4}} \multiput(0,-6)(1,-1){5}{\circle*{.1}}

\put(0,-5){\line(1,-1){5}} \multiput(0,-5)(1,-1){6}{\circle*{.1}}

\put(0,-4){\line(1,-1){6}} \multiput(0,-4)(1,-1){7}{\circle*{.1}}

\put(0,-3){\line(1,-1){7}} \multiput(0,-3)(1,-1){8}{\circle*{.1}}

\multiput(0,-3)(2,-2){4}{\circle{.3}}

\put(1,-2){\vector(-1,-1){1}} \put(2,-1.5){\makebox(0,0){$\overleftarrow{3}$ poids=s(g+1)-2}}

\put(3,-4){\vector(-1,-1){1}} \put(3,-3.5){\makebox(0,0){$\overleftarrow{3}$ poids=s(g+1)-4}}

\put(5,-6){\vector(-1,-1){1}} \put(5,-5.5){\makebox(0,0){$\overleftarrow{3}$ poids=s(g-1)+2}}

\put(7,-8){\vector(-1,-1){1}} \put(7,-7.5){\makebox(0,0){$\overleftarrow{3}$ poids=s(g-1)}}

\end{picture}
 \caption{\label{figure5} $s=4$, $g=2$; $E_2^{p,q}[\Pi^{\oo,v}]$ de la suite spectrale des
cycles évanescents: cas $\Pi_v=\st_4(\pi_v)$. Le dessin du bas représente le terme $E_3$ et donc
l'aboutissement de la suite spectrale.}
\end{figure}

\newpage

\begin{figure}
\centering \setlength{\unitlength}{1cm}
\begin{picture}(9,19)(0,-10)
\linethickness{.1pt}

\put(0,0){\vector(1,0){11}}
\put(0,0){\vector(0,1){9}}
\put(0,0){\circle*{.1}}

\put(0,1){\line(1,-1){1}}
\multiput(0,1)(1,-1){2}{\circle*{.1}}

\put(0,2){\line(1,-1){2}}
\multiput(0,2)(1,-1){3}{\circle*{.1}}

\put(0,3){\line(1,-1){3}}
\multiput(0,3)(1,-1){4}{\circle*{.1}}

\put(0,4){\line(1,-1){4}}
\multiput(0,4)(1,-1){5}{\circle*{.1}}

\put(0,5){\line(1,-1){5}}
\multiput(0,5)(1,-1){6}{\circle*{.1}}

\put(0,6){\line(1,-1){6}}
\multiput(0,6)(1,-1){7}{\circle*{.1}}

\put(0,7){\line(1,-1){7}}
\multiput(0,7)(1,-1){8}{\circle*{.1}}

\put(1,8){\vector(-1,-1){1}}
\put(1,8){\vector(-1,-2){1}}
\put(1,8){\vector(-1,-3){1}}
\put(1,8){\vector(-1,-4){1}}
\put(2,8.5){\makebox(0,0){$h=4 \times 2$}}

\put(0,7){\circle{.2}} \put(-.5,7){\makebox(0,0){$\genfrac{}{}{0pt}{}{\overleftarrow{3}}{10}$}}
\put(0,7){\vector(3,-2){3}}

\put(0,6){\circle{.2}}
\put(-.5,6){\makebox(0,0){$\genfrac{}{}{0pt}{}{[\overleftarrow{2},\overrightarrow{1}]}{8}$}}
\put(0,6){\vector(3,-2){3}}

\put(0,5){\circle{.2}}
\put(-.5,5){\makebox(0,0){$\genfrac{}{}{0pt}{}{[\overleftarrow{1},\overrightarrow{2}]}{6}$}}
\put(0,5){\vector(3,-2){3}}

\put(0,4){\circle{.3}} \put(-.5,4){\makebox(0,0){$\genfrac{}{}{0pt}{}{\overrightarrow{3}}{4}$}}

\put(4,6){\vector(-1,-1){1}} \put(4,6){\vector(-1,-2){1}} \put(4,6){\vector(-1,-3){1}} \put(4,6.5){\makebox(0,0){$h=3
\times 2$}}

\put(3,5){\line(0,-1){2}} \put(3,5){\circle{.2}}
\put(2.5,5){\makebox(0,0){$\genfrac{}{}{0pt}{}{\overleftarrow{2} \overleftarrow{\times}
\overrightarrow{0}}{10}$}} \put(3,5){\vector(3,-2){3}} \put(3,5){\line(1,-1){5}}
\multiput(3,5)(1,-1){6}{\circle*{.1}}

\put(3,4){\circle{.2}} \put(2.5,4){\makebox(0,0){$\genfrac{}{}{0pt}{}{[\overleftarrow{1},\overrightarrow{1}]
\overleftarrow{\times} \overleftarrow{0}}{8}$}} \put(3,4){\vector(3,-2){3}}

\put(3,3){\circle{.3}} \put(2.5,3){\makebox(0,0){$\genfrac{}{}{0pt}{}{\overrightarrow{2}
\overleftarrow{\times} \overleftarrow{0}}{6}$}}

\put(7,4){\vector(-1,-1){1}}
\put(7,4){\vector(-1,-2){1}}
\put(7,4.5){\makebox(0,0){$h=2 \times 2$}}

\put(6,3){\line(0,-1){1}} \put(6,3){\circle{.2}}
\put(6.5,3){\makebox(0,0){$\genfrac{}{}{0pt}{}{\overleftarrow{1} \overleftarrow{\times}
\overrightarrow{1}}{10}$}} \put(6,3){\vector(3,-2){3}} \put(6,3){\line(1,-1){3}}
\multiput(6,3)(1,-1){4}{\circle*{.1}}

\put(6,2){\circle{.3}} \put(6.5,2){\makebox(0,0){$\genfrac{}{}{0pt}{}{\overrightarrow{1}
\overleftarrow{\times} \overrightarrow{1}}{8}$}}

\put(9,2){\vector(0,-1){1}}
\put(9,2.5){\makebox(0,0){$h=1 \times 2$}}

\put(9,1){\circle{.3}} \put(9.5,1){\makebox(0,0){$\genfrac{}{}{0pt}{}{\overleftarrow{0}
\overleftarrow{\times} \overrightarrow{2}}{10}$}} \put(9,1){\line(1,-1){1}}
\multiput(9,1)(1,-1){2}{\circle*{.1}}

\put(0,-10){\vector(1,0){11}} \put(0,-10){\vector(0,1){9}} \put(0,-10){\circle*{.1}}

\put(0,-9){\line(1,-1){1}} \multiput(0,-9)(1,-1){2}{\circle*{.1}}

\put(0,-8){\line(1,-1){2}} \multiput(0,-8)(1,-1){3}{\circle*{.1}}

\put(0,-7){\line(1,-1){3}} \multiput(0,-7)(1,-1){4}{\circle*{.1}}

\put(0,-6){\line(1,-1){4}} \multiput(0,-6)(1,-1){5}{\circle*{.1}}

\put(0,-5){\line(1,-1){5}} \multiput(0,-5)(1,-1){6}{\circle*{.1}}

\put(0,-4){\line(1,-1){6}} \multiput(0,-4)(1,-1){7}{\circle*{.1}}

\put(0,-3){\line(1,-1){7}} \multiput(0,-3)(1,-1){8}{\circle*{.1}}

\put(1,-2){\vector(-1,-4){1}} \put(2,-1.5){\makebox(0,0){$\overrightarrow{3}$ poids=s(g-1)}}

\put(4,-4){\vector(-1,-3){1}} \put(4,-3.5){\makebox(0,0){$\overrightarrow{3}$ poids=s(g-1)+2}}

\multiput(0,-6)(3,-1){4}{\circle{.3}}

\put(3,-5){\line(1,-1){5}} \multiput(3,-5)(1,-1){6}{\circle*{.1}}

\put(7,-6){\vector(-1,-2){1}} \put(7,-5.5){\makebox(0,0){$\overrightarrow{3}$ poids=s(g+1)-4}}

\put(6,-7){\line(1,-1){3}} \multiput(6,-7)(1,-1){4}{\circle*{.1}}

\put(9,-8){\vector(0,-1){1}} \put(9,-7.5){\makebox(0,0){$\overrightarrow{3}$ poids=s(g+1)-2}}

\put(9,-9){\line(1,-1){1}} \multiput(9,-9)(1,-1){2}{\circle*{.1}}

\end{picture}
 \caption{\label{figure6}  $s=4$, $g=2$: $E_2^{p,q}[\Pi^{\oo,v}]$ de la suite spectrale des
cycles évanescents: cas $\Pi_v=\speh_s(\pi_v)$. Le dessin du bas représente le terme $E_3$ et donc
l'aboutissement de la suite spectrale.}
\end{figure}

\newpage

\begin{figure}[b]
\centering \setlength{\unitlength}{1.2cm}
\begin{picture}(9,17)(0,-7)
\linethickness{.1pt}

\put(4,1){\line(0,1){7}} \put(0,5){\line(1,0){8}}

\multiput(4,2)(0,1){7}{\circle*{.1}} \multiput(1,5)(1,0){7}{\circle*{.1}}

\put(1,8){\line(1,-1){6}} \multiput(1,8)(2,-2){4}{\circle{.3}}

\multiput(1,8)(1,-1){7}{\circle*{.1}}

\put(2,6){\line(1,-1){4}} \multiput(2,6)(2,-2){3}{\circle{.3}}

\multiput(2,6)(1,-1){5}{\circle*{.1}}

\put(3,4){\line(1,-1){2}} \multiput(3,4)(2,-2){2}{\circle{.3}}

\multiput(3,4)(1,-1){3}{\circle*{.1}}

\put(4,2){\circle{.3}}

\put(4,2){\vector(1,0){1}} \put(5,2){\vector(1,0){1}} \put(6,2){\vector(1,0){1}}

\put(3,4){\vector(1,0){1}} \put(4,4){\vector(1,0){1}}

\put(2,6){\vector(1,0){1}}

\put(4,7.5){\vector(0,1){1} $j$} \put(7,5){\vector(1,0){1} $i$}

\put(0.5,8){\makebox(0,0){$\overleftarrow{3}$}}

\put(1,6){\makebox(0,0){$\overleftarrow{2} \overrightarrow{\times} \overrightarrow{0}$}}

\put(2.5,4.5){\makebox(0,0){$\overleftarrow{1} \overrightarrow{\times} \overrightarrow{1}$}}

\put(3,2){\makebox(0,0){$\overleftarrow{0} \overrightarrow{\times} \overrightarrow{2}$}}

\multiput(3,6.5)(2,-2){3}{\makebox(0,0){$\overleftarrow{3}$}}

\put(5,2.5){\makebox(0,0){$\overleftarrow{1} \overrightarrow{\times} \overrightarrow{1}$}}

\put(6,2.5){\makebox(0,0){$\overleftarrow{2} \overrightarrow{\times} \overrightarrow{0}$}}

\put(4,4.5){\makebox(0,0){$\overleftarrow{2} \overrightarrow{\times} \overrightarrow{0}$}}

\put(7,1){\vector(-1,1){1}} \put(7,1){\vector(-2,1){2}} \put(7,1){\vector(-3,1){3}}
\put(7,1){\vector(0,1){1}}

\put(7,.7){\makebox(0,0){poids=$-3$}}

\put(3,3){\vector(0,1){1}} \put(3,3){\vector(1,1){1}} \put(3,3){\vector(2,1){2}}

\put(3,2.7){\makebox(0,0){poids=$-1$}}

\put(2,5){\vector(0,1){1}} \put(2,5){\vector(1,1){1}} \put(2,4.9){\makebox(0,0){poids=$1$}}

\put(1,7){\vector(0,1){1}} \put(1,6.9){\makebox(0,0){poids=$3$}}

\put(4,-7){\line(0,1){7}} \put(0,-3){\line(1,0){8}}

\multiput(4,-6)(0,1){7}{\circle*{.1}} \multiput(1,-3)(1,0){7}{\circle*{.1}}

\put(1,0){\line(1,-1){6}} \put(1,0){\circle{.3}}

\multiput(1,0)(1,-1){7}{\circle*{.1}}

\put(2,-2){\line(1,-1){4}} \put(2,-2){\circle{.3}}

\multiput(2,-2)(1,-1){5}{\circle*{.1}}

\put(3,-4){\line(1,-1){2}} \put(3,-4){\circle{.3}}

\multiput(3,-4)(1,-1){3}{\circle*{.1}}

\put(4,-6){\circle{.3}}

\put(4,-.5){\vector(0,1){1} $j$} \put(7,-3){\vector(1,0){1} $i$}

\put(0.5,0){\makebox(0,0){$\overleftarrow{3}$}}

\put(1,-2){\makebox(0,0){$[\overleftarrow{2},\overrightarrow{1}]$}}

\put(2.5,-3.5){\makebox(0,0){$[\overleftarrow{1},\overrightarrow{2}]$}}

\put(3,-6){\makebox(0,0){$\overrightarrow{3}$}}

\end{picture}

\caption{\label{figure7} Termes $MLE_\bullet^{i,j}$ pour $s=4$}

\end{figure}
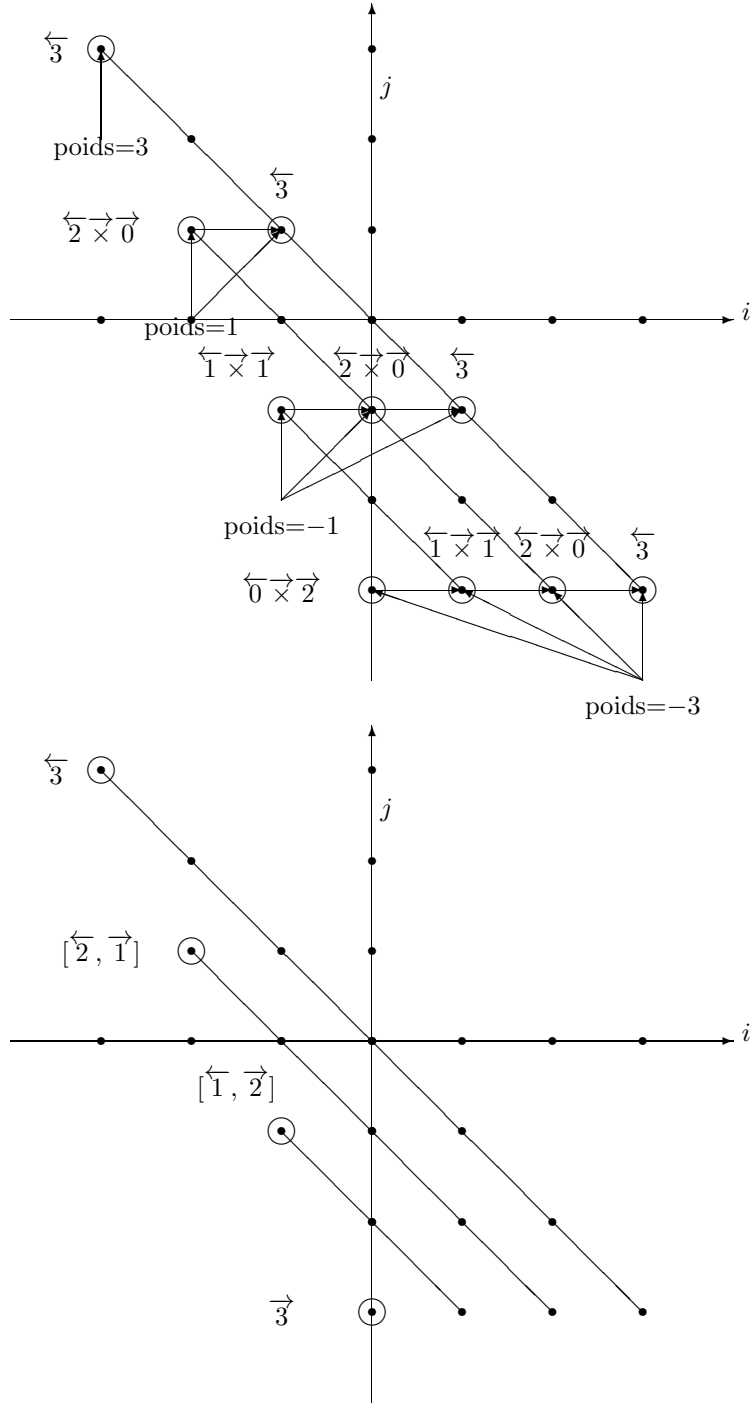

\clearpage


\newpage

\noindent \textit{Aide à la lecture des figures (\ref{figure1}), (\ref{figure2})}: le nombre qui accompagne
les représentations indique le décalage à donner à la représentation de Galois associée; ainsi quand, dans
les $E_2^{p,q}$ des suites spectrales (\ref{ss-dualite}) pour $\PC(g,t,\pi_v,k)$, on écrit
$\genfrac{}{}{0pt}{}{\overleftrightarrow{s-1}}{r}$, il faut lire
$$j^{\geq tg}_! HT(g,t,\pi_v,[\overleftrightarrow{s-1}]_{\pi_v}) \otimes \rec_{F_v}^\vee(\pi_v)
(-\frac{i(g-1)+t-1-k+r}{2}).$$

\begin{figure}[!ht]
\centering \setlength{\unitlength}{.4cm}
\begin{picture}(17,9)(-1,-1)
\linethickness{.1pt} \put(0,7){\circle{.2}}
\put(-.5,7.5){\makebox(0,0){$\genfrac{}{}{0pt}{}{\overleftarrow{3}}{0}$}}

\put(6,7){\circle*{.2}} \put(5.5,6){\makebox(0,0){$\genfrac{}{}{0pt}{}{\overleftarrow{3}
\overrightarrow{\times} \overrightarrow{0}}{0}$}}

\put(7,7){\circle{.2}} \put(7,7){\line(1,-1){7}}
\put(7.5,7.5){\makebox(0,0){$\genfrac{}{}{0pt}{}{\overleftarrow{3} \overleftarrow{\times} \overrightarrow{0}}
{1}$}}

\put(14,1){\circle*{.2}} \put(14.5,1.5){\makebox(0,0){$\genfrac{}{}{0pt}{}{\overleftarrow{3}
\overleftarrow{\times} \overrightarrow{0}}{0}$}}

\put(14,0){\circle*{.1}} \put(14,0){\line(0,1){1}} \put(14,-1){\makebox(0,0){$(0,0)$}}
\put(0,6){\makebox(0,0){$(-14,7)$}} \multiput(0,7)(1,0){8}{\circle*{.1}}
\multiput(7,7)(0,-1){8}{\circle*{.1}} \multiput(7,0)(1,0){8}{\circle*{.1}}

\qbezier(0,7)(3,8)(6,7)

\end{picture}
 \caption{\label{figure1} $g=7$, $s=5$; $E_2^{p,q}$ de la suite spectrale
(\ref{ss-dualite}) pour $t=4$}
\end{figure}

\begin{figure}[!ht]
\centering \setlength{\unitlength}{.45cm}
\begin{picture}(30,16)(-1,1)
\linethickness{.1pt}

\put(0,14){\circle{.2}} \put(-.5,14.5){\makebox(0,0){$\genfrac{}{}{0pt}{}{\overleftarrow{2}}{0}$}}
\qbezier(0,14)(3,15)(6,14)

\put(6,14){\circle*{.2}} \put(5.5,13){\makebox(0,0){$\genfrac{}{}{0pt}{}{0}{\overleftarrow{2}
\overrightarrow{\times} \overrightarrow{0}}$}}

\put(7,14){\circle{.2}} \put(7,14){\line(1,-1){7}}
\put(7,15){\makebox(0,0){$\genfrac{}{}{0pt}{}{\overleftarrow{2} \overleftarrow{\times}
\overrightarrow{0}}{1}$}}

\qbezier(6,14)(9,13)(12,14)

\put(12,14){\circle*{.2}} \put(11.5,13){\makebox(0,0){$\genfrac{}{}{0pt}{}{0}{\overleftarrow{2}
\overrightarrow{\times} \overrightarrow{1}}$}}

\qbezier(7,14)(10,15)(13,14) \put(13,14){\circle*{.2}}
\put(13,15){\makebox(0,0){$\genfrac{}{}{0pt}{}{(\overleftarrow{2} \overleftarrow{\times} \overrightarrow{0})
\overrightarrow{\times} \overrightarrow{0}}{1}$}}

\put(14,14){\circle{.2}} \put(16,14){\makebox(0,0){$\genfrac{}{}{0pt}{}{\overleftarrow{2}
\overleftarrow{\times} \overleftarrow{1}}{2}$}} \put(14,14){\line(1,-1){14}}

\multiput(0,14)(1,0){15}{\circle*{.1}}
\multiput(7,14)(0,-1){8}{\circle*{.1}}
\multiput(7,7)(1,0){15}{\circle*{.1}}

\put(14,8){\circle*{.2}} \put(14,9){\makebox(0,0){$\genfrac{}{}{0pt}{}{\overleftarrow{2}
\overleftarrow{\times} \overrightarrow{0}}{0}$}}

\qbezier(14,8)(17,9)(20,8)

\put(20,8){\circle*{.2}} \put(19.5,7){\makebox(0,0){$\genfrac{}{}{0pt}{}{0}{(\overleftarrow{2}
\overleftarrow{\times} \overrightarrow{0}) \overrightarrow{\times} \overrightarrow{0}}$}}

\put(21,8){\circle{.2}} \put(21.5,9){\makebox(0,0){$\genfrac{}{}{0pt}{}{(\overleftarrow{2}
\overleftarrow{\times} \overleftarrow{0}) \overleftarrow{\times} \overleftarrow{0}}{1}$}}
\put(21,8){\line(1,-1){7}}

\multiput(14,8)(1,0){8}{\circle*{.1}}
\multiput(21,7)(0,-1){8}{\circle*{.1}}
\multiput(21,0)(1,0){8}{\circle*{.1}}

\put(14,7){\line(0,1){1}}
\put(28,0){\line(0,1){2}}

\put(28,2){\circle{.2}} \put(28,3){\makebox(0,0){$\genfrac{}{}{0pt}{}{\overleftarrow{2}
\overleftarrow{\times} \overrightarrow{1}}{0}$}}

\put(28,1){\circle{.1}}

\put(14,6){\makebox(0,0){$(-14,7)$}}
\put(0,13){\makebox(0,0){$(-28,14)$}}
\put(28,-1){\makebox(0,0){$(0,0)$}}

\end{picture}
 \caption{\label{figure2}  $g=7$, $s=5$; $E_2^{p,q}$ de la suite spectrale
(\ref{ss-dualite}) pour $t=3$}
\end{figure}

\clearpage
\newpage

\begin{figure}[!ht]
\centering \setlength{\unitlength}{.48cm}

\begin{picture}(24,25)(-19,-16)
\linethickness{.1pt}

\multiput(-18,0)(1,0){19}{\circle*{.1}} \multiput(-18,0)(5,0){4}{\circle{.3}} \put(-17.7,.2){$5$}

\put(-18,-1){\vector(0,1){1}} \put(-18,-1.5){\makebox(0,0){$(-18,0)$}} \put(-19,1){\vector(1,-1){1}}

\put(-19,1.5){\makebox(0,0){$\overleftarrow{0}$}}

\put(-18,-6){$M^{=6}$} \qbezier(-18,0)(-15.5,-.5)(-13,0)

\put(-13.5,-.5){\line(1,1){2}}

\put(-14,-1.2){$10$}

\put(-13,-6){$M^{=12}$}

\put(-14,1){\vector(1,-1){1}} \put(-14,1.5){\makebox(0,0){$\overleftarrow{0} \overrightarrow{\times}
\overrightarrow{0}$}}

\qbezier(-13,0)(-10.5,.5)(-8,0)

\multiput(-12,1)(1,0){13}{\circle*{.1}}

\multiput(-12,-1)(1,0){13}{\circle*{.1}}

\multiput(-12,1)(5,0){3}{\circle{.3}}

\multiput(-12,-1)(5,0){3}{\circle{.3}}

\put(-12,2){\vector(0,-1){1}}

\put(-12,2.5){\makebox(0,0){$\overleftarrow{1}$}}

\put(-12,-2){\vector(0,1){1}}

\put(-12,-2.5){\makebox(0,0){$\genfrac{}{}{0pt}{}{\overleftarrow{1}}{12}$}}

\qbezier(-12,1)(-9.5,1.5)(-7,1) \qbezier(-12,-1)(-9.5,-.5)(-7,-1)

\put(-8.5,-.5){\line(1,1){3}}
\put(-5.8,1.7){$15$}
\put(-7,-6){$M^{=18}$}
\put(-5.5,0.5){\line(-1,-1){2}}
\put(-5.8,-.3){$17$}
\put(-11,3){\vector(1,-1){3}}

\put(-11,3.5){\makebox(0,0){$\overleftarrow{0} \overrightarrow{\times} \overrightarrow{1}$}}

\put(-8,2){\vector(1,-1){1}}

\put(-8,2.5){\makebox(0,0){$\overleftarrow{1} \overrightarrow{\times} \overrightarrow{0}$}}

\put(-6,3){\vector(0,-1){1}}

\put(-6,3.5){\makebox(0,0){$\overleftarrow{2}$}}

\put(-9,-2){\vector(2,1){2}}

\put(-9,-2.5){\makebox(0,0){$\overleftarrow{1} \overrightarrow{\times} \overrightarrow{0}$}}

\put(-8,-4){\vector(1,2){2}} \put(-8,-4.5){\makebox(0,0){$\overleftarrow{2}$}}

\put(-6,-3){\vector(0,1){1}} \put(-6,-3.5){\makebox(0,0){$\genfrac{}{}{0pt}{}{\overleftarrow{2}}{19}$}}

\multiput(-6,-2)(1,0){7}{\circle*{.1}} \multiput(-6,0)(1,0){7}{\circle*{.1}}
\multiput(-6,2)(1,0){7}{\circle*{.1}} \multiput(-6,-2)(5,0){2}{\circle{.3}}
\multiput(-6,0)(5,0){2}{\circle{.3}} \multiput(-6,2)(5,0){2}{\circle{.3}} \qbezier(-6,-2)(-3.5,-2.5)(-1,-2)
\qbezier(-6,0)(-3.5,.5)(-1,0) \qbezier(-6,2)(-3.5,2.5)(-1,2) \qbezier(-7,-1)(-4.5,-1.5)(-2,-1)
\qbezier(-7,1)(-4.5,1.5)(-2,1) \qbezier(-8,0)(-5.5,-.5)(-3,0)

\put(-3.5,-.5){\line(1,1){4}}
\put(-2.5,-1.5){\line(1,1){3}}
\put(-1.5,-2.5){\line(1,1){2}}
\put(.3,2.5){$20$}
\put(.3,.5){$22$}
\put(.3,-1.5){$24$}
\put(.2,-3){$26$}
\multiput(0,-3)(0,2){4}{\circle{.3}}
\put(-3,-6){$M^{=24}$}
\put(-3,-3){\vector(1,2){1}}

\put(-3,-3.5){\makebox(0,0){$\overleftarrow{1} \overrightarrow{\times} \overrightarrow{1}$}}

\put(-2,-4){\vector(1,2){1}}

\put(-2,-4.5){\makebox(0,0){$\overleftarrow{2} \overrightarrow{\times} \overrightarrow{0}$}}

\put(0,-5){\vector(0,1){2}}

\put(0,-5.5){\makebox(0,0){$\overleftarrow{3}$}}

\put(1,-1){\vector(-1,0){2}}

\put(1.5,-1){$\overleftarrow{3}$}

\put(1,1){\vector(-1,0){2}}

\put(1.5,1){$\overleftarrow{3}$}

\put(1,3){\vector(-1,0){2}}

\put(1.5,3){$\overleftarrow{3}$}

\put(1,0){\vector(-1,0){2}}

\put(1.5,0){$\overleftarrow{2} \overrightarrow{\times} \overrightarrow{0}$}

\put(-1,5){\vector(0,-1){3}}

\put(-1,5.5){\makebox(0,0){$\overleftarrow{2} \overrightarrow{\times} \overrightarrow{0}$}}


\put(-2,3){\vector(0,-1){2}}

\put(-2,3.5){\makebox(0,0){$\overleftarrow{1} \overrightarrow{\times} \overrightarrow{1}$}}

\put(-4,4){\vector(1,-4){1}}

\put(-3,-13){\vector(1,2){1}}

\put(-4,4.5){\makebox(0,0){$\overleftarrow{0} \overrightarrow{\times} \overrightarrow{2}$}}

\put(2,-2){\vector(-1,1){2}}

\put(2.5,-2.5){$(0,0)$}

\multiput(-18,-10)(1,0){19}{\circle*{.1}} \multiput(-18,-10)(5,0){4}{\circle{.3}}

\put(-18,-11){\vector(0,1){1}} \put(-18,-11.5){\makebox(0,0){$(-18,0)$}} \put(-19,-9){\vector(1,-1){1}}

\put(-19,-8.5){\makebox(0,0){$\overleftarrow{0}$}}

\put(-14,-9){\vector(1,-1){1}} \put(-14,-8.5){\makebox(0,0){$\overrightarrow{1}$}}

\multiput(-12,-9)(1,0){13}{\circle*{.1}}

\multiput(-12,-11)(1,0){13}{\circle*{.1}}

\multiput(-12,-11)(5,0){3}{\circle{.3}}

\put(-11,-6.5){\makebox(0,0){$\overrightarrow{2}$}}

\put(-12,-12){\vector(0,1){1}}

\put(-12,-12.5){\makebox(0,0){$\overleftarrow{1}$}}

\put(-11,-6.5){\makebox(0,0){$\overrightarrow{2}$}}

\put(-9,-12){\vector(2,1){2}}

\put(-9,-12.5){\makebox(0,0){$[\overleftarrow{1},\overrightarrow{1}]$}}

\put(-11,-7){\vector(1,-1){3}}

\put(-6,-13){\vector(0,1){1}} \put(-6,-13.5){\makebox(0,0){$\overleftarrow{2}$}}

\multiput(-6,-2)(1,0){7}{\circle*{.1}} \multiput(-6,0)(1,0){7}{\circle*{.1}}
\multiput(-6,2)(1,0){7}{\circle*{.1}}

\multiput(-6,-12)(5,0){2}{\circle{.3}}

\put(0,-13){\circle{.3}}

\put(-3,-13.5){\makebox(0,0){$[\overleftarrow{1},\overrightarrow{2}]$}}

\put(-2,-14){\vector(1,2){1}}

\put(-2,-14.5){\makebox(0,0){$[\overleftarrow{2},\overrightarrow{1}]$}}

\put(0,-15){\vector(0,1){2}}

\put(0,-15.5){\makebox(0,0){$\overleftarrow{3}$}}

\put(-4,-8){\vector(1,-2){1}}

\put(-4,-7.5){\makebox(0,0){$\overrightarrow{3}$}}

\put(2,-12){\vector(-1,1){2}}

\put(2.5,-12.5){$(0,0)$}

\end{picture}
\caption{\label{fig0} La figure illustre les théorèmes (\ref{theo-global1}) et (\ref{theo-global2}) dans le
cas $g=6$ et $s=4$, où on représente en $(i,k)$ les $h^i gr_k$ en indiquant le support de la forme $\bar
X_{U^p,m}^{(24-6r)}$ avec $1 \leq r \leq 4$, la représentation associée et le poids. Concrètement quand on
écrit en $(i,k)$, $\pi$, il faut lire $j^{\geq 6r}_! HT(6,r,\pi_v)\pi) \otimes \rec_{F_v}^\vee(\pi_v)
|-|^{-(23+i+k)/2}$. Le dessin du bas est le terme $E_2$ et donc l'aboutissement de la suite spectrale
(\ref{suite-spectrale}).}
\end{figure}
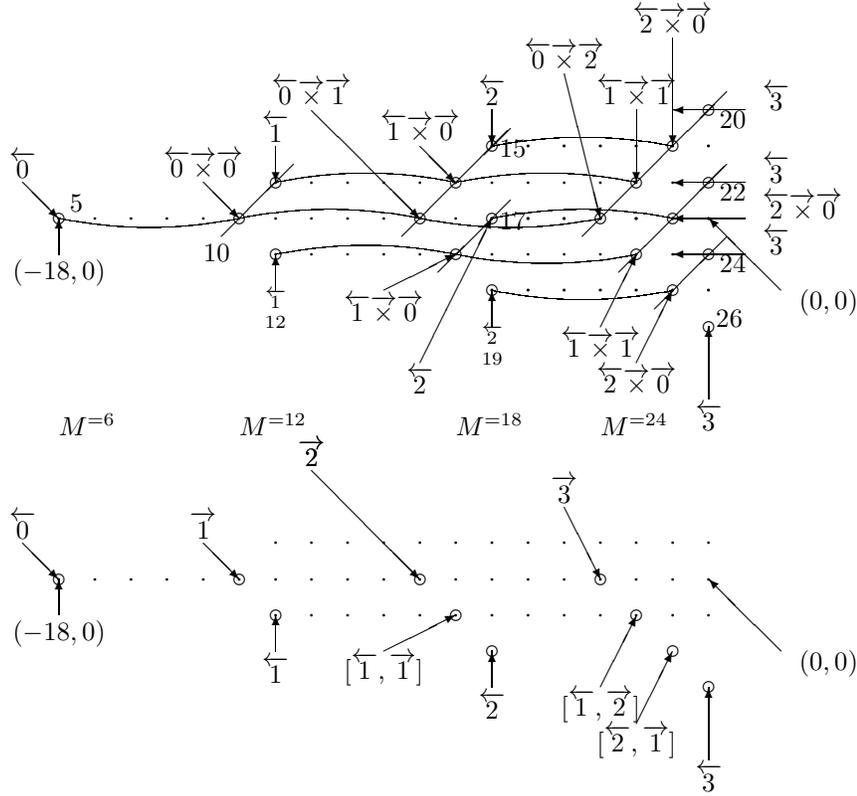


\clearpage
\newpage

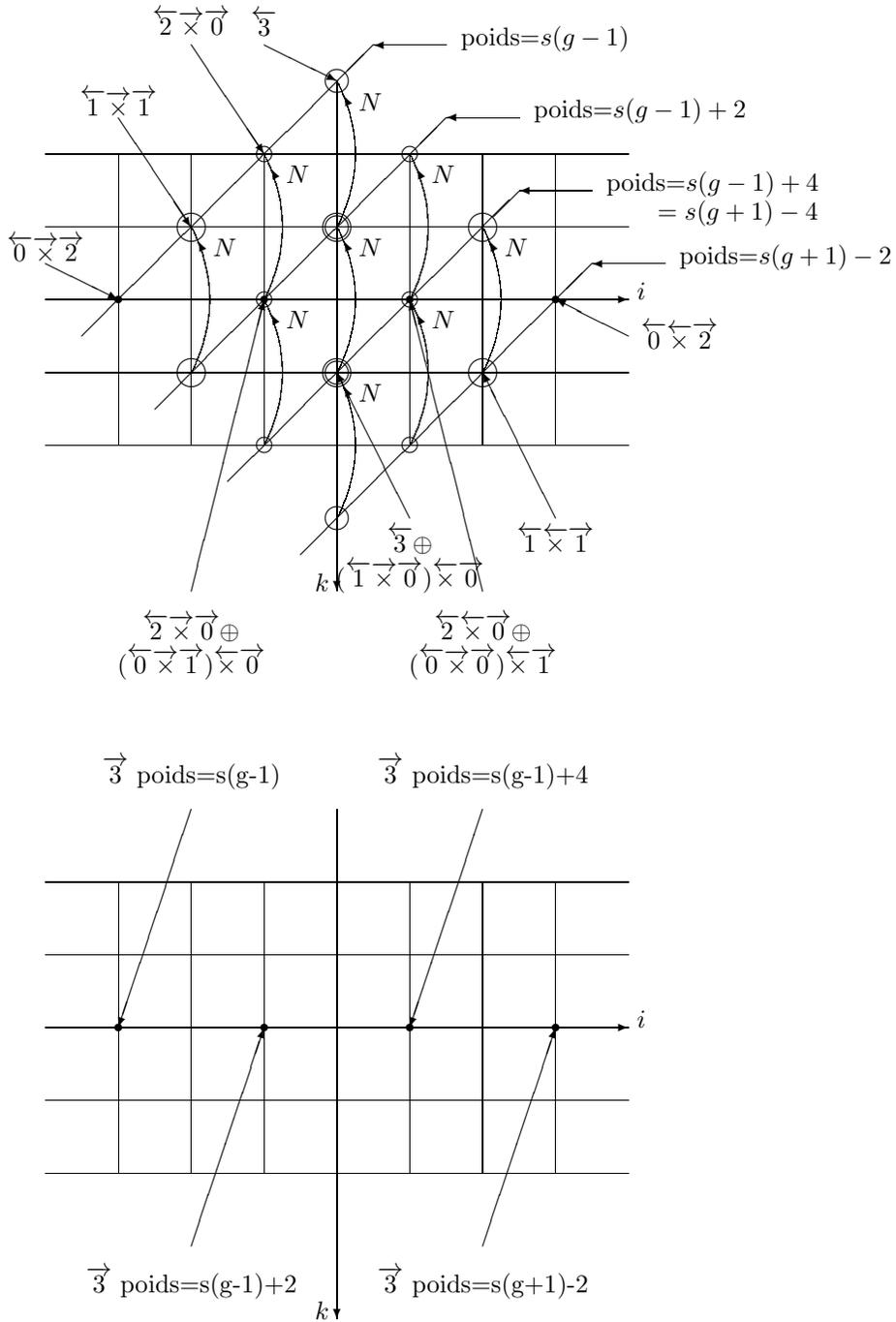
\begin{figure}[!ht]

\centering \setlength{\unitlength}{1cm}
\begin{picture}(15,20)(0,-10)
\linethickness{.1pt}

\multiput(1,3)(0,1){5}{\line(1,0){8}} \multiput(2,3)(1,0){7}{\line(0,1){4}}

\multiput(5,2)(0,2){4}{\circle{.3}}

\multiput(4,3)(0,2){3}{\circle{.2}} \multiput(6,3)(0,2){3}{\circle{.2}}

\multiput(3,4)(2,0){3}{\circle{.4}} \multiput(3,6)(2,0){3}{\circle{.4}}

\multiput(2,5)(2,0){4}{\circle*{.1}}

\put(8,5){\vector(1,0){1} $i$} \put(5,8){\vector(0,-1){7}} \put(4.7,1){$k$}

\multiput(4.5,1.5)(-1,1){4}{\line(1,1){4}}

\multiput(6.5,8.5)(1,-1){4}{\vector(-1,0){1}}

\put(6.7,8.5){poids=$s(g-1)$} \put(7.7,7.5){poids=$s(g-1)+2$} \put(8.7,6.5){poids=$s(g-1)+4$}
\put(9.4,6.1){$=s(g+1)-4$} \put(9.7,5.5){poids=$s(g+1)-2$}

\qbezier(3,4)(3.5,5)(3,6) \put(3.1,5.8){\vector(-1,2){0}} \qbezier(7,4)(7.5,5)(7,6)
\put(7.1,5.8){\vector(-1,2){0}}

\multiput(3.3,5.6)(1,1){3}{$N$} \multiput(4.3,4.6)(1,1){3}{$N$} \multiput(5.3,3.6)(1,1){3}{$N$}

\qbezier(4,3)(4.5,4)(4,5) \put(4.1,4.8){\vector(-1,2){0}} \qbezier(4,5)(4.5,6)(4,7)
\put(4.1,6.8){\vector(-1,2){0}} \qbezier(6,3)(6.5,4)(6,5) \put(6.1,4.8){\vector(-1,2){0}}
\qbezier(6,5)(6.5,6)(6,7) \put(6.1,6.8){\vector(-1,2){0}}

\qbezier(5,2)(5.5,3)(5,4) \put(5.1,3.8){\vector(-1,2){0}}\qbezier(5,4)(5.5,5)(5,6)
\put(5.1,5.8){\vector(-1,2){0}}\qbezier(5,6)(5.5,7)(5,8) \put(5.1,7.8){\vector(-1,2){0}}

\put(4,8.5){\vector(2,-1){1}} \put(4,8.8){\makebox(0,0){$\overleftarrow{3}$}}

\put(3,8.5){\vector(2,-3){1}} \put(3,8.8){\makebox(0,0){$\overleftarrow{2} \overrightarrow{\times}
\overrightarrow{0}$}}

\put(2,7.5){\vector(2,-3){1}} \put(2,7.7){\makebox(0,0){$\overleftarrow{1} \overrightarrow{\times}
\overrightarrow{1}$}}

\put(1,5.5){\vector(2,-1){1}} \put(1,5.7){\makebox(0,0){$\overleftarrow{0} \overrightarrow{\times}
\overrightarrow{2}$}}

\put(3,1){\vector(1,4){1}} \put(3,.5){\makebox(0,0){$\overleftarrow{2} \overrightarrow{\times}
\overrightarrow{0} \oplus$}} \put(3,0){\makebox(0,0){$(\overleftarrow{0} \overrightarrow{\times}
\overrightarrow{1}) \overleftarrow{\times} \overrightarrow{0}$}}

\put(7,1){\vector(-1,4){1}} \put(7,.5){\makebox(0,0){$\overleftarrow{2} \overleftarrow{\times}
\overrightarrow{0} \oplus$}} \put(7,0){\makebox(0,0){$(\overleftarrow{0} \overrightarrow{\times}
\overrightarrow{0}) \overleftarrow{\times} \overrightarrow{1}$}}

\put(6,2){\vector(-1,2){1}} \put(6,1.7){\makebox(0,0){$\overleftarrow{3} \oplus$}}
\put(6,1.2){\makebox(0,0){$(\overleftarrow{1} \overrightarrow{\times} \overrightarrow{0})
\overleftarrow{\times} \overrightarrow{0}$}}

\put(8,2){\vector(-1,2){1}} \put(8,1.7){\makebox(0,0){$\overleftarrow{1} \overleftarrow{\times}
\overrightarrow{1}$}}

\put(9,4.5){\vector(-2,1){1}} \put(9.7,4.5){\makebox(0,0){$\overleftarrow{0} \overleftarrow{\times}
\overrightarrow{2}$}}

\multiput(1,-7)(0,1){5}{\line(1,0){8}} \multiput(2,-7)(1,0){7}{\line(0,1){4}}

\put(3,-2){\vector(-1,-3){1}} \put(3,-1.5){\makebox(0,0){$\overrightarrow{3}$ poids=s(g-1)}}

\put(3,-8){\vector(1,3){1}} \put(3,-8.5){\makebox(0,0){$\overrightarrow{3}$ poids=s(g-1)+2}}

\put(7,-2){\vector(-1,-3){1}} \put(7,-1.5){\makebox(0,0){$\overrightarrow{3}$ poids=s(g-1)+4}}

\put(7,-8){\vector(1,3){1}} \put(7,-8.5){\makebox(0,0){$\overrightarrow{3}$ poids=s(g+1)-2}}

\multiput(2,-5)(2,0){4}{\circle*{.1}}

\put(8,-5){\vector(1,0){1} $i$} \put(5,-2){\vector(0,-1){7}} \put(4.7,-9){$k$}

\end{picture}

\caption{\label{figure8} $H^i(gr_{k,\r_\oo})[\Pi^{\oo,v}]$ pour $\Pi_v \simeq \speh_4(\pi_v)$ avec $\pi_v$
une représentation cuspidale de $GL_g(F_v)$. Le dessin du bas représente le terme $E_2$ et donc
l'aboutissement de la suite spectrale.}

\end{figure}


\clearpage
\newpage

\bibliographystyle{plain}
\bibliography{bib-ok}

\def\cprime{$'$}
\begin{thebibliography}{10}

\bibitem{ast}
A.~A. Beilinson, J.~Bernstein, and P.~Deligne.
\newblock Faisceaux pervers.
\newblock In {\em Analysis and topology on singular spaces, I (Luminy, 1981)},
  volume 100 of {\em Ast\'erisque}, pages 5--171. Soc. Math. France, Paris,
  1982.

\bibitem{boy}
P~Boyer.
\newblock Mauvaise réduction des variétés de {D}rinfeld et correspondance de
  {L}anglands locale.
\newblock {\em Invent. Math.}, 138(3):573--629, 1999.

\bibitem{ca}
H.~Carayol.
\newblock Nonabelian {L}ubin-{T}ate theory.
\newblock In {\em Automorphic forms, Shimura varieties, and {L}-functions},
  volume~11 of {\em Perspect. Math.}, pages 15--39. Academic Press, Boston, MA,
  1990.

\bibitem{dr1}
V.~G. Drinfel'd.
\newblock Elliptic modules.
\newblock {\em Mat. Sb. (N.S.)}, 136(94):594--627, 1974.

\bibitem{falt}
Gerd Faltings.
\newblock A relation between two moduli spaces studied by {V}. {G}. {D}rinfeld.
\newblock In {\em Algebraic number theory and algebraic geometry}, volume 300
  of {\em Contemp. Math.}, pages 115--129. Amer. Math. Soc., Providence, RI,
  2002.

\bibitem{h-t}
M.~Harris, R.~Taylor.
\newblock {\em The geometry and cohomology of some simple {S}himura varieties},
  volume 151 of {\em Annals of Mathematics Studies}.
\newblock Princeton University Press, Princeton, NJ, 2001.

\bibitem{H-L}
Michael Harris and Jean-Pierre Labesse.
\newblock Conditional base change for unitary groups.
\newblock {\em Asian J. Math.}, 8(4):653--683, 2004.

\bibitem{he}
Guy Henniart.
\newblock Sur la conjecture de {L}anglands locale pour {${\rm GL}\sb n$}.
\newblock {\em J. Th\'eor. Nombres Bordeaux}, 13(1):167--187, 2001.
\newblock 21st Journ\'ees Arithm\'etiques (Rome, 2001).

\bibitem{ito}
Tetsushi Ito.
\newblock Weight-monodromy conjecture for {$p$}-adically uniformized varieties.
\newblock {\em Invent. Math.}, 159(3):607--656, 2005.

\bibitem{Ko7}
R.E. Kottwitz.
\newblock Stable trace formula: elliptic singular terms.
\newblock {\em Math. Ann.}, 275(3):365--399, 1986.

\bibitem{Ko1}
R.E. Kottwitz.
\newblock Points on some {\Large s}himura varieties over finite fields.
\newblock {\em J. Amer. Math. Soc.}, 5(2):373--444, 1992.

\bibitem{lrs}
G.~Laumon, M.~Rapoport, and U.~Stuhler.
\newblock {${\mathcal D}$}-elliptic sheaves and the {L}anglands correspondence.
\newblock {\em Invent. Math.}, 113(2):217--338, 1993.

\bibitem{Lubin1}
J.~Lubin.
\newblock Canonical subgroups of formal groups.
\newblock {\em Trans. Amer. Math. Soc.}, 251:103--127, (1979).

\bibitem{s-s}
P.~Schneider and U.~Stuhler.
\newblock The cohomology of {$p$}-adic symmetric spaces.
\newblock {\em Invent. Math.}, 105(1):47--122, 1991.

\bibitem{se}
Jean-Pierre Serre.
\newblock {\em Corps locaux}.
\newblock Hermann, Paris, 1968.
\newblock Deuxi\`eme \'edition, Publications de l'Universit\'e de Nancago, No.
  VIII.

\bibitem{y-t}
Yoshida~T. Taylor~R.
\newblock Compatibility of local and global langlands correspondences.
\newblock {\em Preprint}, 2005.

\bibitem{ze}
A.~V. Zelevinsky.
\newblock Induced representations of reductive {${p}$}-adic groups. {II}. {O}n
  irreducible representations of {${\rm GL}(n)$}.
\newblock {\em Ann. Sci. \'Ecole Norm. Sup. (4)}, 13(2):165--210, 1980.

\end{thebibliography}

\end{document}